\documentclass{amsbook}
\usepackage{amscd}
\usepackage{amsmath,amsthm}     
\usepackage{amssymb}
\usepackage{euscript}           
\usepackage[all,2cell,ps]{xy} \UseAllTwocells \SilentMatrices
\CompileMatrices		
\UseComputerModernTips		

\newcommand{\comment}[1]{}

\newcommand{\cat}[1]{\EuScript {#1}}

\newcommand{\lra}{\longrightarrow}
\newcommand{\lla}{\longleftarrow}
\newcommand{\ra}{\rightarrow}
\newcommand{\la}{\leftarrow}

\newcommand{\dbra}{\rightrightarrows}           
\newcommand{\Ra}{\Rightarrow}
\newcommand{\La}{\Leftarrow}

\newcommand{\dra}{\xymatrix@C=10pt{ \ar[r]|\circ &}}
\newcommand{\dla}{\xymatrix@C=10pt{ & \ar[l]|\circ }}
\newcommand{\dlra}{\xymatrix@C=20pt{ \ar[r]|\circ &}}
\newcommand{\dlla}{\xymatrix@C=20pt{ & \ar[l]|\circ }}

\newcommand{\tdot}{^{^{_{_\bullet}}}} 	

\newcommand{\sdot}{_{_{^{^\bullet}}}} 	

\newcommand{\comp}{\circ}

\newcommand{\eps}{\epsilon}

\newdir{ >}{{}*!/-5pt/\dir{>}}                  
\newdir{ |}{{}*!/-5pt/\dir{|}}                  

\newcommand{\Sum} {\sqcup}

\DeclareMathOperator*{\colim}{colim\,}
\DeclareMathOperator*{\wcolim}{wcolim\,}
\DeclareMathOperator*{\wlim}{wlim\,}
\DeclareMathOperator*{\initial}{{\bf 0}}
\DeclareMathOperator*{\final}{{\bf 1}}
\DeclareMathOperator*{\terminal}{\final}
\DeclareMathOperator*{\initialcat}{\emptyset}
\DeclareMathOperator*{\finalcat}{{\it e}}
\DeclareMathOperator*{\terminalcat}{\finalcat}

\newcommand{\corollaryref}[1] {Cor. \ref{#1}}
\newcommand{\definitionref}[1] {Def. \ref{#1}}

\newcommand{\lemmaref}[1] {Lemma \ref{#1}}
\newcommand{\propositionref}[1] {Prop. \ref{#1}}
\newcommand{\theoremref}[1] {Thm. \ref{#1}}

\newcommand{\remarkref}[1] {Remark \ref{#1}}
\newcommand{\chapterref}[1] {Chap. \ref{#1}}
\newcommand{\sectionref}[1] {Section \ref{#1}}

\newcommand{\equationref}[1] {(\ref{#1})}

\makeatletter
\newcommand{\interitemtext}[1]{
\begin{list}{}
{\itemindent=0mm\labelsep=0mm
\labelwidth=0mm\leftmargin=0mm
\addtolength{\leftmargin}{-\@totalleftmargin}}
\item #1
\end{list}}
\makeatother

\newcounter{theoremcounter}
\numberwithin{theoremcounter}{section}
\numberwithin{section}{chapter}
\numberwithin{equation}{chapter}

\theoremstyle{plain} 

\newtheorem{thm}[theoremcounter]{Theorem}
\newtheorem{cor}[theoremcounter]{Corollary}
\newtheorem{lem}[theoremcounter]{Lemma}
\newtheorem{prop}[theoremcounter]{Proposition}

\theoremstyle{definition}
\newtheorem{defn}[theoremcounter]{Definition}
\theoremstyle{remark}
\newtheorem{rem}[theoremcounter]{Remark}

\makeindex

\begin{document}
\frontmatter
\title[Cofibrations in Homotopy Theory]{Cofibrations in Homotopy Theory}
\author{Andrei R\u{a}dulescu-Banu}
\address{86 Cedar St \\ Lexington, MA 02421 USA}
\email{andrei@alum.mit.edu}
\dedicatory{Pentru Bianca, Cristina \c{s}i Daria}

\date{\today}


\begin{abstract}
We define Anderson-Brown-Cisinski (ABC) cofibration categories, and construct homotopy colimits of diagrams of objects in ABC cofibration categories. Homotopy colimits for Quillen model categories are obtained as a particular case. We attach to each ABC cofibration category a left Heller derivator. A dual theory is developed for homotopy limits in ABC fibration categories and for right Heller derivators. These constructions provide a natural framework for 'doing homotopy theory' in ABC (co)fibration categories.
\end{abstract}

\maketitle

\setcounter{page}{4}
\tableofcontents

\chapter*{Preface}
\label{chap:Introduction}

Model categories, introduced by Daniel Quillen \cite{Quillen1}, are a remarkably succesful framework for expressing homotopy theoretic ideas in an axiomatic way. A Quillen model category (${\cat M}$, ${\cat W}$, ${\cat Cof}, {\cat Fib}$) consists of a category ${\cat M}$, and three distinguished classes of maps: the weak equivalences ${\cat W}$, the cofibrations ${\cat Cof}$ and the fibrations ${\cat Fib}$, subject to a list of axioms (\definitionref{defn:modelcat}). 

Let us fix some notation. If ${\cat D}$ is a small category we denote ${\cat M}^{\cat D}$ the category of ${\cat D}$-diagrams in ${\cat M}$. For a functor $u \colon {\cat D}_1 \ra {\cat D}_2$ of small categories we denote $u^* \colon {\cat M}^{{\cat D}_2} \ra {\cat M}^{{\cat D}_1}$ the functor defined by $(u^*X)_{d_1} = X_{ud_1}$ for objects $d_1 \eps {\cat D}_1$.  For a category with a class of weak equivalences $({\cat M}, {\cat W})$ we denote ${\bf ho}{\cat M} = {\cat M}[{\cat W}^{-1}]$ its homotopy category, and we will always consider the weak equivalences on ${\cat M}^{\cat D}$ to be {\it pointwise}, i.e. $f$ is a weak equivalence in ${\cat M}^{\cat D}$ if $f_d$ is a weak equivalence in ${\cat M}$ for all objects $d \eps {\cat D}$.

A Quillen model category ${\cat M}$ admits homotopy pushouts and homotopy pullbacks. These are the total left (resp. right) derived functors of the pushout (resp. pullback) functor in ${\cat M}$. If ${\cat M}$ is pointed one can construct the homotopy cofiber and fiber of a map in ${\bf ho}{\cat M}$, and the suspension and loop space of an object in ${\bf ho}{\cat M}$. One can then form the cofibration sequence and the fibration sequence of a map in ${\bf ho}{\cat M}$.

Furthermore, a Quillen model category ${\cat M}$ admits all small homotopy colimits and limits (also called the homotopy left and right Kan extensions). For a functor $u \colon {\cat D}_1 \ra {\cat D}_2$ of small categories, the homotopy colimit ${\bf L}\colim^u \colon {\bf ho}({\cat M}^{{\cat D}_1}) \ra {\bf ho}({\cat M}^{{\cat D}_2})$ is the total left derived functor of the colimit functor $\colim^u$, and is left adjoint to ${\bf ho}u^* \colon {\bf ho}({\cat M}^{{\cat D}_2}) \ra {\bf ho}({\cat M}^{{\cat D}_1})$. Dually, the homotopy limit ${\bf R}\lim^u$ is the total right derived functor of $\lim^u$, and is right adjoint to ${\bf ho}u^*$.

For the general construction of homotopy colimits, we recommend the work of Dwyer, Kan, Hirschhorn and Smith \cite{Dwyer-Kan-Hirschhorn-Smith}, Hirschhorn \cite{Hirschhorn}, Weibel \cite{Weibel} (after Thomason's unpublished notes), Chach\'{o}lski and Scherer \cite{ChacScherer} and Cisinski \cite{Cisinski4}. Each of these references presents a different perspective on homotopy (co)limits. It is apparent from these sources that in a Quillen model category the homotopy colimits can actually be constructed just manipulating weak equivalences and cofibrations. 

\bigskip

To summarize, homotopy colimits are the total left derived functors of the colimit, so their {\it definition} requires just the presence of weak equivalences. The {\it existence} of homotopy colimits can be proved in presence of the Quillen model category axioms, where their {\it construction} can be performed using just cofibrations and weak equivalences.

It is therefore natural to ask if one can simplify Quillen's axioms, and separate a minimal set of axioms required by cofibrations and weak equivalences in order to still be able to prove existence of homotopy colimits. Along the way, this leads us to rethink the role of cofibrations in abstract homotopy theory.

\bigskip

In a very influential paper, Ken Brown \cite{Brown} formalized some of these obsevations by defining categories of fibrant objects and working out in detail their properties. Reversing arrows, one defines categories of cofibrant objects, and Brown's work carries over by duality to categories of cofibrant objects. 

A category of cofibrant objects (${\cat M}$, ${\cat W}$, ${\cat Cof}$) consists of a category ${\cat M}$, the class of weak equivalences ${\cat W}$ and the class of cofibrations ${\cat Cof}$, subject to a list of axioms (the duals of the axioms of \cite{Brown}). These axioms require in particular all objects to be cofibrant.

For a pointed category of cofibrant objects, Brown was able to construct homotopy cofibers of maps, suspensions of objects and the cofibration sequence of a map. These constructions exist in dual form for categories of fibrant objects.

\bigskip

Building on Brown's work, Don Anderson \cite{Anderson1} extended Brown's axioms for a category of cofibrant objects by dropping the requirement that all objects be cofibrant. Anderson called the categories defined by his new axioms {\it left homotopical}; our text changes terminology and calls them {\it Anderson-Brown-Cisinski cofibration categories} (or just cofibration categories for simplicity). The cofibration category axioms we use are slightly more general than Anderson's.

Anderson's main observation was that the cofibration category axioms on ${\cat M}$ suffice for the construction of a left adjoint of ${\bf ho}u^*$, for any functor $u$ of small categories. It is implicit in his work that the left adjoint of ${\bf ho}u^*$ is a left derived of $\colim^u$.

Unfortunately for the history of this subject, Anderson's paper \cite{Anderson1} contains statements but omits proofs, and has a title (``Fibrations and Geometric Realizations``) that does not reflect the generality of his work. Also, Anderson quit mathematics shortly after his paper was published, the proofs of \cite{Anderson1} got lost and as a result his whole theory laid dormant for twenty five years.

\bigskip

We can be grateful to Denis-Charles Cisinski \cite{Cisinski2}, \cite{Cisinski4} for bringing back to light Brown and Anderson's ideas. Cisinski simplifies Anderson's arguments, and provides for cofibration categories a complete construction of homotopy colimits along functors $u \colon {\cat D}_1 \ra {\cat D}_2$ with ${\cat D}_1$ finite and direct. 

Cisinski has also worked out the construction of homotopy colimits along arbitrary functors $u$ of small categories, as well as the end result regarding the derivability of cofibration categories (our \chapterref{chap:derivators}). While this part of his work remains unpublished, he was kind enough to share with me its outline. I would like to thank him for suggesting the correct formulation of axioms CF5-CF6, and for patiently explaining to me the finer points of excision. 

\bigskip

The goal of these notes is then to work out an account of homotopy colimits from the axioms of an Anderson-Brown-Cisinski cofibration category, and show that they satisfy the axioms of a left Heller derivator. There are a number of properties of homotopy colimits that are a formal consequence of the left Heller derivator axioms, but they are outside of the scope of our text. We will instead try to investigate the relation with the better-known Quillen model categories, and compare with other axiomatizations that have been proposed for cofibration categories. 

While some of the proofs we propose may be new, the credit for this theory should go entirely to Brown, Anderson and Cisinski. It was our choice in this text to make use of approximation functors and abstract Quillen equivalences, and for that we were influenced by the work of Dwyer, Kan, Hirschhorn and Smith \cite{Dwyer-Kan-Hirschhorn-Smith}. Our treatment of direct categories bears the influence of Daniel Kan's theory of Reedy categories outlined in \cite{Dwyer-Kan-Hirschhorn}, \cite{Hovey} and \cite{Hirschhorn}.

\bigskip

I would like to thank Haynes Miller and Daniel Kan, my mentors in abstract homotopy, for their gracious support and encouragement. I am grateful to Denis-Charles Cisinski, Philip Hirschhorn and Haynes Miller for the conversations we had on the subject of this text.

\subsection*{Comparing ABC and Quillen model categories}
\label{subsec:compareabcquillen}

Anderson-Brown-Cisinski (ABC) cofibration categories (${\cat M}$, ${\cat W}$, ${\cat Cof}$) are defined in \sectionref{sec:CofibrationCategoriesFibrationCategoriesAxioms}. We will make a distinction between ABC precofibration categories (satisfying axioms CF1-CF4) and ABC cofibration categories (satisfying the full set of axioms CF1-CF6). 

ABC fibration categories (${\cat M}$, ${\cat W}$, ${\cat Fib}$) are defined by duality, and ABC model categories (${\cat M}$, ${\cat W}$, ${\cat Cof}$, ${\cat Fib}$) by definition carry an ABC cofibration and fibration category structure.

Any Quillen model category is an ABC model category, but the class of ABC cofibration categories behaves differently in some respects:
\begin{enumerate}
\item
If (${\cat M}$, ${\cat W}$) admits a Quillen model category structure, the choice of cofibrations determines the fibrations, and viceversa. This is not true of an ABC cofibration category structure, where the choice of cofibrations is completely independent of the choice of fibrations.
\item
If (${\cat M}$, ${\cat W}$, ${\cat Cof}$) is an ABC cofibration category, then so is (${\cat M}^{\cat D}$, ${\cat W}^{\cat D}$, ${\cat Cof}^{\cat D}$) for a small category ${\cat D}$. Quillen model categories have this property only in particular cases, for example if they are cofibrantly generated.
\item
If ${\cat M}$ is a locally small Quillen model category, then ${\bf ho}{\cat M}$ is also locally small. ABC cofibration categories do not have this property.
\item
If (${\cat M}$, ${\cat W}$, ${\cat Cof}$, ${\cat Fib}$) is a Quillen model category, then ${\cat W}$ is saturated ${\cat W} = \overline{\cat W}$. Similarly, if (${\cat M}$, ${\cat W}$, ${\cat Cof}$) is an ABC cofibration category then ${\cat W}$ is saturated. If (${\cat M}$, ${\cat W}$, ${\cat Cof}$) is only a precofibration category, then ${\cat W}$ is not necessarily saturated, but (${\cat M}$, $\overline{\cat W}$, ${\cat Cof}$) is still an ABC precofibration category.
\item
If (${\cat M}$, ${\cat W}$, ${\cat Cof}$) is a left proper ABC cofibration category, one has a {\it maximal} left proper CF1-CF4 cofibration category structure on (${\cat M}$, ${\cat W}$), taking as cofibrations the left proper maps ${\cat PrCof}$. In this context, the weak equivalences therefore determine the left proper maps as a preferred class of CF1-CF4 cofibrations.
\end{enumerate}

\subsection*{Outline of the text}
\label{subsec:outline}

We start \chapterref{chap:CofibrationAndFibrationCategories} with the ABC (pre)cofibration category axioms, we proceed with their elementary properties, and end the chapter with a discussion of proper ABC cofibration categories. For a cofibration category ${\cat M}$, we will denote ${\cat M}_{cof}$ its full subcategory of cofibrant objects.

Going to \chapterref{chap:modcatcatcofobj}, we compare ABC cofibration categories with other axiomatic systems that have been proposed for categories with cofibrations. Most notably, we show that any Quillen model category is an ABC model category.

Topological spaces, with homotopy equivalences as weak equivalences, Hurewicz cofibrations as cofibrations, and Hurewicz fibrations as fibrations, form an ABC model category. Complexes of objects in an abelian category ${\cat A}$, with quasi-isomorphisms as weak equivaleces, and monics as cofibrations form an ABC precofibration category, which satisfies CF5 if ${\cat A}$ satisfies the Grothendieck axiom AB4, and satisfies CF6 if ${\cat A}$ satisfies the Grothendieck axiom AB5. These examples are explained in \chapterref{chap:examples}.

In \chapterref{chap:kanextensions} we recall the language of 2-categories and that of Kan extensions. This is a prerequisite for \chapterref{chap:catweq}, where we introduce axiomatically the {\it left approximation} functors $t \colon ({\cat M}^{'}, {\cat W}^{'}) \ra ({\cat M}, {\cat W})$ between two category pairs. The Approximation \theoremref{thm:approxthm} states that left approximation functors induce an equivalence of categories ${\bf ho}t \colon {\bf ho}{\cat M}^{'} \ra {\bf ho}{\cat M}$. We prove an existence result for total derived functors (\theoremref{thm:generalexistenceleftderived}), and a Quillen-type adjunction property between the total derived functors of an adjoint pair of functors (\theoremref{thm:generalexistencetotalderivedadjoint2}). 

In \chapterref{chap:homotcat}, we take as model the inclusion ${\cat M}_{cof} \ra {\cat M}$ and define {\it cofibrant approximation} functors $t \colon {\cat M}^{'} \ra {\cat M}$ between precofibration categories. Cofibrant approximation functors are in particular left approximations, and by the Approximation \theoremref{thm:approxthm} they induce an equivalence of categories ${\bf ho}t$. In particular, this shows that ${\bf ho}{\cat M}_{cof} \ra {\bf ho}{\cat M}$ is an equivalence of categories.

In the second half of \chapterref{chap:homotcat} we then recall the theory of cylinders and left homotopic maps in a precofibration category, and develop the formalism of homotopy calculus of fractions. Going to \chapterref{chap:homcalsappl}, we give some applications of the formalism of homotopy calculus of fractions.

In \chapterref{chap:constrcat}, we recall over and under categories, and elementary properties of limits and colimits.

In \chapterref{chap:limcolimcofinal}, we reach our main objective. We define the homotopy colimit ${\bf L}\colim^u \colon {\bf ho}({\cat M}^{{\cat D}_1})$ $\ra$ ${\bf ho}({\cat M}^{{\cat D}_2})$ as a left derived functor, and show that it exists and it is left adjoint to ${\bf ho}u^* \colon {\bf ho}({\cat M}^{{\cat D}_2})$ $\ra$ ${\bf ho}({\cat M}^{{\cat D}_1})$. We show that if ${\cat D}$ is a small category and (${\cat M}$, ${\cat W}$, ${\cat Cof}$) is a cofibration category, then the diagram category ${\cat M}^{\cat D}$ is again a cofibration category, with pointwise weak equivalences and pointwise cofibrations.

Finally, in \chapterref{chap:derivators} we recall the notion of a {\it left Heller derivator}, and show that the homotopy colimits in a cofibration category satisfy the axioms of a left Heller derivator. The purpose of the last chapter is simply to assert the results of \chapterref{chap:limcolimcofinal} within the axiomatic language of derivators.

\mainmatter

\chapter{Cofibration categories}
\label{chap:CofibrationAndFibrationCategories}

This chapter defines Anderson-Brown-Cisinski (or ABC) cofibration, fibration and model categories. For simplicity, we refer to ABC cofibration (fibration, model) categories as just cofibration (fibration, model) categories, when no confusion with Quillen model categories or Baues cofibration categories is possible.

What is an ABC cofibration category? It is a category ${\cat M}$ with two distinguished classes of maps, the weak equivalences and the cofibrations, satisfying a set of axioms denoted CF1-CF6. An ABC fibration category is a category ${\cat M}$ with weak equivalences and fibrations, satisfying the dual axioms F1-F6. An ABC model category is a category ${\cat M}$ with weak equivalences, cofibrations and fibrations that is at the same time a cofibration and a fibration category.

Any Quillen model category is an ABC model category (\propositionref{prop:relclosedmodcat}). Diagrams indexed by a small category in an ABC cofibration category form again a cofibration category (\theoremref{thm:generalpointwisecofstructure}), a property not enjoyed in general by Quillen model categories.

The ultimate goal using the cofibration category axioms CF1-CF6 is to construct (in \chapterref{chap:limcolimcofinal}) homotopy colimits in ${\cat M}$ indexed by small categories ${\cat D}$, and more generally to construct 'relative' homotopy colimits along functors $u \colon {\cat D}_1 \ra {\cat D}_2$ of small catgories.

The category ${\cat M}$ will not be assumed in general to be cocomplete. Under the simplifying assumption that ${\cat M}$ is cocomplete, however, the homotopy colimit along a functor $u \colon {\cat D}_1 \ra {\cat D}_2$ is the total left derived functor of $\colim^u \colon {\cat M}^{{\cat D}_1} \ra {\cat M}^{{\cat D}_2}$. It is a known fact that a colimit indexed by a small category ${\cat D}$ can be constructed in terms of pushouts of small sums of objects in ${\cat M}$. Furthermore, a 'relative' colimit $\colim^u$ can be described in terms of absolute colimits indexed by the over categories $(u \downarrow d_2)$ for $d_2 \eps {\cat D}_2$ (see \lemmaref{lem::computerelcolim}).

It is perhaps not surprising then that the cofibration category axioms specify an approximation property of maps by cofibrations (axiom CF4), as well as the behaviour of cofibrations under pushouts (axiom CF3) and under small sums (axiom CF5). The axiom CF6 has a technical rather than conceptual motivation - we simply need it to make the rest of our arguments work (but see \sectionref{sec:oddsandends}).

A large part of the theory can be developed actually from a subset of the axioms, namely the axioms CF1-CF4. The theory of homotopic maps, of homotopy calculus of fractions and of cofibrant approxiomation functors of \chapterref{chap:homotcat} only requires this smaller set of axioms. A category with weak equivalences and cofibrations satisfying the axioms CF1-CF4 will be called a {\it precofibration} category. From the precofibration category axioms, it turns out that one can construct all the homotopy colimits indexed by {\it finite, direct} categories \cite{Cisinski2}. 

One idea worth repeating is that while we need both the cofibrations and the weak equivalences to {\it construct} homotopy colimits, the homotopy colimits are {\it characterized} in the end just by the weak equivalences. So when working with a cofibration category with a fixed class of weak equivalences, it would be desirable to have a class of cofibrations as large as possible. This is where the concept of {\it left proper} maps becomes useful (see \sectionref{sec:properness}). In a {\it left proper} precofibration category $({\cat M}$, ${\cat W}$, ${\cat Cof})$, any cofibration $A \ra B$ with $A$ cofibrant is a left proper map, and the class of left proper maps denoted ${\cat PrCof}$ yields again a precofibration category structure $({\cat M}$, ${\cat W}$, ${\cat PrCof})$. But if ${\cat M}$ is a left proper CF1-CF6 cofibration category, then $({\cat M}$, ${\cat W}$, ${\cat PrCof})$ may not necessarily satisfy the axioms CF5 and CF6. Dual results hold for right proper prefibration categories.

\section{The axioms}
\label{sec:CofibrationCategoriesFibrationCategoriesAxioms}

\index{ABC, or Anderson-Brown-Cisinski!(co)fibration category}
\begin{defn}[Anderson-Brown-Cisinski cofibration categories]
\label{defn:cofcat}
\mbox{} \\
An ABC cofibration category (${\cat M}$, ${\cat W}$, ${\cat Cof}$) consists of a category ${\cat M}$, and two distinguished classes of maps of ${\cat M}$ - the weak equivalences (or trivial maps) ${\cat W}$ and the cofibrations ${\cat Cof}$, subject to the axioms CF1-CF6. The initial object $\initial$ of ${\cat M}$ exists by axiom CF1, and an object $A$ is called {\it cofibrant} if the map $\initial \ra A$ is a cofibration.

The axioms are:
\begin{description}
\item[CF1]
${\cat M}$ has an initial object $\initial$, which is cofibrant. Cofibrations are stable under composition. All isomorphisms of ${\cat M}$ are weak equivalences, and all isomorphisms with the domain a cofibrant object are trivial cofibrations. 

\item[CF2]
\index{axiom!2 out of 3 axiom}(Two out of three axiom) Suppose $f$ and $g$ are maps such that $gf$ is defined. If two of $f$, $g$, $gf$ are weak equivalences, then so is the third.

\item[CF3]
\index{axiom!pushout axiom}(Pushout axiom) Given a solid diagram in ${\cat M}$, with $i$ a cofibration and  $A$, $C$ cofibrant
\begin{center}
$ \xymatrix {
	A \ar[r] \ar@{>->}[d]_i &
	C \ar@{>-->}[d]^j \\
	B \ar@{-->}[r] &
	D
	}$
\end{center}
then
\begin{enumerate}
\item
the pushout exists in ${\cat M}$ and $j$ is a cofibration, and
\item
if additionally $i$ is a trivial cofibration, then so is $j$.
\end{enumerate}

\item[CF4]
\index{axiom!factorization axiom}(Factorization axiom) Any map $f \colon A \ra B$ in ${\cat M}$ with $A$ cofibrant factors as $f = rf'$, with $f'$ a cofibration and $r$  a weak equivalence

\item[CF5]
If $f_i \colon A_i \ra B_i$ for $i \eps I$ is a set of cofibrations with $A_i$ cofibrant, then
\begin{enumerate}
\item
$\Sum A_i$, $\Sum B_i$ exist and are cofibrant, and $\Sum f_i$ is a cofibration.
\item
if additionally all $f_i$ are trivial cofibrations, then so is $\Sum f_i$.
\end{enumerate}

\item[CF6]
For any countable direct sequence of cofibrations with $A_0$ cofibrant
\begin{center}
$ \xymatrix {
        A_0 \ar@{>->}[r]^{a_{0}} &
        A_1 \ar@{>->}[r]^{a_{1}} &
        A_2 \ar@{>->}[r]^{a_{2}} &
	... \\
	}$
\end{center}
\begin{enumerate}
\item
the colimit object $\colim A_n$ exists and the transfinite composition $A_0 \ra \colim A_n$ is a cofibration.
\item
if additionally all $a_i$ are trivial cofibrations, then so is $A_0 \ra \colim A_n$.
\end{enumerate}
\end{description}

If (${\cat M}$, ${\cat W}$, ${\cat Cof}$) only satisfies the axioms CF1-CF4, it is called a precofibration category.
\end{defn}

Pushouts are defined by an universal property, and are only defined up to an unique isomorphism. Since all isomorphisms with cofibrant domain are trivial cofibrations, it does not matter which isomorphic representative of the pushout we choose in CF3.

In the axiom CF4, the map $f^{'}$ is called a {\it cofibrant replacement} of $f$. If $r \colon A^{'} \ra A$ is a weak equivalence with $A^{'}$ cofibrant, the object $A^{'}$ is called a {\it cofibrant replacement} of $A$.

If $A \ra B$ is a cofibration with $A$ cofibrant then $B$ is cofibrant. But there may exist cofibrations $A \ra B$ with $A$ not cofibrant. If we denote ${\cat Cof}^{'}$ the class of cofibrations $A \ra B$ with $A$ cofibrant, then (${\cat M}$, ${\cat W}$, ${\cat Cof}^{'}$) is again a cofibration category. 

We will sometimes refer to a cofibration category as just ${\cat M}$. We will also denote ${\cat M}_{cof}$ \index{${\cat M}_{cof}$, ${\cat M}_{fib}$} the full subcategory of cofibrant objects of ${\cat M}$. 

The category ${\cat M}_{cof}$ is a cofibration category, and in fact so is any full subcategory ${\cat M}^{'}$ of ${\cat M}$ that includes ${\cat M}_{cof}$, with the induced structure (${\cat M}^{'}$, ${\cat W} \cap {\cat M}^{'}$, ${\cat Cof} \cap {\cat M}^{'}$). 

If ${\cat M}$ is a precofibration category, then (${\cat M}$, ${\cat W}$, ${\cat Cof}^{'}$), ${\cat M}_{cof}$ and any ${\cat M}^{'}$ as above are precofibration categories. We will also refer to precofibration categories as CF1-CF4 cofibration categories.

\begin{defn}[Anderson-Brown-Cisinski fibration categories]
\label{defn:fibcat}
\mbox{} \\
An ABC fibration category (${\cat M}$, ${\cat W}$, ${\cat Fib}$) consists of a category ${\cat M}$, and two distinguished classes of maps - the weak equivalences ${\cat W}$ and the fibrations ${\cat Fib}$, subject to the axioms F1-F6. The terminal object $\terminal$ of ${\cat M}$ exists by axiom F1, and an object $A$ of ${\cat Fib}$ is called {\it fibrant} if the map $A \ra \terminal$ is a fibration. 

The axioms are:
\begin{description}
\item[F1]
${\cat M}$ has a final object $\final$, which is fibrant. Fibrations are stable under composition. All isomorphisms of ${\cat M}$ are weak equivalences, and all isomorphisms with fibrant codomain are trivial fibrations. 

\item[F2]
(Two out of three axiom) Suppose $f$ and $g$ are maps such that $gf$ is defined. If two of $f$, $g$, $gf$ are weak equivalences, then so is the third.

\item[F3]
(Pullback axiom) Given a solid diagram in ${\cat M}$, with $p$ a fibration and $A$, $C$ fibrant,
\begin{center}
$ \xymatrix {
	D \ar@{-->}[r] \ar@{-->>}[d]_q &
	B \ar@{->>}[d]^p \\
	C \ar[r] &
	A
	}$
\end{center}
then
\begin{enumerate}
\item
the pullback exists in ${\cat M}$ and $q$ is a fibration, and
\item
if additionally $p$ is a trivial fibration, then so is $q$.
\end{enumerate}

\item[F4]
(Factorization axiom) Any map $f \colon A \ra B$ in ${\cat M}$ with $B$ fibrant factors as $f = f's$, with $s$ a weak equivalence and $f'$  a fibration.

\item[F5]
If $f_i \colon A_i \ra B_i$ for $i \eps I$ is a set of fibrations with $B_i$ fibrant, then 
\begin{enumerate}
\item
$\times A_i$, $\times B_i$ exist and are fibrant, and $\times f_i$ is a fibration
\item
if additionally all $f_i$ are trivial fibrations, then so is $\times f_i$.
\end{enumerate}

\item[F6]
For any countable inverse sequence of fibrations with $A_0$ fibrant
\begin{center}
$ \xymatrix {
    ... \ar@{->>}[r]^{a_{2}} &
    A_2 \ar@{->>}[r]^{a_{1}} &
    A_1 \ar@{->>}[r]^{a_{0}} &
    A_0
  }$
\end{center}
\begin{enumerate}
\item
the limit object $\lim A_i$ exists and the transfinite composition $\lim A_n \ra A_0$ is a fibration
\item
if additionally all $a_i$ are trivial fibrations, then so is $\lim A_n \ra A_0$.
\end{enumerate}
\end{description}

If (${\cat M}$, ${\cat W}$, ${\cat Fib}$) only satisfies the axioms F1-F4, it is called a prefibration category.
\end{defn}

The axioms are dual in the sense that (${\cat M}$, ${\cat W}$, ${\cat Cof}$) is a cofibration category if and only if (${\cat M}^{op}, {\cat W}^{op}, {\cat Cof}^{op}$) is a fibration category. 

In the axiom F4, the map $f^{'}$ is called a {\it fibrant replacement} of $f$. If $r \colon A \ra A^{'}$ is a weak equivalence with $A^{'}$ fibrant, the object $A^{'}$ is called a {\it fibrant replacement} of $A$.

If we denote ${\cat Fib}^{'}$ the class of fibrations $A \ra B$ with $B$ fibrant then (${\cat M}$, ${\cat W}$, ${\cat Fib}^{'}$) again is a fibration category.

We will denote ${\cat M}_{fib}$ to be the full subcategory of fibrant objects of a fibration category ${\cat M}$.  The category ${\cat M}_{fib}$ as well as any full subcategory ${\cat M}^{'}$ of ${\cat M}$ that includes ${\cat M}_{fib}$ satisfy again the axioms of a fibration category. 

If ${\cat M}$ is a prefibration category, then so are (${\cat M}$, ${\cat W}$, ${\cat Fib}^{'}$), ${\cat M}_{fib}$ and any ${\cat M}^{'}$ as above.

\begin{defn}[Anderson-Brown-Cisinski model categories]
\label{defn:abcmodel}
\mbox{} \\
An ABC model category \index{ABC, or Anderson-Brown-Cisinski!model category} $({\cat M}$, ${\cat W}$, ${\cat Cof}, {\cat Fib})$ consists of a category ${\cat M}$ and three distinguished classes of maps ${\cat W}, {\cat Cof}, {\cat Fib}$ with the property that $({\cat M}$, ${\cat W}$, ${\cat Cof})$ is an ABC cofibration category and that $({\cat M}$, ${\cat W}$, ${\cat Fib})$ is an ABC fibration category.
\end{defn}

We will say that ${\cat M}$ is an ABC premodel category if it is only an ABC precofibration and prefibration category.

In all sections of this chapter except \sectionref{sec:transfincomp}, we will do our work assuming only the precofibration category axioms CF1-CF4 (and dually the prefibration category axioms F1-F4). In \sectionref{sec:transfincomp}, we will assume that the full set of axioms is verified.

\section{Sums and products of objects}
\label{sec:sumsproducts}

In general, the objects of a precofibration category are not closed under finite sums. But finite sums of cofibrant objects exist and are cofibrant. Dually, in a prefibration category finite products of fibrant objects exist and are fibrant. In fact we can prove the slightly more general statement:

\begin{lem}
\label{lem:cofsums}
\mbox{}
\begin{enumerate}
\item
Suppose that ${\cat M}$ is a precofibration category. If $f_i \colon A_i \ra B_i$ for $i = 0, ..., n$ are cofibrations with $A_i$ cofibrant, then $\Sum A_i$, $\Sum B_i$ exist and are cofibrant, and $\Sum f_i$ is a cofibration which is trivial if all $f_i$ are trivial.
\item
Suppose that ${\cat M}$ is a prefibration category. If $f_i \colon A_i \ra B_i$ for $i = 0, ..., n$ are fibrations with $B_i$ fibrant, then $\times A_i$, $\times B_i$ exist and are fibrant, and $\times f_i$ is a fibration which is trivial if all $f_i$ are trivial.
\end{enumerate}
\end{lem}

\begin{proof}
We will prove (1), and observe that statement (2) is dual to (1). Using induction on $n$, we can reduce the problem to two maps $f_0 \colon A_0 \ra B_0$ and $f_1 ; A_1 \ra B_1$. If we prove the statement for $f_0, 1_{A_1}$ and $1_{B_0}, f_1$ then the statement follows for $f_0, f_1$. So it suffices to show that if $f \colon A \ra B$ is a (trivial) cofibration and $A$, $C$ are cofibrant, then $A \Sum C$, $B \Sum C$ exist and are cofibrant and $f \Sum 1_C$ is a (trivial) cofibration. From axiom CF3 (1) applied to 
\begin{center}
$ \xymatrix {
	\initial \ar@{>->}[r] \ar@{>->}[d] &
	C \ar@{>->}[d] \\
	A \ar@{>->}[r] &
	A \Sum C
	}$
\end{center}
we see that $A \Sum C$ exists and is cofibrant, and $A \ra A \Sum C$ is a cofibration. Similarly, $B \Sum C$ exists and is cofibrant, and from CF3 applied to
\begin{center}
$ \xymatrix {
	A \ar@{>->}[r] \ar@{>->}[d]_{f} &
	A \Sum C \ar@{>->}[d]^{f \Sum 1_C} \\
	B \ar@{>->}[r] & 
	B \Sum C
	}$
\end{center}
we see that $f \Sum 1_C$ is a cofibration which is trivial if $f$ is trivial.
\end{proof}

\section{Factorization lemmas}
\label{sec:factlemmas}
The Brown Factorization Lemma is an improvement of the factorization axiom CF4 for maps between cofibrant objects.

\index{lemma!Brown factorization}
\begin{lem}[Brown factorization, \cite{Brown}]
\label{lem:brownfact}
\mbox{}
\begin{enumerate}
\item
Let ${\cat M}$ be a precofibration category, and $f \colon A \ra B$ be a map between cofibrant objects. Then $f$ factors as $f = rf'$, where $f'$ is a cofibration and $r$ is a left inverse to a trivial cofibration. 
\item
Let ${\cat M}$ be a prefibration category, and $f \colon A \ra B$ be a map between fibrant objects. Then $f$ factors as $f = f's$, where $f'$ is a fibration and $s$ is a right inverse to a trivial fibration.
\end{enumerate}
\end{lem}

\begin{proof}
The statements are dual, so it suffices to prove (1). We need to construct $f'$, $r$ and $s$ with $f = rf'$ and $rs = 1_B$. 

If we apply the factorization axiom to $f + 1_B$, we get a diagram
\begin{center}
$ \xymatrix {
        A \Sum B \ar[rr]^{f + 1_{B}} \ar@{>->}[rd]_{f^{'} + s} &  & B & \\
        & B^{'} \ar[ru]_{r}^\sim & &
	}$
\end{center}
Since $f^{'} + s$ is a cofibration and $A, B$ are cofibrant, the maps $f^{'}$  and $s$ are cofibrations. The map $r$ is a weak equivalence, and from the commutativity of the diagram we have $rs = 1_B$, therefore $s$ is also a weak equivalence.
\end{proof}

\begin{rem}
\label{rem:brownfact}
We have in fact proved a stronger statement. We have shown that any map $f \colon A \ra B$ between cofibrant objects in a precofibration category factors as $f = rf'$, with $rs = 1_B$ where $f'$, $f^{'} + s$ are cofibrations and $s$ is a trivial cofibration.

Dually, any map $f \colon A \ra B$ between fibrant objects in a prefibration category factors as $f = f's$, with $rs = 1_A$ where $f'$ and $(f^{'}, r)$ are fibrations and $r$ is a trivial fibration.
\end{rem}

Next lemma is a relative version of the factorization axiom CF4 (resp. F4). 

\begin{lem}[Relative factorization of maps]
\label{lem:relativefactorization}
\mbox{}
\begin{enumerate}
\item
Let ${\cat M}$ be a precofibration category, and let
\begin{center}
$\xymatrix {
    A_1 \ar[r]^{f_1} \ar[d]_a & B_1 \ar[d]^b \\
    A_2 \ar[r]^{f_2} & B_2
  }$
\end{center}
be a commutative diagram with $A_1, A_2$ cofibrant. Suppose that $f_1 = r_1 f^{'}_1$ is a factorization of $f_1$ as a cofibration followed by a weak equivalence. Then there exists a commutative diagram
\begin{center}
$\xymatrix{     
    A_1 \ar@{>->}[r]^{f^{'}_1} \ar[d]_a & A^{'}_1 \ar[r]^{r_1}_\sim \ar[d]_(.44){a^{'}} & B_1 \ar[d]^b \\
    A_2 \ar@{>->}[r]^{f^{'}_2} & A^{'}_2 \ar[r]^{r_2}_\sim & B_2
}$ 
\end{center}
where $r_2 f^{'}_2$ is a factorization of $f_2$ as a cofibration followed by a weak equivalence and such that $A_2 \Sum_{A_1} A^{'}_1 \ra A^{'}_2$ is a cofibration.
\item
The dual of (1) holds for prefibration categories.
\end{enumerate}
\end{lem}

\begin{proof}
To prove (1), in the commutative diagram
\begin{center}
$\xymatrix {
    A_1 \ar[d]_{a} \ar@{>->}[r]^{f^{'}_1} &
    A^{'}_1 \ar[d] \ar[rr]_{\sim}^{r_1} & &
    B_1 \ar[d]^b \\
    A_2 \ar@{>->}[r] &
    A_2 \Sum_{A_1} A^{'}_1 \ar@{>->}[r]^-s &
    A^{'}_2 \ar[r]_{\sim}^{r_2} &
    B_2
  }$
\end{center}
the pushout $A_2 \Sum_{A_1} A^{'}_1$ exists by CF3, and we construct the cofibration $s$ and the weak equivalence $r_2$ using the factorization axiom CF4 applied to $A_2 \Sum_{A_1} A^{'}_1 \lra B_2$. 
\end{proof}

The Brown Factorization Lemma has the following relative version:

\begin{lem}[Relative Brown factorization]
\label{lem:relbrownfact}
\mbox{}
\begin{enumerate}
\item
Suppose that ${\cat M}$ is a precofibration category, and that
\begin{center}
$ \xymatrix {
        A_1 \ar[r]^{f_1} \ar[d]_a &
	B_1 \ar[d]^b \\
	A_2 \ar[r]^{f_2} &
	B_2
	}$
\end{center}
is a commutative diagram with cofibrant objects. Suppose that $f_1 = r_1f^{'}_1$, $r_1 s_1 = 1$ is a Brown factorization of $f_1$, with $f^{'}_1$, $f^{'}_1 + s_1$ cofibrations and $s_1$ a trivial cofibration. Then there exists a Brown factorization $f_2 = r_2f^{'}_2$, $r_2 s_2 = 1 $ with $f^{'}_2$, $f^{'}_2 + s_2$ cofibrations and $s_2$ a trivial cofibration and a map $b^{'}$ such that in the diagram
\begin{center}
$ \xymatrix {
        A_1 \ar@{>->}[r]^{f^{'}_1} \ar[d]_a &
        B^{'}_1 \ar@<1ex>[r]^{r_1} \ar[d]^{b^{'}} &
	B_1 \ar[d]^b \ar@<1ex>@{>->}[l]^{s_1}_{\sim} \\
	A_2 \ar@{>->}[r]^{f^{'}_2} &
	B^{'}_2 \ar@<1ex>[r]^{r_2} &
	B_2 \ar@<1ex>@{>->}[l]^{s_2}_{\sim}
	}$
\end{center}
we have that $b^{'}f^{'}_1 = f^{'}_2 a$, $br_1 = r_2 b^{'}$, $b^{'}s_1 = s_2 b$,  and that
\begin{center}
$A_2 \Sum_{A_1} B^{'}_1 \ra B^{'}_2$ and $B_2 \Sum_{B_1} B^{'}_1 \ra B^{'}_2$
\end{center}
are a cofibration (resp. a trivial cofibration).
\item
The dual of (1) holds for prefibration categories.
\end{enumerate}
\end{lem}

\begin{proof}
To prove (1), denote $A_3 = A_2 \Sum_{A_1} B^{'}_1$ and $B_3 = B_2 \Sum_{B_1} B^{'}_1$. 
\begin{center}
$ \xymatrix {
        A_1 \ar@{>->}[rr]^{f^{'}_1} \ar[d]_a &&
        B^{'}_1 \ar@<1ex>[rr]^{r_1} \ar[dr] \ar[dl] &&
	B_1 \ar[d]^b \ar@<1ex>@{>->}[ll]^{s_1}_{\sim} \\
	A_2 \ar@{>->}[r]^{f^{'}} &
	A_3 \ar[rr]^{f_3} && B_3 \ar@<1ex>[r]^{r} &
	B_2 \ar@<1ex>@{>->}[l]^{s}_{\sim}
	}$
\end{center}
We apply \lemmaref{lem:relativefactorization} to the commutative diagram
\begin{center}
$\xymatrix {
    A_1 \Sum B_1 \ar@{>->}[r]^-{f^{'}_1 + s_1} \ar[d]_{f^{'} a \Sum s b} & B^{'}_1 \ar[r]^{r_1}_\sim & B_1 \ar[d]^{sb} \\
    A_3 \Sum B_3 \ar[rr]^-{f_3 + 1_{B_3}} & & B_3
  }$
\end{center}
and we construct a commutative diagram
\begin{center}
$ \xymatrix {
        A_1 \Sum B_1 \ar[rr]^{f_1 + 1_{B_1}} \ar[dd]_{f^{'} a \Sum s b} \ar@{>->}[rd]_{f^{'}_1 + s_1} &  & B_1 \ar[dd]^{sb} & \\
        & B^{'}_1 \ar[ru]_{r_1}^\sim \ar[dd]^(.3){b^{'}} & & \\
        A_3 \Sum B_3 \ar'[r][rr] \ar@{>->}[rd]_{f^{'}_3 + s_3} & & B_3 & \\
        & B^{'}_2 \ar[ru]_{r_3}^\sim & &
	}$
\end{center}
We now set $f^{'}_2 = f^{'}_3f^{'}$, $r_2 = rr_3$ and $s_2 = s_3s$.
\end{proof}

\section{Extension lemmas}
\label{sec:extlemmas}
The Gluing Lemma describes the behavior of cofibrations and weak equivalences under pushouts, and is one of the basic building blocks we will employ in the construction of homotopy colimits.

\index{lemma!gluing lemma}
\begin{lem}[Gluing Lemma]
\label{lem:gluinglemma1}
\mbox{}
\begin{enumerate}
\item
Let ${\cat M}$ be a precofibration category. In the diagram
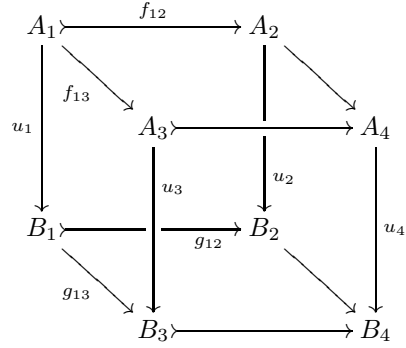
\begin{center}
$ \xymatrix {
        A_1 \ar@{>->}[rr]^{f_{12}} \ar[dd]_{u_1} \ar[rd]_{f_{13}} &  & A_2 \ar'[d][dd]^{u_2} \ar[rd] & \\
        & A_3 \ar@{>->}[rr] \ar[dd]^(.3){u_3} & & A_4 \ar[dd]^{u_4} \\
        B_1 \ar@{>->}'[r][rr]_{g_{12}} \ar[rd]_{g_{13}} & & B_2 \ar[rd] & \\
        & B_3 \ar@{>->}[rr] & & B_4
	}$
\end{center}
suppose that $A_1$, $A_3$, $B_1$, $B_3$ are cofibrant, that $f_{12}, g_{12}$ are cofibrations, and that the top and bottom faces are pushouts. 
\begin{enumerate}
\item
If $u_1, u_3$ are cofibrations and the natural map \\ $B_1 \Sum_{A_1} A_2 \lra B_2$ is a cofibration, then $u_2, u_4$ and the natural map $B_3 \Sum_{A_3} A_4 \lra B_4$ are cofibrations.
\item
If $u_1, u_2, u_3$ are weak equivalences, then $u_4$ is a weak equivalence.
\end{enumerate}
\item
Let ${\cat M}$ be a prefibration category. In the diagram
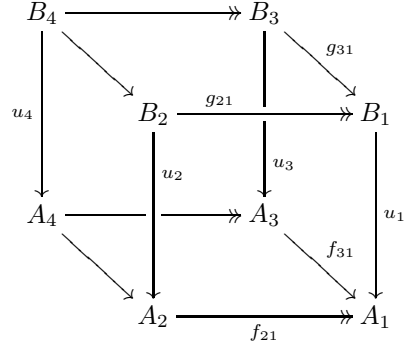
\begin{center}
$ \xymatrix {
        B_4 \ar@{->>}[rr] \ar[dd]_{u_4} \ar[rd] &  & B_3 \ar'[d][dd]^{u_3} \ar[rd]^{g_{31}} & \\
        & B_2 \ar@{->>}[rr]^(.3){g_{21}} \ar[dd]^(.3){u_2} & & B_1 \ar[dd]^{u_1} \\
        A_4 \ar@{->>}'[r][rr] \ar[rd] & & A_3 \ar[rd]^{f_{31}} & \\
        & A_2 \ar@{->>}_{f_{21}}[rr] & & A_1
	}$
\end{center}
suppose that $A_1$, $A_3$, $B_1$, $B_3$ are fibrant, that $f_{21}, g_{21}$ are fibrations, and that the top and bottom faces are pullbacks. 
\begin{enumerate}
\item
If $u_1, u_3$ are fibrations and the natural map \\ $B_2 \lra B_1 \times_{A_1} A_2 $ is a fibration, then $u_2, u_4$ and the natural map $B_4 \lra B_3 \times_{A_3} A_4$ are fibrations.
\item
If $u_1, u_2, u_3$ are weak equivalences, then $u_4$ is a weak equivalence.
\end{enumerate}
\end{enumerate}
\end{lem}

\begin{proof}
The statements are dual, and we will prove only (1). Let $B^{'}_2 = B_1 \Sum_{A_1} A_2$ and $B^{'}_4 = B_3 \Sum_{A_3} A_4$ be the pushout of the front and back faces of the diagram of (1). These pushouts exist because of the pushout axiom CF3.
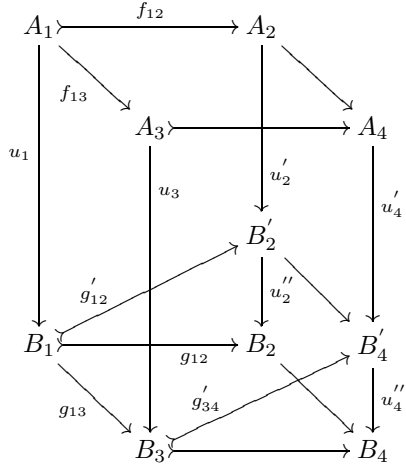
\begin{center}
$ \xymatrix {
        A_1 \ar@{>->}[rr]^{f_{12}} \ar[ddd]_(.4){u_1} \ar[rd]_{f_{13}} & & A_2 \ar[dd]^(.7){u^{'}_2} \ar[rd] & \\
        & A_3 \ar@{>->}[rr] \ar[ddd]^(.2){u_3} & & A_4 \ar[dd]^(.35){u^{'}_4} \\
        & & B^{'}_2 \ar[d]^{u^{''}_2} \ar[rd] & \\
        B_1 \ar@{>->}[rr]_(.7){g_{12}} \ar[rd]_{g_{13}} \ar@{>->}[rru]^(.3){g^{'}_{12}} & & B_2 \ar[rd] & B^{'}_4 \ar[d]^{u^{''}_4} \\
        & B_3 \ar@{>->}[rr] \ar@{>->}[rru]^(.3){g^{'}_{34}} & & B_4
	}$
\end{center}
The maps $g^{'}_{12}$ and $g^{'}_{34}$ are pushouts of cofibrations, therefore cofibrations. Furthermore, we observe that $B^{'}_4 = B_3 \Sum_{B_1} B^{'}_2$, and therefore $u^{''}_4$ is the pushout of $u^{''}_2$ along $B^{'}_2 \ra B^{'}_4$. 

Let's prove (a). If $u_1, u_3$ and $u^{''}_2$ are cofibrations, then by the pushout axiom $u^{'}_2$ and $u^{'}_4$ are cofibrations. $u^{''}_4$ is a pushout of $u^{''}_2$, therefore also a cofibration. It follows that $u_2, u_4$ are cofibrations.

Let's now prove (b), first under the assumption that
\begin{equation}
\label{eqn:gluinglemmasimplifying}
u_1, u_3 \;\; \mathrm{ and } \;\; u^{''}_2 \;\; \mathrm{ are} \;\; \mathrm{ cofibrations }
\end{equation}
If they are, since $u_1$, $u_2$, $u_3$ are weak equivalences we see that $u_1$, $u_3$ and their pushouts $u^{'}_2$, $u^{'}_4$ must be trivial cofibrations. From the two out of three axiom, $u^{''}_2$ is a weak equivalence, therefore a trivial cofibration, and so its pushout $u^{''}_4$ also is a trivial cofibration, which shows that $u_4 = u^{''}_4u^{'}_4$ is a weak equivalence.

For general weak equivalences $u_1$, $u_2$, $u_3$, we use the relative Brown factorization lemma to construct the diagram
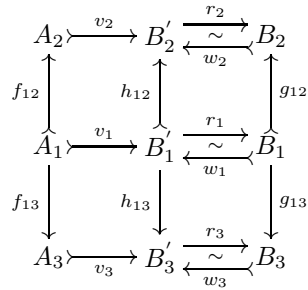
\begin{center}
$ \xymatrix {
    A_2 \ar@{>->}[r]^{v_2} & B^{'}_2 \ar@<1ex>[r]^{r_2}_\sim & B_2 \ar@{>->}@<1ex>[l]^{w_2} \\
    A_1 \ar@{>->}[u]^{f_{12}} \ar[d]_{f_{13}} \ar@{>->}[r]^{v_1} & 
    B^{'}_1 \ar@{>->}[u]^{h_{12}} \ar[d]_{h_{13}} \ar@<1ex>[r]^{r_1}_\sim & 
    B_1 \ar@{>->}[u]_{g_{12}} \ar[d]^{g_{13}} \ar@{>->}@<1ex>[l]^{w_1} \\
    A_3 \ar@{>->}[r]_{v_3} & B^{'}_3 \ar@<1ex>[r]^{r_3}_\sim & B_3 \ar@{>->}@<1ex>[l]^{w_3} \\
	}$
\end{center}
In this diagram:
\begin{enumerate}
\item
$v_1$, $w_1$, $r_1$ are constructed as a Brown factorization of $u_1$ as in \remarkref{rem:brownfact}
\item
$v_i$, $w_i$, $r_i$ for $i = 2, 3$ are constructed as relative Brown factorizations of $u_2$ resp. $u_3$ over the Brown factorization $v_1$, $w_1$, $r_1$.
\end{enumerate}

The maps ($w_1$, $w_2$, $w_3$) are trivial cofibrations and ($u_1$, $u_2$, $u_3$) are weak equivalences, so ($v_1$, $v_2$, $v_3$) are trivial cofibrations.

Statement (b) is true for 
\begin{enumerate}
\item[-]
($v_1$, $v_2$, $v_3$) resp. ($w_1$, $w_2$, $w_3$) because they satisfy property (\ref{eqn:gluinglemmasimplifying})
\item[-]
therefore true for ($r_1$, $r_2$, $r_3$) as a left inverse to ($w_1, w_2, w_3$)
\item[-]
therefore true for ($u_1$, $u_2$, $u_3$) as the composition of ($v_1$, $v_2$, $v_3$) and ($r_1$, $r_2$, $r_3$).
\end{enumerate}
\end{proof}

As a corollary we have

\index{excision}
\index{lemma!excision lemma}
\begin{lem}[Excision]
\label{lem:excision}
\mbox{}
\begin{enumerate}
\item
Let ${\cat M}$ be a precofibration category. In the diagram below
\begin{center}
$ \xymatrix {
        A \ar[r]^f_{\sim} \ar@{>->}[d]_i & C \\
	B & 
	}$
\end{center}
suppose that $A$, $C$ are cofibrant, $i$ is a cofibration and $f$ is a weak equivalence. Then the pushout of $f$ along $i$ is again a weak equivalence.
\item
Let ${\cat M}$ be a prefibration category. In the diagram below
\begin{center}
$ \xymatrix {
        & B \ar@{->>}[d]^p \\
        C \ar[r]^f_{\sim} & A
	}$
\end{center}
suppose that $A$, $C$ are fibrant, $p$ is a fibration and $f$ is a weak equivalence. Then the pullback of $f$ along $p$ is again a weak equivalence.
\end{enumerate}
\end{lem}

\begin{proof}
Part (1) is a particular case of the Gluing Lemma (1) for $f_{12} = g_{12} = i$, $u_1 = f_{13} = 1_A$, $g_{13} = u_3 = f$ and $u_2 = 1_B$. Part (2) is dual.
\end{proof}

It is now easy to see that in the presence of the rest of the axioms, the axiom CF3 (2) is equivalent to the Gluing Lemma (1) (b) and to excision. A dual statement holds for the fibration axioms.

\begin{lem}[Equivalent formulation of CF3]
\label{lem:equivaxiomcf3}
\mbox{}
\begin{enumerate}
\item
If $({\cat M}$, ${\cat W}$, ${\cat Cof})$ satisfies the axioms CF1-CF2, CF3 (1) and CF4, then the following are equivalent:
\begin{enumerate}
\item
It satisfies CF3 (2)
\item
It satisfies the Gluing Lemma (1) (b)
\item
It satisfies the Excision \lemmaref{lem:excision} (1).
\end{enumerate}
\item
If $({\cat M}$, ${\cat W}$, ${\cat Fib})$ satisfies the axioms F1-F2, F3 (1) and F4, then the following are equivalent:
\begin{enumerate}
\item
It satisfies F3 (2)
\item
It satisfies the Gluing Lemma (2) (b)
\item
It satisfies the Excision \lemmaref{lem:excision} (2).
\end{enumerate}
\end{enumerate}
\end{lem}

\begin{proof}
It suffices to prove (1). We have proved (a) $\Ra$ (b) as \lemmaref{lem:gluinglemma1}, and we have seen (b) $\Ra$ (c) in the proof of \lemmaref{lem:excision}.

For (c) $\Ra$ (a), suppose that $A, C$ are cofibrant, that $i \colon A \rightarrowtail B$ is a trivial cofibration and that $f \colon A \ra C$ is a map. Factor $f$ as a cofibration $f^{'}$ followed by a weak equivalence $r$, and using axiom CF3 (1) construct the pushouts $g^{'}, s$ of $f^{'}, r$
\begin{center}
$ \xymatrix {
    A \ar@{>->}[r]^{f^{'}} \ar@{>->}[d]_{i}^{\sim} & C^{'} \ar[r]^r_\sim \ar@{>->}[d]^{i^{'}}_\sim
 & C \ar@{>->}[d]^j \\
    B \ar@{>->}[r]^{g^{'}} & B^{'} \ar[r]^{s}_\sim & D
  }$
\end{center}
The maps $g^{'}, i^{'}$ and $j$ are cofibrations by axiom CF3 (1). Using excision, $i^{'}$ and $s$ are weak equivalences, so by the 2 out of 3 axiom CF2 the map $j$ is also a weak equivalence.
\end{proof}

\section{Cylinder and path objects}
\label{sec:cylobj}
We next define cylinder objects in a precofibration category, and show that cylinder objects exist. Dually, we define and prove existence of path objects in a prefibration category. 

\index{cylinder object}
\index{path object}
\begin{defn}[Cylinder and path objects]
\label{defn:cylobj}
\mbox{}
\begin{enumerate}
\item
Let ${\cat M}$ be a precofibration category, and $A$ a cofibrant object of ${\cat M}$. A cylinder object for $A$ consists of an object $IA$ and a factorization of the codiagonal $ \nabla \colon \xymatrix { A \Sum A \ar@{>->}[r]^-{i_0+i_1} & IA \ar[r]^{p}_{\sim} & A } $, with $i_0+i_1$ a cofibration and $p$ a weak equivalence.

\item
Let ${\cat M}$ be a prefibration category, and $A$ a fibrant object of ${\cat M}$. A path object for $A$ consists of an object $A^I$ and a factorization of the diagonal $ \Delta \colon \xymatrix { A \ar[r]^{i}_{\sim} & A^I \ar@{->>}[r]^-{(p_0,p_1)} & A \times A } $, with $(p_0,p_1)$ a fibration and $i$ a weak equivalence.
\end{enumerate}
\end{defn}

\begin{lem}[Existence of cylinder and path objects]
\label{lem:existencecyl}
\mbox{}
\begin{enumerate}
\item
Let ${\cat M}$ be a precofibration category, and $A$ a cofibrant object of ${\cat M}$. Then $A$ admits a (non-functorial) cylinder object.

\item
Let ${\cat M}$ be a prefibration category, and $A$ a fibrant object of ${\cat M}$. Then $A$ admits a (non-functorial) path object.
\end{enumerate}
\end{lem}

\begin{proof}
To prove (1), observe that if $A$ is cofibrant then the sum $A \Sum A$ exists and is cofibrant by \lemmaref{lem:cofsums}, and we can then use the factorization axiom CF4 to construct a cylinder object $ \xymatrix { A \Sum A \ar@{>->}[r]^-{i_0+i_1} & IA \ar[r]^{p}_{\sim} & A } $. The statement (2) follows from duality.
\end{proof}

Observe that for cylinder objects, the inclusion maps $i_0,i_1 \colon A \lra IA$ are trivial cofibrations. For path objects, the projection maps $p_0,p_1 \colon A^I \lra A$ are trivial fibrations.

\begin{lem}[Relative cylinder and path objects]
\label{lem:existencerelcyl}
\mbox{}
\begin{enumerate}
\item
Let ${\cat M}$ be a precofibration category, and $f \colon A \lra B$ a map with $A, B$ cofibrant objects. Let $IA$ be a cylinder of $A$. Then there exists a cylinder $IB$ and a commutative diagram
\begin{center}
$\xymatrix {
    A \Sum A \ar[d]_{f \Sum f} \ar@{>->}[r] &
    IA \ar[d]_{If} \ar[r]^{\sim} &
    A \ar[d]^f \\
    B \Sum B \ar@{>->}[r] &
    IB \ar[r]^{\sim} &
    B
  }$
\end{center}
with $ \xymatrix { (B \Sum B) \Sum_{A \Sum A} IA \ar@{>->}[r] & IB } $ a cofibration.
\item
Let ${\cat M}$ be a prefibration category,  and $f \colon A \lra B$ a map with $A, B$ fibrant objects. Let $B^I$ be a path object for of $B$. Then there exists a path object $A^I$ and a commutative diagram
\begin{center}
$ \xymatrix {
	A \ar[d]_f \ar[r]^{\sim} & 
	A^I \ar[d]_{f^I} \ar@{->>}[r] &
        A \times A \ar[d]^{f \times f} \\
	B \ar[r]^{\sim} &
        B^I \ar@{->>}[r] &
        B \times B
	}$
\end{center}
with $ \xymatrix { A^I  \ar@{->>}[r] & B^I \times_{B \times B} (A \times A) } $ a fibration.
\end{enumerate}
\end{lem}

\begin{proof}
To prove (1), apply \lemmaref{lem:relativefactorization} to the diagram 
\begin{center}
$\xymatrix {
    A \Sum A \ar[d]_{f \Sum f} \ar@{>->}[r] &
    IA  \ar[r]^{\sim} &
    A \ar[d]^f \\
    B \Sum B \ar[rr]^{\nabla} && B
  }$
\end{center}
Statement (2) is dual to (1).
\end{proof}

\section{Elementary consequences of CF5 and CF6}
\label{sec:transfincomp}
In the previous sections we have proved a number of elementary lemmas that are consequences of the precofibration category axioms CF1-CF4. In this section, we will do the same bringing in one by one the cofibration category axioms CF5 and CF6.

A word on the motivation behind the two additional axioms CF5-CF6. In the construction of homotopy colimits indexed by small diagrams in a cofibration category, it turns out that the role of small {\it direct categories} (\definitionref{def:directinvcat}) is essential, because an arbitrary small diagram can be approximated by a diagram indexed by a small direct category (\sectionref{sec:colimarbitrarycat}). 

For a small direct category ${\cat D}$, its degreewise filtration can be used to show that colimits indexed by ${\cat D}$ may be constructed using small sums, pushouts and countable direct transfinite compositions (at least if the base category is cocomplete). To put things in perspective, the axiom CF3 is a property of pushouts, the axiom CF5 is a property of small sums of maps and the axiom CF6 is a property of countable direct transfinite compositions of maps.

Let us clarify for a moment what we mean in CF6 and F6 by transfinite direct and inverse compositions of maps.
\begin{defn}
\label{defn:transcomp}
Let ${\cat M}$ be a category, and let $k$ be an ordinal. 
\begin{enumerate}
\item
A {\it direct $k$-sequence} of maps (or a {\it direct sequence of length $k$})
\begin{center}
$ \xymatrix {
        A_0 \ar[r]^{a_{01}} &
        A_1 \ar[r]^{a_{12}} &
	... \ar[r] &
        A_i \ar[r] &
        ... &
        (i < k)
	}$
\end{center}
consists of a collection of objects $A_i$ for $i < k$ and maps $a_{i_1i_2} \colon A_{i_1} \ra A_{i_2}$ for $i_1 < i_2 < k$, such that $a_{i_2i_3}a_{i_1i_2} = a_{i_1i_3}$ for all $i_1 < i_2 < i_3 < k$. The map $A_0 \ra \colim^{i < k} A_i$, if the colimit exists, is called the {\it transfinite composition} of the direct $k$-sequence. 
\item
An {\it inverse $k$-sequence} of maps (or an {\it inverse sequence of maps of length $k$})
\begin{center}
$ \xymatrix {
    ... \ar[r] &
    A_i \ar[r] &
    ... \ar[r]^{a_{21}} &
    A_1 \ar[r]^{a_{10}} &
    A_0 &
    (i < l)
  }$
\end{center}
consists of a collection of objects $A_i$  for $i < k$ and maps $a_{i_2i_1} \colon A_{i_2} \ra A_{i_1}$ for $i_1 < i_2 < k$, such that $a_{i_2i_1}a_{i_3i_2} = a_{i_3i_1}$ for all $i_1 < i_2 < i_3 < k$. The map $\lim^{i < k} A_i \ra A_0$, if the limit exists, is called the transfinite composition of the inverse $k$-sequence.
\end{enumerate}
\end{defn}

A direct $k$-sequence of maps is nothing but a diagram in ${\cat M}$ indexed by ordinals $< k$. A map of direct $k$-sequences is a map of such diagrams. 
 
If ${\cat M}$ is a precofibration category, a direct $k$-sequence of (trivial) cofibrations is a direct $k$-sequence in which all maps $a_{ij}$, $i < j < k$ are (trivial) cofibrations.  If ${\cat M}$ is a fibration category, an inverse $k$-sequence of (trivial) fibrations is an inverse $k$-sequence in which all maps $a_{ji}$, $i < j < k$ are (trivial) fibrations.

A sufficient condition for the additional axioms CF5-CF6 to be satisfied is that cofibrations and trivial cofibrations with cofibrant domain are stable under all {\it small} transfinite direct compositions.

\begin{lem}
\label{lem:transfinitecompositionstability}
\mbox{}
\begin{enumerate}
\item
If $({\cat M}, {\cat W}, {\cat Cof})$ satisfies axioms CF1-CF4, and if cofibrations (resp. trivial cofibrations) with cofibrant domain are stable under transfinite compositions of direct $k$-sequences for any small ordinal $k$, then $({\cat M}, {\cat W}, {\cat Cof})$ also satisfies axioms CF5-CF6.
\item
If $({\cat M}, {\cat W}, {\cat Fib})$ satisfies axioms F1-F4, and if fibrations (resp. trivial fibrations) with fibrant codomain are stable under transfinite compositions of inverse $k$-sequences for any small ordinal $k$, then $({\cat M}, {\cat W}, {\cat Fib})$ also satisfies axioms F5-F6.
\end{enumerate}
\end{lem}

\begin{proof}
We only prove (1). Axiom CF6 is clearly verified for ${\cat M}$, since it states that cofibrations (resp. trivial cofibrations) with cofibrant domain are stable under {\it countable} transfinite direct compositions. For the axiom CF5, given a set of (trivial) cofibrations $f_i \colon A_i \ra B_i$ for $i \eps I$ with $A_i$ cofibrant, we choose a well ordering of $I$. Denote $I^+$ the well ordered set $I$ with a maximal element adjoined. $I^+$ can be viewed as the succesor ordinal of $I$, and all elements $i \eps I$ can be viewed as the ordinals smaller than $I$.

We show that for any $i \eps I^+$, we have that $\Sum_{_{k < i}} f_k \colon \Sum_{_{k < i}} A_k \ra \Sum_{_{k < i}} B_k$ is well defined and is a (trivial) cofibration with a cofibrant domain. We use transfinite induction, and the initial step is trivial. Suppose the statement is true for all elements $<i$, and let's prove it for $i$.

If $i$ is a succesor ordinal, the statement for $i$ follows from \lemmaref{lem:cofsums}. Suppose that $i$ is a limit ordinal. 

For any $i^{''} < i^{'} < i$, the inclusion
\[
\Sum_{_{k < i^{''}}} A_k \ra \Sum_{_{k < i^{'}}} A_k
\]
is a cofibration, using the inductive hypothesis. The transfinite composition of these cofibrations defines $\Sum_{_{k < i}} A_k$, which therefore exists and is cofibrant. Similarly, $\Sum_{_{k < i}} B_k$ exists and is cofibrant.

For any $i^{''} < i^{'} < i$, the map
\[
\Sum_{_{k < i^{''}}} B_k \Sum_{_{i^{''} \le k < i}} A_k \ra \Sum_{_{k < i^{'}}} B_k \Sum_{_{i^{'} \le k < i}} A_k
\]
given by
\[
(\Sum_{_{k < i^{''}}} 1_{B_k}) \Sum (\Sum_{_{i^{''} \le k < i^{'}}} f_k) \Sum (\Sum_{_{i^{'} \le k < i}} 1_{A_k}) 
\]
is a well defined (trivial) cofibration with cofibrant domain, using the inductive hypothesis. The transfinite composition of these (trivial) cofibrations with cofibrant domain defines $\Sum_{_{k < i}} f_k$, which is therefore a (trivial) cofibration. The statement of our lemma now follows if we take $i$ to be the maximal element of $I^+$.
\end{proof}

The next two lemmas describe properties of the additional axiom CF5 (resp. F5).

\begin{lem}
\label{lem:smallcofsums}
\mbox{}
\begin{enumerate}
\item
Suppose that ${\cat M}$ is a CF1-CF4 (resp. CF1-CF5) cofibration category. If $f_i \colon A_i \ra B_i$ for $i \eps I$ is a finite (resp. small) set of weak equivalences between cofibrant objects, then $\Sum f_i$ is a weak equivalence.
\item
Suppose that ${\cat M}$ is an F1-F4 (resp. F1-F5) fibration category. If $f_i \colon A_i \ra B_i$ for $i \eps I$ is a finite (resp. small) set of weak equivalences between fibrant objects, then $\times f_i$ is a weak equivalence.
\end{enumerate}
\end{lem}

\begin{proof}
We only prove (1). Using the Brown Factorization Lemma, write $f_i = r_if^{'}_i $, where $f^{'}_i$ is a trivial cofibration and $r_i$ is a left inverse to a trivial cofibration $s_i$. Under both alternative hypotheses, the maps $\Sum f^{'}_i$ and $\Sum s_i$ are trivial cofibrations, so $\Sum r_i$ and therefore $\Sum f_i$ are weak equivalences.
\end{proof}

\begin{lem}[Equivalent formulation of CF5]
\label{lem:equivformaxiomcf5}
\mbox{}
\begin{enumerate}
\item
Suppose that $({\cat M}$, ${\cat W}$, ${\cat Cof})$ satisfies axioms CF1-CF4 and CF5 (1). Then the following are equivalent:
\begin{enumerate}
\item
It satisfies axiom CF5 (2)
\item
The class of weak equivalences between cofibrant objects is stable under small sums.
\end{enumerate}
\item
Suppose that $({\cat M}$, ${\cat W}$, ${\cat Fib})$ satisfies axioms F1-F4 and F5 (1). Then the following are equivalent:
\begin{enumerate}
\item
It satisfies axiom F5 (2)
\item
The class of weak equivalences between fibrant objects is stable under small products.
\end{enumerate}
\end{enumerate}
\end{lem}

\begin{proof}
The implication $(a) \Ra (b)$ is a consequence of \lemmaref{lem:smallcofsums}, and $(b) \Ra (a)$ is trivial.
\end{proof}

\comment{
We have seen from the CF1-CF4 axioms that (trivial) cofibrations with cofibrant domain are stable under {\it finite} sums (this is \lemmaref{lem:cofsums}). We can improve on that using the axioms CF1-CF5 as follows.

\begin{lem}
\label{lem:partialCF6}
\mbox{}
\begin{enumerate}
\item
Let ${\cat M}$ be a CF1-CF5 cofibration category. Then (trivial) cofibrations with cofibrant domain are stable under {\it countable} sums.
\item
Let ${\cat M}$ be an F1-F5 fibration category. Then (trivial) fibrations with cofibrant domain are stable under {\it countable} products.
\end{enumerate}
\end{lem}

\begin{proof}
The proof is very similar to that of \lemmaref{lem:transfinitecompositionstability}, and is omitted.
\end{proof}
} 

Here is a consequence of the axiom CF6 (resp. F6).

\begin{lem}
\label{lem:transcompequiv}
\mbox{}
\begin{enumerate}
\item
Let ${\cat M}$ be a CF1-CF4 cofibration category satisfying CF6. For any map of countable direct sequences of cofibrations with $A_0$, $B_0$ cofibrant and all $f_n$ weak equivalences
\begin{center}
$ \xymatrix {
        A_0 \ar@{>->}[r]^{a_{0}} \ar[d]_{f_0}^{\sim} &
        A_1 \ar@{>->}[r]^{a_{1}} \ar[d]_{f_1}^{\sim} &
        A_2 \ar@{>->}[r]^{a_{2}} \ar[d]_{f_2}^{\sim} &
	... \\
        B_0 \ar@{>->}[r]^{b_{0}} &
        B_1 \ar@{>->}[r]^{b_{1}} &
        B_2 \ar@{>->}[r]^{b_{2}} &
	...
	}$
\end{center}
the colimit map $\colim f_n \colon \colim A_n \ra \colim B_n$ is a weak equivalence between cofibrant objects.
\item
Let ${\cat M}$ be an F1-F4 fibration category satisfying CF6. For any map of countable inverse sequences of fibrations with $A_0$, $B_0$ fibrant and all $f_n$ weak equivalences
\begin{center}
$ \xymatrix {
    ... \ar@{->>}[r]^{a_{2}} & \
    A_2 \ar@{->>}[r]^{a_{1}} \ar[d]^{f_2}_{\sim} &
    A_1 \ar@{->>}[r]^{a_{0}} \ar[d]^{f_1}_{\sim} &
    A_0 \ar[d]^{f_0}_{\sim} \\
    ... \ar@{->>}[r]^{b_{2}} &
    B_2 \ar@{->>}[r]^{b_{1}} &
    B_1 \ar@{->>}[r]^{b_{0}} &
    B_0
  }$
\end{center}
the limit map $\lim f_n \colon \lim A_n \ra \lim B_n$ is a weak equivalence between fibrant objects.
\end{enumerate}
\end{lem}

\begin{proof}
To prove (1), observe that using Rem. \ref{rem:brownfact} and \lemmaref{lem:relbrownfact}, we can inductively construct Brown factorizations of $f_n$ that make the diagram
\begin{center}
$ \xymatrix {
        A_0 \ar@{>->}[r]^{a_{0}} \ar@{>->}[d]_{f_0^{'}}^{\sim} &
        A_1 \ar@{>->}[r]^{a_{1}} \ar@{>->}[d]_{f_1^{'}}^{\sim} &
        A_2 \ar@{>->}[r]^{a_{2}} \ar@{>->}[d]_{f_2^{'}}^{\sim} &
	... \\
        B_0^{'} \ar@{>->}[r]^{b^{'}_{0}} \ar@<1ex>[d]^{r_0} &
        B_1^{'} \ar@{>->}[r]^{b^{'}_{1}} \ar@<1ex>[d]^{r_1} &
        B_2^{'} \ar@{>->}[r]^{b^{'}_{2}} \ar@<1ex>[d]^{r_2} &
	... \\
        B_0 \ar@{>->}[r]^{b_{0}} \ar@{>->}@<1ex>[u]^{s_0}_{\sim} &
        B_1 \ar@{>->}[r]^{b_{1}} \ar@{>->}@<1ex>[u]^{s_1}_{\sim} &
        B_2 \ar@{>->}[r]^{b_{2}} \ar@{>->}@<1ex>[u]^{s_2}_{\sim} &
	...
	}$
\end{center}
commutative, such that additionally all maps
\begin{center}
$B^{'}_{n-1} \Sum_{A_{n-1}} A_n \ra B^{'}_{n}$ and $B^{'}_{n-1} \Sum_{B_{n-1}} B_n \ra B^{'}_n$
\end{center}
are trivial cofibrations. If we show that $\colim f^{'}_n$ and $\colim s_n$ are weak equivalences, it will follow that $\colim f_n$ is a weak equivalence. 

It suffices therefore to prove that $\colim f_n$ is a weak equivalence under the additional assumption that $f_0$ and all $B_{n-1} \Sum_{A_{n-1}} A_n \ra B_{n}$ are trivial cofibrations. The objects denoted $A = \colim A_n$ and $B = \colim B_n$ are cofibrant by CF6 (1). 

\begin{center}
$\xymatrix{
    A_0 \ar@{>->}[r]^{a_{0}} \ar@{>->}[d]_{f_0} & A_1 \ar@{>->}[r]^{a_{1}} \ar@{>->}[d] & ... \ar@{>->}[r] & A \ar@{>->}[d] \\
    B_0 \ar@{>->}[rd]_{b_{0}} \ar@{>->}[r] & B_0 \Sum_{A_0} A_1 \ar@{>->}[r] \ar@{>->}[d] & ... \ar@{>->}[r] & B_0 \Sum_{A_0} A \ar@{>->}[d] \\
    & B_1 \ar@{>->}[r] \ar@{>->}[rd]_{b_1} & ... \ar@{>->}[r] & B_1 \Sum_{A_1} A \ar@{>->}[d] \\
    & & ... \ar@{>->}[rd] & ... \ar@{>->}[d] \\
    &&& B
    }$
\end{center}
The map $A \ra B$ factors as the composition of the direct sequence of maps $A \ra B_0 \Sum_{A_0} A$ followed by $B_{n-1} \Sum_{A_{n-1}} A \ra B_{n} \Sum_{A_{n}} A$ for $n \ge 1$, and each map in the sequence is a trivial cofibration as the pushout of the trivial cofibrations $f_0$ resp. $B_{n-1} \Sum_{A_{n-1}} A_{n} \ra B_{n}$, so by CF6 (2) the map $A \ra B$ is a trivial cofibration.

The proof of statement (2) is dual to the proof of (1). 
\end{proof}

\begin{lem}[Equivalent formulation of CF6]
\label{lem:equivformaxiomcf6}
\mbox{}
\begin{enumerate}
\item
Suppose that $({\cat M}$, ${\cat W}$, ${\cat Cof})$ satisfies axioms CF1-CF4 and CF6 (1). Then the following are equivalent:
\begin{enumerate}
\item
It satisfies axiom CF6 (2)
\item
It satisfies the conclusion of \lemmaref{lem:transcompequiv} (1)
\item
It satisfies the conclusion of \lemmaref{lem:transcompequiv} (1) for any map $f_n$ of countable direct sequences of cofibrations, with each $f_n$ a trivial cofibration.
\end{enumerate}
\item
Suppose that $({\cat M}$, ${\cat W}$, ${\cat Fib})$ satisfies axioms F1-F4 and F6 (1). Then the following are equivalent:
\begin{enumerate}
\item
It satisfies axiom F6 (2)
\item
It satisfies the conclusion of \lemmaref{lem:transcompequiv} (2)
\item
It satisfies the conclusion of \lemmaref{lem:transcompequiv} (2) for any map $f_n$ of countable inverse sequences of fibrations, with each $f_n$ a trivial fibration.
\end{enumerate}
\end{enumerate}
\end{lem}

\begin{proof}
We only prove (1). The implication (a) $\Ra$ (b) is proved by \lemmaref{lem:transcompequiv}, and (b) $\Ra$ (c) is trivial. 

Let us prove (c) $\Ra$ (a). Suppose we have a countable direct sequence of trivial cofibrations with $A_0$ cofibrant
\begin{center}
$ \xymatrix {
        A_0 \ar@{>->}[r]^{a_{0}}_\sim &
        A_1 \ar@{>->}[r]^{a_{1}}_\sim &
        A_2 \ar@{>->}[r]^{a_{2}}_\sim &
	... \\
	}$
\end{center}
Then $\colim A_n$ exists and $A_0 \ra \colim A_n$ is a cofibration by axiom CF6 (1). If we view $A_0$ as a constant, countable direct sequence of identity maps, we get a map $f_n = a_{n-1} ... a_0$ of countable direct sequences of cofibrations, with each $f_n$ a trivial cofibration. From the conclusion of \lemmaref{lem:transcompequiv} (1) we see that $A_0 \ra \colim A_n$ is a weak equivalence.
\end{proof}

\section{Over and under categories}
\label{sec:catoverunder}
\index{category!over, under category}
If $u \colon {\cat A} \ra {\cat B}$ is a functor and $B$ is an object of ${\cat B}$, the over category $(u \downarrow B)$ by definition has:
\begin{enumerate}
\item
as objects, pairs $(A, g)$ of an object $A \eps {\cat A}$ and a map $g \colon uA \ra B$
\item
as maps $(A_1, g_1) \ra (A_2, g_2)$, the maps $f \colon A_1 \ra A_2$ such that $g_2 \circ uf = g_1$.
\end{enumerate}
The under category $(B \downarrow u)$ by definition has:
\begin{enumerate}
\item
as objects, pairs $(A, g)$ with $A \eps {\cat A}$ and $g \colon b \ra uA$
\item
as maps $(A_1, g_1) \ra (A_2, g_2)$, the maps $f \colon A_1 \ra A_2$ such that $uf \circ g_1 = g_2$.
\end{enumerate}

The two definitions are dual in the sense that $(B \downarrow u) \cong (u^{op} \downarrow B)^{op}$.

We have a canonical functor $i_{u, B} \colon (u \downarrow B) \ra {\cat A}$ that sends an object $(A, g: ua \ra B)$ to $A \eps {\cat A}$ and a map $(A_1, g_1) \ra (A_2, g_2)$ to the component map $A_1 \ra A_2$. Dually, we have a canonical functor $i_{B, u} \colon (B \downarrow u) \ra {\cat A}$ defined by $i_{B, u} = (i_{u^{op}, B})^{op}$.

If $A \eps {\cat A}$ is an object, for simplicity we denote $({\cat A} \downarrow A)$ for $(1_{\cat A} \downarrow A)$ and $(A \downarrow {\cat A})$ for $(A \downarrow 1_{\cat A})$. 

Suppose that ${\cat M}$ is a category and ${\cat F}$ is a class of maps of ${\cat M}$. For example, ${\cat M}$ could be a cofibration category, and ${\cat F}$ could be ${\cat W}$ or ${\cat Cof}$. For an object $A \eps {\cat M}$, we say that $i_{1_{\cat M}, A}^{-1} {\cat F}$ is the class of maps of $({\cat M} \downarrow A)$ {\it induced} by ${\cat F}$. Dually, $i_{A, 1_{\cat M}}^{-1} {\cat F}$ is the class of maps induced by ${\cat F}$ on $(A \downarrow {\cat M})$.

If $({\cat M}, {\cat W}, {\cat Cof})$ is a precofibration category, we therefore obtain an induced class of weak equivalences and cofibrations on $({\cat M} \downarrow A)$ and on $(A \downarrow {\cat M})$. Dually if $({\cat M}, {\cat W}, {\cat Fib})$ is a prefibration category we obtain an induced class of weak equivalences and fibrations on $({\cat M} \downarrow A)$ and on $(A \downarrow {\cat M})$.

The following result is a simple consequence of the definitions.

\begin{prop}
\label{prop:overcofcat}
\mbox{}
\begin{enumerate}
\item
Suppose that ${\cat M}$ is a (pre)cofibration category, and $A \eps {\cat M}$ is an object. Then:
\begin{enumerate}
\item
$({\cat M} \downarrow A)$ is a (pre)cofibration category
\item
If $A$ is cofibrant, $(A \downarrow {\cat M})$ is a (pre)cofibration category
\end{enumerate}
with respect to the induced weak equivalences and cofibrations.
\item
Suppose that ${\cat M}$ is a (pre)fibration category, and $A \eps {\cat M}$ is an object. Then:
\begin{enumerate}
\item
$(A \downarrow {\cat M})$ is a (pre)fibration category
\item
If $A$ is fibrant, $({\cat M} \downarrow A)$ is a (pre)fibration category 
\end{enumerate}
with respect to the induced weak equivalences and fibrations. $\square$
\end{enumerate}
\end{prop}

\section{Properness}
\label{sec:properness}
Sometimes, a precofibration category $({\cat M}, {\cat W}, {\cat Cof})$ admits more than one precofibration structure with weak equivalences ${\cat W}$. Under an additional condition (left properness) one can show that $({\cat M}, {\cat W})$ admits an intrinsic structure of a CF1-CF4 precofibration category, larger than ${\cat Cof}$, defined in terms of what we will call {\it left proper} maps. 

\begin{defn}
\label{defn:propercofcat}
\mbox{}
\begin{enumerate}
\item
A precofibration category $({\cat M}, {\cat W}, {\cat Cof})$ is {\it left proper} if it satisfies
\begin{description}
\item[PCF] 
Given a solid diagram in ${\cat M}$ with $i$ a cofibration and $A$ cofibrant,
\begin{center}
$ \xymatrix {
	A \ar[r]^r \ar@{>->}[d]_i &
	C \ar@{>-->}[d] \\
	B \ar@{-->}[r]^{r^{'}} &
	D
	}$
\end{center}
then the pushout exists in ${\cat M}$. Moreover, if $r$ is a weak equivalence then so is $r^{'}$.
\end{description}
\item
A prefibration category $({\cat M}, {\cat W}, {\cat Fib})$ is {\it right proper} if 
\begin{description}
\item[PF] 
Given a solid diagram in ${\cat M}$ with $p$ a fibration and $A$ fibrant,
\begin{center}
$ \xymatrix {
        D \ar@{-->}[r]^{r^{'}} \ar@{-->}[d] & B \ar@{->>}[d]^p \\
        C \ar[r]^r_{\sim} & A
	}$
\end{center}
then the pullback exists in ${\cat M}$. Moreover, if $r$ is a weak equivalence then so is $r^{'}$.
\end{description}
\item
An ABC model category is proper if its underlying precofibration and prefibration categories are left resp. right proper.
\end{enumerate}
\end{defn}

From the Excision Lemma, a precofibration category with all objects cofibrant is left proper. A prefibration category with all objects fibrant is right proper.

We will define proper maps in the context of what we call {\it category pairs}.

\index{category!with weak equivalences}
\begin{defn}
\label{defn:catweakequiv}
A category pair $({\cat M}$, ${\cat W})$ consists of a category ${\cat M}$ with a class of weak equivalence maps ${\cat W}$, where ${\cat W}$  is stable under composition and includes the identity maps of ${\cat M}$. We may view ${\cat W}$ as defining a subcategory with the same objects as ${\cat M}$. 
\end{defn}

\begin{defn}
\label{defn:propercofibration}
Suppose that $({\cat M}, {\cat W})$ is a category pair.
\begin{enumerate}
\item
A map $f \colon A \ra B$ is called {\it left proper} if for any diagram of full maps with $r$ a weak equivalence
\begin{center}
$\xymatrix {
    A \ar[r] \ar[d]_f & C_1 \ar[r]^r_{\sim} \ar@{-->}[d] & C_2 \ar@{-->}[d] \\
    B \ar@{-->}[r] & D_1 \ar@{-->}[r]^{r^{'}} & D_2
  }$
\end{center}
the pushouts exist, and the map $r^{'}$ is again a weak equivalence. 
\item
A map $f \colon B \ra A$ is called {\it right proper} if for any diagram of full maps with $r$ a weak equivalence
\begin{center}
$\xymatrix {
    D_2 \ar@{-->}[r]^{r^{'}} \ar@{-->}[d] & D_1 \ar@{-->}[r] \ar@{-->}[d] & B \ar[d]^f \\
    C_2 \ar[r]^r_{\sim} & C_1 \ar[r] & A
  }$
\end{center}
the pullbacks exist, and the map $r^{'}$ is again a weak equivalence. 
\end{enumerate}
\end{defn}

We say that an object $A$ is {\it left proper} if the map $\initial \ra A$ is left proper. An object $A$ is {\it right proper} if the map $A \ra \terminal$ is right proper.

The class of left proper maps of $({\cat M}, {\cat W})$ will be denoted ${\cat PrCof}$. The class of right proper maps will be denoted ${\cat PrFib}$. Observe that the left proper maps are stable under composition and under pushout. The right proper maps are stable under composition and under pullback. All isomorphisms are left and right proper.

The left proper weak equivalences will be called {\it trivial left proper} maps, and the right proper weak equivalences will be called {\it trivial right proper} maps.

\begin{thm}
\label{thm:leftproperstruct}
\mbox{}
\begin{enumerate}
\item
If $({\cat M}, {\cat W}, {\cat Cof})$ is a left proper precofibration category, then 
\begin{enumerate}
\item
Any cofibration $A \ra B$ with $A$ cofibrant is left proper.
\item
Any map of ${\cat M}$ factors as a left proper map followed by a weak equivalence.
\item
Trivial left proper maps are stable under pushout.
\item
$({\cat M}, {\cat W}, {\cat PrCof})$ is a precofibration category.
\end{enumerate}
\item
If $({\cat M}, {\cat W}, {\cat Fib})$ is a right proper prefibration category, then 
\begin{enumerate}
\item
Any fibration $A \ra B$ with $B$ fibrant is right proper.
\item
Any map of ${\cat M}$ factors as a weak equivalence followed by a right proper map.
\item
Trivial right proper maps are stable under pushout.
\item
$({\cat M}, {\cat W}, {\cat PrFib})$ is a prefibration category.
\end{enumerate}
\end{enumerate}
\end{thm}

\begin{proof}
We only prove (1). For (1) (a), suppose we have a diagram
\begin{center}
$\xymatrix {
    A \ar[r]^f \ar@{>->}[d]_i & C_1 \ar[r]^r_{\sim} \ar[d] & C_2 \ar[d] \\
    B \ar[r] & D_1 \ar[r]^{r^{'}} & D_2
  }$
\end{center}
with $A$ cofibrant, $i$ a cofibration and $r$ a weak equivalence, such that both squares are pushouts. We'd like to show that $r^{'}$ is a weak equivalence. 

Denote $f = f^{'}s$ a factorization of $f$ as a cofibration $f^{'}$ followed by a weak equivalence $s$. 
\begin{center}
$\xymatrix {
    A \ar@{>->}[r]^{f_1} \ar@{>->}[d]_i & C^{'}_1 \ar[r]^s_\sim \ar@{>->}[d]^{i^{'}} & C_1 \ar[r]^r_{\sim} \ar[d] & C_2 \ar[d] \\
    B \ar@{>->}[r] & D^{'}_1 \ar[r]^{s^{'}} & D_1 \ar[r]^{r^{'}} & D_2
  }$
\end{center}
Denote $i^{'}$ the pushout of $i$, and $s^{'}$ the pushout of $s$. Since $i^{'}$ is a cofibration and $s$, $rs$ are weak equivalences, from the left properness of ${\cat M}$ we see that $s^{'}$, $r^{'}s^{'}$ and therefore $r^{'}$ are weak equivalences.

For (1) (b), suppose $f \colon A \ra B$ is a map in ${\cat M}$. We'd like to construct a factorization $f = rf^{'}$ with $f^{'}$ left proper and $r$ a weak equivalence.

Let $a \colon A^{'} \ra A$ be a cofibrant replacement of $A$, and $fa = r^{'}_1f_1$ factorization of $fa$ as a cofibration $f_1$ followed by a weak equivalence $r_1$.
\begin{center}
$\xymatrix{ 
    &&& B \\
    A \ar[rr]_{f^{'}} \ar[urrr]^f && Y^{'} \ar[ur]^(.3)r & \\
    A^{'} \ar[u]_\sim^{a^{'}} \ar@{>->}[rr]_{f_1} && A_1 \ar[u]^{b^{'}}_\sim \ar[uur]_{r_1^{'}}^\sim &
  }$
\end{center}
Define $f^{'}$ as the pushout of $f_1$. The map $b^{'}$ is a weak equivalence as the pushout of $a^{'}$, since ${\cat M}$ is left proper. The map $r$ is a weak equivalence by the two out of three axiom, and the map $f^{'}$ is left proper as the pushout of $f_1$ which is a cofibration with cofibrant domain (therefore left proper).

For (1) (c), let $i \colon A \ra B$ be a trivial left proper map and let $f \colon A \ra C$ be a map. In the diagram
\begin{center}
$\xymatrix {
    A \ar[r]^{f^{'}} \ar[d]_i^\sim & C^{'} \ar[r]^r_{\sim} \ar[d]^{i^{'}} & C \ar[d]^{j} \\
    B \ar[r] & D^{'} \ar[r]^{r^{'}} & D
  }$
\end{center}
using part (b) we have factored $f$ as $rf^{'}$, with $f^{'}$ left proper and $r$ a weak equivalence. We define $i^{'}$ and $j$ as the pushouts of $i$. We therefore have that $i^{'}$ and $j$ are left proper. The map $i^{'}$ is a weak equivalence since $f^{'}$ is left proper. The map $r^{'}$ is a weak equivalence since $i$ is proper. From the two out of three axiom, the map $j$ is also a weak equivalence, therefore a trivial left proper map.

For (1) (d), the axioms CF1, CF2 and CF3 (1) are trivially verified. Part (c) proves axiom CF3 (2), and part (b) proves axiom CF4.
\end{proof}

It does not appear to be the case that $({\cat M}, {\cat W}, {\cat PrCof})$ necessarily satisfies CF5 or CF6 if $({\cat M}, {\cat W}, {\cat Cof})$ does.

In the rest of the section, we will review briefly a number of elementary properties of proper maps.

\begin{prop}
\label{prop:propermaps1}
\mbox{}
\begin{enumerate}
\item
Suppose that $({\cat M}, {\cat W}, {\cat Cof})$ is a left proper precofibration category. Then a weak equivalence is a trivial left proper map iff all its pushouts exist and are weak equivalences.
\item
Suppose that $({\cat M}, {\cat W}, {\cat Fib})$ is a right proper prefibration category. Then a weak equivalence is a trivial right proper map iff all its pullbacks exist and are weak equivalences.
\end{enumerate}
\end{prop}

\begin{proof}
We only prove (1). Implication $\Ra$ follows from \theoremref{thm:leftproperstruct} (1) (c). 

For $\La$, suppose that $i \colon A \ra B$ is a weak equivalence whose pushouts remain weak equivalences. We'd like to show that $i$ is left proper. For any map $f$ and weak equivalence $r$, we construct the diagram with pushout squares
\begin{center}
$\xymatrix {
    A \ar[r]^f \ar[d]_i^\sim & C_1 \ar[r]^r_{\sim} \ar[d]^{i_1}_\sim & C_2 \ar[d]^{i_2}_\sim \\
    B \ar[r] & D_1 \ar[r]^{r^{'}} & D_2
  }$
\end{center}
The maps $i_1$, $i_2$ and are weak equivalences as pushouts of $i$. From the 2 out of 3 axiom, the map $r^{'}$ is a weak equivalence, which shows that $i$ is left proper.
\end{proof}

\begin{prop}
\label{prop:propermaps2}
\mbox{}
\begin{enumerate}
\item
Suppose that $({\cat M}, {\cat W}, {\cat Cof})$ is a cocomplete, left proper precofibration category. If $f, g$ are two composable maps such that $gf$ is left proper and $g$ is trivial left proper. Then $f$ is left proper.
\item
Suppose that $({\cat M}, {\cat W}, {\cat Fib})$ is a complete, right proper prefibration category. If $f, g$ are two composable maps such that $gf$ is left proper and $g$ is trivial left proper. Then $f$ is left proper.
\end{enumerate}
\end{prop}

\begin{proof}
We only prove (1). In the diagram
\begin{center}
$\xymatrix {
    A_1 \ar[r] \ar[d]_f & B_1 \ar[r]^r_{\sim} \ar[d] & C_1 \ar[d] \\
    A_2 \ar[r] \ar[d]_g^\sim & B_2 \ar[r]^{r^{'}} \ar[d]^{g^{'}} & C_2 \ar[d]^{g^{''}} \\
    A_3 \ar[r] & B_3 \ar[r]^{r^{''}} & C_3
  }$
\end{center}
$r$ is a weak equivalence and all squares are pushouts. To prove that $f$ is left proper, we need to show that $r^{'}$ is a weak equivalence. But $r^{''}$ is a weak equivalence since $gf$ is proper. $g^{'}$ and $g^{''}$ are trivial left proper as pushouts of $g$. It follows from the two out of three axiom that $r^{'}$ is a weak equivalence. 
\end{proof}

Recall the definition of the retract of a map
\begin{defn}
\label{defn:mapretract}
\index{retract}
A map $f \colon A \ra B$ in a category ${\cat M}$ is a retract of $g \colon C \ra D$ if there exists a commutative diagram
\begin{center}
$\xymatrix{
	A \ar[r] \ar[d]_f \ar@/^1pc/[rr]^{1_A} &
	C \ar[d]^g \ar[r] & A \ar[d]^f \\
	B \ar[r] \ar@/_1pc/[rr]_{1_B} &
	D \ar[r] & B
	}$
\end{center}
\end{defn}

Note that the saturation $\overline{\cat W}$ of the class of weak equivalences ${\cat W}$ is closed under retracts.

\begin{prop}
\label{prop:propermaps3}
Suppose that $({\cat M}, {\cat W})$ is a category pair and that ${\cat W}$ is closed under retracts.
\begin{enumerate}
\item
Assume that ${\cat M}$ is cocomplete. Then the class of left proper maps and that of trivial left proper maps are both closed under retracts.
\item
Assume that ${\cat M}$ is complete. Then the class of right proper maps and that of trivial right proper maps are both closed under retracts.
\end{enumerate}
\end{prop}

\begin{proof}
Follows directly from the definitions.
\end{proof}

\chapter{Relation with other axiomatic systems}
\label{chap:modcatcatcofobj}

We would like to describe in this chapter how ABC cofibration categories relate to Brown's categories of cofibrant objects, to Quillen model categories and to other axiomatizations that have been proposed for categories with cofibrations. 

Aside from the goal of bringing together and comparing various axiomatizations that have been proposed for (co)fibrations and weak equivalences, this allows us to tap into a large class of examples of ABC model categories. 

For example, simplicial sets $sSets$ form a Quillen model category \cite{Quillen1}, with inclusions as cofibrations, with maps satisfying the Kan extension property as fibrations and with maps whose geometric realization is a topological homotopy equivalence as weak equivalences. Any Quillen model category is an ABC model category, and therefore $sSets$ is an ABC model category. By \theoremref{thm:generalpointwisecofstructure}, so is any diagram category $sSets^{\cat D}$ for a small category ${\cat D}$, and by \definitionref{defn::restricteddiagrams} and \theoremref{thm:generalpointwisecofibstruct2} so are the ${\cat D}_2$-reduced ${\cat D}_1$-diagrams $sSets^{({\cat D}_1, {\cat D}_2)}$ for a small category pair $({\cat D}_1, {\cat D}_2)$.

Other basic examples of ABC model categories are explained in \chapterref{chap:examples}.

The list of alternative cofibration category axiomatizations discussed in this chapter is by no means exhaustive. 

We have omitted Alex Heller's notion of h-c categories \cite{HellerHC1}, \cite{HellerHC2}, \cite{HellerHC3}. As we have seen, the ABC cofibration category axioms are written in terms of {\it cofibrations} and {\it weak equivalences}.  In contrast, Heller's h-c category axioms are written in terms of cofibrations and a homotopy relation $\simeq$ on maps. A map $f \colon A \ra B$ in a h-c category by definition is a weak equivalence if there exists $g \colon B \ra A$ with $gf \simeq 1_A$, $fg \simeq 1_B$.

We have also omitted the categories with a natural cylinder \cite{Kamps}, \cite{Shitanda} and \cite{Kamps-Porter}, which are a surprising elaboration of ideas surrounding the Kan extension property for cubical sets \cite{Kan1}, \cite{Kan2}.


\section{Brown's categories of cofibrant objects}
In his paper \cite{Brown}, Brown defines categories of fibrant objects (and dually categories of cofibrant objects). We list below Brown's axioms, stated in the cofibration setting, slightly modified but equivalent to the actual axioms of \cite{Brown}.

\begin{defn}[Categories of cofibrant objects]
\label{defn:catcofobj}
\index{category!of (co)fibrant objects}
A category of cofibrant objects $({\cat M}$, ${\cat W}$, ${\cat Cof})$ consists of a category ${\cat M}$ and two distinguished classes of maps ${\cat W}, {\cat Cof}$ - the weak equivalences and respectively cofibrations of ${\cat M}$, subject to the axioms below:
\begin{description}
\item[CFObj1] 
All isomophisms of ${\cat M}$ are trivial cofibrations. ${\cat M}$ has an initial object $\initial$, and all objects of ${\cat M}$ are cofibrant. Cofibrations are stable under composition. 
\item[CFObj2] 
(Two out of three axiom) If $f, g$ are maps of ${\cat M}$ such that $gf$ is defined, and if two of $f, g, gf$ are weak equivalences, then so is the third.
\item[CFObj3] 
(Pushout axiom) Given a solid diagram in ${\cat M}$, with $i$ a cofibration,
\begin{center}
$ \xymatrix {
	A \ar[r] \ar@{>->}[d]_i &
	C \ar@{>->}[d]^j \\
	B \ar@{-->}[r] &
	D
	}$
\end{center}
then the pushout exists in ${\cat M}$ and $j$ is a cofibration. If additionally $i$ is a trivial cofibration, then $j$ is a trivial cofibration.
\item[CFObj4]
(Cylinder axiom) For any object $A$ of ${\cat M}$, the codiagonal $\nabla \colon A \Sum A \ra A$ admits a factorization as a cofibration followed by a weak equivalence.
\end{description}
\end{defn}

The axioms for categories of fibrant objects are dual to those of \definitionref{defn:catcofobj}, and are denoted FObj1-FObj4.

The precofibration categories (satisfying the minimal axioms CF1-CF4 but not the additional axioms CF5-CF6) are essentialy a modification of Brown's categories of cofibrant objects - in the sense that we allow objects to be non-cofibrant. The following lemma explains the precise relationship between precofibration categories and categories of cofibrant objects.

\begin{prop}
\label{prop:relcatcofobj}
\mbox{}
\begin{enumerate}
\item
Any category of cofibrant objects is a precofibration category. Conversely, if ${\cat M}$ is a precofibration category, then ${\cat M}_{cof}$ is a category of cofibrant objects. 
\item
Any category of fibrant objects is a prefibration category. Conversely, if ${\cat M}$ is a prefibration category, then ${\cat M}_{fib}$ is a category of fibrant objects. 
\end{enumerate}
\end{prop}

\begin{proof}
We only prove (1). Implication $\Leftarrow$ is an easy consequence of the axioms. For the other direction $\Rightarrow$, the only axiom that needs to be proved is the factorization axiom CF4.

For that, it suffices to prove Brown's factorization lemma \ref{lem:brownfact} in the context of the axioms CFObj1-CFObj4. We want to show that any map $f \colon A \ra B$ factors as $f = rf'$, where $f'$ is a cofibration and $rs = 1_B$ for a trivial cofibration $s$.

Choose a cylinder $IA$, and construct $s$ as the pushout of the trivial cofibration $i_0$.

\begin{center}
$ \xymatrix {
        A \ar[r]^f \ar@{>->}[d]_{i_0}^{\sim} &
	B \ar@{>->}[d]^s_{\sim} & \\
	IA \ar[r]^F &
	B'
	}$
\end{center}

Notice that $i_0$ has $p \colon IA \lra A$ as a left inverse, and it follows that $s$ has a left inverse $r$. Let $f'$ be $Fi_1$, which satisfies $f = rf'$, and to complete the proof it remains to prove that $f'$ is a cofibration. 

We notice that $B'$ is also the pushout of the diagram below
\begin{center}
$ \xymatrix {
        A \Sum A \ar[r]^{f \Sum 1} \ar@{>->}[d]_{i_0 + i_1} &
	B \Sum A \ar@{>->}[d]^{s + f'} & \\
	IA \ar[r]^F &
	B'
	}$
\end{center}
so $f'$ is $ \xymatrix { A \ar@{>->}[r] & B \Sum A \ar@{>->}[r]^-{s + f'} & B' } $, therefore a cofibration.
\end{proof}

\section{Quillen model categories}
\index{Quillen!model category}
Quillen's model categories involve both cofibrations and fibrations, and come with built-in Eckman-Hilton duality between cofibrations and fibrations. 

To start, recall the definition of the left (and right) lifting property of maps.

\begin{defn}
\label{defn:llp}
\index{left lifting property LLP}
\index{right lifting property RLP}
For a solid commutative diagram in a category ${\cat M}$
\begin{center}
$ \xymatrix {
	A \ar[r] \ar@{->}[d]_i &
	C \ar@{->}[d]^p \\
	B \ar[r] \ar@{-->}[ur] &
	D
	}$
\end{center}
if a dotted arrow exists making the diagram commutative we say that $i$ has the LLP (left lifting property) with respect to $p$, and that $p$ has the RLP (right lifting property) with respect to $i$.
\end{defn}

We will use the (closed) Quillen model category axiom formulation of \cite{Hirschhorn}, except that we do not require functorial factorization of maps in the axiom M5. 

\begin{defn}[Quillen model categories]
\label{defn:modelcat}
A Quillen model category (${\cat M}$, ${\cat W}$, ${\cat Cof}$, ${\cat Fib}$) consists of a category ${\cat M}$ and three distinguished classes of maps ${\cat W}$, ${\cat Cof}$, ${\cat Fib}$ - the weak equivalences, the cofibrations and respectively the fibrations of ${\cat M}$, subject to the axioms below:
\begin{description}
\item[M1] 
${\cat M}$ is complete and cocomplete. 
\item[M2] 
(Two out of three axiom) If $f, g$ are maps of ${\cat M}$ such that $gf$ is defined, and if two of $f, g, gf$ are weak equivalences, then so is the third.
\item[M3] 
(Retract axiom) Weak equivalences, cofibrations and fibrations are closed under retracts.
\item[M4]
(Lifting axiom) Cofibrations have the LLP with respect to trivial fibrations, and trivial cofibrations have the LLP with respect to fibrations.
\item[M5]
(Factorization axiom) Any map of ${\cat M}$ admits a factorization as a cofibration followed by a trivial fibration, and a factorization as a trivial cofibration followed by a fibration.
\end{description}
\end{defn}

The axiom M3 states that given a commutative diagram
\begin{center}
$\xymatrix{
	A \ar[r] \ar[d]_f \ar@/^1pc/[rr]^{1_A} &
	C \ar[d]^g \ar[r] & A \ar[d]^f \\
	B \ar[r] \ar@/_1pc/[rr]_{1_B} &
	D \ar[r] & B
	}$
\end{center}
if $g$ is a weak equivalence (resp. cofibration, resp. fibration) then so is its retract $f$.

\begin{prop}
\label{prop:llpcharact}
In a Quillen model category any two of the following classes of maps of ${\cat M}$ - the cofibrations, the trivial cofibrations, the fibrations and the trivial fibrations - determine each other by the following rules: a map is
\begin{enumerate}
\item
A cofibration $\Leftrightarrow$ it has the LLP with respect to all trivial fibrations
\item
A trivial cofibration $\Leftrightarrow$ it has the LLP with respect to all fibrations
\item
A fibration $\Leftrightarrow$ it has the RLP with respect to all trivial cofibrations
\item
A trivial fibration $\Leftrightarrow$ it has the RLP with respect to all cofibrations $\square$
\end{enumerate}
\end{prop}

\begin{proof}
This is a direct consequence of the axioms M1-M5.
\end{proof}

\begin{prop}
\label{prop:relclosedmodcat}
Any Quillen model category is an ABC model category. 
\end{prop}

\begin{proof} A Quillen model category trivially satisfies the axioms CF2 and CF4. The axioms CF1, CF3, CF5 and CF6 are satisfied as a consequence of \propositionref{prop:llpcharact}. A dual argument shows that a Quillen model category satisfies the axioms F1-F6.
\end{proof}

If ${\cat M}$ is a Quillen model category and ${\cat D}$ is a small category, then the category of diagrams ${\cat M}^{\cat D}$ does not generally form a Quillen model category (except in important particular cases, for example when ${\cat M}$ is cofibrantly generated or when ${\cat D}$ is a Reedy category). But we will see further down (\theoremref{thm:generalpointwisecofstructure}) that ${\cat M}^{\cat D}$ carries an ABC model category structure, and in that sense one can always 'do homotopy theory' on ${\cat M}^{\cat D}$.

A Quillen model category ${\cat M}$ is called {\it left proper} if for any pushout diagram
\begin{center}
$ \xymatrix {
	A \ar[r]^j \ar@{>->}[d]_i &
	C \ar@{>-->}[d] \\
	B \ar@{-->}[r]^{j^{'}} &
	D
	}$
\end{center}
with $i$ a cofibration and $j$ a weak equivalence, the map $j^{'}$ is a weak equivalence. ${\cat M}$ is called {\it right proper} if for any pullback diagram
\begin{center}
$ \xymatrix {
	D \ar@{-->}[r]^{q^{'}} \ar@{-->>}[d] &
	B \ar@{->>}[d]^p \\
	C \ar[r]^q &
	A
	}$
\end{center}
with $p$ a fibration and $q$ a weak equivalence, the map $q^{'}$ is a weak equivalence. ${\cat M}$ is called {\it proper} if it is left and right proper.

From this definition and from \propositionref{prop:relclosedmodcat} we immediately get
\begin{prop}
\label{prop:relproperclosedmodcat}
Any proper Quillen model category is a proper ABC model category. $\square$
\end{prop}

\section{Baues cofibration categories}
We next turn our attention to the notion of (co)fibration category as defined by Baues. We state below the axioms of a Baues cofibration category, in a slightly modified but equivalent form to \cite{Baues1}, Sec. 1.1.

\begin{defn}
\label{defn:bauescofcat}
A Baues cofibration category $({\cat M}$, ${\cat W}$, ${\cat Cof})$ consists of a category ${\cat M}$ and two distinguished classes of maps ${\cat W}$ and ${\cat Cof}$ - the weak equivalences and the cofibrations of ${\cat M}$ - subject to the axioms below:
\begin{description}
\item[BCF1] 
All isomophisms of ${\cat M}$ are trivial cofibrations. Cofibrations are stable under composition. 
\item[BCF2] 
(Two out of three axiom) If $f, g$ are maps of ${\cat M}$ such that $gf$ is defined, and if two of $f, g, gf$ are weak equivalences, then so is the third.
\item[BCF3] 
(Pushout and excision axiom) Given a solid diagram in ${\cat M}$, with $i$ a cofibration,
\begin{center}
$ \xymatrix {
	A \ar[r]^f \ar@{>->}[d]_i &
	C \ar@{>-->}[d]^j \\
	B \ar@{-->}[r]^g &
	D
	}$
\end{center}
then the pushout exists in ${\cat M}$ and $j$ is a cofibration. Moreover:
\begin{enumerate}
\item[(a)]
If $i$ is a trivial cofibration, then $j$ is a trivial cofibration
\item[(b)]
If $f$ is a weak equivalence, then $g$ is a weak equivalence.
\end{enumerate}
\item[BCF4]
(Factorization axiom) Any map of ${\cat M}$ admits a factorization as a cofibration followed by a weak equivalence.
\item[BCF6]
(Axiom on fibrant models) For each object $A$ of ${\cat M}$ there is a trivial cofibration $A \ra B$, with $B$ satisfying the property that each trivial cofibration $C \ra B$ admits a left inverse.
\end{description}
\end{defn}

A Baues cofibration category may not have  an initial object, but if it does then by BCF1 the initial object is cofibrant. We will not state the axioms for a Baues fibration category - they are dual to the above axioms.

\begin{prop}[Relation with Baues cofibration categories]
\label{prop:relcatcofcatbaues}
\mbox{}
\begin{enumerate}
\item
Any Baues cofibration category with an initial object is a left proper precofibration category.
\item
Any Baues fibration category with a terminal object is a right proper prefibration category.
\end{enumerate}
\end{prop}

\begin{proof}
Easy consequence of the axioms and of \lemmaref{lem:equivaxiomcf3}.
\end{proof}

\section{Waldhausen categories}
We list below the axioms we'll use for a Waldhausen cofibration category. These axioms are equivalent to the axioms Cof1-Cof3, Weq1 and Weq2 of \cite{Waldhausen1}. Axioms for a Waldhausen fibration category will of couse be dual to the axioms below. 
\begin{defn}
\label{defn:waldcofcat}
A Waldhausen cofibration category $({\cat M}$, ${\cat W}$, ${\cat Cof})$ consists of a category ${\cat M}$ and two distinguished classes of maps ${\cat W}$ and ${\cat Cof}$ - the weak equivalences and the cofibrations of ${\cat M}$, subject to the axioms below:
\begin{description}
\item[WCF1] 
${\cat M}$ is pointed. All isomophisms of ${\cat M}$ are trivial cofibrations. All objects of ${\cat M}$ are cofibrant. Cofibrations are stable under composition. 
\item[WCF2] 
(Pushout axiom) Given a solid diagram in ${\cat M}$, with $i$ a cofibration,
\begin{center}
$ \xymatrix {
	A \ar[r] \ar@{>->}[d]_i &
	C \ar@{>-->}[d]^j \\
	B \ar@{-->}[r] &
	D
	}$
\end{center}
then the pushout exists in ${\cat M}$ and $j$ is a cofibration
\item[WCF3]
(Gluing axiom) In the diagram below
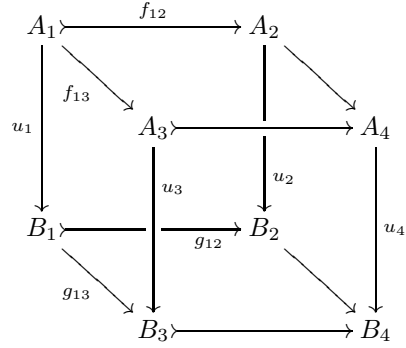
\begin{center}
$ \xymatrix {
        A_1 \ar@{>->}[rr]^{f_{12}} \ar[dd]_{u_1} \ar[rd]_{f_{13}} &  & A_2 \ar'[d][dd]^{u_2} \ar[rd] & \\
        & A_3 \ar@{>->}[rr] \ar[dd]^(.3){u_3} & & A_4 \ar[dd]^{u_4} \\
        B_1 \ar@{>->}'[r][rr]_{g_{12}} \ar[rd]_{g_{13}} & & B_2 \ar[rd] & \\
        & B_3 \ar@{>->}[rr] & & B_4
	}$
\end{center}
if $f_{12}, g_{12}$ are cofibrations, $u_1, u_2, u_3$ are weak equivalences and the top and bottom faces are pushouts, then $u_4$ is a weak equivalence.
\end{description}
\end{defn}

We have the following

\begin{prop}[Relation with Waldhausen cofibration categories]
\label{prop:relcatcofcatwaldhausen}
\mbox{}
\begin{enumerate}
\item
\begin{enumerate}
\item
If ${\cat M}$ is a pointed precofibration category, then ${\cat M}_{cof}$ is a Waldhausen cofibration category.
\item
If ${\cat M}$ is a Waldhausen cofibration category satisfying the 2 out of 3 axiom CF2 and the cylinder axiom CFObj4, then it is a precofibration category.
\end{enumerate}
\item
\begin{enumerate}
\item
If ${\cat M}$ is a pointed prefibration category, then ${\cat M}_{fib}$ is a Waldhausen fibration category.
\item
If ${\cat M}$ is a Waldhausen fibration category satisfying the 2 out of 3 axiom F2 and the path object axiom FObj4, then it is an prefibration category.
\end{enumerate}
\end{enumerate}
\end{prop}

\begin{proof}
Part (1) (a) follows from the Gluing Lemma \ref{lem:gluinglemma1}, part (1) (b) from \lemmaref{lem:equivaxiomcf3} and \propositionref{prop:relcatcofobj}, and part (2) is dual to the above.
\end{proof}

Waldhausen categories have been further studied by Thomason-Trobaugh \cite{Thomason-Trobaugh}, Weiss and Williams \cite{Weiss1}, \cite{Weiss2}, \cite{Weiss3}.
\chapter{Examples of ABC and Quillen model categories}
\label{chap:examples}

\section{Topological spaces}
\label{sec:exampletopspaces}
\index{$Top$}
In the category of topological spaces $Top$, we denote $B^A$ for the space of continuous maps $A \ra B$, with the compact-open topology. The exponential map $f \mapsto (a \mapsto (b \mapsto f(a, b)))$ defines a bijection
\begin{equation}
\label{eqn:exp}
Top(A \times B, C) \cong Top(A, C^B)
\end{equation}
for any topological spaces $A, B, C$ with $B$ locally compact and Hausdorff separated. 


Denote $I$ the unit interval in $Top$. For any space $A$, we denote $i_0, i_1 \colon A \ra I \times A$ the inclusions $i_k(a) = (k, a)$ for $k = 0, 1$, and denote $p \colon I \times A \ra A$ the second factor projection. We also denote $p_0, p_1 \colon A^I \ra A$ the evaluation maps $p_k(f) = f(k)$, and $i \colon A \ra A^I$ the constant-value map $i(a)(t) = a$. 

\begin{defn}
\label{defn:homotopy}
Suppose that $f_0, f_1 \colon A \ra B$ is a pair of maps in $Top$. A homotopy $h \colon f_0 \simeq f_1$ (also denoted $f_0 \xymatrix{\ar[r]|h &} f_1$) is a map $h \colon I \times A \ra B$ with $hi_k = f_k$. Using the bijection \eqref{eqn:exp}, a homotopy $h \colon I \times A \ra B$ amounts to a map $h^{'} \colon A \ra B^I$ with $p_kh^{'} = f_k$.
\end{defn}

\begin{defn}
\label{defn:homotopyundera}
Suppose that $f_0, f_1$ are two maps under $A$, meaning that $f_kb=c$ for $k = 0, 1$.
\[
\xymatrix{
  & A \ar[dl]_b \ar[dr]^{c} & \\
  B \ar@<.5ex>[rr]^{f_0} \ar@<-.5ex>[rr]_{f_1} & & C
}
\]
A homotopy under $A$, denoted $h \colon f_0 \overset{A}{\simeq} f_1$ is a homotopy $h \colon f_0 \simeq f_1$ which is {\it constant} when restricted to $I \times A$, meaning that $h \comp (I \times b) = cp$, where $p \colon I \times A \ra A$ is the projection.
\end{defn}

\begin{defn}
\label{defn:homotopyovera}
Suppose that $f_0, f_1$ are two maps over $A$, meaning that $cf_k=b$ for $k = 0, 1$.
\[
\xymatrix{
  B \ar@<.5ex>[rr]^{f_0} \ar@<-.5ex>[rr]_{f_1} \ar[dr]_b & & C \ar[dl]^c \\
  & A &
}
\]
A homotopy over $A$, denoted $h \colon f_0 \underset{A}{\simeq} f_1$ is a homotopy $h \colon f_0 \simeq f_1$ which is {\it constant} when corestricted to $C$, meaning that $ch = bp$, where $p \colon I \times B \ra B$ is the projection.
\end{defn}

Homotopy $\simeq$ and its relative counterparts $\overset{A}{\simeq}$, $\underset{A}{\simeq}$ are equivalence relations.

\begin{defn}
\label{defn:homotequiv}
A map $f \colon A \ra B$ is a {\it homotopy equivalence} if it admits a map $g \colon B \ra A$ with $fg \simeq 1_B$, $gf \simeq 1_A$.
\end{defn}

Homotopy equivalences satisfy the 2 out of 3 axiom, and are closed under small sums and small products.

\index{Hurewicz (co)fibration}
\begin{defn}
\label{defn:hep}
\mbox{}
\begin{enumerate}
\item
A map $i \colon A \ra B$ is a {\it Hurewicz cofibration} if it has the left lifting property with respect to all maps of the form $p_0 \colon C^I \ra C$.
\item
A map $p \colon A \ra B$ is a {\it Hurewicz fibration}  if it has the right lifting property with respect to all maps of the form $i_0 \colon C \ra I \times C$.
\end{enumerate}
\end{defn}

By definition, a Hurewicz cofibration (or fibration) is {\it trivial} if it is also a homotopy equivalence. The following result is immediate.

\begin{lem}
\label{lem:hurcomp}
\mbox{}
\begin{enumerate}
\item
(Trivial) Hurewicz cofibrations are closed under compositions and small sums.
\item
(Trivial) Hurewicz fibrations are closed under compositions and small products. $\square$
\end{enumerate}
\end{lem}

Using the bijection \eqref{eqn:exp} we see that $i \colon A \ra B$ is a Hurewicz cofibration iff it has the {\it homotopy extension property (HEP)}, meaning that for any maps $f, h$ with $fi = hi_0$ 
\begin{equation}
\label{eqn:hep}
\xymatrix{
  A \ar[r]^-{i_0} \ar[d]_i & I \times A \ar[d]_{I \times i} \ar[rdd]^h & \\
  B \ar[r]^-{i_0} \ar[drr]_f & I \times B \ar@{-->}[rd]^(.3){h^{'}} & \\
  && C
}
\end{equation}
there exists a homotopy $h^{'}$ keeping the diagram commutative. Dually, a map $p \colon A \ra B$ is a Hurewicz fibration iff it has the {\it homotopy lifting property (HLP)}, meaning that for any maps $f, h$ with $pf = p_0h$ 
\begin{equation}
\label{eqn:hlp}
\xymatrix{
  C \ar[rrd]^f \ar[rdd]_h \ar@{-->}[rd]^(.9){h^{'}} && \\
  & A^I \ar[r]_{p_0} \ar[d]^{p^I} & A \ar[d]^p \\
  & B^I \ar[r]_{p_0} & B
}
\end{equation}
there exists a homotopy $h^{'}$ keeping the diagram commutative.

\begin{lem}
\label{lem:hurewiczreplace}
\mbox{}
\begin{enumerate}
\item
Suppose that $i$ is a Hurewicz cofibration and that the diagram 
\[
\xymatrix{
  & A \ar@{>->}[dl]_i \ar[dr]^{g} & \\
  B \ar[rr]^{f} & & C
}
\]
is homotopy commutative $fi \simeq g$. Then there exists a map $f^{'} \simeq f$ with $f^{'}i = g$.
\item
Suppose that $p$ is a Hurewicz fibration and that the diagram
\[
\xymatrix{
  A \ar[rr]^{f} \ar[dr]_g & & B \ar@{->>}[dl]^{p} \\
  & C &
}
\]
is homotopy commutative $pf \simeq g$. Then there exists a map $f^{'} \simeq f$ with $pf^{'} = g$.
\end{enumerate}
\end{lem}

\begin{proof}
We only prove (1). Denote $h \colon fi \simeq g$ a homotopy. Using HEP for $i$, we extend $h$ to a homotopy $h^{'} \colon I \times B \ra C$ with $h^{'} \comp (I \times i) = h$ and $h^{'}i_0 = f$. Then $f^{'} = h^{'}i_1 \colon B \ra C$ has the desired properties.
\end{proof}

\index{Lemma!Dold lema}
\begin{lem}[A. Dold]
\label{lem:dold}
\mbox{}
\begin{enumerate}
\item
Suppose that $i$, $i^{'}$ are Hurewicz cofibrations and $fi = i^{'}$
\[
\xymatrix{
  & A \ar@{>->}[dl]_i \ar@{>->}[dr]^{i^{'}} & \\
  B \ar[rr]^{f} & & C
}
\]
If $f$ is a homotopy equivalence, then it is a homotopy equivalence {\it under} $A$.
\item
Suppose that $p$, $p^{'}$ are Hurewicz fibrations and $p^{'}f = p$
\[
\xymatrix{
  A \ar@{->>}[rr]^{p} \ar[dr]_g & & B \ar@{->>}[dl]^{p^{'}} \\
  & C &
}
\]
If $f$ is a homotopy equivalence, then it is a homotopy equivalence {\it over} $C$.
\end{enumerate}
\end{lem}

\begin{proof}
We only prove (1), and we do that in several steps.

Step 1. It suffices to show that $f$ admits a right homotopy inverse $f^{'}$ under $A$. For in that case, by the same reasoning $f^{'}$ has a right homotopy inverse $f^{''}$ under $A$, which has itself a right homotopy inverse $f^{'''}$ under $A$. We therefore have $f \overset{A}{\simeq} f^{''}$ and $f^{'} \overset{A}{\simeq} f^{'''}$, from which we see that $f^{'}$ is also a left homotopy inverse of $f$ under $A$.

Step 2. To show that $f$ admits a right homotopy inverse $f^{'}$ under $A$, it suffices to assume that $i = i^{'}$ and $f \simeq 1_B$. Indeed, suppose that $f^{'} \colon C \ra B$ is a homotopy inverse of $f$ (so $f^{'}f \simeq 1_{B}$). By \lemmaref{lem:hurewiczreplace}, we can find $f^{''} \simeq f^{'}$ with $f^{''}i^{'} = i$. Then $f^{''}f \simeq 1_B$, so by our assumption $f^{''}f \overset{A}{\simeq} 1_B$, and $f^{''}$ is a right homotopy inverse under $A$ to $f$.

Step 3. In the commutative diagram with $f \simeq 1_B$
\[
\xymatrix{
  & A \ar@{>->}[dl]_i \ar@{>->}[dr]^{i} & \\
  B \ar[rr]^{f} & & B
}
\]
we'd like to show that $f \overset{A}{\simeq} 1_B$. Denote $f \xymatrix{\ar[r]|h &} 1_{B}$ the homotopy $h \colon f \simeq 1_B$. Using HEP for $i$, the homotopy $i \xymatrix{\ar[rr]|{h\comp{(I \times i)}} &&} i$ extends to a homotopy $1_{B} \xymatrix{\ar[r]|{h^{'}} &} f^{'}$, for some map $f^{'} \colon B \ra B$. We will show that $f^{'}f \overset{A}{\simeq} 1_B$, which completes the proof of Step3, and with it the proof of our Lemma.

The composite $f^{'}f \xymatrix{&& \ar[ll]|{h^{'} \comp (I \times f)}} f \xymatrix{\ar[r]|{h} &} 1_{B}$ defines a homotopy $f^{'}f \simeq 1_{B}$. This homotopy restricts via $i$ to the homotopy $i \xymatrix{&& \ar[ll]|{h \comp (I \times i)}} i \xymatrix{\ar[rr]|{h \comp (I \times i)} &&} i$. The square
\[
\xymatrix{
  i \ar[d]|{id} & \; \ar@{}[drr]|{\alpha} & i \ar[ll]|{h \comp (I \times i)} \ar[rr]|{h \comp (I \times i)} && i \ar[d]|{id} \\
  i \ar[rrrr]|{id} &&& \; & i
}
\]
can be filled with a homotopy $\alpha \colon I^2 \times A \ra B$. We apply HEP to the Hurewicz cofibration $I \times A \ra I \times B$ with respect to the top edge of the previous diagram. This constructs a map $\beta \colon I^2 \times B \ra B$
\[
\xymatrix{
  f^{'}f \ar[d] & \; \ar@{}[drr]|{\beta} & f \ar[ll]|{h^{'} \comp (I \times f)} \ar[rr]|{h} && 1_{B} \ar[d] \\
  \, \ar[rrrr] &&& \; & \, 
}
\]
Tracing the left, bottom and right edges above yields the desired homotopy $f^{'}f \overset{A}{\simeq} 1_B$.
\end{proof}

\begin{lem}
\label{lem:trivhurcof}
\mbox{}
\begin{enumerate}
\item
Any trivial Hurewicz cofibration $i \colon A \ra B$ admits a strong deformation retract, i.e. a map $r \colon B \ra A$ such that $ri = 1_A$ and $ir \overset{A}{\simeq} 1_B$.
\item
Any trivial Hurewicz fibration $p \colon A \ra B$ admits a strong deformation section, i.e. a map $s \colon B \ra A$ such that $ps = 1_B$ and $sp \underset{A}{\simeq} 1_A$.
\end{enumerate}
\end{lem}

\begin{proof}
We only prove (1). Denote $r^{'}$ a homotopy inverse of $i$. By \lemmaref{lem:hurewiczreplace} applied to the homotopy commutative diagram
\[
\xymatrix{
  & A \ar@{>->}[dl]_i \ar[dr]^{1_A} & \\
  B \ar[rr]^{r^{'}} & & A
}
\]
there exists $r \colon B \ra A$, with $r \simeq r^{'}$ and $ri = 1_A$. From $ir \simeq 1_B$ in the commutative diagram
\[
\xymatrix{
  & A \ar@{>->}[dl]_i \ar@{>->}[dr]^{i} & \\
  B \ar[rr]^{ir} & & B
}
\]
by Dold's \lemmaref{lem:dold} we have $ir \overset{A}{\simeq} 1_B$.
\end{proof}

\begin{lem}
\label{lem:hurpush}
\mbox{}
\begin{enumerate}
\item
(Trivial) Hurewicz cofibrations are stable under pushouts.
\item
(Trivial) Hurewicz fibrations are stable under pullbacks.
\end{enumerate}
\end{lem}

\begin{proof}
We only prove (1). Hurewicz cofibrations are defined as having the left lifting property with respect to all $p_0 \colon B^I \ra B$, so Hurewicz cofibrations are stable under pushouts. Strong deformation retracts are also stable under pushouts, and in view of \lemmaref{lem:trivhurcof} so are trivial Hurewicz cofibrations.
\end{proof}

\begin{lem}
\label{lem:hurcolim}
\mbox{}
\begin{enumerate}
\item
If $A_0 \rightarrowtail A_1  \rightarrowtail ... \rightarrowtail A_n  \rightarrowtail ...$ is a sequence of (trivial) Hurewicz cofibrations, then its composition $A_0 \ra \colim^n A_n$ is again a (trivial) Hurewicz cofibration.
\item
If $... \ra A_n \ra ... \ra A_1 \ra A_0$ is a sequence of (trivial) Hurewicz fibrations, then $\lim^n A_n \ra A$ is again a (trivial) Hurewicz fibration.
\end{enumerate}
\end{lem}

\begin{proof}
We only prove (1). Denote $i_n \colon A_n \rightarrowtail A_{n+1}$, and $i \colon A_0 \ra \colim^n A_n = A$. If all $i_n$ are Hurewicz cofibrations, then $i$ has the left lifting property with respect to all $p_0 \colon B^I \ra B$, so $i$ is also a Hurewicz cofibration.

Suppose now that each $i_n$ is trivial. By \lemmaref{lem:trivhurcof}, there exist strong deformation retracts $r_n \colon A_{n+1} \ra A_n$ with $r_ni_n = 1_{A_n}$, and with relative homotopies $1_{A_{n+1}} \overset{A_n}{\simeq} i_nr_n$.

Define $r \colon A \ra A_0$ the colimit of the compositions $r_0r_1 ... r_n$. It is a retract of $i$, and we need to show that $1_A \simeq ir$.

The homotopy $1_{A_1} \overset{A_0}{\simeq} i_0r_0$ extends by HEP for $A_1 \rightarrowtail A$ to a homotopy $h_0 \colon 1_{A} \overset{A_0}{\simeq} s_0$, for some map $s_0 \colon A \ra A$ with $s_0 |_{A_1} = i_0r_0$.

The composite homotopy $s_0 |_{A_2} \simeq 1_{A_2} \simeq i_1i_0r_0r_1$, by Dold's \lemmaref{lem:dold} yields a homotopy $s_0 |_{A_2} \overset{A_1}{\simeq} i_1i_0r_0r_1$, which extends by HEP for $A_2 \rightarrowtail A$ to a homotopy $h_1 \colon s_0 \overset{A_1}{\simeq} s_1$, for some map $s_1 \colon A \ra A$ with $s_1 |_{A_2} = i_1i_0r_0r_1$.

By induction, we construct maps $s_n \colon A \ra A$ with $s_n |_{A_{n+1}} = i_n...i_0$$r_0...r_n$, and homotopies $h_n \colon s_{n-1} \overset{A_{n}}{\simeq} s_n$. Stitching together all the $h_n$, we obtain the desired homotopy $h \colon 1_A \simeq ir$ as
\[
h(1 - \frac{1}{2^n} + \frac{t}{2^{n+1}}, a) = h_n(t, a) \;\;\; (n \ge 0, \; 0 \le t \le 1, \; a \eps A)
\]
\end{proof}

\begin{thm}[The Hurewicz model structure]
\label{thm:hurtop}
$Top$ is an ABC model category, with all objects cofibrant and fibrant, with:
\begin{enumerate}
\item
Homotopy equivalences as weak equivalences
\item
Hurewicz fibrations as fibrations
\item
Hurewicz cofibrations as cofibrations.
\end{enumerate}
\end{thm}

\begin{proof}
The axioms CF1, CF2, CF5 are straightforward. The axiom CF4 is given by the classic mapping cylinder construction in $Top$. The axiom CF3 is proved by \lemmaref{lem:hurpush}, and CF6 by \lemmaref{lem:hurcolim}. All maps $p_0 \colon B^I \ra B$ admit the constant map $B \ra B^I$ as a section, from which any space is Hurewicz cofibrant. A similar proof shows that any space is also Hurewicz fibrant.
\end{proof}

We will mention without proof two more results regarding model structures on $Top$.

\begin{thm}[Strom, \cite{Strom1}]
\label{thm:stromtop}
$Top$ is a Quillen model category, with:
\begin{enumerate}
\item
Homotopy equivalences as weak equivalences
\item
Hurewicz fibrations as fibrations
\item
Hurewicz cofibrations that are closed maps, as cofibrations.
\end{enumerate}
\end{thm}

\index{Serre fibration}
\begin{defn}
\label{defn:serrefib}
\mbox{}
A map $p \colon A \ra B$ is a {\it Serre fibration}  if it has the right lifting property with respect to all face maps  $I^{n-1} \ra I^n$.
\end{defn}

\begin{thm}[Quillen, \cite{Quillen1}]
\label{thm:quillentop}
$Top$ is a Quillen model category, with:
\begin{enumerate}
\item
Weak homotopy equivalences as weak equivalences
\item
Serre fibrations as fibrations
\item
Maps with the left lifting property with respect to trivial Serre fibrations as cofibrations.
\end{enumerate}
\end{thm}

\section{Abelian categories}
\label{sec:exampleabcat}

We recall the following basic results regarding abelian categories:

\comment{
Suppose that ${\cat A}$ is an abelian category. This means that 

\begin{description}
\item[A0]
${\cat A}$ is pointed
\item[A1]
${\cat A}$ is closed under finite sums and products
\item[A3]
Every map has a kernel and a cokernel
\item[A4]
Every monic is a kernel of a map. Every epic is a cokernel of a map.
\end{description}

From these axioms, one shows \cite{Freyd1} that:
\begin{enumerate}
\item
${\cat A}$ is closed under finite limits and colimits
\item
For any pair of objects, the natural map $A \Sum B \ra A \times B$ is an isomorphism, whose common value is denoted $A \oplus B$
\item
${\cat A}$ is additive, i.e. ${\cat A}(A, B)$ is a natural abelian group with composition induced by the diagonal $B \oplus B \ra B$
\item
A map $f \colon A \ra B$ is monic iff $Ker \; f = \initial$, and is epic iff $Coker \; f = \initial$
\end{enumerate}
} 

\begin{lem}
\label{lem:pushabelian1}
In an abelian category, given a square
\[
\xymatrix{
  A \ar[r]^f \ar[d]_g & B \ar[d]^{g^{'}} \\
  C \ar[r]^{f^{'}} & D
}
\]
consider the sequence $\xymatrix{A \ar[r]^-{(f, g)} & B \oplus C \ar[r]^-{f^{'} - g^{'}} & D}$. Then:
\begin{enumerate}
\item
$A \ra B \oplus C \ra D$ is zero iff the square commutes
\item
$A \ra B \oplus C \ra D \ra 0$ is exact iff the square is a pushout
\item
$0 \ra A \ra B \oplus C \ra D $ is exact iff the square is a pullback
\item
$0 \ra A \ra B \oplus C \ra D \ra 0$ is exact iff the square is a pushout and a pullback
\item
If the square is a pushout with $f$ or $g$ monic, then it is also a pullback
\item
If the square is a pullback with $f^{'}$ or $g^{'}$ epic, then it is also a pushout $\square$
\end{enumerate}
\end{lem}

\begin{proof}
Left to the reader (who may wish consult for example \cite{Freyd1}).
\end{proof}

\begin{lem}
\label{lem:pushabelian2}
In an abelian category:
\begin{enumerate}
\item
The pushout of a map $f$ is epic iff $f$ is epic.
\item
The pushout of a monic is monic.
\item
The pullback of a map $f^{'}$ is monic iff $f^{'}$ is monic.
\item
The pullback of an epic is epic $\square$
\end{enumerate}
\end{lem}

\begin{proof}
This is a ready consequence of \lemmaref{lem:pushabelian1}.
\end{proof}

In fact, any small abelian category admits, by a theorem of Mitchell, a fully faithful exact functor to a category of modules. As a consequence, to prove any statement involving small diagrams, (co)kernels, (co)images, monics, epics, pushouts, pullbacks and finite sums (like \lemmaref{lem:pushabelian1} and \lemmaref{lem:pushabelian2}) in an arbitrary abelian category ${\cat A}$ it suffices to prove that same statement for categories of modules.

\subsection*{Complexes of objects}
\index{$C({\cat A}), C^+({\cat A}), C^-({\cat A}), C^b({\cat A})$}
Denote $C({\cat A})$ the category of $\mathbb{Z}$-graded complexes of objects in ${\cat A}$ 
\begin{center}
$... \lra A_{n+1} \overset{a_{n+1}}{\lra} A_n \overset{a_{n}}{\lra} A_{n-1} \lra ...$

$a_n a_{n+1} = 0$ for any $n \eps {\mathbb Z}$
\end{center}
Denote $C^+({\cat A})$, $C^-({\cat A})$, $C^b({\cat A})$ the full subcategories of $C({\cat A})$ having as objects the bounded below, bounded above resp. bounded complexes. All four are abelian categories, with (co)kernels and (co)images computed degreewise. We recall the classic

\begin{lem}[Snake lemma]
\label{lem:snake}
\index{lemma!snake lemma}
Any short exact sequence $0 \ra A \sdot \overset{f}{\ra} B \sdot \overset{g}{\ra} C \sdot \ra 0$ in $C({\cat A})$ induces a natural long exact sequence in homology
\[
... \lra H_n A \overset{H_nf}{\lra} H_n B \overset{H_ng}{\lra} H_n C \overset{\delta_n}{\lra} H_{n-1} A \lra ...
\]
By Mitchell's Theorem, it suffices to define the connecting homomorphism $\delta_n$ for the case when ${\cat A}$ is a category of modules, in which case in the next diagram
\[
\xymatrix{
  0 \ar[r] & A_n \ar[d]_{a_n} \ar[r]^{f_n} & B_n \ar[d]^{b_n} \ar[r]^{g_n} & C_n \ar[d]^{c_n} \ar[r] & 0 \\
  0 \ar[r] & A_{n-1} \ar[r]^{f_{n-1}} & B_{n-1} \ar[r]^{g_{n-1}} & C_{n-1} \ar[r] & 0
}
\]
$\delta_n$ sends the representative of an element $x_n \eps Ker \, c_n$ to the representative of $y_{n-1} \eps Ker \, a_{n-1}$, where $y_{n-1}$  is constructed observing that there exists $z_n \eps B_n$ with $g_nz_n = x_n$, so $g_{n-1}b_n z_n = 0$ which allows us to define $y_{n-1} = f^{-1}_{n-1}b_nz_n$. Furthermore, the definition of $\delta_n$ outlined above does not depend on the choices involved. $\square$
\end{lem}

\index{quasi-isomorphism}
A map of complexes $A \sdot \ra B \sdot$ is called a {\it quasi-isomorphism} if it induces an isomorphism in homology. The homotopy category of $C({\cat A})$ (resp. $C^{+,-,b}({\cat A})$) with respect to quasi-isomorphisms is denoted $D({\cat A})$ (resp. $D^{+,-,b}({\cat A})$), and is called the {\it derived} category of ${\cat A}$.

\begin{thm}
\label{thm:abcatisabccat1}
For an abelian category ${\cat A}$, the categories of chain complexes $C({\cat A})$ (resp. $C^{+,-,b}({\cat A})$), with quasi-isomorphisms as weak equivalences, with monics as cofibrations and with epics as fibrations are all pointed ABC premodel categories, with all objects at once fibrant and cofibrant.
\end{thm}

\begin{proof}
Axiom CF3 (1) is given by \lemmaref{lem:pushabelian2} (2). Axiom CF3 (2) follows from the Snake lemma. To prove CF4, we construct the factorization of a map $f \colon A \sdot \ra B \sdot$ as a monic $f^{'}$ followed by a quasi-isomorphism $r$, as follows:
\[
\xymatrix {
  A \sdot \ar@{>->}[r]^-{f^{'}} & A \sdot \oplus A \sdot[1] \oplus B \sdot \ar[r]^-r_-\sim & B \sdot
}
\]
In the middle complex, the boundary map
\[
A_n \oplus A_{n-1} \oplus B_n \ra A_{n-1} \oplus A_{n-2} \oplus B_{n-1}
\]
is given by the matrix
\[
\begin{pmatrix}
  a_n & 0 & 0 \\
  id & -a_{n-1} & 0  \\
  0 & 0 & b_n 
\end{pmatrix}
\]

The ABC prefibration category axioms are proved in a dual fashion.
\end{proof}

\subsection*{Exact categories}
For an abelian category ${\cat A}$, suppose that ${\cat E} \subset {\cat A}$ is an exact subcategory, i.e. a full subcategory with the property that for any exact sequence
\[
0 \ra A \ra B \ra C \ra 0
\]
if two out of $A, B, C$ are in ${\cat E}$, then so is the third. Denote $C({\cat E})$ (resp. $C^{+, -, b}({\cat E})$) the category of complexes (resp. bounded below, bounded above, bounded complexes) of maps in ${\cat E}$. 

We say that a monic (resp. epic) map $f$ of ${\cat A}$ is {\it admissible}, if it is also a map of ${\cat E}$. We can now state the following stronger variant of \theoremref{thm:abcatisabccat1}:

\begin{thm}
\label{thm:abcatisabccat2}
Suppose that ${\cat A}$ is an abelian category and that ${\cat E} \subset {\cat A}$ is an exact subcategory. Then $C({\cat E})$ (resp. $C^{+,-,b}({\cat E})$), with quasi-isomorphisms as weak equivalences, with admissible monics as cofibrations and admissible epics as fibrations are all pointed ABC premodel categories, with all objects at once fibrant and cofibrant.
\end{thm}

\begin{proof}
Axiom CF3 (1) follows from \lemmaref{lem:pushabelian2} (2), observing that the push in ${\cat E}$ of an admissible monic exists and is again an admissible monic. Axioms CF3 (2) and CF4 are proved the same way as for \theoremref{thm:abcatisabccat1}.
\end{proof}

\subsection*{The Grothendieck axioms AB4, AB5} 
These are axioms that specify exactness properties for small sums, resp. for filtered colimits. Our reference here is Grothendieck's {\it Toh\^{o}ku} paper \cite{Grothendieck4}, also explained at length in \cite{Popescu1}. We start by recalling a classic adjunction lemma for abelian categories. 

\begin{lem}
\label{lem:adjabelian}
Suppose that $u \dashv v$, $u \colon {\cat A} \rightleftarrows {\cat B} \colon v$ is a pair of adjoint functors between two abelian categories. Then $u$ is right exact, and $v$ is left exact.
\end{lem}

\begin{proof}
Suppose that $A_0 \overset{f}{\ra} A_1 \overset{g}{\ra} A_2 \ra 0$ is an exact sequence in ${\cat A}$. Suppose that $B \eps {\cat B}$ is an object and $h \colon uA_1 \ra B$ is a map with $h \comp uf = 0$. 
\begin{equation}
\label{eqn:adjabelian1}
\xymatrix{
  uA_0 \ar[r]^{uf} & uA_1 \ar[r]^{ug} \ar[rd]_{h} & uA_2 \ar@{-->}[d]^{k^{'}} \ar[r] & 0 \\
  & & B &
}
\end{equation}
Denote $h^{'} \colon A_1 \ra vB$ the map adjoint to $h$. We have $h^{'}f = 0$, so there exists an unique map $k$ making the next diagram commutative,
\begin{equation}
\label{eqn:adjabelian2}
\xymatrix{
  A_0 \ar[r]^f & A_1 \ar[r]^g \ar[rd]_{h^{'}} & A_2 \ar@{-->}[d]^k \ar[r] & 0 \\
  & & vB &
}
\end{equation}
hence there exists a unique map $k^{'}$ (the adjoint of $k$) making the diagram \eqref{eqn:adjabelian1}  commutative. This shows that $ug = Coker(uf)$, so the top row of \eqref{eqn:adjabelian1} is exact.
\end{proof}

If ${\cat D}$ is a small category and ${\cat A}$ is abelian, then the category of ${\cat D}$-diagrams in ${\cat A}$, denoted ${\cat A}^{\cat D}$, is again abelian, with (co)kernels and (co)images computed pointwise. If ${\cat A}$ is closed under colimits indexed by ${\cat D}$, then the functor 
\[
\colim\nolimits^{\cat D} \colon {\cat A}^{\cat D} \ra {\cat A}
\]
is right exact by \lemmaref{lem:adjabelian}.

Grothendieck's axiom AB4 states that ${\cat A}$ is closed under small sums, and that for any set ${\cat D}$ the sum functor $\oplus^{\cat D}$ is left exact. 

An abelian category is always closed under pushouts; if it is closed under small sums then it is closed under all small colimits. Grothendieck's axiom AB5 states that ${\cat A}$ is closed under small sums (therefore under small colimits), and that for any {\it filtered} small category ${\cat D}$ the functor $\colim\nolimits^{\cat D}$ is left exact. Axiom AB5 implies axiom AB4.

One can formulate dual axioms AB4$^*$, AB5$^*$ for products and filtered limits, however, as explained in loc. cit., if an abelian category satisfies AB4 and AB4$^*$ then it is trivial in the sense that all its objects are isomorphic with $0$.

The next theorem is now immediate:

\begin{thm}
\label{thm:abcatisabccat3}
If an abelian category ${\cat A}$ satisfies AB4 (resp. AB5), then $C({\cat A})$, with quasi-isomorphisms as weak equivalences and monics as cofibrations satisfies CF5 (resp CF6). $\square$
\end{thm}

We will mention without proof the following result of Tibor Beke \cite{Beke1}:
\begin{thm}
\label{thm:abcatisabccat4}
If an abelian category ${\cat A}$ satisfies AB5 and has a generator (i.e. an object $A$ with the property that ${\cat A}(A, -)$ is faithful), then $C({\cat A})$, with quasi-isomorphisms as weak equivalences, with monics as cofibrations, and maps with the RLP with respect to monics as fibrations, forms a Quillen model category.
\end{thm}

\chapter{Kan extensions}
\label{chap:kanextensions}

The purpose of this chapter is to introduce the language of 2-categories and the apparatus of Kan extensions.

\section{The language of 2-categories}
\label{sec:functorsnaturammaps}

\index{2-category}
Recall that a 2-category ${\cat C}$ is a category enriched over categories. This means by definition that for each two objects $A, B$ of ${\cat C}$ the maps $A \ra B$ form the objects of a category ${\cat C}(A, B)$. The the composition functor $c_{ABC} \colon {\cat C}(A, B) \times {\cat C}(B, C) \ra {\cat C}(A, C)$ is required to be associative and to have $1_A \colon * \ra {\cat C}(A, A)$ as a left and right unit, where $*$ denotes the point-category. 

The objects of a 2-category ${\cat C}$ are called 0-cells, the objects of ${\cat C}(A, B)$ are called 1-cells and the morphisms of ${\cat C}(A, B)$ are called 2-cells. A good introduction to 2-categories can be found in Kelly-Street \cite{Kelly-Street1} or Borceux \cite{Borceux}.

This section describes the notation we use for compositions of 1-cells and 2-cells in a 2-category. Each notation has a full form and a simple form. The simple form of the notation is ambiguous, and is only used if it is clear from the context which functor or natural map operation we refer to. 

We denote as usual 1-cells $f: A \ra B$ with a single arrow. Between 1-cells $f, g \colon A \ra B$, we denote 2-cells as $\alpha \colon f \Ra g$, or just $\alpha \colon f \ra g$ if no confusion can occur. 

The composition of 1-cells $f, g$ 
\begin{equation}
\label{eqn:comp1}
\xymatrix{
    A \ar[r]^f &
    B \ar[r]^g &
    C
}
\end{equation}
is just the composition at the level of unenriched hom-sets ${\cat C}(-, -)$ and is denoted $g \comp f \colon A \ra C$, or in simple form $gf$.

The composition of 2-cells $\alpha, \beta$ 
\begin{equation}
\label{eqn:comp2}
\xymatrix{
    A \ar@<4ex>[rr]^(.3)f _\,="f" \ar[rr]^(.3)g ^\,="g1" _\,="g2"  \ar@<-4ex>[rr] ^(.3)h^\,="h" & &
    B
    \ar@{=>} "f";"g1" ^{\alpha}
    \ar@{=>} "g2";"h" ^{\beta}
}
\end{equation}
is composition at the level of ${\cat C}(-, -)$ and is denoted $\beta \odot \alpha \colon f \Ra h$, or in simple form $\beta \alpha$.

The composition of a 2-cell $\alpha$ with a 1-cell $f$
\begin{equation}
\label{eqn:comp3}
\xymatrix{
    A \ar[rr]^f & &
    B \ar@<2ex>[rr]^(.3)g _\,="g" \ar@<-2ex>[rr] ^(.3)h^\,="h" & & 
    C
    \ar@{=>} "g";"h" ^{\alpha}
}
\end{equation}
is just $c_{ABC}(1_f, \alpha)$, and we denote it $\alpha . f \colon gf \Ra hf$ or in simple form $\alpha f$. The composition in the other direction
\begin{equation}
\label{eqn:comp4}
\xymatrix{
    A \ar@<2ex>[rr]^(.3)f _\,="f" \ar@<-2ex>[rr] ^(.3)g^\,="g" & &
    B \ar[rr]^h & & 
    C
    \ar@{=>} "f";"g" ^{\alpha}
  }
\end{equation}
is $c_{ABC}(\alpha, 1_h)$ and we denote it $h . \alpha \colon hf \Ra hg$, or in simple form $h \alpha$.

The notations we have established up until now can be used to completely describe compositions of 1- and 2-cells. However it is convenient to introduce the $\star$ notation to denote the composition of adjacent of 2-cells of planar diagrams.

We will denote the composition of 2-cells $\alpha \colon jf \Ra g$ and $\beta \colon h \Ra ij$
\begin{equation}
\label{eqn:comp5}
\xymatrix {
    & B \ar[dr]^h \ar[dd]^j & \\
    A \ar[ur]^f \ar[dr]_g \ar@{}[rr]|(.3){\big\Downarrow \alpha}|(.7){\big\Downarrow \beta} & & D \\
    & C \ar[ur]_i &
  }
\end{equation}
as $\beta \star \alpha = (i . \alpha) \odot (\beta . f)$.

We use the same $\star$ notation to denote the composition of 2-cells $\alpha \colon f \Ra jg$ and $\beta \colon hj \Ra i$
\begin{equation}
\label{eqn:comp6}
\xymatrix {
    & B \ar[dr]^h & \\
    A \ar[ur]^f \ar[dr]_g \ar@{}[rr]|(.3){\big\Downarrow \alpha}|(.7){\big\Downarrow \beta} & & D  \\
    & C \ar[ur]_i \ar[uu]_j &
  }
\end{equation}
as $\beta \star \alpha = (\beta . g) \odot (h . \alpha)$.

In particular, taking $j$ to be an identity map in (\ref{eqn:comp5}) or (\ref{eqn:comp6}) we get the composition of 2-cells $\beta \star \alpha = c_{ABC}(\alpha, \beta)$
\begin{equation}
\label{eqn:comp7}
\xymatrix{
    A \ar@<2ex>[rr]^(.3)f _\,="f" \ar@<-2ex>[rr] ^(.3)g^\,="g" & &
    B \ar@<2ex>[rr]^(.3)h _\,="h" \ar@<-2ex>[rr] ^(.3)i^\,="i" & & 
    C
    \ar@{=>} "f";"g" ^{\alpha}
    \ar@{=>} "h";"i" ^{\beta}
  }
\end{equation}
Since $c_{ABC}$ is a functor, $(i . \alpha) \odot (\beta . f) = (\beta . g) \odot (h . \alpha)$ and the $\star$ notation is consistent in (\ref{eqn:comp7}) with (\ref{eqn:comp5}) and (\ref{eqn:comp6}). 

In simple form, if no confusion is possible we denote $\beta \alpha$ for $\beta \star \alpha$.

There is a general theorem called the Pasting Theorem regarding compositions of 2-cells of planar diagrams, for which we refer to \cite{Power} . 

\index{2-category!opposites ${\cat C}^{1-op}$, ${\cat C}^{2-op}$, ${\cat C}^{1,2-op}$}
A 2-category ${\cat C}$ admits three flavours of opposites, denoted respectively:
\begin{enumerate}
\item
${\cat C}^{1-op}$ (also denoted ${\cat C}^{co}$), reverting the direction of 1-cells
\item
${\cat C}^{2-op}$ (also denoted ${\cat C}^{op}$), reverting the direction of 2-cells
\item
${\cat C}^{1,2-op}$ (also denoted ${\cat C}^{coop}$), reverting the direction of both 1- and 2-cells
\end{enumerate}

\index{2-category!2-functor}
\index{2-category!2-functor!lax or right weak}
\index{2-category!2-functor!op-lax or left weak}
\index{2-category!2-functor!pseudo 2-functor}
\index{2-category!2-functor!strict}
Let us also recall the various flavours of 2-functors $F \colon {\cat C}_1 \ra {\cat C}_2$. All flavours send the 0, 1 and 2-cells of ${\cat C}_1$ respectively to 0, 1 and 2-cells of ${\cat C}_2$, and preserve the composition and units of 2-cells on the nose. The functor $F$ is called

\begin{enumerate}
\item
{\it lax} or {\it right weak} if it preserves the composition and units of 1-cells up to canonical 2-cells $\Gamma_{f, g} \colon  F_gF_f \Ra F_{gf}$ for $A \overset{f}{\ra} B \overset{g}{\ra} C$ in ${\cat C}_1$ and $\Delta_A \colon 1_{F_A} \Ra F_{1_A}$ for $A$ an object of ${\cat C}_1$, so that if $c \overset{h}{\ra} d$ in ${\cat C}_1$ then the diagrams below are commutative
\[
\xymatrix{
  F_hF_gF_f \ar[d]_{F_h \star \Gamma_{f, g}} \ar[rr]^{\Gamma_{g, h} \star F_f} && F_{hg}F_f \ar[d]^{\Gamma_{f, hg}} \\
  F_hF_{gf} \ar[rr]_{\Gamma_{gf, h}} && F_{hgf}
}
\]
\[
\xymatrix{
  F_f \comp 1_{F_A} \ar[dr]_{=} \ar[r]^{F_f \star \Delta_A} & F_f \comp F_{1_A} \ar[d]^{\Gamma_{1_A, f}} && F_{1_B} \comp F_{f} \ar[d]_{\Gamma_{f, 1_B}} & 1_{F_B} \comp F_f \ar[l]_{\Delta_B \star F_f} \ar[ld]^{=} \\
  & F_{f \comp 1_A} && F_{1_B \comp f} &
}
\]
$\Gamma_{f, g}$ are called {\it composition} 2-cells and $\Delta_A$ are called {\it unit} 2-cells.
\item
{\it op-lax} or {\it left weak} if composition and unit of 1-cells are functorial up to natural 2-cells $\Gamma^{'}_{f, g}, \Delta^{'}_A$ going in the opposite direction, making commutative three diagrams dual to the ones listed above
\item
a {\it pseudo} functor if it is lax and all $\Gamma_{f, g}$, $\Delta_A$ are {\it isomorphisms}. Notice that this is equivalent to saying that $F$ is op-lax and  all $\Gamma^{'}_{f, g}$, $\Delta^{'}_A$ are isomorphisms. A pseudo-functor thus preserves composition and units of 1-cells up to canonical isomorphisms $\Gamma_{f, g}$, $\Delta_A$.
\item
{\it strict} if it is a pseudo-functor whose all $\Gamma_{f, g}$, $\Delta_A$ are {\it identities}. A strict 2-functor thus preserves composition and units of 1-cells on the nose.
\end{enumerate}

If $F \colon {\cat C}_1 \ra {\cat C}_2$ is lax or op-lax, we will also say that $F$ is
\begin{enumerate}
\item
{\it pseudo-unital} if all its unit 2-cells $\Delta_a \colon 1_{Fa} \Ra F_{1a}$ (or $\Delta^{'}_a$) are {\it isomorphisms}
\item
{\it strictly unital}, if all its unit 2-cells $\Delta_a \colon 1_{Fa} \Ra F_{1a}$ (or $\Delta^{'}_a$) are {\it identities}
\end{enumerate} 

A pseudo-functor is pseudo-unital, and a strict 2-functor is strictly unital\footnote{In fact, any lax pseudo-unital (or op-lax pseudo-unital, or pseudo-) functor can be ``rectified'' to a strictly unital lax (resp. strictly unital op-lax, strictly unital pseudo-) functor. We are purposefully vague on this point, to avoid having to talk about natural transformations of 2-functors.}.

A 2-subcategory ${\cat C}^{'}$ of a 2-category ${\cat C}$ consists of a subclass of objects $Ob {\cat C}^{'} \subset Ob{\cat C}$ along with subcategories ${\cat C}^{'}(A, B) \subset {\cat C}(A, B)$ for objects $A, B \eps Ob {\cat C}^{'}$ that are stable under the composition rule $c_{ABC}$ for $A, B, C \eps Ob{\cat C}^{'}$ and include the image of the unit $1_A \colon * \ra {\cat C}(A, A)$ for $A \eps Ob{\cat C}^{'}$. A 2-subcategory is 2-{\it full} if ${\cat C}^{'}(A, B) = {\cat C}(A, B)$ for any $A, B \eps Ob {\cat C}^{'}$.

The category of categories $CAT$ forms a 2-category, denoted $2CAT$, with categories as 0-cells, functors as 1-cells and natural maps as 2-cells. We will use the notation introduced in this section to denote compositions of functors and natural maps. In particular, natural maps will be denoted with a double arrow $\Ra$ (or simply with $\ra$ if no confusion is possible), and composition of natural maps viewed as adjacent 2-cells will be denoted with $\star$. 

The category of small categories forms a 2-subcategory of $2CAT$, denoted $2{\cat Cat}$.

Any category ${\cat C}$ is in a canonical way a 2-category, having only identity 2-cells. In particular, we can talk of lax, op-lax, pseudo and strict 2-functors ${\cat C} \ra 2CAT$. A strict 2-functor ${\cat C} \ra 2CAT$ is the same as a functor ${\cat C} \ra CAT$.

For the rest of the chapter, we limit ourselves to the case of $2CAT$. As an exercise, the reader may enjoy reformulating the definitions and the results that follow so that they make sense in an arbitrary 2-category (cf. \cite{Kelly-Street1}, \cite{Street1}).

\section{Adjoint functors}
\label{sec:adjointfunctors}
\index{functor!adjoint}
Recall that an adjunction $u_1 \dashv u_2$ between two functors $u_1 : {\cat A} \rightleftarrows {\cat B} : u_2$ is a bijection of sets
\begin{center}
$\zeta \colon {\cat B}(u_1A, B) \cong {\cat A}(A, u_2B)$ 
\end{center}
natural in objects $A \eps {\cat A}$ and $B \eps {\cat B}$. For example, if $u_1$, $u_2$ are equivalences of categories, then $u_1 \dashv u_2$ (and $u_2 \dashv u_1$).

Whenever we say that $u_1, u_2$ is an adjoint pair we refer to a particular bijection $\zeta$. The following proposition encodes the basic properties of adjoint functors that we will need.

\begin{prop}
\label{prop:adjfunctors}
Suppose that $u_1 : {\cat A} \rightleftarrows {\cat B} : u_2$ is a pair of functors.
\begin{enumerate}
\item
An adjunction $u_1 \dashv u_2$ is uniquely determined by natural maps $\phi \colon 1_{\cat A} \Ra u_2u_1$ (the {\it unit}) and $\psi \colon u_1u_2 \Ra 1_{\cat B}$ (the {\it counit} of the adjunction) with the property that both the following compositions are identities
\begin{center}
$\xymatrix {
    {u_1} \ar@{=>}[r]^-{{u_1} {\phi}} & {u_1}{u_2}{u_1} \ar@{=>}[r]^-{{\psi} {u_1}} & {u_1} & {u_2} \ar@{=>}[r]^-{{\phi} {u_2}} & {u_2}{u_1}{u_2} \ar@{=>}[r]^-{{u_2} {\psi}} & {u_2} 
  }
$
\end{center}
\item
If $u_1 \dashv u_2$ is an adjunction with unit $\phi$ and counit $\psi$, then
\begin{enumerate}
\item
$u_1$ (resp. $u_2$) is fully faithful iff $\phi$ (resp. $\psi$) is a natural isomorphism.
\item
$u_1$ and $u_2$ are inverse equivalences of categories iff both $\phi$ and $\psi$ are natural isomorphisms.
\end{enumerate}
\end{enumerate}
\end{prop}

\begin{proof}
See for example Mac Lane \cite{MacLane}.
\end{proof}

Here is another way to state part (2) of the previous proposition. The proof is left to the reader.
\begin{prop}
\label{prop:adjfunctors2}
Suppose that $u_1 \dashv u_2$ is an adjoint pair of functors as above. Then the following statements are equivalent:
\begin{enumerate}
\item[(1)]
(resp. (1r), resp. (1l)). 
For any objects $A \eps {\cat A}$, $B \eps {\cat B}$, a map $u_1A \ra B$ is an isomorphism iff (resp. if, resp. only if) its adjoint $A \ra u_2B$ is an isomorphism
\item[(2)]
(resp. (2r), resp. (2l)). 
$u_1$ and $u_2$ are inverse equivalences of categories (resp. $u_2$ is fully faithful, resp $u_1$ is fully faithful).
\end{enumerate}
\end{prop}

\section{Kan extensions}
\label{sec:leftrightkanextension}

\index{functor!derived}
\index{functor!kan extension}
\begin{defn}
\label{defn::def:kanext}
Consider two functors $u \colon {\cat A} \ra {\cat B}$ and $\gamma \colon {\cat A} \ra {\cat A}^{'}$
\begin{enumerate}
\item
A left Kan extension of $u$ along $\gamma$ is a pair $({\bf L}_\gamma u, \eps)$ where ${\bf L}_\gamma u \colon {\cat A}^{'} \ra {\cat B}$ is a functor and $\eps \colon {\bf L}_\gamma u \comp \gamma \Ra u$ is a natural map
\begin{center}
$\xymatrix {
    {\cat A} \ar[rrd]^u \ar[d]_{\gamma} &  \\ 
    {\cat A}^{'} \ar[rr]_{{\bf L}_\gamma u} \ar@{}[ru]|(.42){\;\; \eps \; \Uparrow } && {\cat B}
  }$
\end{center}
satisfying the following universal property: if $(\Lambda, \lambda)$ is another pair of a functor $\Lambda \colon {\cat A}^{'} \ra {\cat B}$ and natural map $\lambda \colon \Lambda \comp \gamma \Ra u$
\begin{center}
$\xymatrix {
    {\cat A} \ar@<1ex>[rrd]^u \ar[d]_{\gamma} &&  \\ 
    {\cat A}^{'} \ar@/^1pc/[rr]_(.6){{\bf L}_\gamma u}_(.4)\,="d" \ar@/_1pc/[rr]_(.6){\Lambda} ^(.4)\,="c" \ar@{}[ru]|(.45){\;\; \eps \; \Uparrow } && {\cat B}
    \ar@{} "c";"d" ^{\delta \; \Uparrow}
  }$
\end{center}
then there exists a unique natural map $\delta \colon \Lambda \Ra {\bf L}_\gamma u$ with $\eps \star \delta = \lambda$
\item
A right Kan extension of $u$ along $\gamma$ is a pair $({\bf R}_\gamma u, \nu)$ where ${\bf R}_\gamma u \colon {\cat A}^{'} \ra {\cat B}$ is a functor and $\nu \colon u \Ra {\bf R}_\gamma u \comp \gamma$ is a natural map
\begin{center}
$\xymatrix {
    {\cat A} \ar[rrd]^u \ar[d]_{\gamma} &&  \\ 
    {\cat A}^{'} \ar[rr]_{{\bf R}_\gamma u} \ar@{}[ru]|(.42){\;\; \nu \; \Downarrow } & & {\cat B}
  }$
\end{center}
satisfying the following universal property: if $(\Lambda, \lambda)$ is another pair of a functor $\Lambda \colon {\cat A}^{'} \ra {\cat B}$ and a natural map $\lambda \colon u \Ra \Lambda \comp \gamma $
\begin{center}
$\xymatrix {
    {\cat A} \ar@<1ex>[rrd]^u \ar[d]_{\gamma} &&  \\ 
    {\cat A}^{'} \ar@/^1pc/[rr]_(.6){{\bf R}_\gamma u}_(.4)\,="d" \ar@/_1pc/[rr]_(.6){\Lambda} ^(.4)\,="c" \ar@{}[ru]|(.45){\;\; \nu \; \Downarrow } && {\cat B}
    \ar@{} "c";"d" ^{\delta \; \Downarrow}
  }$
\end{center}
then there exists a unique natural map $\delta \colon {\bf R}_\gamma u \Ra \Lambda$ with $\delta \star \nu = \lambda$
\end{enumerate}
\end{defn}

The left Kan extension $({\bf L}_\gamma u, \eps)$ is also called in some references the {\it left derived} of $u$ along $\gamma$ and the right Kan extension $({\bf R}_\gamma u, \eps)$ the {\it right derived} of $u$ along $\gamma$. Since it is defined by an universal property, if the left (or right) Kan extension exists then it is unique up to a unique isomorphism. 

The next proposition is an existence criterion for Kan extensions: if $\gamma$ admits a left (resp. right) adjoint, then $u$ admits a left (resp. right) Kan extension along $\gamma$.

\begin{prop}
\label{prop:existencekanextension}
Consider two functors $u \colon {\cat A} \ra {\cat B}$ and $\gamma \colon {\cat A} \ra {\cat A}^{'}$.
\begin{enumerate}
\item
If $\gamma$ admits a left adjoint $\gamma^{'}$ with adjunction unit $\phi \colon 1_{{\cat A}^{'}} \Ra \gamma\gamma^{'}$ and counit $\psi \colon \gamma^{'}\gamma \Ra 1_{{\cat A}}$, then $(u\gamma^{'}, u \psi)$ is a left Kan extension of $u$ along $\gamma$
\item
If $\gamma$ admits a right adjoint $\gamma^{'}$ with adjunction unit $\phi \colon 1_{\cat A} \Ra \gamma^{'}\gamma$ and counit $\psi \colon \gamma\gamma^{'} \Ra 1_{{\cat A}^{'}}$, then $(u\gamma^{'}, u \phi)$ is a right Kan extension of $u$ along $\gamma$
\end{enumerate}
\end{prop}

\begin{proof}
We only prove (1). 
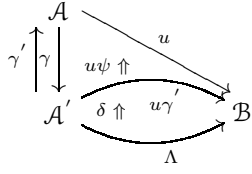
\begin{center}
$\xymatrix {
    {\cat A} \ar@<1ex>[rrd]^u \ar[d]_{\gamma} &&  \\ 
    {\cat A}^{'} \ar@/^1pc/[rr]_(.6){u \gamma^{'}}_(.4)\,="d" \ar@/_1pc/[rr]_(.6){\Lambda} ^(.4)\,="c" \ar@{}[ru]|(.45){\;\; u\psi \; \Uparrow } \ar@<2ex>[u]^(.55){\gamma^{'}} && {\cat B}
    \ar@{} "c";"d" ^{\delta \; \Uparrow}
  }$
\end{center}

For any pair $(\Lambda, \lambda)$ with $\Lambda \colon {\cat A}^{'} \ra {\cat B}$ and $\lambda \colon \Lambda \gamma \Ra u$, we'd like to show that there exists a unique natural map $\delta \colon \Lambda \Ra u \gamma^{'}$ with
\begin{equation}
\label{eqn:existencekanextension1}
\lambda = (u \psi) \comp (\delta \gamma)
\end{equation}
To show existence, we define $\delta \colon \Lambda \Ra \Lambda \gamma \gamma^{'} \Ra u \gamma^{'}$ as the composition 
\begin{equation}
\label{eqn:existencekanextension2}
\delta = (\lambda \gamma^{'}) \comp (\Lambda \phi)
\end{equation}
$\delta$ defined by (\ref{eqn:existencekanextension2}) satisfies (\ref{eqn:existencekanextension1}), using the commutative diagram
\begin{center}
$\xymatrix{
    \Lambda \gamma \ar@{=>}[r]^-{\Lambda \phi \gamma} \ar@{=>}[rd]_{1_{\Lambda \gamma}} & \Lambda \gamma \gamma^{'} \gamma \ar@{=>}[r]^{\lambda \gamma^{'} \gamma} \ar@{=>}[d]_{\Lambda \gamma \psi} & u \gamma^{'} \gamma \ar@{=>}[d]^{u\psi} \\
    & \Lambda \gamma \ar@{=>}[r]^{\lambda} & u
    }$
\end{center} 
To show uniqueness, if a map $\delta$ satisfies (\ref{eqn:existencekanextension1}) then $\delta$ satisfies (\ref{eqn:existencekanextension2}) because of the commutative diagram
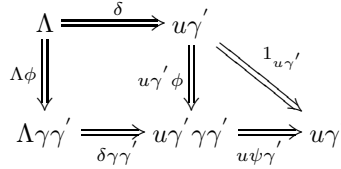
\begin{center}
$\xymatrix{
    \Lambda  \ar@{=>}[d]_-{\Lambda \phi} \ar@{=>}[r]^{\delta} & u \gamma^{'} \ar@{=>}[d]_{u \gamma^{'} \phi} \ar@{=>}[rd]^{1_{u \gamma^{'}}} &  \\
     \Lambda \gamma \gamma^{'} \ar@{=>}[r]_{\delta \gamma \gamma^{'}} & u \gamma^{'} \gamma \gamma^{'} \ar@{=>}[r]_{ u \psi \gamma^{'}} & u \gamma^{'}
    }$
\end{center} 
\end{proof}

Kan extensions commute with composition {\it along the base functor} $\gamma$, as seen in the next propoition.

\begin{prop}
\label{prop:basecompositionkanextension}
Consider the functors ${\cat A} \overset{u}{\ra} {\cat B}$, ${\cat A} \overset{\gamma}{\ra} {\cat A}^{'} \overset{\gamma^{'}}{\ra} {\cat A}^{''}$.
\begin{enumerate}
\item
Suppose that $({\bf L}_\gamma u, \eps)$ exists. Then $({\bf L}_{\gamma^{'}} {\bf L}_\gamma u, \eps^{'})$ exists iff $({\bf L}_{\gamma^{'}\gamma} u, \eps^{''})$ exists. If they both exist, then the latter is isomorphic to $({\bf L}_{\gamma^{'}}{\bf L}_\gamma u, \eps \star \eps^{'})$.
\item
Suppose that $({\bf R}_\gamma u, \nu)$ exists. Then $({\bf R}_{\gamma^{'}} {\bf R}_\gamma u, \nu^{'})$ exists iff $({\bf R}_{\gamma^{'}\gamma} u, \nu^{''})$ exists, in which case the latter is isomorphic to $({\bf R}_{\gamma^{'}}{\bf R}_\gamma u, \nu^{'} \star \nu)$.
\end{enumerate}
\end{prop}

\begin{proof}
Immediate using the universal property of the left (resp. right) Kan extensions. 
\end{proof}

We state a corollary needed in the proof of \theoremref{thm:genholim}.
\begin{cor}
\label{cor:neededbygenhomlim}
Consider two functors $u \colon {\cat A} \ra {\cat B}$ and $\gamma \colon {\cat A} \ra {\cat A}^{'}$, and a pair of inverse equivalences of categories $\gamma^{'} : {\cat A}^{'} \leftrightarrows {\cat A}^{''} : \gamma^{''}$, with natural isomorphisms $\phi \colon 1_{{\cat A}^{'}} \Ra \gamma^{''}\gamma^{'}$ and $\psi \colon 1_{{\cat A}^{''}} \Ra \gamma^{'}\gamma^{''}$.
\begin{enumerate}
\item
The left Kan extension $({\bf L}_\gamma u, \eps)$ exists iff the left Kan extension  $({\bf L}_{{\gamma}^{'}\gamma} u, \eps^{'})$ exists. If they both exist, then the latter is isomorphic to $(({\bf L}_\gamma u) \gamma^{''}$, $\eps \star \psi)$.
\item
The right Kan extension $({\bf R}_\gamma u, \nu)$ exists iff the right Kan extension $({\bf R}_{{\gamma}^{'}\gamma} u$, $\nu^{'})$ exists, in which case the latter is isomorphic to $(({\bf R}_\gamma u) \gamma^{''}$, $\phi \star \nu)$.
\end{enumerate}
\end{cor}

\begin{proof}
Consequence of \propositionref{prop:existencekanextension} and \propositionref{prop:basecompositionkanextension}.
\end{proof}

\chapter{Categories with weak equivalences}
\label{chap:catweq}

ABC cofibration categories can be essentially thought of as 'nicely behaved' categories with weak equivalences. But how much can we say about categories with weak equivalences without bringing in the cofibrations? We will try to find an answer in this chapter.

We will denote ${\bf ho}{\cat M}$ the homotopy category of a category ${\cat M}$ with a class of weak equivalences ${\cat W}$. The homotopy category ${\bf ho}{\cat M}$ is defined by a universal property, but admits a description in terms of generators and relations, starting with the objects and maps of ${\cat M}$ and formally adding inverses of the maps in ${\cat W}$.

The total left (resp. right) derived of a functor $u \colon {\cat M}_1 \ra {\cat M}_2$ between two category pairs $({\cat M}_i, {\cat W}_i)$ for $i = 1, 2$ with localization functors $\gamma_i \colon {\cat M}_i \ra {\bf ho}{\cat M}_i$ is defined as the left (resp. right) Kan extension of $\gamma_2 u$ along $\gamma_1$.

We will use left, resp. right {\it approximation functors} (\definitionref{defn:wleftapproximation}) as a tool for an existence theorem for total derived functors (\theoremref{thm:generalexistenceleftderived}), and a rather technical adjointness property of total derived functors (\theoremref{thm:generalexistencetotalderivedadjoint2}).

The left, resp. right approximation functors $f \colon ({\cat M}^{'}, {\cat W}^{'}) \ra ({\cat M}, {\cat W})$ among other things have the property that they induce an equivalence ${\bf ho}t \colon {\bf ho}{\cat M}^{'} \ra {\bf ho}{\cat M}$. This is proved by the Approximation \theoremref{thm:approxthm}. 

\section{Universes and smallness}
If ${\cat U} \subset {\cat U}^{'}$ are two universes, a $({\cat U}^{'}, {\cat U})$-category ${\cat C}$ has by definition a ${\cat U}^{'}$-small set of objects $Ob{\cat C}$ and ${\cat U}$-small Hom-sets. A category ${\cat C}$ is ${\cat U}$-small if it is a $({\cat U}, {\cat U})$-category, and is {\it locally} ${\cat U}$-small if it has ${\cat U}$-small Hom-sets.

For a fixed universe pair ${\cat U} \subset {\cat U}^{'}$, the ${\cat U}^{'}$-small sets are also referred to as {\it classes}, while the ${\cat U}$-small sets are referred to simply as {\it sets}, or {\it small sets}. We will denote ${\cat Cat}$, resp. $CAT$ the category of ${\cat U}$-small categories, resp. $({\cat U}^{'}, {\cat U})$ categories, with functors as 1-cells. Both ${\cat Cat}$ and $CAT$ actually carry 2-category structures, with natural transformations as 2-cells.

\section{The homotopy category}
\label{sec:whomotopycat}

\begin{defn}
\label{defn:hcat}
Suppose we have a $({\cat U}^{'}, {\cat U})$-category ${\cat M}$ with a class of weak equivalences ${\cat W}$. The homotopy category ${\bf ho}{\cat M}$ by definition is a ${\cat U}^{'}$-small category equipped with a localization functor $\gamma_{\cat M} \colon {\cat M}$ $\ra$ ${\bf ho}{\cat M}$, with the properties that 
\begin{enumerate}
\item
$\gamma_{\cat M}$ sends weak equivalences to isomorphisms
\item
\label{defn:hcat:2}
for any other such functor $\gamma^{'} \colon {\cat M}$ $\ra$ ${\cat M}^{'}$ that sends weak equivalences to isomorphisms, there exists a unique functor $\delta \colon {\bf ho}{\cat M}$ $\ra$ ${\cat M}^{'}$ such that $\delta \gamma_{\cat M} = \gamma^{'}$
\end{enumerate}
\end{defn}

From its universal property, the homotopy category is uniquely defined up to an {\it isomorphism} of categories.

The homotopy category always exists, and is constructed in the next theorem.

\begin{thm}[Gabriel-Zisman]
\label{thm:gz}
For a category with a class of weak equivalences $({\cat M}, {\cat W})$, the category 
\begin{enumerate}
\item
With the same objects as ${\cat M}$ 
\item
With maps between $X$ and $Y$ the equivalence classes of zig-zags 
\[
\xymatrix{ X = A_0 \ar@{-}[r]^-{f_1} &  A_1 \ar@{-}[r]^-{f_2} & ... \ar@{-}[r]^-{f_n} & A_n = Y}
\]
where $f_i$ are maps in ${\cat M}$ going either forward or backward, and all the maps going backward are in ${\cat W}$, where
\begin{enumerate}
\item
Two zig-zags are equivalent if they can be obtained one from another by a finite number of the following three operations and their inverses:
\begin{enumerate}
\item
Skipping elements $A \overset{1_{A}}{\ra} A$ or $A \overset{1_{A}}{\la} A$
\item
Replacing $A \overset{f}{\ra} B \overset{g}{\ra} C$ with $A \overset{gf}{\ra} C$ or $A \overset{f}{\la} B \overset{g}{\la} C$ with $A \overset{fg}{\la} C$
\item
Skipping elements $A \overset{f}{\ra} B \overset{f}{\la} A$ or $A \overset{f}{\la} B \overset{f}{\ra} A$
\end{enumerate}
\item
Composition of maps is induced by the concatenation of zig-zags
\end{enumerate}
\end{enumerate}
is a homotopy category with a localization functor that preserves the objects, and sends a map $X \overset{f}{\ra}Y$ to itself viewed as a zig-zag.
\end{thm}

\begin{proof}
See \cite{Gabriel-Zisman}.
\end{proof}

The homotopy category is sometimes also denoted ${\cat M}[{\cat W}^{-1}]$. 

\index{${\cat W}, \overline{\cat W}$}
\index{weak equivalences ${\cat W}$}
\index{weak equivalences ${\cat W}$!saturation $\overline{\cat W}$}
The saturation $\overline{\cat W}$ of ${\cat W}$ by definition is the class of maps of ${\cat M}$ that become isomorphisms in ${\bf ho}{\cat M}$. The saturation $\overline{\cat W}$ is closed under composition, includes the isomorphisms of ${\cat M}$, and ${\cat M}[\overline{W}^{-1}]$ is isomorphic to ${\cat M}[{\cat W}^{-1}]$.

If ${\cat M}$ is a locally ${\cat U}$-small category with a class of weak equivalences ${\cat W}$, then ${\bf ho}{\cat M} = {\cat M}[{\cat W}^{-1}]$ is not necessarily locally ${\cat U}$-small. For example, if ${\cat M}$ is the $({\cat U}^{'}, {\cat U})$-category of sets ${\cat Set}$, then ${\cat W} = {\cat Cof} = {\cat Fib} = {\cat Set}$  gives ${\cat M}$ a structure of ABC model category, but ${\cat Set}[{\cat Set}^{-1}]$ is only locally ${\cat U}^{'}$-small. 

On the other hand if a locally ${\cat U}$-small category ${\cat M}$ satisfies the hypothesis of \corollaryref{cor:smallhomotopycat} (for example if it is a Quillen model category), then ${\bf ho}{\cat M}$ is again locally ${\cat U}$-small.

One benefit of the intrinsic description of the homotopy category in terms of zig-zags of maps is that it shows that the homotopy category remains the same {\it independent of the universe pair} $({\cat U}^{'}, {\cat U})$ that we start with.

\section{Homotopic maps}
\label{sec:generalhomotopicmaps}

\begin{defn}
\label{defn:homotopicmaps}
For a category ${\cat M}$ with a class of weak equivalences ${\cat W}$, say that two maps $f, g \colon A \ra B$ in ${\cat M}$ are {\it homotopic} (write $f \simeq g$) if $\gamma_{\cat M} f = \gamma_{\cat M} g$.
\end{defn}

The homotopy relation $\simeq$ is an equivalence relation on $Hom_{\cat M}(A, B)$, and the quotient set $Hom_{\cat M}(A, B) / _\simeq$ is classically denoted with the bracket notation $[A, B]$. 

Suppose that $f \simeq g \colon A \ra B$. If $u \colon A^{'} \ra A$, then $fu \simeq gu$. If $u \colon B \ra B^{'}$, then $uf \simeq ug$. 

The quotient category ${\cat M}/_\simeq$, also denoted $\pi {\cat M}$, can therefore be constructed with the same objects as ${\cat M}$ and with $Hom_{\pi {\cat M}} (A, B) = [A, B]$. The induced functor $\pi {\cat M} \ra {\bf ho}{\cat M}$ is an inclusion, and the localization functor $\gamma_{\cat M}$ then factors as the composition ${\cat M} \ra \pi {\cat M} \ra {\bf ho}{\cat M}$.

\section{Left and right approximation functors}
\label{sec:approximationfunctors}

For the rest of the text, it is convenient to restrict ourselves to the case when ${\cat W}$ is closed under composition and includes the identity maps of ${\cat M}$. In the language of \definitionref{defn:catweakequiv}, these are the {\it category pairs} $({\cat M}$, ${\cat W})$. 

A category pair {\it functor} $t \colon ({\cat M}^{'}$, ${\cat W}^{'})$ $\ra$ $({\cat M}$, ${\cat W})$ by definition will be a functor $t \colon {\cat M}^{'}$ $\ra$ ${\cat M}$ that preserves the weak equivalences. A category pair functor $t$ as above induces a functor at the level of the homotopy categories, denoted ${\bf ho}t \colon {\bf ho}{\cat M}^{'}$ $\ra$ ${\bf ho}{\cat M}$.

Our next goal is to find sufficient conditions for a category pair functor to induce an {\it equivalence} at the level of homotopy categories.

\index{approximation!left, right}
\begin{defn}
\label{defn:wleftapproximation}
A category pair functor $t \colon ({\cat M}^{'}$, ${\cat W}^{'})$ $\ra$ $({\cat M}$, ${\cat W})$ is a {\it left approximation} if
\begin{description}
\item[LAP1]
For any object $A \eps {\cat M}$ there exists an object $A^{'} \eps {\cat M}^{'}$ and a weak equivalence map $tA^{'} \overset{\sim}{\ra} A$
\item[LAP2]
For any maps $f \eps {\cat M}$ and $p \eps {\cat W}$ as below
\[
\xymatrix@R=15pt@C=15pt { 
&&&& B \\
&&&& \\
tA^{'} \ar[rrrruu]^{f} \ar@{-->}[rr]_{ti^{'}} && tD^{'} \ar@{-->}[rruu]^{\sigma} && tB^{'} \ar@{-->}[ll]^{tj^{'}}_\sim \ar[uu]_{p}^\sim
}
\]
\begin{enumerate}
\item
There exist maps $A^{'} \overset{i^{'}}{\lra} D^{'} \underset{\sim}{\overset{j^{'}}{\lla}} B^{'}$ with $j^{'} \eps {\cat W}^{'}$ and a weak equivalence $\sigma \eps {\cat W}$ making the diagram commutative.
\item
For any other such $i^{''}, j^{''}, \sigma^{''}$ the equality $j^{'-1}i^{'} = j^{''-1}i^{''}$ holds in ${\bf ho}{\cat M}^{'}$.
\end{enumerate}
\end{description}

Dually, the functor $t$ is a {\it right approximation} if
\begin{description}
\item[RAP1]
For any object $A \eps {\cat M}$ there exists an object $A^{'} \eps {\cat M}^{'}$ and a weak equivalence map $A \overset{\sim}{\ra} tA^{'}$
\item[RAP2]
For any maps $f \eps {\cat M}$ and $i \eps {\cat W}$ as below
\[
\xymatrix@R=15pt@C=15pt { 
A \ar[dd]_{i}^\sim \ar@{-->}[rrdd]^{\tau} \ar[rrrrdd]^{f}  &&&& \\
&&&& \\
tA^{'}  && tS^{'} \ar@{-->}[ll]^{tq^{'}}_\sim \ar@{-->}[rr]_{tp^{'}} && tB^{'}  
}
\]
\begin{enumerate}
\item
There exist maps $A^{'} \underset{\sim}{\overset{q^{'}}{\lla}} S^{'} \overset{p^{'}}{\lra} B^{'}$ with $q^{'} \eps {\cat W}^{'}$ and a weak equivalence $\tau \eps {\cat W}$ making the diagram commutative.
\item
For any other such $p^{''}, q^{''}, \tau^{''}$ the equality $p^{'}q^{'-1} = p^{''}q^{''-1}$ holds in ${\bf ho}{\cat M}^{'}$.
\end{enumerate}
\end{description}
\end{defn}


\begin{lem}
\label{lem:sufficientcondequiv}
Suppose that $t \colon ({\cat M}^{'}, {\cat W}^{'}) \ra ({\cat M}, {\cat W})$ is a category pair functor, and $f \colon A^{'} \ra C^{'}$ is a map in ${\cat M}^{'}$. Suppose that either
\begin{enumerate}
\item
$t$ is a left approximation, and there exist $p, q \eps {\cat W}$ making commutative the diagram
\[
\xymatrix{
  & B & \\
  tA^{'} \ar[ru]^p_\sim \ar[rr]_{tf^{'}} && tC^{'} \ar[lu]^\sim_q
}
\]
\item
$t$ is a right approximation, and there exist $i, j \eps {\cat W}$ making commutative the diagram
\[
\xymatrix{
  & B \ar[ld]_i^\sim \ar[rd]^j_\sim & \\
  tA^{'} \ar[rr]_{tf^{'}} && tC^{'} 
}
\]
\end{enumerate}
Then the image of $f^{'}$ in ${\bf ho}{\cat M}^{'}$ is an isomorphism, i.e. $f \eps \overline{\cat W}^{'}$.
\end{lem}

\begin{proof}
We only prove (1). Using LAP2 (1), we can construct $A^{'} \underset{\sim}{\overset{i^{'}}{\lra}} B^{'} \overset{j^{'}}{\lla} C^{'}$ with $i^{'} \eps {\cat W}^{'}$ and a weak equivalence $\sigma \eps {\cat W}$ with $\sigma \comp ti^{'} = p$ and $\sigma \comp tj^{'} = q$. 
\[
\xymatrix {
& B & \\ 
tA^{'} \ar[rr]_(.3){tf^{'}}|{\;\;\;} \ar[ru]^p \ar@{-->}[rd]_{ti^{'}} && tC^{'} \ar@{-->}[ld]^{tj^{'}} \ar[lu]_q \\
& tB^{'} \ar@{-->}[uu]^(.7){\sigma} & 
} 
\]
From LAP2 (2), $j^{'}f^{'}$ is invertible in ${\bf ho}{\cat M}^{'}$. Similarly we can construct $j^{''} \eps {\cat M}^{'}$ such that $j^{''}j^{'}$ is invertible in ${\bf ho}{\cat M}^{'}$. We deduce that $j^{'}$ and therefore $f^{'}$ are invertible in ${\bf ho}{\cat M}^{'}$.
\end{proof}

\index{cleavage!left (right) cleavage}
\begin{defn}
\label{defn:wleftcleavage}
Let $t \colon {\cat M}^{'} \ra {\cat M}$ be a left approximation of category pairs $({\cat M}^{'}$, ${\cat W}^{'})$ and $({\cat M}$, ${\cat W})$. A {\it left cleavage} ${\cat C} = (C, D, p, i, j, \sigma)$ along $t$ consists of the following data:
\begin{enumerate}
\item
For any object $A$ of ${\cat M}$, an object $C(A)$ of ${\cat M}^{'}$ and a weak equivalence $p_A \colon tC(A) \ra A$
\item
For any map $f \colon A \ra B$, maps $C(A) \overset{i(f)}{\lra} D(f) \underset{\sim}{\overset{j(f)}{\lla}} C(B)$ and a commutative diagram
\[
\xymatrix{
  A \ar[rr]^{f} && B \\
  tC(A) \ar[u]^{p_A}_\sim \ar[r]_{ti(f)} & tD(f) \ar[ur]^{\sigma(f)} & tC(B) \ar[l]^{tj(f)}_\sim \ar[u]^\sim_{p_B}
}
\]
with $j(f)$, $\sigma(f)$ weak equivalences.
\end{enumerate}

A {\it left cleavage map} ${\cat C}_1 \overset{\cat F}{\ra} {\cat C}_2$ from ${\cat C}_1$ $=$ ($C_1$, $D_1$, $p_1$, $i_1$, $j_1$, $\sigma_1$) to ${\cat C}_2$ $=$ ($C_2$, $D_2$, $p_2$, $i_2$, $j_2$, $\sigma_2$) consists of a family of maps $c(A) \colon C_1(A) \ra C_2(A)$ and $d(f) \colon D_1(f) \ra D_2(f)$, satisfying $p_{1A} = p_{2A} \comp tc(A)$, $d(f)i_1(f) = i_2(f)c(A)$ and $d(f)j_1(f) = j_2(f)c(B)$. 

Notice that by \lemmaref{lem:sufficientcondequiv}, we must have $c(A), d(f) \eps \overline{\cat W}^{'}$.

Composition of left cleavage maps is given by componentwise composition of $c(A), d(f)$. Left cleavages and left cleavage maps along $t$ form therefore a category.
\end{defn}

The dual definition for right approximations is

\begin{defn}
\label{defn:wrightcleavage}
Let $t \colon {\cat M}^{'} \ra {\cat M}$ be a right approximation of category pairs $({\cat M}^{'}$, ${\cat W}^{'})$ and $({\cat M}$, ${\cat W})$. A {\it right cleavage} ${\cat R} = (R, S, i, p^{'}, q^{'}, \tau)$ along $t$ consists of the following data:
\begin{enumerate}
\item
For any object $A$ of ${\cat M}$, an object $R(A)$ of ${\cat M}^{'}$ and a weak equivalence $i_A \colon A \ra tR(A)$
\item
For any map $f \colon A \ra B$, maps $R(A) \underset{\sim}{\overset{q(f)}{\lla}} S(f) \overset{p(f)}{\lra} R(B)$ and a commutative diagram
\[
\xymatrix{
  A \ar[rr]^{f} \ar[d]_{i_A}^\sim \ar[dr]^{\tau(f)} && B \ar[d]_\sim^{i_B}   \\
  tR(A)  & tS(f)  \ar[l]^{tq(f)}_\sim \ar[r]_{tp(f)} & tR(B) 
}
\]
with $q(f)$, $\tau(f)$ weak equivalences.
\end{enumerate}

A {\it right cleavage map} ${\cat R}_1 \overset{\cat F}{\ra} {\cat R}_2$ from ${\cat R}_1$ $=$ ($R_1$, $S_1$, $i_1$, $p^{'}_1$, $q^{'}_1$, $\tau_1$) to ${\cat R}_2$ = ($R_2$, $S_2$, $i_2$, $p^{'}_2$, $q^{'}_2$, $\tau_2$) consists of a family of maps $r(A) \colon R_1(A) \ra R_2(A)$ and $s(f) \colon S_1(f) \ra S_2(f)$, satisfying $r(A) \comp i_{1A} = i_{2A}$, $r(B)p_1(f) = p_2(f)s(f)$ and $r(A)q_1(f) = q_2(f)s(f)$. 
\end{defn}

\begin{lem}
\label{lem:cleavagemap}
Suppose that $t \colon ({\cat M}^{'}, {\cat W}^{'}) \ra ({\cat M}, {\cat W})$ is a category pair functor.
\mbox{}
\begin{enumerate}
\item
If $t$ is a left approximation and ${\cat C}_1, {\cat C}_2$ are left cleavages, there exist a left cleavage ${\cat C}_3$ and left cleavage maps ${\cat C}_1 \ra {\cat C}_3 \la {\cat C}_2$.
\item
If $t$ is a right approximation and ${\cat R}_1, {\cat R}_2$ are right cleavages, there exist a right cleavage ${\cat R}_3$ and right cleavage maps ${\cat R}_1 \la {\cat R}_3 \ra {\cat R}_2$.
\end{enumerate}
\end{lem}

\begin{proof}
Apply repeatedly the axioms LAP2 (1), resp. RAP2 (1).
\end{proof}

\index{cleavage!normalized}
\index{cleavage!regular object}
\index{cleavage!regular}
\begin{defn}
\label{defn:wregularcleavageobj}
Suppose that $({\cat M}^{'}$, ${\cat W}^{'})$ and $({\cat M}$, ${\cat W})$ are two category pairs.
\begin{enumerate}
\item
Let $t \colon {\cat M}^{'} \ra {\cat M}$ be a left approximation functor, and ${\cat C} = (C, D, p, i, j, \sigma)$ be a left cleavage. 
\begin{enumerate}
\item
The left cleavage ${\cat C}$ is {\it normalized} if for any object $A$ of ${\cat M}$ we have $D(1_A) = C(A)$, $i(1_A) = j(1_A) = 1_{C(A)}$ and $\sigma(1_A) = p_A$.
\item
An object $A^{'} \eps {\cat M}^{'}$ is {\it regular} with respect to ${\cat C}$ (referred to simply as {\it regular} if no confusion is possible) if $C(tA^{'}) = A^{'}$ and $p_{tA^{'}} = 1_{tA^{'}}$.
\item
The left cleavage ${\cat C}$ is {\it regular} if it is normalized, and any object $A^{'} \eps {\cat M}^{'}$ admits a regular object $A^{''} \eps {\cat M}^{'}$ isomorphic to $A^{'}$ in ${\bf ho}{\cat M}^{'}$.
\end{enumerate}
\item
Let $t \colon {\cat M}^{'} \ra {\cat M}$ be a right approximation functor, and ${\cat R} = (R, S, i, p^{'}, q^{'}, \tau)$ be a right cleavage along $t$. 
\begin{enumerate}
\item
The right cleavage ${\cat R}$ is {\it normalized} if for any object $A$ of ${\cat M}$ we have $S(1_A) = R(A)$, $p(1_A) = q(1_A) = 1_{R(A)}$ and $\tau(1_A) = i_A$.
\item
An object $A^{'} \eps {\cat M}^{'}$ is {\it regular} with respect to a right cleavage along $t$ if $R(tA^{'}) = A^{'}$ and $i_{tA^{'}} = 1_{tA^{'}}$.
\item
The right cleavage ${\cat R}$ is {\it regular} if it is normalized, and any object $A^{'} \eps {\cat M}^{'}$ admits a regular object $A^{''} \eps {\cat M}^{'}$ isomorphic to $A^{'}$ in ${\bf ho}{\cat M}^{'}$.
\end{enumerate}
\end{enumerate}
\end{defn}

\begin{lem}
\label{lem:admitframing}
\mbox{}
\begin{enumerate}
\item
Any left approximation admits a regular left cleavage.
\item
Any right approximation admits a regular right cleavage. 
\end{enumerate}
\end{lem}

\begin{proof}
We only prove (1), and perform the construction in a number of steps.

(i) For any object $A \eps {\cat M}$ of the form $A = tA^{'}$, pick exactly one such $A^{'}$ and define $C(tA^{'}) = A^{'}$, $p_{tA^{'}} = 1_{tA^{'}}$.

(ii). For any $A \eps {\cat M}$ for which this data is not defined, use LAP1 to define the object $C(A)$ and the weak equivalence $p_A$.

(iii). For any $A \eps {\cat M}$, define $D(1_A) = C(A)$, $i(1_A) = j(1_A) = 1_{C(A)}$ and $\sigma(1_A) = p_A$.

(iv). For any map $f \colon A \ra B$ in ${\cat M}$ for which this data is not defined, use LAP2 (1) to construct the desired objects and maps $D(f), i(f) \eps {\cat M^{'}}$, $j(f) \eps {\cat W}^{'}$ and $\sigma(f) \eps {\cat W}$.

This constructs a left cleavage, which is normalised because of step (iii). 

Let us check condition (1) of \definitionref{defn:wregularcleavageobj}. For $A^{'} \eps {\cat M}^{'}$, we have $tC(tA^{'}) = tA^{'}$. The object $C(tA^{'})$ is regular, and by \lemmaref{lem:sufficientcondequiv} it is isomorphic to $A^{'}$ in ${\bf ho}{\cat M}$. In conclusion, we have constructed a regular left cleavage.
\end{proof}

If a left (resp. right) approximation $t \colon {\cat M}^{'} \ra {\cat M}$ is {\it injective on objects}, we could actually construct a normalized left (resp. right) cleavage where all objects of ${\cat M}^{'}$ were regular. 

\section{The approximation theorem}

Recall that a functor $u \colon {\cat M}_1 \ra {\cat M}_2$ is called
\index{functor!full}
\index{functor!faithful}
\index{functor!essentially surjective}
\begin{enumerate}
\item
{\it essentially surjective} if any object of ${\cat M}_2$ is isomorphic to an object in the image of $u$
\item
{\it full} if any map in $Hom_{{\cat M}_2}(uA, uB)$ is in the image of $u$, for all objects $A, B$ of ${\cat M}_1$
\item
{\it faithful} if $u$ is injective on $Hom_{{\cat M}_1}(A, B)$ for all objects $A, B$ of ${\cat M}_1$
\end{enumerate}

The functor $u$ is an equivalence of categories if and only if it is essentially surjective, full and faithful.

\begin{thm}[The approximation theorem]
\label{thm:approxthm}
\mbox{} \\
A left (or right) approximation functor $t \colon ({\cat M}^{'}$, ${\cat W}^{'})$ $\ra$ $({\cat M}$, ${\cat W})$ induces an equivalence of homotopy categories ${\bf ho} t \colon{\bf ho}{\cat M}^{'} \ra {\bf ho}{\cat M}$.
\end{thm}

\begin{proof}
Suppose that $t$ is a left approximation. For a left cleavage ${\cat C} = (C, D, p, i, j, \sigma)$ along $t$, construct a functor ${\bf s}_{\cat C} \colon {\bf ho}{\cat M}^{'} \ra {\bf ho}{\cat M}$, as follows:
\begin{enumerate}
\item
On objects $A \eps {\cat M}$, define ${\bf s}_{\cat C}A  = C(A)$
\item
For a zig-zag ${\cat Z}$ in ${\cat M}$
\begin{equation}
\label{eqn:approxarehtpyequiv}
\xymatrix{ A = A_0 \ar@{-}[r]^-{f_1} &  A_1 \ar@{-}[r]^-{f_2} & ... \ar@{-}[r]^-{f_n} & A_n = B}
\end{equation}
where the maps $f_k$ go either forward or backward, and all the maps going backward are in ${\cat W}$, denote $A^{'}_{k} = C(A_k)$ and $p_k = p_{A_k}$.
\begin{enumerate}
\item
If an $f_k$ goes forward, then construct the commutative diagram
\begin{equation}
\label{eqn:approxarehtpyequiv1}
\xymatrix{
  A_{k-1} \ar[rr]^{f_k} && A_k \\
  tA^{'}_{k-1} \ar[u]^{p_{k-1}}_\sim \ar[r]_{ti^{'}_k} & tD^{'}_k \ar[ur]^{\sigma_k} & tA^{'}_k \ar[l]^{tj^{'}_k}_\sim \ar[u]^\sim_{p_k}
}
\end{equation}
\item
If an $f_k$ goes backward, then construct the commutative diagram
\begin{equation}
\label{eqn:approxarehtpyequiv2}
\xymatrix{
  A_{k-1}  && A_k \ar[ll]^\sim_{f_k} \\
  tA^{'}_{k-1} \ar[u]^{p_{A_{k-1}}}_\sim \ar[r]_{ti^{'}_k} & tD^{'}_k \ar[ul]_{\sigma_k} & tA^{'}_k \ar[l]^{tj^{'}_k}_\sim \ar[u]^\sim_{p_k}
}
\end{equation}
\end{enumerate}
where in both cases we denoted $D^{'}_k = D(f_k)$, $i^{'}_k = i(f_k)$, $j^{'}_k = j(f_k)$ and $\sigma_k = \sigma(f_k)$.
\end{enumerate}

The maps $i^{'}_k, j^{'}_k$ collect together to a zig-zag denoted ${\cat Z}^{'}$ in ${\cat M}^{'}$, from $A^{'}_0$ to $A^{'}_n$, with all backward going maps being weak equivalences. We define ${\bf s}_{\cat C}$(image of ${\cat Z}$ in ${\bf ho}{\cat M}$) = image of ${\cat Z}^{'}$ in ${\bf ho}{\cat M}^{'}$.

Let's show that the definition of ${\bf s}_{\cat C}$ on maps does not depend on the choices involved. Let ${\cat Z}_1$, ${\cat Z}_2$ by two zig-zags of the form \equationref{eqn:approxarehtpyequiv}, and denote ${\cat Z}^{'}_1$, ${\cat Z}^{'}_2$ the associated zig-zags constructed in ${\cat M}^{'}$.

(i). If ${\cat Z}_2$ is obtained from ${\cat Z}_1$ by inserting an element $A_k \overset{1_{A_k}}{\lra} A_k$, then ${\cat Z}^{'}_2$ is obtained from ${\cat Z}^{'}_1$ by inserting an element 
\[
A_k^{'} \overset{i^{'}}{\lra} D^{'} \underset{\sim}{\overset{j^{'}}{\lla}} A^{'}_k
\]
with the property that there exists a commutative diagram
\[
\xymatrix{
  A_k \ar[rr]^{1_{A_k}} && A_k \\
  tA^{'}_k \ar[u]^{p_k}_\sim \ar[r]_{ti^{'}} & tD^{'} \ar[ur]^{\sigma} & tA^{'}_k \ar[l]^{tj^{'}}_\sim \ar[u]^\sim_{p_k}
}
\]
with $\sigma \eps {\cat W}$. Using LAP2 (2) we see that ${\cat Z}^{'}_1, {\cat Z}^{'}_2$ define the same element in ${\bf ho}{\cat M}^{'}$.

(ii). If ${\cat Z}_2$ is obtained from ${\cat Z}_1$ by inserting an element $A_k \overset{1_{A_k}}{\lla} A_k$, then similar to (i) we can show that ${\cat Z}^{'}_1, {\cat Z}^{'}_2$ define the same element in ${\bf ho}{\cat M}^{'}$.

(iii).  If ${\cat Z}_2$ is obtained from ${\cat Z}_1$ by replacing an element $A_{k-1} \overset{f_k}{\lra} A_{k} \overset{f_{k+1}}{\lra} A_{k+1}$ with $A_{k-1} \overset{f_{k+1}f_k}{\lra} A_{k+2}$ , then ${\cat Z}^{'}_2$ is obtained from ${\cat Z}^{'}_1$ by replacing the element 
\[
A_{k-1}^{'} \overset{i^{'}_{k}}{\lra} D^{'}_{k} \underset{\sim}{\overset{j^{'}_{k}}{\lla}} A^{'}_{k} \overset{i^{'}_{k+1}}{\lra} D^{'}_{k+1} \underset{\sim}{\overset{j^{'}_{k+1}}{\lla}} A^{'}_{k+1}
\]
with an element
\[
A_{k-1}^{'} \overset{i^{'}}{\lra} D^{'} \underset{\sim}{\overset{j^{'}}{\lla}} A^{'}_{k+1}
\]
with the property that there exists a commutative diagram
\[
\xymatrix{
  A_{k-1} \ar[rr]^{f_{k+1}f_k} && A_{k+1} \\
  tA^{'}_{k-1} \ar[u]^{p_{k-1}}_\sim \ar[r]_{ti^{'}} & tD^{'} \ar[ur]^{\sigma} & tA^{'}_{k+1} \ar[l]^{tj^{'}}_\sim \ar[u]^\sim_{p_{k+2}}
}
\]
with $\sigma \eps {\cat W}$. Using LAP2 (1) construct $D^{''}, i^{''} \eps {\cat M}^{'}$, $j^{''} \eps {\cat W}^{''}$ and $\sigma^{''} \eps {\cat W}$ making commutative the diagram
\[
\xymatrix {
& A_{k+1} & \\ 
\; tD^{'}_{k} \; \ar[ru]^{f_{k+1}\sigma_{k}} \ar@{-->}[rd]_{ti^{''}} && tD^{'}_{k+1} \ar@{-->}[ld]^{tj^{''}}_\sim \ar[lu]_{\sigma_{k+1}}^\sim \\
& tD^{''} \ar@{-->}[uu]^{\sigma^{''}} & 
} 
\]

In the diagram
\[
\xymatrix {
  && tA^{'}_k \ar[ld]_{tj^{'}_k}^\sim \ar[rd]^{ti^{'}_{k+1}} &&& \\
  tA^{'}_{k-1} \ar[r]^{ti^{'}_{k}} \ar[rrdd]_(.3){ti^{'}} & tD^{'}_{k} \ar[rd]_{ti^{''}} && tD^{'}_{k+1} \ar[ld]^{tj^{''}}_\sim & tA^{'}_{k+1} \ar[l]_{tj^{'}_{k+1}}^\sim \ar[lldd]^(.3){tj^{'}}_(.3)\sim & \\
  && tD^{''} \ar@{-->}[rrrd]^(.7){\sigma^{''}} &&&  \\
  && tD^{'} \ar@{-->}[rrr]_(.7){\sigma} &&& A_{k+1}
} 
\]
we have $\sigma^{''} \comp t(i^{''}j^{'}_k) = \sigma^{''} \comp t(j^{''}i^{'}_{k+1})$, $\sigma^{''} \comp t(i^{''}i^{'}_k) = \sigma \comp ti^{'}$ and $\sigma^{''} \comp t(j^{''}j^{'}_{k+1}) = \sigma \comp tj^{'}$. Using LAP2 (2) we see that ${\cat Z}^{'}_1, {\cat Z}^{'}_2$ define the same element in ${\bf ho}{\cat M}^{'}$.

(iv).  If ${\cat Z}_2$ is obtained from ${\cat Z}_1$ by replacing weak equivalences $A_{k-1} \underset{\sim}{\overset{f_k}{\lla}} A_{k} \underset{\sim}{\overset{f_{k+1}}{\lla}} A_{k+1}$ with $A_{k-1} \underset{\sim}{\overset{f_{k+1}f_k}{\lla}} A_{k+1}$ , then similar to (iii) one shows that ${\cat Z}^{'}_1, {\cat Z}^{'}_2$ define the same element in ${\bf ho}{\cat M}^{'}$.

(v). If ${\cat Z}_2$ is obtained from ${\cat Z}_1$ by inserting weak equivalences $A_{k} \underset{\sim}{\overset{f}{\lra}} A \underset{\sim}{\overset{f}{\lla}} A_{k}$, then ${\cat Z}^{'}_2$ is obtained from ${\cat Z}^{'}_1$ by inserting an element
\[
A_{k}^{'} \overset{i^{'}}{\lra} D^{'} \underset{\sim}{\overset{j^{'}}{\lla}} A^{'} \underset{\sim}{\overset{j^{'}}{\lra}} D^{'} \overset{i^{'}}{\lla}  A^{'}_k
\]
with the property that there exists a commutative diagram
\[
\xymatrix{
  A_{k}  \ar[rr]^f_\sim && A && A_{k} \ar[ll]_{f}^\sim \\
  tA^{'}_{k} \ar[u]^{p_{k}}_\sim \ar[r]_{ti^{'}} & tD^{'} \ar[ur]^{\sigma} & tA^{'} \ar[l]^{tj^{'}}_\sim \ar[r]_{tj^{'}}^\sim  \ar[u]^\sim_{p} & tD^{'} \ar[lu]_\sigma & tA^{'}_{k} \ar[u]_{p_{k}}^\sim \ar[l]^{ti^{'}} 
}
\]
with $\sigma \eps {\cat W}$. Using \lemmaref{lem:sufficientcondequiv}, we see that $i^{'} \eps \overline{\cat W}^{'}$, and so ${\cat Z}^{'}_1, {\cat Z}^{'}_2$ define the same element in ${\bf ho}{\cat M}^{'}$.

(vi). If ${\cat Z}_2$ is obtained from ${\cat Z}_1$ by replacing weak equivalences $A_k$ with $A_{k} \underset{\sim}{\overset{f}{\lla}} A \underset{\sim}{\overset{f}{\lra}} A_{k}$, similar to (v) one can show that ${\cat Z}^{'}_1, {\cat Z}^{'}_2$ define the same element in ${\bf ho}{\cat M}^{'}$.

From (i)-(vi) we conclude that ${\bf s}_{\cat C}$ is well defined on the maps of ${\bf ho}{\cat M}$.

(vii). If the end of ${\cat Z}_1$ and the beginning of ${\cat Z}_2$ coincide, then so do the end of ${\cat Z}^{'}_1$ and the beginning of ${\cat Z}^{'}_2$. From this observation, we see that ${\bf s}_{\cat C}$ preserves composition of maps.

(viii). A proof similar to (i) shows that ${\bf s}_{\cat C}$ preserves unit morphisms.

Using (vii) and (viii), we see that ${\bf s}_{\cat C} \colon {\bf ho}{\cat M} \ra {\bf ho}{\cat M}^{'}$ is a functor.

The functor ${\bf s}_{\cat C}$ is essentially surjective on objects - this can be verified using LAP1, LAP2 and \lemmaref{lem:sufficientcondequiv}.

From the commutativity of \equationref{eqn:approxarehtpyequiv1} and \equationref{eqn:approxarehtpyequiv2}, we get a natural isomorphism $({\bf ho} t) \comp {\bf s}_{\cat C} \cong 1_{{\bf ho}{\cat M}}$. This shows that ${\bf s}_{\cat C}$ is faithful. 

The functor ${\bf s}_{\cat C}$ is also full - to see that, pick a zig-zag ${\cat Z}^{''}$ in ${\cat M}^{'}$
\[
\xymatrix{ B^{'}_0 \ar@{-}[r]^-{g^{'}_1} &  B^{'}_1 \ar@{-}[r]^-{g^{'}_2} & ... \ar@{-}[r]^-{g^{'}n} & B^{'}_n}
\]
and denote $A_k = tB^{'}_k$, $f_k = tg^{'}_k$. The equivalence class of the zig-zag ${\cat Z}$ of \equationref{eqn:approxarehtpyequiv} is mapped by ${\bf s}_{\cat C}$ to the equivalence class of the zig-zag ${\cat Z}^{'}$ described earlier. For each $k$, using LAP2 (1) we can construct objects $A^{''}_k \eps {\cat M}^{'}$, maps $p^{''}_k \colon A^{''}_k \ra  B^{'}_k$, $q^{''}_k \colon A^{''}_k \ra  A^{'}_k$ and $p^{'}_k \colon tA^{''}_k \ra A_k$  such that $q^{''}_k \eps {\cat W}^{'}$, $p^{'}_k \eps {\cat W}$, $p^{'}_k \comp tp^{''}_k = 1_{A_k}$, $p^{'}_k \comp tq^{''}_k = p_k$. The maps $p^{''}_k$ are in $\overline{\cat W}^{'}$ by \lemmaref{lem:sufficientcondequiv}. 

Using LAP2 (1) twice in a row we can then construct objects $D^{''}_k \eps {\cat M}^{'}$ and maps $i^{''}_k \colon A^{''}_{k-1} \ra D^{''}_k$, $j^{''} \colon A^{''}_k \ra D^{''}_k$, $\nu^{'}_k \colon D^{'}_k \ra D^{''}_k$ and $\tau_k \colon tD^{''}_k \ra (\mathrm{\; target \; of \; }f_k)$ such that $j^{''}_k \eps {\cat W}^{'}$, $\nu^{'}_k \eps \overline{\cat W}^{'}$ with $\nu^{'}_k i^{'}_k q^{''}_k = i^{''}_k$, $\nu^{'}_k j^{'}_k q^{''}_k = j^{''}_k$ and $\tau_k \comp t \nu^{'}_k = \sigma_k$. 

The maps $p^{''}_k, q^{''}_k, i^{''}_k, j^{''}_k $ yield a zig-zag ${\cat Z}^{'''}$. Using LAP2 (2), we see that ${\cat Z}^{'}$ and ${\cat Z}^{'''}$, resp. ${\cat Z}^{''}$ and ${\cat Z}^{'''}$ define isomorphic maps in ${\bf ho}{\cat M}^{'}$. This concludes the proof that ${\bf s}_{\cat C}$ is full.

In conclusion, ${\bf s}_{\cat C}$ is an equivalence of categories, and therefore so is its quasi-inverse ${\bf ho}t$.
\end{proof}

\begin{rem}
Since ${\bf ho} t$ is an equivalence, for any other left cleavage ${\cat C}^{'}$ we have a canonical isomorphism ${\bf s}_{\cat C} \cong {\bf s}_{{\cat C}^{'}}$. By \lemmaref{lem:cleavagemap}, there exist a left cleavage ${\cat C}^{''}$ and a zig-zag of left cleavage maps ${\cat C} \ra {\cat C}^{''} \la {\cat C}^{'}$. It is straightforward to see that any such zig-zag induces the canonical isomorphisms ${\bf s}_{\cat C} \cong {\bf s}_{{\cat C}^{''}} \cong {\bf s}_{{\cat C}^{'}}$.
\end{rem}

\begin{cor}
\label{cor:equivleftapprox}
Suppose that $t \colon ({\cat M}^{'}$, ${\cat W}^{'})$ $\ra$ $({\cat M}$, ${\cat W})$ is a category pair functor.
\begin{enumerate}
\item
The following statements are equivalent:
\begin{enumerate}
\item
$t$ is a left approximation
\item
$t$ satisfies LAP1, LAP2 (1) and induces an equivalence of categories ${\bf ho}t$
\end{enumerate}
\item
The following statements are equivalent:
\begin{enumerate}
\item
$t$ is a right approximation
\item
$t$ satisfies RAP1, RAP2 (1) and induces an equivalence of categories ${\bf ho}t$
\end{enumerate}
\end{enumerate}
\end{cor}

\begin{proof}
(b) $\Ra$ (a) is immediate, and (a) $\Ra$ (b) follows from the Approximation \theoremref{thm:approxthm}.
\end{proof}

\begin{cor}
\label{cor:otherleftapprox}
Left (resp. right) approximation functors are closed under composition.
\end{cor}

\begin{proof}
Category pair functors satisfying LAP1 and LAP2 (1), resp. RAP1 and RAP2 (1) are closed under composition. The corollary now is a consequence of \corollaryref{cor:equivleftapprox}.
\end{proof}

\section{Total derived functors}
\label{sec:totalleftrightderived}

A functor $u \colon {\cat M}_1 \ra {\cat M}_2$ between two category pairs $({\cat M}_1, {\cat W}_1), ({\cat M}_2, {\cat W}_2)$ descends to a functor ${\bf ho} u \colon {\bf ho} {\cat M}_1 \ra {\bf ho} {\cat M}_2$ if and only if $u({\cat W}_1) \subset \overline{\cat W}_2$, where $\overline{\cat W}_2$ denotes the saturation of ${\cat W}_2$ in ${\cat M}_2$.

In the general case however ${\bf ho} u$ does not exist, and the best we can hope for is the existence of a left (or a right) Kan extension of $\gamma_{{\cat M}_2}u$ along $\gamma_{{\cat M}_1}$, also called the {\it total left} (resp. right) {\it derived functors} of $u$.

\index{functor!total derived}
\begin{defn}
\label{defn::def:totalleftkanextensions}
Suppose that $({\cat M}_1, {\cat W}_1), ({\cat M}_2, {\cat W}_2)$ are two categories with weak equivalences, with localization functors denoted $\gamma_{{\cat M}_1}$ and respectively $\gamma_{{\cat M}_1}$, and suppose that $u \colon {\cat M}_1 \ra {\cat M}_2$ is a functor.
\begin{enumerate}
\item
The {\it total left derived} functor of $u$, denoted $({\bf L}u, \eps)$ is the left Kan extension $({\bf L}_{\gamma_{{\cat M}_1}} (\gamma_{{\cat M}_2} u), \eps)$ of $\gamma_{{\cat M}_2} u$ along $\gamma_{{\cat M}_1}$
\begin{center}
$\xymatrix {
    {\cat M}_1 \ar[r]^u \ar[d]_{\gamma_{_{{\cat M}_1}}} & 
    {\cat M}_2 \ar[d]^{\gamma_{_{{\cat M}_2}}} \\
    {\bf ho}{\cat M}_1 \ar[r]_{{\bf L}u}  & 
    {\bf ho}{\cat M}_2 \ultwocell<\omit>{<0>\eps}
  }$
\end{center}
\item
The {\it total right derived} functor of $u$, denoted $({\bf R}u, \nu)$, is the right Kan extension $({\bf R}_{\gamma_{{\cat M}_1}} (\gamma_{{\cat M}_2} u), \nu)$ of $\gamma_{{\cat M}_2} u$ along $\gamma_{{\cat M}_1}$
\begin{center}
$\xymatrix {
    {\cat M}_1 \ar[r]^u \ar[d]_{\gamma_{_{{\cat M}_1}}} \drtwocell<\omit>{<0>\eta} & 
    {\cat M}_2 \ar[d]^{\gamma_{_{{\cat M}_2}}} \\
    {\bf ho}{\cat M}_1 \ar[r]_{{\bf R}u} & 
    {\bf ho}{\cat M}_2
  }$
\end{center}
\end{enumerate}
\end{defn}

The total left and derived functors $({\bf L}u, \eps), ({\bf R}u, \nu)$ are defined in terms of the localization functors $\gamma_{{\cat M}_1}, \gamma_{{\cat M}_2}$ and therefore will not change if we replace in the definition ${\cat W}_1$, ${\cat W}_2$ with their saturations $\overline{\cat W}_1$, $\overline{\cat W}_2$.

Note that if $u({\cat W}_1) \subset \overline{\cat W}_2$ then ${\bf L}u = {\bf R}u = {\bf ho}u$.

\section{An existence criterion for total derived functors}

\begin{thm}
\label{thm:generalexistenceleftderived}
Given three categories with weak equivalences $({\cat M}^{'}_1$, ${\cat W}^{'}_1)$, $({\cat M}_1$, ${\cat W}_1)$ and $({\cat M}_2$, ${\cat W}_2)$ and two functors $\xymatrix {{\cat M}^{'}_1 \ar[r]^t & {\cat M}_1 \ar[r]^u & {\cat M}_2}$.
\begin{enumerate}
\item
If $t$ is a left approximation and $ut$ preserves weak equivalences, then $u$ admits a total left derived functor $({\bf L}u, \eps)$. The natural map $\eps \colon ({\bf L}u)(tA^{'}) \Ra utA^{'}$ is an isomorphism in objects $A^{'}$ of ${\cat M}^{'}_1$.
\item
If $t$ is a right approximation and $ut$ preserves weak equivalences, then $u$ admits a total right derived functor $({\bf R}u, \eta)$. The natural map $\eta \colon utA^{'} \Ra ({\bf R}u)(tA^{'})$ is an isomorphism in $A^{'}$.
\end{enumerate}
\end{thm}

\begin{proof}
We only prove (1). Pick a regular left cleavage ${\cat C} = (C, D, p, i, j, \sigma)$ along $t$. 
\begin{equation}
\label{eqn:leftcleavage}
\xymatrix { 
A \ar[rr]^f && B \\
tC(A) \ar[u]^{p_A}_\sim \ar[r]_{ti(f)} & 
tD(f) \ar[ur]^{\sigma(f)} &
tC(B) \ar[l]_\sim^{tj(f)} \ar[u]_{p_B}^\sim
}
\end{equation}

Recall that the functor ${\bf s}_{\cat C}$ in the proof of the Approximation \theoremref{thm:approxthm} is a quasi-inverse of ${\bf ho}t$, and define ${\bf L}u = {\bf ho}(ut) \comp {\bf s}_{\cat C}$.

Let us spell out in detail the functor ${\bf L}u$. For an object $A$ of ${\cat M}_1$, we have $({\bf L}u) (A) = utC(A)$. On maps $f \colon A \ra B$ of ${\cat M}_1$, we have ${\bf L}u(f) = (utj(f))^{-1} uti(f)$.

We define $\eps \colon ({\bf L}u)(A) \ra uA$ as $u(p_A)$. The commutativity of the diagram (\ref{eqn:leftcleavage}) implies that $\eps \colon ({\bf L}u) \gamma_{{\cat M}_1} \Ra \gamma_{{\cat M}_2} u$ is a natural map.

We need to show that the pair $({\bf L}u, \eps)$ is {\it terminal} among pairs $(\Lambda, \lambda)$ where $\Lambda \colon {\bf ho} {\cat M}_1 \ra {\bf ho} {\cat M}_2$ is a functor and $\lambda \colon \Lambda \gamma_{{\cat M}_1} \Ra \gamma_{{\cat M}_2} u$ is a natural transformation. For any object $A$ of ${\cat M}_1$ the sequence of full maps in ${\bf ho}{\cat M}_2$ and their inverses
\begin{equation}
\label{eqn:terminalproperty}
\xymatrix { 
    \Lambda(A) \ar@{-->}[rr]^{\lambda(A)} & & uA \\
    \Lambda tC(A)  \ar[u]^-{\Lambda (p_A)}_-\sim \ar[rr]^-{\lambda (tC(A))} & & utC(A) = {\bf L}u(A) \ar@{-->}[u]_{u(p_A) = \eps(A)}
  }
\end{equation}
defines a map $\delta \colon \Lambda(A) \ra {\bf L}u(A)$. For maps $f \colon A \ra B$ of ${\cat M}_1$ we have a commutative diagram in ${\bf ho}{\cat M}_2$
\begin{equation}
\label{eqn:terminalproperty2}
\xymatrix { 
    \Lambda(A) \ar[rr]^{\Lambda(f)} && \Lambda(B) \\
    \Lambda tC(A) \ar[u]^-{\Lambda (p_A)}_-\sim \ar[r]^{\Lambda ti(f)} \ar[d]_-{\lambda (tC(A))} &
    \Lambda tD(f) \ar[ru]^-{\Lambda \sigma(f)} \ar[d]_-{\lambda (tD(f))} &
    \Lambda tC(B) \ar[u]_-{\Lambda (p_B)}^-\sim \ar[l]_{\Lambda tj(f)}^\sim \ar[d]^-{\lambda (tC(B))} \\
    utC(A) \ar[r]_{uti(f)} & utD(f) & utC(B) \ar[l]^{utj(f)}_\sim    
  }
\end{equation}
where $utj(f)$ is a weak equivalence since $j(f)$ is. The commutativity of $(\ref{eqn:terminalproperty2})$ implies that $\delta \colon \Lambda \Ra {\bf L}u$ is natural in maps of ${\cat M}_1$, therefore natural in maps of ${\bf ho}{\cat M}_1$. 

The commutativity of diagram (\ref{eqn:terminalproperty}) shows that we have $\eps \star \delta = \lambda$. 

The natural map $\eps: ({\bf L}u)(tA^{'}) \Ra utA^{'}$ is an isomorphism for {\it regular} objects $A^{'} \eps {\cat M}^{'}_1$, therefore for any object $A^{'}  \eps {\cat M}^{'}_1$. To see that, note that for $A^{'}$ regular the map $\eps: ({\bf L}u)(tA^{'}) \Ra utA^{'}$ can be identified with $1_{utA^{'}}$. Furthermore, our left cleavage was assumed to be regular so any object in ${\cat M}^{'}_1$ is isomorphic in ${\bf ho}{\cat M}^{'}_1$ to a regular object $A^{'}$.
\end{proof}

\begin{cor}
\label{cor:generalexistenceleftderived}
\mbox{}
\begin{enumerate}
\item
If $t, u$ are as in \theoremref{thm:generalexistenceleftderived} (1) and ${\bf s}$ denotes any quasi-inverse of ${\bf ho}t$, then ${\bf ho}(ut) \, {\bf s}$ is naturally isomorphic to ${\bf L}u$.
\item
If $t, u$ are as in \theoremref{thm:generalexistenceleftderived} (2) and ${\bf s}$ denotes any quasi-inverse of ${\bf ho}t$, then ${\bf ho}(ut) \, {\bf s}$ is naturally isomorphic to ${\bf R}u$. $\square$
\end{enumerate}
\end{cor}

\section{The abstract Quillen adjunction property}

We next state an adjunction property of total derived functors. 
\begin{thm}[Abstract Quillen adjunction]
\label{thm:generalexistencetotalderivedadjoint}
Given four categories with weak equivalences $({\cat M}^{'}_1$, ${\cat W}^{'}_1)$, $({\cat M}_1$, ${\cat W}_1)$, $({\cat M}^{'}_2$, ${\cat W}^{'}_2)$, $({\cat M}_2$, ${\cat W}_2)$ and four functors
\begin{center}
$\xymatrix {{\cat M}^{'}_1 \ar[r]^{t_1} & {\cat M}_1 \ar@<2pt>[r]^{u_1} & {\cat M}_2 \ar@<2pt>[l]^{u_2} & {\cat M}^{'}_2 \ar[l]_{t_2} }$
\end{center}
where 
\begin{enumerate}
\item
$u_1 \dashv u_2$ is an adjoint pair
\item
$t_1$ is a left approximation, $t_2$ is a right approximation
\item
$u_1t_1$ and $u_2t_2$ preserve weak equivalences
\end{enumerate}
then ${\bf L}u_1 \dashv {\bf R}u_2$ is a naturally adjoint pair 
\begin{center}
$\xymatrix {{\bf ho}{\cat M}_1 \ar@<2pt>[r]^{{\bf L}u_1} & {\bf ho}{\cat M}_2 \ar@<2pt>[l]^{{\bf R}u_2} }$
\end{center}
If additionally
\begin{enumerate}
\item[(4)]
(resp. (4r), resp. (4l)).
For any objects $A^{'} \eps {\cat M}^{'}_1$, $B^{'} \eps {\cat M}^{'}_2$, a map $u_1t_1A^{'} \ra t_2B^{'}$ is a weak equivalence iff (resp. if, resp. only if) its adjoint $t_1A^{'} \ra u_2t_2B^{'}$ is a weak equivalence
\end{enumerate}
then ${\bf L}u_1$ and ${\bf R}u_2$ are inverse equivalences of categories (resp. ${\bf R}u_2$ is fully faithful, resp ${\bf L}u_1$ is fully faithful).
\end{thm}

This theorem suggests the following
\index{Quillen!abstract adjunction}
\index{Quillen!abstract equivalence}
\begin{defn}
\label{defn:abstractQuillenadjoint}
We will call the functors $u_1, u_2$ satisfying the properties (1), (2), (3) of \theoremref{thm:generalexistencetotalderivedadjoint} an {\it abstract Quillen adjoint} pair with respect to $t_1, t_2$. If the additional property (4) is satisfied, we will call $u_1, u_2$ an {\it abstract Quillen equivalence} pair with respect to $t_1, t_2$.
\end{defn}

There is a very nice, conceptual proof due to Georges Maltsiniotis \cite{Maltsiniotis3} of \theoremref{thm:generalexistencetotalderivedadjoint}, using the universal property in the definition of total derived functors and \theoremref{thm:generalexistenceleftderived}. 

In this text however, we will let \theoremref{thm:generalexistencetotalderivedadjoint} be a consequence of the theorem below. In preparation, let us introduce more notations and definitions. Given a diagram of functors
\begin{equation}
\label{eqn:partialadj}
\xymatrix{
    {\cat A} \ar[r]^{v_1} \ar[d]_{t_1} & {\cat B} \\
    {\cat C} & {\cat D} \ar[l]^{v_2} \ar[u]_{t_2}
  }
\end{equation}
a {\it partial adjunction} \index{functor!partially adjoint} between $v_1, v_2$ with respect to $t_1, t_2$ is a bijection
\begin{center}
$\zeta \colon Hom_{{\cat B}}(v_1A, t_2D) \cong Hom_{{\cat C}}(t_1A, v_2D)$ 
\end{center}
natural in objects $A \eps {\cat A}_1, D \eps {\cat D}$. Whenever we say that $v_1, v_2$ is an adjoint pair with respect to $t_1, t_2$ we refer to a particular bijection $\zeta$. Note that $\zeta^{-1}$ defines a partial adjunction between $t_1, t_2$ with respect to $v_1, v_2$. 

If in addition we have adjoint pairs $t_1 \dashv s_1$, $s_2 \dashv t_2$ 
\begin{equation}
\label{eqn:partialadj2}
\xymatrix{
    {\cat A} \ar[r]^{v_1} \ar@<-2pt>[d]_{t_1} & {\cat B} \ar@<-2pt>[d]_{s_2} \\
    {\cat C} \ar@<-2pt>[u]_{s_1} & {\cat D} \ar[l]^{v_2} \ar@<-2pt>[u]_{t_2}
  }
\end{equation}
then the partial adjunctions of $v_1, v_2$ with respect to $t_1, t_2$ are in a one to one correspondence with adjunctions $s_2v_1 \dashv s_1v_2$. In the diagram (\ref{eqn:partialadj2}), assuming that $s_1, t_1$ and $s_2, t_2$ are inverse equivalences of categories, then the partial adjuctions of $v_1, v_2$ with respect to $t_1, t_2$ are in a one to one correspondence with adjunctions $s_2v_1 \dashv s_1v_2$, and furthermore in a one to one correspondence with adjunctions $v_1s_1 \dashv v_2 s_2$.

In a diagram of the form
\begin{equation}
\label{eqn:partialadj3}
\xymatrix{
    {\cat A}^{'} \ar[r]^{a} & {\cat A} \ar[r]^{v_1} \ar[d]_{t_1} & {\cat B} & \\
    & {\cat C} & {\cat D} \ar[l]^{v_2} \ar[u]_{t_2} & {\cat D}^{'} \ar[l]_{d}
  }
\end{equation}
a partial adjunction between $v_1, v_2$ with respect to $t_1, t_2$ will induce a partial adjunction between $v_1a, v_2d$ with respect to $t_1a, t_2d$, given by 
\begin{center}
$\zeta \colon Hom_{{\cat B}}(v_1aA^{'}, t_2dD^{'}) \cong Hom_{{\cat C}}(t_1aA^{'}, v_2dD^{'})$ 
\end{center}

We can now state
\begin{thm}[Abstract Quillen partial adjunction]
\label{thm:generalexistencetotalderivedadjoint2}
Given four categories with weak equivalences $({\cat M}^{'}_1$, ${\cat W}^{'}_1)$, $({\cat M}_1$, ${\cat W}_1)$, $({\cat M}^{'}_2$, ${\cat W}^{'}_2)$, $({\cat M}_2$, ${\cat W}_2)$ and four functors $t_1, t_2, v_1, v_2$
\begin{center}
$\xymatrix{
    {\cat M}^{'}_1 \ar[r]^{v_1} \ar[d]_{t_1} & {\cat M}_2 \\
    {\cat M}_1 & {\cat M}^{'}_2 \ar[l]_{v_2} \ar[u]_{t_2} &
    }$
\end{center}
such that:
\begin{enumerate}
\item
$v_1, v_2$ are partially adjoint with respect to $t_1, t_2$, meaning that there exists a bijection
\begin{center}
$\zeta \colon Hom_{{\cat M}_2}(v_1A^{'}, t_2B^{'}) \cong Hom_{{\cat M}_1}(t_1A^{'}, v_2B^{'})$
\end{center}
natural in $A^{'} \eps {\cat M}^{'}_1$ and $B^{'} \eps {\cat M}^{'}_2$
\item
$t_1$ is a left approximation, $t_2$ is a right approximation
\item
$v_1$ and $v_2$ preserve weak equivalences
\end{enumerate}

Then ${\bf ho}v_1$, ${\bf ho}v_2$ are naturally partial adjoint with respect to ${\bf ho}t_1$, ${\bf ho}t_2$. Equivalently, denote ${\bf s}_i$ a quasi-inverse of ${\bf ho}t_i$, and let ${\bf V}_i = {\bf ho}(v_i) {\bf s}_i$ for $i = 1, 2$. Then ${\bf V}_1 \dashv {\bf V}_2$
\begin{center}
$\xymatrix {{\bf ho}{\cat M}_1 \ar@<2pt>[r]^{{\bf V}_1} & {\bf ho}{\cat M}_2 \ar@<2pt>[l]^{{\bf V}_2}}$
\end{center}
is a naturally adjoint pair.

If additionally
\begin{enumerate}
\item[(4)]
(resp. (4r), resp. (4l)).
For any objects $A^{'} \eps {\cat M}^{'}_1$, $B^{'} \eps {\cat M}^{'}_2$, a map $v_1A^{'} \ra t_2B^{'}$ is a weak equivalence iff (resp. if, resp. only if) its partial adjoint $t_1A^{'} \ra v_2B^{'}$ is a weak equivalence
\end{enumerate}
then ${\bf V}_1$ and ${\bf V}_2$ are inverse equivalences of categories (resp. ${\bf V}_2$ is fully faithful, resp ${\bf V}_1$ is fully faithful).
\end{thm}

The following definition is suggested:
\index{Quillen!abstract partial adjunction}
\index{Quillen!abstract partial equivalence}
\begin{defn}
\label{defn:abstractQuillenadjoint2}
We will call the functors $v_1, v_2$ satisfying the properties (1), (2), (3) of \theoremref{thm:generalexistencetotalderivedadjoint2} an {\it abstract Quillen partially adjoint} pair with respect to $t_1, t_2$. If the additional property (4) is satisfied, we will call $v_1, v_2$ an {\it abstract Quillen partial equivalence} pair with respect to $t_1, t_2$.
\end{defn}

\begin{proof}[Proof of \theoremref{thm:generalexistencetotalderivedadjoint} assuming \theoremref{thm:generalexistencetotalderivedadjoint2}]
Since $u_1 \dashv u_2$ is an adjoint pair, we see that $v_1 = u_1t_1$, $v_2 = u_2t_2$ is partially adjoint with respect to $t_1, t_2$. From \corollaryref{cor:generalexistenceleftderived}, we have natural isomorphisms ${\bf L}u_1 \cong {\bf ho}(u_1t_1) \, {\bf s}_1 = {\bf V}_1$ and ${\bf R}u_2 \cong {\bf ho}(u_2t_2) \, {\bf s}_2 = {\bf V}_2$. The statement now follows.
\end{proof}

\begin{proof}[Proof of \theoremref{thm:generalexistencetotalderivedadjoint2}]
\mbox{} 
If we can prove the conclusion for a particular choice of ${\bf s}_1 \colon {\bf ho}{\cat M}_1 \ra {\bf ho}{\cat M}^{'}_1$ and ${\bf s}_2 \colon {\bf ho}{\cat M}_2 \ra {\bf ho}{\cat M}^{'}_2$, then the conclusion follows for any ${\bf s}_1$ and ${\bf s}_2$.

We pick a regular left cleavage ${\cat C} = (C, D, p, i, j, \sigma)$ along $t_1$, and a regular right cleavage ${\cat R} = (R, S, i, p^{'}, q^{'}, \tau)$ along $t_2$. As in the proof of the Approximation \theoremref{thm:approxthm}, we get a quasi-inverse ${\bf s}_1 = {\bf s}_{\cat C}$ to ${\bf ho}t_1$ and a quasi-inverse ${\bf s}_2 = {\bf s}_{\cat R}$ to ${\bf ho}t_2$. We will work with these particular choices ${\bf s}_1$ and ${\bf s}_2$.

Denote ${\bf V}_1 = {\bf ho}(v_1) {\bf s}_1$ and ${\bf V}_2 = {\bf ho}(v_2) {\bf s}_2$. We will construct natural maps ${\bf \Phi} \colon 1_{{\bf ho}{\cat M}_1} \Ra {\bf V}_2{\bf V}_1$ and ${\bf \Psi} \colon {\bf V}_1{\bf V}_2 \Ra 1_{{\bf ho}{\cat M}_2}$, and show that they are the unit and counit of an adjunction between ${\bf V}_1$ and ${\bf V}_2$. We start by constructing a natural map
\begin{center}
$\overline{\bf \Phi} \colon {\bf ho}t_1 \Ra {\bf ho}(v_2) \, {\bf s}_2 \, {\bf ho}(v_1)$
\end{center}
by defining $\overline{\bf \Phi} (A^{'}) = \zeta i_{v_1A^{'}}$ for any object $A^{'}$ of ${\bf ho}{\cat M}^{'}_1$, where $i_{v_1A^{'}}$ and $\zeta i_{v_1A^{'}}$ are the maps 
\begin{center}
$\xymatrix{
    v_1A^{'} \ar[r]^-{i_{v_1A^{'}}}_-\sim & t_2R(v_1A^{'}) & t_1A^{'} \ar[rrr]^-{\zeta i_{v_1A^{'}}} &&& v_2R(v_1A^{'})
  }$
\end{center}

Given a map $f \colon A^{'} \ra B^{'}$ in ${\cat M}^{'}_1$ we get a commutative diagram
\[
\xymatrix {
  v_1A^{'} \ar[rr]^{v_1f} \ar[d]_-{i_{v_1A^{'}}}^-\sim \ar[dr]^{\tau(v_1f)} && 
  v_1B^{'} \ar[d]^-{i_{v_1B^{'}}}_-\sim \\
  t_2R(v_1A^{'}) &
  t_2S(v_1f) \ar[l]^{t_2q(v_1f)}_\sim \ar[r]_{t_2p(v_1f)} &
  t_2R(v_1B^{'}) 
}
\]
in ${\cat M}_2$, where $q(v_1f)$ is a weak equivalence. Applying the natural bijection $\zeta$ we get a commutative diagram in ${\cat M}_1$
\[
\xymatrix {
  t_1A^{'} \ar[rr]^{t_1f} \ar[d]_-{\zeta i_{v_1A^{'}}} \ar[dr]^{\zeta \tau(v_1f)} &&
  t_1B^{'} \ar[d]^-{\zeta i_{v_1B^{'}}} \\
  v_2R(v_1A^{'}) &
  v_2S(v_1f) \ar[l]^{v_2q(v_1f)}_\sim \ar[r]_{v_2p(v_1f)} &
  v_2R(v_1B^{'})
}
\]
where $v_2q(v_1f)$ is a weak equivalence since $q(v_1f)$ is. The commutativity of the second diagram shows that $\overline{\bf \Phi}$ is natural in maps of ${\cat M}^{'}_1$, and therefore natural in maps of ${\bf ho}{\cat M}^{'}_1$.

Since ${\bf ho}t_1$ and ${\bf s}_1$ are quasi-inverses of each other, the natural map $\overline{\bf \Phi} \colon {\bf ho}t_1 \Ra {\bf ho}(v_2) \, {\bf s}_2 \, {\bf ho}(v_1)$ yields the desired natural map 
\begin{center}
${\bf \Phi} \colon 1_{{\bf ho}{\cat M}_1} \Ra  {\bf ho}(v_2) \, {\bf s}_2 \, {\bf ho}(v_1) \, {\bf s}_1 = {\bf V}_2{\bf V}_1$
\end{center}

We dually construct a natural map
\begin{center}
$\overline{\bf \Psi} \colon {\bf ho}(v_1) \, {\bf s}_1 \, {\bf ho}(v_2) \Ra {\bf ho}t_2$ 
\end{center}
by defining $\overline{\bf \Psi} (B^{'}) = \zeta^{-1} p_{v_2B^{'}}$ for any object $B^{'}$ of ${\bf ho}{\cat M}^{'}_2$, where the maps $p_{v_2B^{'}}$ and $\zeta^{-1} p_{v_2B^{'}}$ are the maps
\begin{center}
$\xymatrix {
    t_1C(v_2B^{'}) \ar[r]^-{p_{v_2B^{'}}}_-\sim & v_2B^{'} &  & v_1C(v_2B^{'}) \ar[rr]^-{\zeta^{-1} p_{v_2B^{'}}} && t_2B^{'}
  }$
\end{center}
The proof that $\overline{\bf \Psi}$ is a natural map is dual to the proof that $\overline{\bf \Phi}$ is a natural map. Since ${\bf ho}t_2$ and ${\bf s}_2$ are quasi-inverses of each other, the natural map $\overline{\bf \Psi}$ yields the desired natural map 
\begin{center}
${\bf \Psi} \colon {\bf V}_1{\bf V}_2 = {\bf ho}(v_1) \, {\bf s}_1 \, {\bf ho}(v_2) \, {\bf s}_2 \Ra 1_{{\bf ho}{\cat M}_2} $
\end{center}

It remains to show that the natural maps ${\bf \Phi}, {\bf \Psi}$ are the unit resp. the counit of an adjunction between the functors ${\bf V}_1, {\bf V}_2$. In other words, we need to prove that the following composites are identities.
\begin{center}
\begin{equation}
\label{eqn:generaladj11}
\xymatrix {
    {\bf V}_1 \ar@{=>}[r]^-{{\bf V}_1 {\bf \Phi}} & {\bf V}_1{\bf V}_2{\bf V}_1 \ar@{=>}[r]^-{{\bf \Psi} {\bf V}_1} & {\bf V}_1
  }
\end{equation}
\begin{equation}
\label{eqn:generaladj12}
\xymatrix {
    {\bf V}_2 \ar@{=>}[r]^-{{\bf \Phi} {\bf V}_2} & {\bf V}_2{\bf V}_1{\bf V}_2 \ar@{=>}[r]^-{{\bf V}_2 {\bf \Psi}} & {\bf V}_2 
  }
\end{equation}
\end{center}

We only prove that (\ref{eqn:generaladj11}) is an identity, since the proof for (\ref{eqn:generaladj12}) is dual. It suffices to prove that (\ref{eqn:generaladj11}) is an identity on objects of the form $A = t_1A^{'}$, with $A^{'} \eps {\cat M}^{'}_1$ regular with respect to the left cleavage along $t_1$.

If $A^{'}$ is a regular object of ${\cat M}^{'}_1$, we have that $C(t_1A^{'}) = A^{'}$ and $p_{tA^{'}} = 1_{t_1A^{'}} \colon t_1C(t_1A^{'}) \ra t_1A^{'}$.

Denote $h = \zeta i_{v_1A^{'}}$. The left cleavage diagram associated to $h$ yields a commutative diagram in ${\cat M}_1$
\[
\xymatrix { 
  t_1A^{'} \ar[rr]^{h = \zeta i_{v_1A^{'}}} &&
  v_2R(v_1A^{'}) \\
  t_1C(t_1A^{'}) = t_1 A^{'} \ar[u]^-{p_{t_1A^{'}} = 1_{t_1A^{'}}} \ar[r]_-{t_1i(h)} &
  t_1D(h) \ar[ur]^{\sigma(h)} &
  t_1C(v_2R(v_1A^{'})) \ar[l]_-\sim^-{t_1j(h)} \ar[u]_-{p_{v_2R(v_1A^{'})}}^\sim 
}
\]

Applying the natural bijection $\zeta^{-1}$ to the diagram we get a commutative diagram in ${\cat M}_2$
\[
\xymatrix { 
  v_1A^{'} \ar[rr]^{\zeta^{-1}h = i_{v_1A^{'}}}_\sim &&
  t_2R(v_1A^{'}) \\
  v_1A^{'} \ar[u]^{1_{v_1A^{'}}} \ar[r]_-{v_1i(h)} &
  v_1D(h) \ar[ur]^{\zeta^{-1}\sigma(h)} &
  v_1C(v_2R(v_1A^{'})) \ar[l]_-\sim^-{v_1j(h)} \ar[u]_-{\zeta^{-1}p_{v_2R(v_1A^{'})}} 
}
\]
In this diagram, $j(h)$ and therefore $v_1j(h)$ are weak equivalences. 

Since $A = t_1A^{'}$, in ${\bf ho}{\cat M}_1$ we can identify $v_1A^{'}$ with ${\bf V}_1A$ and $v_1C(v_2R(v_1A^{'}))$ with ${\bf V}_1{\bf V}_2{\bf V}_1A$. Under this identification the composition $(v_1j(h))^{-1} \comp v_1i(h)$ becomes ${\bf V}_1{\bf \Phi}(A)$, and the composition $(i_{v_1A^{'}})^{-1}\comp (\zeta^{-1}p_{v_2R(v_1A^{'})})$ becomes ${\bf \Psi}{\bf V}_1(A)$. It follows that the composition (\ref{eqn:generaladj11}) is an identity. Dually, (\ref{eqn:generaladj12}) is an identity and we have proved that ${\bf V}_1 \dashv {\bf V}_2$ are adjoint with adjunction unit ${\bf \Phi}$ and counit ${\bf \Psi}$.

For the second part of the theorem, we assume hypothesis (4l) and we will show that ${\bf \Phi}$ is a natural isomorphisms. From \propositionref{prop:adjfunctors}, this will imply that ${\bf V}_1$ is fully faithful. 

For any object $A^{'}$ of ${\bf ho}{\cat M}^{'}_1$, the map
\begin{center}
$\xymatrix{
    v_1A^{'} \ar[r]^-{i_{v_1A^{'}}}_-\sim & t_2R(v_1A^{'}) 
  }$
\end{center}
and therefore from hypothesis (4) the map
\begin{center}
$\xymatrix{
    t_1A^{'} \ar[rrr]^-{\zeta i_{v_1A^{'}}}_-\sim &&& v_2R(v_1A^{'})
  }$
\end{center}
are weak equivalences. It follows that the natural maps $\overline{\bf \Phi}$ and therefore ${\bf \Phi}$ are isomorphisms. A dual proof shows that hypothesis (4l) implies that ${\bf \Psi}$ is a natural isomorphism, and therefore ${\bf V}_2$ is fully faithful. 

If hypothesis (4) is satisfied, then both ${\bf \Phi}$ and ${\bf \Psi}$ are isomorphisms, therefore ${\bf V}_1, {\bf V}_2$ are inverse equivalences of categories (cf. \propositionref{prop:adjfunctors}).
\end{proof}

\chapter{The homotopy category of a cofibration category}
\label{chap:homotcat}

Our goal in this chapter is to describe the homotopy category of a cofibration category. All the definitions and the results of this chapter actually only require the smaller set of precofibration category axioms CF1-CF4. 

For a precofibration category $({\cat M}$, ${\cat W}$, ${\cat Cof})$, recall that we have denoted ${\cat M}_{cof}$ the full subcategory of cofibrant objects of ${\cat M}$. We will show that the functor ${\bf ho}{\cat M}_{cof} \ra {\bf ho} {\cat M}$ is an equivalence of categories (Anderson, \cite{Anderson1}). In fact, we will develop an axiomatic description of {\it cofibrant approximation} functors $t \colon {\cat M}^{'} \ra {\cat M}$, which are modelled on the properties of the inclusion ${\cat M}_{cof} \ra {\cat M}$. Cofibrant approximation functors are left approximations in the sense of \definitionref{defn:wleftapproximation}, and by the Approximation \theoremref{thm:approxthm} they induce an isomorphism at the level of the homotopy category. 

Cofibrant aproximation functors will resurface later in \sectionref{sec:colimarbitrarycat}, when we will reduce the construction of homotopy colimits indexed by arbitrary small diagrams to the construction of homotopy colimits indexed by small direct categories.

We then turn to the study of homotopic maps in a precofibration category. In ${\cat M}_{cof}$ we will define the {\it left} homotopy relation $\simeq_l$ on maps, and show that $f \simeq_l g$ iff $f \simeq g$. The localization of ${\cat M}_{cof}$ modulo homotopy is denoted $\pi {\cat M}_{cof}$. 

We show that the class of weak equivalences between cofibrant objects admits a calculus of fractions in $\pi {\cat M}_{cof}$. As a consequence we obtain a description of ${\bf ho} {\cat M}_{cof}$ whereby any map in ${\bf ho} {\cat M}_{cof}$ can be written (up to homotopy!) as a 'left fraction' $f t^{-1}$, with $t$ a weak equivalence. The theory of homotopic maps and calculus of fractions up to homotopy as described here is due to Brown \cite{Brown}.

\section{Fibrant and cofibrant approximations}
\label{sec:cofapproximations}
We are interested in a class of precofibration category functors which are left approximations (\definitionref{defn:wleftapproximation}), and therefore
\begin{enumerate}
\item
induce isomorphisms when passed to the homotopy category (consequence of the Approximation \theoremref{thm:approxthm}), and 
\item
serve as resolutions for computing total left derived functors (\theoremref{thm:generalexistenceleftderived}). 
\end{enumerate}

These are the {\it cofibrant approximation} functors, defined below. The cofibrant approximation functors should be thought of as an axiomatization of the inclusion ${\cat M}_{cof}$ $\ra$ ${\cat M}$, where ${\cat M}$ is a precofibration category.

\index{approximation!(co)fibrant}
\begin{defn}(Cofibrant approximation functors)
\label{defn::cofapproximation}
\mbox{}
Let ${\cat M}$ be a precofibration category. A functor $t \colon {\cat M}^{'} \ra {\cat M}$ is a {\it cofibrant approximation} if ${\cat M}^{'}$ is a precofibration category with all objects cofibrant and
\begin{description}
\item[CFA1]
$t$ preserves the initial object and cofibrations
\item[CFA2]
A map $f$ of ${\cat M}^{'}$ is a weak equivalence if and only if $tf$ is a weak equivalence
\item[CFA3]
If $A \ra B$, $A \ra C$ are cofibrations of ${\cat M}^{'}$ then the natural map $tB \Sum_{tA} tC \ra t(B \Sum_{A} C)$ is an isomorphism
\item[CFA4]
Any map $f \colon tA \ra B$ factors as $f = r \comp tf^{'}$ with $f^{'}$ a cofibration of ${\cat M}^{'}$ and $r$ a weak equivalence of ${\cat M}$.
\end{description}
\end{defn}

A cofibrant approximation functor in particular sends any object to a cofibrant object, and sends trivial cofibrations to trivial cofibrations. If ${\cat M}$ is a precofibration category, then the inclusion ${\cat M}_{cof} \ra {\cat M}$ is a cofibrant approximation.

The dual definition for prefibration categories is

\begin{defn}(Fibrant approximation functors)
\label{defn::fibapproximation}
\mbox{}
Let ${\cat M}$ be a prefibration category. A functor $t \colon {\cat M}^{'} \ra {\cat M}$ is a {\it fibrant approximation} if ${\cat M}^{'}$ is a prefibration category with all objects fibrant and
\begin{description}
\item[FA1]
$t$ preserves the final object and fibrations
\item[FA2]
A map $f$ of ${\cat M}^{'}$ is a weak equivalence if and only if $tf$ is a weak equivalence
\item[FA3]
If $B \ra A$, $C \ra A$ are fibrations of ${\cat M}^{'}$ then the natural map $ t(B \times_{A} C) \ra tB \times_{tA} tC$ is an isomorphism
\item[FA4]
Any map $f \colon A \ra tB$ factors as $f = tf^{'} \comp s$ with $f^{'}$ a fibration of ${\cat M}^{'}$ and $s$ a weak equivalence of ${\cat M}$.
\end{description}
\end{defn}

We will need the lemmas below. Recall that two maps $f, g \colon A\ra B$ in a caetgry pair $({\cat M}, {\cat W})$ are homotopic $f \simeq g$ by definition if they have the same image in ${\bf ho}{\cat M}$. The prototypic example of homotopic maps in  a precofibration category are the cylinder structure maps $i_0, i_1$ for any cylinder $IA$ on a cofibrant object $A$
\[
\xymatrix { A \Sum A \ar@{>->}[r]^-{i_0+i_1} & IA \ar[r]^{p}_{\sim} & A }
\]
for $p \eps {\cat W}$ implies that $i_0, i_1$ have the same image in ${\bf ho}{\cat M}$.

\begin{lem}
\label{lem::weakequivcancelation}
\mbox{}
\begin{enumerate}
\item
Let $t \colon {\cat M}^{'} \ra {\cat M}$ be a cofibrant approximation of a precofibration category. For any maps $f, g \colon A \ra B$ of ${\cat M}^{'}$ and weak equivalence $b \colon tB \ra B^{'}$ of ${\cat M}$ with $b \comp tf = b \comp tg$, we have that $f \simeq g$.
\item
Let $t \colon {\cat M}^{'} \ra {\cat M}$ be a fibrant approximation of a prefibration category. For any maps $f, g \colon A \ra B$ of ${\cat M}^{'}$ and weak equivalence $a \colon A^{'} \ra tA$ of ${\cat M}$ with $tf \comp a = tg \comp a$, we have that $f \simeq g$.
\end{enumerate}
\end{lem}

\begin{proof}
We only prove (1). We may assume that $f + g \colon A \Sum A \ra B$ is a cofibration. Indeed, for general $f, g$ we factor $f + g$ as a cofibration $f^{'} + g^{'}$ followed by a weak equivalence $r$. The map $tr$ is a weak equivalence, and so is $b \comp tr$. If we proved that $f^{'} \simeq g^{'}$ then it would follow that $f \simeq g$.

So assume that $f + g \colon A \Sum A \ra B$ is a cofibration. Pick a cylinder $IA$, and construct step by step the commutative diagram below
\begin{center}
$\xymatrix { 
    tA \Sum tA \ar@{>->}[r]^-{ti_0 + ti_1} \ar@{>->}[d]_{tf + tg} & tIA \ar[rr]^{tp}_\sim \ar@{>->}[d]^{th} & & A \ar[d]^{b \comp tf = b \comp tg} \\
    tB \ar@{>->}[r]^{tb_1} & tB_1 \ar@{>->}[r]^{tb_2} & tB_2 \ar[r]^{b_3}_\sim & B^{'}    
} $
\end{center}
In this diagram, the bottom horizontal composition is $b \colon tB \ra B^{'}$. We construct $B_1 = B \Sum_{A \Sum A} IA$ with component maps $b_1$ and $h$. $tB_1$ is the pushout of the left square of our diagram, and using the pushout property we construct the map $tB_1 \ra B^{'}$. We then construct $b_3 \comp tb_2$ as the CFA4 factorization of $tB_1 \ra B^{'}$. 

The maps $b_1, b_2, h$ are cofibrations. By the 2 out of 3 axiom CF2, the maps $tb_2 \comp tb_1$ and therefore $b_2b_1$ are weak equivalences. Since $i_0 \simeq i_1$, we get $f \simeq g$.
\end{proof}

\begin{lem}
\label{lem:cleavagechoice}
\mbox{}
\begin{enumerate}
\item
Let $t \colon {\cat M}^{'} \ra {\cat M}$ be a cofibrant approximation of a precofibration category. For any commutative diagram
\begin{center}
$\xymatrix { 
    & tA^{'} \ar[dl]_{ti^{'}} \ar[dr]^{ti^{''}} & & \\
    tD^{'} \ar[drrr]^(.3){\sigma^{'}}  & & tD^{''} \ar[dr]^{\sigma^{''}} & \\
    & tB^{'} \ar[ul]^(.3){tj^{'}} \ar[ur]_(.3){tj^{''}} & & B     
} $
\end{center}
with $j^{'}, j^{''}$ weak equivalences of ${\cat M}^{'}$ and $\sigma^{'}, \sigma^{''}$ weak equivalences of ${\cat M}$ we have that $j^{'-1} i^{'} = j^{''-1} i^{''}$ in ${\bf ho}{\cat M}^{'}$
\item
Let $t \colon {\cat M}^{'} \ra {\cat M}$ be a fibrant approximation of a prefibration category. For any commutative diagram
\begin{center}
$\xymatrix { 
    & tA^{'}  & & A \ar[dlll]_(.25){\tau^{'}}  \ar[dl]^(.3){\tau^{''}} \\
    tD^{'} \ar[ur]^(.7){tq^{'}} \ar[dr]_{tp^{'}}  & & tD^{''} \ar[ul]_(.7){tq^{''}} \ar[dl]^{tp^{''}} & \\
    & tB^{'} & &     
} $
\end{center}
with $q^{'}, q^{''}$ weak equivalences of ${\cat M}^{'}$ and $\tau^{'}, \tau^{''}$ weak equivalences of ${\cat M}$ we have that $ p^{'} q^{'-1} = p^{''} q^{''-1} $ in ${\bf ho}{\cat M}^{'}$
\end{enumerate}
\end{lem}

\begin{proof}
We only prove (1). We may assume that $j^{'}, j^{''}$ are trivial cofibrations. Indeed, in the case of $j^{'}$ suppose that $j^{'} = p^{'} j^{'}_1$ is a Brown factorization with $s^{'}$ a right inverse of $p^{'}$. Then we may replace $i^{'}$, $j^{'}$, $\sigma^{'}$ with $s^{'}i^{'}$, $j^{'}_1$ resp. $\sigma^{'} p^{'}$.

Suppose now that $j^{'}, j^{''}$ are trivial cofibrations. Construct the sum $D^{'} \Sum_{B^{'}} D^{''}$ with component maps the trivial cofibrations $h^{'}$ and $h^{''}$. In the diagram below
\begin{center}
$\xymatrix { 
    & tA^{'} \ar@{>->}[dl]_{ti^{'}} \ar@{>->}[dr]^{ti^{''}} & & \\
    tD^{'} \ar@{>->}[r]^(.4){th^{'}}  & t(D^{'} \Sum_{B^{'}} D^{''}) \ar[drr]^(.7){\sigma} & tD^{''} \ar@{>->}[l]_(.4){th^{''}} & \\
    & tB^{'} \ar@{>->}[ul]^(.3){tj^{'}} \ar@{>->}[ur]_(.3){tj^{''}} & & B     
} $
\end{center}
the bottom triangle is a pushout by CFA3. The map $\sigma$ exists by the universal property of the pushout, since $\sigma^{'} \comp tj^{'} = \sigma^{''} \comp tj^{''}$, and is a weak equivalence by the 2 out of 3 axiom. \lemmaref{lem::weakequivcancelation} (1) applied to $h^{'}i^{'}, h^{''}i^{''}, \sigma$ implies that $h^{'}i^{'} \simeq h^{''}i^{''}$ in ${\cat M}^{'}$, and we conclude that $j^{'-1} i^{'} = j^{''-1} i^{''}$ in ${\bf ho}{\cat M}^{'}$.
\end{proof}

\begin{thm}
\label{thm:cofapproxareleftapprox}
\mbox{}
\begin{enumerate}
\item
Cofibrant approximation functors are left approximations.
\item
Fibrant approximation functors are right approximations.
\end{enumerate}
\end{thm}

\begin{proof}
We only prove (1). If $t \colon {\cat M}^{'} \ra {\cat M}$ is a cofibrant approximation of a precofibration category ${\cat M}$, then $t$ sends weak equivalences to weak equivalences by CFA2. Axiom LAP1 folows from CFA4. To prove the axiom LAP2, , use CFA4 to construct a factorization $f + p \colon tA^{'} \Sum tB^{'} \ra B$ as a cofibration $ti^{'} + tj^{'}$ followed by a weak equivalence $\sigma$. Since $\sigma \comp tj^{'} = p$, the maps $tj^{'}$ and therefore $j^{'}$ are weak equivalences. Axiom LAP2 (2) is proved by \lemmaref{lem:cleavagechoice}.
\end{proof}

As corollaries of the Approximation \theoremref{thm:approxthm}, we note:
\begin{thm}[Anderson]
\label{thm::catcofobjequiv}
\mbox{}
\begin{enumerate}
\item
Given a precofibration category ${\cat M}$, the inclusion $i_{\cat M} \colon {\cat M}_{cof} \ra {\cat M}$ induces an equivalence of categories ${\bf ho}i_{\cat M} \colon {\bf ho}{\cat M}_{cof} \ra {\bf ho}{\cat M}$
\item
Given a prefibration category ${\cat M}$, the inclusion $j_{\cat M} \colon {\cat M}_{fib} \ra {\cat M}$ induces an equivalence of categories ${\bf ho}j_{\cat M} \colon {\bf ho}{\cat M}_{fib} \ra {\bf ho}{\cat M}$
\end{enumerate}
\end{thm}

More generally we have:

\begin{thm}
\label{thm::homotopycatequiv}
\mbox{}
\begin{enumerate}
\item
A cofibrant approximation $t \colon {\cat M}^{'} \ra {\cat M}$ of a precofibration category induces an equivalence of categories ${\bf ho}{\cat M}^{'} \ra {\bf ho}{\cat M}$.
\item
A fibrant approximation $t \colon {\cat M}^{'} \ra {\cat M}$ of a prefibration category induces an equivalence of categories ${\bf ho}{\cat M}^{'} \ra {\bf ho}{\cat M}$.
\end{enumerate}
\end{thm}

It should be noted that the last theorem is actually a particular case of an even more general result of Cisinski, for which we refer the reader to \cite{Cisinski2}, 3.12.

\begin{rem}
\label{rem:intermediatecofcat}
If $t \colon {\cat M}^{'} \ra {\cat M}$ is a cofibrant approximation of a precofibration category ${\cat M}$, suppose that $i \colon {\cat M}_1 \hookrightarrow {\cat M}$ is a subcategory that includes the image of $t$, with weak equivalences and cofibrations induced from ${\cat M}$. Then $i$ as well as the corestriction $t_1 \colon {\cat M}^{'} \ra {\cat M}_1$ of $t$ both define cofibrant approximations. By \theoremref{thm::homotopycatequiv} both functors ${\bf ho}i$, ${\bf ho}t_1$ are equivalences of categories.
\end{rem}

\begin{rem}
\label{rem:intermediatecofcat2}
Suppose that $t \colon {\cat M}^{'} \ra {\cat M}$ is a functor between precofibration categories such that $t$ restricted to ${\cat M}^{'}_{cof}$ is a cofibrant approximation. In view of \theoremref{thm::catcofobjequiv} and \theoremref{thm::homotopycatequiv}, it is not hard to see that $t$ induces a composite equivalence of categories ${\bf ho}{\cat M}^{'} \la {\bf ho}{\cat M}^{'}_{cof} \ra {\bf ho}{\cat M}$. The proper way to formulate this is to say that the {\it total left derived} functor of t is an equivalence, which we will prove as \theoremref{thm:existenceleftderivedequiv} in the next section.
\end{rem}

\section{Total derived functors in cofibration categories}
\label{sec:leftrightderivedcofcat}

The following result describes a sufficient condition for the existence of a total left resp. right derived functor:

\begin{thm}
\label{thm:existenceleftderived}
Let $({\cat M}_2, {\cat W}_2)$ be a category with weak equivalences.
\begin{enumerate}
\item
If ${\cat M}_1$ is a precofibration category and $u \colon {\cat M}_1 \ra {\cat M}_2$ is a functor that sends trivial cofibrations between cofibrant objects to weak equivalences, then $u$ admits a total left derived functor $({\bf L}u, \eps)$. The natural map $\eps \colon ({\bf L}u)(A) \Ra uA$ is an isomorphism for $A$ cofibrant.
\item
If ${\cat M}_1$ is a prefibration category and $u \colon {\cat M}_1 \ra {\cat M}_2$ is a functor that sends trivial fibrations between fibrant objects to weak equivalences, then $u$ admits a total right derived functor $({\bf R}u, \eta)$. The natural map $\eta \colon uA \Ra ({\bf R}u)(A)$ is an isomorphism for $A$ fibrant.
\end{enumerate}
\end{thm}

More generally:
\begin{thm}
\label{thm:existenceleftderivedapprox}
Let $({\cat M}_2, {\cat W}_2)$ be a category with weak equivalences.
\begin{enumerate}
\item
If $t \colon {\cat M}^{'}_1 \ra {\cat M}_1$ is a cofibrant approximation of a precofibration category ${\cat M}_1$ and $u \colon {\cat M}_1 \ra {\cat M}_2$ is a functor such that $ut$ sends trivial cofibrations to weak equivalences, then $u$ admits a total left derived functor $({\bf L}u, \eps)$. The natural map $\eps \colon ({\bf L}u)(tA) \Ra utA$ is an isomorphism for objects $A$ of ${\cat M}^{'}_1$.
\item
If $t \colon {\cat M}^{'}_1 \ra {\cat M}_1$ is a fibrant approximation of a prefibration category ${\cat M}_1$ and $u \colon {\cat M}_1 \ra {\cat M}_2$ is a functor that sends trivial fibrations to weak equivalences, then $u$ admits a total right derived functor $({\bf R}u, \eta)$. The natural map $\eta \colon utA \Ra ({\bf R}u)(tA)$ is an isomorphism for objects $A$ of ${\cat M}^{'}_1$.
\end{enumerate}
\end{thm}

\begin{proof}
We only prove (1). We may assume that ${\cat W}_2$ is saturated. The composition $ut$ sends trivial cofibrations to weak equivalences, therefore using the Brown Factorization \lemmaref{lem:brownfact} $ut$ sends weak equivalences to weak equivalences since ${\cat W}_2$ is saturated. The result follows from \theoremref{thm:generalexistenceleftderived} applied to the cofibrant approximation $t$ and the functor $u$.
\end{proof}

And the next result describes a sufficient condition for the total left (resp. right) derived functor to be an equivalence of categories:

\begin{thm}
\label{thm:existenceleftderivedequiv}
\mbox{}
\begin{enumerate}
\item
If $t \colon {\cat M}_1 \ra {\cat M}_2$ is a functor between precofibration categories such that its restriction $({\cat M}_1)_{cof} \ra {\cat M}_2$ is a cofibrant approximation, then $t$ admits a total left derived functor $({\bf L}t, \eps)$ and ${\bf L}t$ is an equivalence of categories. The natural map $\eps \colon ({\bf L}t)(A) \Ra tA$ is an isomorphism for $A$ cofibrant.
\item
If $t \colon {\cat M}_1 \ra {\cat M}_2$ is a functor between prefibration categories such that its restriction $({\cat M}_1)_{fib} \ra {\cat M}_2$ is a fibrant approximation, then $t$ admits a total right derived functor $({\bf R}t, \eps)$ and ${\bf R}t$ is an equivalence of categories. The natural map $\eta \colon tA \Ra ({\bf R}t)(A)$ is an isomorphism $A$ fibrant.
\end{enumerate}
\end{thm}

\begin{proof}
For part (1), denote $i_{{\cat M}_1} \colon ({\cat M}_1)_{cof} \ra {\cat M}_1$ the inclusion. The functors $i_{{\cat M}_1}$, $ti_{{\cat M}_1}$ are cofibrant approximations and induce equivalences of categories ${\bf ho}i_{{\cat M}_1}$, ${\bf ho}(ti_{{\cat M}_1})$ by \theoremref{thm::homotopycatequiv}. 

From \theoremref{thm:existenceleftderivedapprox} applied to the cofibrant approximation $i_{{\cat M}_1}$ followed by $t$ we see that $t$ admits a total left derived functor $({\bf L}t, \eps)$. Furthermore, $\eps \colon ({\bf L}t)(A) \Ra tA$ is an isomorphism for $A$ cofibrant, therefore ${\bf L}t {\bf ho}i_{{\cat M}_1} = {\bf ho}(ti_{{\cat M}_1})$ and ${\bf L}t$ is an equivalence of categories since ${\bf ho}i_{{\cat M}_1}$ and ${\bf ho}(ti_{{\cat M}_1})$ are equivalences.
\end{proof}

We will also introduce in the context of (co)fibration categories the notion of Quillen adjoint functors and of Quillen equivalences.

\index{Quillen!adjunction}
\index{Quillen!equivalence}
\begin{defn}
\label{defn:quillenadjoint}
Consider four functors
\begin{center}
$\xymatrix {{\cat M}^{'}_1 \ar[r]^{t_1} & {\cat M}_1 \ar@<2pt>[r]^{u_1} & {\cat M}_2 \ar@<2pt>[l]^{u_2} & {\cat M}^{'}_2 \ar[l]_{t_2} }$
\end{center}
where 
\begin{enumerate}
\item
$u_1 \dashv u_2$ is an adjoint pair.
\item
$t_1$ is a cofibrant approximation of a precofibration category ${\cat M}_1$ and $t_2$ is a fibrant approximation of a prefibration category ${\cat M}_2$
\item
$u_1t_1$ sends trivial cofibrations to weak equivalences and $u_2t_2$ sends trivial fibrations to weak equivalences
\end{enumerate}

We then say that $u_1, u_2$ is a {\it Quillen adjoint} pair with respect to $t_1$ and $t_2$. If additionally 
\begin{enumerate}
\item[(4)]
For any objects $A^{'} \eps {\cat M}^{'}_1$, $B^{'} \eps {\cat M}^{'}_2$, a map $u_1t_1A^{'} \ra t_2B^{'}$ is a weak equivalence iff its adjoint $t_1A^{'} \ra u_2t_2B^{'}$ is a weak equivalence
\end{enumerate}
we say that $u_1, u_2$ is a {\it Quillen equivalence} pair with respect to $t_1$ and $t_2$.
\end{defn}

If the functors $t_1, t_2$ are implied by the context, we may refer to $u_1, u_2$ as a Quillen pair of adjoint functors (resp. equivalences) without direct reference to $t_1$ and $t_2$.

\begin{thm}[Quillen adjunction]
\label{thm:quillenadjoint}
\mbox{}
\begin{enumerate}
\item
A Quillen adjoint pair $u_1 \dashv u_2$ with respect to $t_1, t_2$
\begin{center}
$\xymatrix {{\cat M}^{'}_1 \ar[r]^{t_1} & {\cat M}_1 \ar@<2pt>[r]^{u_1} & {\cat M}_2 \ar@<2pt>[l]^{u_2} & {\cat M}^{'}_2 \ar[l]_{t_2} }$
\end{center}
induces a pair of adjoint functors ${\bf L}u_1 \dashv {\bf R}u_2$
\begin{center}
${\bf L}u_1 : {\bf ho}{\cat M}_1 \rightleftarrows {\bf ho}{\cat M}_2 : {\bf R}u_2$
\end{center}
\item
If additionally $u_1, u_2$ satisfy (4l) (resp. (4r)) of \theoremref{thm:generalexistencetotalderivedadjoint} with respect to $t_1, t_2$, then ${\bf L}u_1$ (resp. ${\bf R}u_2$) are fully faithful.
\item
If $u_1, u_2$ is a pair of Quillen equivalences with respect to $t_1, t_2$ then
\begin{center}
${\bf L}u_1 : {\bf ho}{\cat M}_1 \rightleftarrows {\bf ho}{\cat M}_2 : {\bf R}u_2$
\end{center}
is a pair of equivalences of categories.
\end{enumerate}
\end{thm}

\begin{proof}
This is a corollary of \theoremref{thm:generalexistencetotalderivedadjoint}.
\end{proof}

We leave it to the reader to formulate the definition of abstract Quillen partially adjoint functors in the context of cofibration and fibration categories, and to state the analogue in this context of \theoremref{thm:generalexistencetotalderivedadjoint2} 


\section{Homotopic maps}
\label{chap:homotopicmaps}
\index{homotopic maps}
We now turn to a more detailed study of homotopic maps in precofibration categories.

Start with ${\cat M}$ a precofibration category, and let $f,g \colon A \ra B$ be two maps with $A, B$ cofibrant. A {\it left homotopy} from $f$ to $g$ is a commutative diagram
\begin{equation}
\label{eqn:lefthomotopy}
\xymatrix{
    A \Sum A \ar[r]^-{f + g} \ar@{>->}[d]_{i_0 + i_1} & B \ar@{>->}[d]^b_{\sim} \\
    IA \ar[r]^H & B^{'}
    }
\end{equation}
with $IA$ a cylinder of $A$ and $b$ a trivial cofibration. We thus have that $Hi_0 = bf$ and $Hi_1 = bg$. The map $H$ is called the {\it left homtopy map} between $f$ and $g$, and we say that the left homotopy goes through the cylinder $IA$ and through the trivial cofibration $b$ . We write $f \simeq_l g$ to say that $f, g$ are left homotopic.

A left homotopy $f \simeq_l g$ with $B = B^{'}$ and $b = 1_B$ is called {\it strict}. We should be careful to point out that Brown uses the notation $\simeq$ differently - to denote what we call strict left homotopy.

Clearly if $f \simeq_l g$ then $f \simeq g$. Our goal will be to show that the notions of homotopy and left homotopy coincide (\theoremref{thm:homotmaps}).

Here is the dual setup for a prefibration category ${\cat M}$. Suppose that $f,g \colon B \ra A$ be two maps with $A, B$ fibrant. A {\it right homotopy} $f \simeq_r g$ is a commutative diagram
\begin{equation}
\label{eqn:righthomotopy}
\xymatrix{
    B^{'} \ar[r]^H \ar@{->>}[d]_b^{\sim} & A^I \ar@{->>}[d]^{(p_0,\,p_1)} \\
    B \ar[r]^-{(f,\,g)} & A \times A
    }
\end{equation}
with $A^I$ a path object of $A$ and $b$ a trivial fibration. A {\it strict} right homotopy additionally satisfies $B = B^{'}$ and $b = 1_B$.

\begin{thm}[Brown, \cite{Brown}]
\label{thm:homotmaps}
\mbox{}
\begin{enumerate}
\item
In a precofibration category ${\cat M}$, two maps $f, g \colon A \ra B$ with $A, B$ cofibrant are left homotopic iff they are homotopic.
\item
In a prefibration category ${\cat M}$, two maps $f, g \colon A \ra B$ with $A, B$ fibrant are right homotopic iff they are homotopic. 
\end{enumerate}
\end{thm}

\begin{proof}
We only prove part (1). Given a left homotopy \equationref{eqn:lefthomotopy}, since $i_0 \simeq i_1$ we have $Hi_0 \simeq Hi_1$, therefore $bf \simeq bg$ and so $f \simeq g$ since $b$ is a weak equivalence.

If $f \simeq g$ are homotopic in ${\cat M}$ on the other hand, by \theoremref{thm::catcofobjequiv} we see that they are homotopic inside ${\cat M}_{cof}$. From  \theoremref{thm::descrhomot} (c) below it will follow that $f \simeq_l g$.
\end{proof}

Let us now work our way to complete the proof of \theoremref{thm:homotmaps}.

Suppose that $A$ is a cofibrant object in the precofibration category ${\cat M}$. If $IA$ is a cylinder of $A$, a {\it refinement} of $IA$ consists of a cylinder $I^{'}A$ and a trivial cofibration $j \colon IA \ra I^{'}A$ such that the diagram below commutes
\begin{center}
$\xymatrix{
    & A \Sum A \ar@{>->}[dl]_{i_0 + i_1} \ar@{>->}[dr]^{i^{'}_0 + i^{'}_1} & \\
    IA \ar[dr]_p^\sim \ar@{>->}[rr]^j & & I^{'}A \ar[dl]^{p^{'}}_\sim \\
    & A &
  }$
\end{center}

If $f \simeq_l g \colon A \ra B$ through $I^{'}A$ with homotopy map $H^{'}$, then $H^{'}j$ defines a homotopy through $IA$.

Note that given any two cylinders $IA$, $I^{'}A$, we can construct a common refinement $I^{''}A$ by factoring $IA \Sum_{A \Sum A} I^{'}A \ra A$ as a cofibration $IA \Sum_{A \Sum A} I^{'}A \ra I^{''}A$ followed by a weak equivalence $I^{''}A \ra A$.

This allows us to prove the following lemma:

\begin{lem}
\label{lem::choosecyl}
\mbox{}
\begin{enumerate}
\item
Let ${\cat M}$ be a precofibration category, and $f \simeq_l g \colon A \ra B$ with $A, B$ cofibrant. Let $IA$ be a cylinder of $A$. Then one can construct a homotopy $f \simeq_l g$ through $IA$.
\item
Let ${\cat M}$ be a prefibration category, and $f \simeq_r g \colon B \ra A$ with $A, B$ fibrant. Let $A^I$ be a path object for $A$. Then one can construct a homotopy $f \simeq_r g$ through $A^I$.
\end{enumerate}
\end{lem}

\begin{proof}
We only prove (1). Assume that there exists a homotopy $f \simeq_l g$ through another cylinder $I^{'}A$. Construct a common refinement $I^{''}A$ of $IA$ and $I^{'}A$. To prove (1), it suffices to construct a homotopy through the refinement $I^{''}A$. In the commutative diagram below
\begin{center}
$\xymatrix{
    A \Sum A \ar[r]^{f + g} \ar@{>->}[d]_{i^{'}_0 + i^{'}_1} & B \ar@{>->}[d]^{b}_{\sim} \\
    I^{'}A \ar[r]^{H^{'}} \ar@{>->}[d]_j^\sim & B^{'} \ar@{>->}[d]^{b^{'}}_\sim \\
    I^{''}A \ar[r]^{H^{''}} & B^{''}
    }$
\end{center}
$H^{'}, B^{'}$ and $b$ define the homotopy $f \simeq_l g$, and $j$ is the cylinder refinement map from $I^{'}A$ to $I^{''}A$. We construct $B^{''}$ as the pushout of $j$ and $H^{'}$, and we have the desired homotopy $H^{''}, B^{''}$ and $b^{'}b$ through $I^{''}A$.
\end{proof}

We can now prove
\begin{thm}
\label{thm::homotopyequivrel}
\mbox{}
\begin{enumerate}
\item
If ${\cat M}$ is a precofibration category, then $\simeq_l$ is an equivalence relation in ${\cat M}_{cof}$. Furthermore if $f \simeq_l g \colon A \ra B$ with $A, B$ cofibrant then
\begin{enumerate}
\item
If $h \colon B \ra C$ with $C$ cofibrant then $hf \simeq_l hg$
\item
If $h \colon C \ra A$ with $C$ cofibrant then $fh \simeq_l gh$
\end{enumerate}
\item
If ${\cat M}$ is a prefibration category, then $\simeq_r$ is an equivalence relation in ${\cat M}_{fib}$. Furthermore if $f \simeq_r g \colon A \ra B$ with $A, B$ fibrant then
\begin{enumerate}
\item
If $h \colon B \ra C$ with $C$ fibrant then $hf \simeq_r hg$
\item
If $h \colon C \ra A$ with $C$ fibrant then $fh \simeq_r gh$
\end{enumerate}
\end{enumerate}
\end{thm}

\begin{proof}
We only prove (1). Clearly $\simeq_l$ is symmetric and reflexive. 

To see that $\simeq_l$ is transitive, assume $f \simeq_l g \simeq_l h \colon A \ra B$. From the previous lemma, we may assume that both homotopies go through the same cylinder $IA$. Denote these homotopies $H_1, B_1^{'}, b_1$ and $H_2, B_2^{'}, b_2$. Taking the pushout of $b_1$ and $b_2$, we obtain homotopies $H_1^{'}, B^{'}, b$ and $H_2^{'}, B^{'}, b$. 

Notice that both diagrams below are pushouts
\begin{center}
$\xymatrix{
    A \ar@{>->}[r]_\sim^{i_0} \ar@{>->}[d]^\sim_{i_1} & IA \ar@{>->}[d]^\sim & & A \Sum A \Sum A \,\,\, \ar@{>->}[rr]^-{A \Sum (i_0 + i_1)} \ar@{>->}[d]_-{(i_0 + i_1) \Sum A} & & A \Sum IA \ar@{>->}[d]^-{i_0 \Sum IA} \\
    IA \ar@{>->}[r]^-\sim & \,\,\, IA \Sum_A IA & & IA \Sum A \,\,\, \ar@{>->}[rr]^-{IA \Sum i_1} & & IA \Sum_A IA
    }$
\end{center}
The factorization $\nabla \colon \xymatrix{A \Sum A \ar@{>->}[r]^-{i_0 \Sum i_1} & IA \Sum_A IA \ar[r]^-{p + p}_-\sim & A}$ is a cylinder: the map $p+p$ is a weak equivalence because of the first diagram, and the map $i_0 \Sum i_1$ is a cofibration as seen if we precompose the second diagram with the cofibration $(i_0, i_2) \colon A \Sum A \ra  A \Sum A \Sum A$. The commutative diagram below then defines a homotopy from $f$ to $h$.
\begin{center}
$\xymatrix{
    A \Sum A \ar[r]^{f + h} \ar@{>->}[d]_{i_0 \Sum i_1} & B \ar@{>->}[d]^{b}_{\sim} \\
    IA \Sum_A IA \ar[r]^-{H^{'}_1 + H^{'}_2} & B^{'}
    }$
\end{center}

To prove (a), let $IA, H, B^{'}, b$ define a homotopy $f \simeq_l g$. In the diagram below let $C^{'}$ be the pushout of $b, h$.
\begin{center}
$\xymatrix{
    A \Sum A \ar[r]^{f + g} \ar@{>->}[d]_{i_0 \Sum i_1} & B \ar[r]^h \ar@{>->}[d]^b_\sim & C \ar@{>->}[d]^c_\sim \\
    IA \ar[r]^H & B^{'} \ar[r]^{h^{'}} & C^{'}
    }$
\end{center}
The outer rectangle defines a homotopy $hf \simeq_l hg$.

For (b), use \lemmaref{lem:existencerelcyl} to construct relative cylinders $IC, IA$ along $h$. From \lemmaref{lem::choosecyl}, we can construct a homotopy $f \simeq_l g$ through $IA$. Precomposing this homotopy with $IC \ra IA$ yields a homotopy $fh \simeq_l gh$.
\end{proof}

For a precofibration category ${\cat M}$, if we factor ${\cat M}_{cof}$ modulo $\simeq_l$ we obtain a category $\pi_l {\cat M}_{cof}$, with same objects as ${\cat M}_{cof}$. By \theoremref{thm::homotopyequivrel}, the morphisms of $\pi_l {\cat M}_{cof}$ are given by $Hom_{\pi_l {\cat M}_{cof}}(A,B) = Hom_{{\cat M}_{cof}}(A,B)_{/ \simeq_l}$. We define weak equivalences in $\pi_l {\cat M}_{cof}$ to be homotopy classes of maps that have one (and hence all) representatives weak equivalece maps of ${\cat M}$.

Of course, in view of \theoremref{thm:homotmaps} ultimately $\pi_l {\cat M}_{cof} \cong \pi {\cat M}_{cof}$. 

For a prefibration category, $\pi_r {\cat M}_{fib}$ denotes the factorization of ${\cat M}_{fib}$ modulo $\simeq_r$. Weak equivalences in $\pi_r {\cat M}_{fib}$ are by definition homotopy classes of maps that have one (and hence all) representatives weak equivalence maps of ${\cat M}$.

\section{Homotopy calculus of fractions}
\label{chap:homotopycalcoffractions}

We will show that for a precofibration category ${\cat M}$, the category $\pi_l {\cat M}_{cof}$ admits a calculus of left fractions in the sense of Gabriel-Zisman with respect to weak equivalences. A nice way to phrase this is to say that ${\cat M}_{cof}$ admits a {\it homotopy} calculus of left fractions. Dually, given a prefibration category ${\cat M}$, its category of fibrant objects ${\cat M}_{fib}$ admits a homotopy calculus of right fractions.

We say that a category pair $({\cat M}$, ${\cat W})$ satisfies the 2 out of 3 axiom provided that for any composable morphisms $f, g$ of ${\cat M}$, if two of $f, g, gf$ are in ${\cat W}$ then so is the third. Weak equivalences in a precofibration category ${\cat M}$ satisfy the 2 out of 3 axiom, and so do weak equivalences in $\pi_l {\cat M}_{cof}$. 

For category pairs $({\cat M}$, ${\cat W})$ satisfying the 2 out of 3 axiom, the Gabriel-Zisman calculus of fractions takes a simplified form that we recall below. The general case of calculus of fractions - when weak equivalences do not necessarily satisfy the 2 out of 3 axiom - is described in \cite{Gabriel-Zisman} at pag. 12.

\index{calculus of fractions}
\begin{thm}[Simplified calculus of left fractions]
\label{thm::gabrielzisman}
Suppose that $({\cat M}, {\cat W})$ is a category pair satisfying the following:
\begin{enumerate}
\item[(a)]
The 2 out of 3 axiom
\item[(b)]
Any diagram of solid maps with $a \eps {\cat W}$
\begin{equation}
\label{eqn:gabrielzismanhypothesis1}
\xymatrix{
  A \ar[r] \ar[d]_a^\sim & B \ar@{-->}[d]^b_\sim \\
  A^{'} \ar@{-->}[r] & B^{'}
}
\end{equation}
extends to a commutative diagram with $b \eps {\cat W}$
\item[(c)]
For any maps $A^{'} \overset{s}{\underset{\sim}{\ra}} A \overset{f}{\underset{g}{\dbra}} B$ with $fs = gs$ and $s \eps {\cat W}$, there exists $B \overset{s^{'}}{\underset{\sim}{\ra}} B^{'}$ in ${\cat W}$ with $s^{'}f = s^{'}g$.
\end{enumerate}
Then:
\begin{enumerate}
\item
Each map in $Hom_{{\bf ho}{\cat M}} (A, B)$ can be written as a left fraction $s^{-1}f$
\begin{center}
$\xymatrix{ A \ar[r]^f & B^{'} & B \ar[l]_-s^-\sim
  }$
\end{center}
with $s$ a weak equivalence. 
\item
Two fractions $s^{-1}f$, $t^{-1}g$ are equal in ${\bf ho}{\cat M}$ if and only if there exist weak equivalences $s^{'}, t^{'}$ as in the diagram below
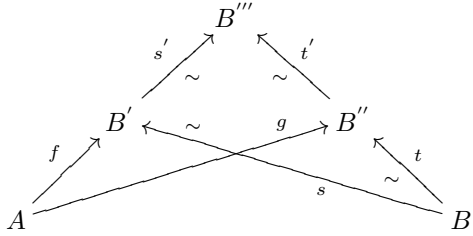
\begin{equation}
\label{eqn:descrhomotopy1}
\xymatrix{ 
    & & B^{'''} & & \\
    & B^{'} \ar[ur]_\sim^{s^{'}} & &  B^{''} \ar[ul]^\sim_{t^{'}} & \\
    A \ar[ur]^f \ar[urrr]^(.8)g & & & & B \ar[ulll]^(.4)s_(.8)\sim \ar[ul]^\sim_t
  }
\end{equation}
so that $s^{'}s = t^{'}t$ and $s^{'}f = t^{'}g$.
\end{enumerate}
If furthermore weak equivalences are left cancellable, in the sense that for any pair of maps $f, g \colon A \ra B$ and weak equivalence $h \colon B \ra B^{'}$ with $hf = hg$ we have $f = g$, then
\begin{enumerate}
\item[(3)]
Two maps $f, g \colon A \ra B$ are equal in ${\bf ho}{\cat M}$ if and only if $f = g$.
\end{enumerate}
\end{thm}

The dual result for right fractions is
\begin{thm}[Simplified calculus of right fractions]
\label{thm::gabrielzisman2}
Suppose that $({\cat M}, {\cat W})$ is a category pair satisfying the folowing:
\begin{enumerate}
\item[(a)]
The 2 out of 3 axiom
\item[(b)]
Any diagram of solid maps with $a \eps {\cat W}$
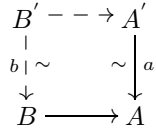
\begin{equation}
\label{eqn:gabrielzismanhypothesis2}
\xymatrix{
  B^{'} \ar@{-->}[r] \ar@{-->}[d]_b^\sim & A^{'} \ar[d]^a_\sim  \\
  B \ar[r] & A
}
\end{equation}
extends to a commutative diagram with $b \eps {\cat W}$
\item[(c)]
For any maps $A \overset{f}{\underset{g}{\dbra}} B \overset{t}{\underset{\sim}{\ra}} B^{'}$ with $tf = tg$ and $t \eps {\cat W}$, there exists $A^{'} \overset{t^{'}}{\underset{\sim}{\ra}} A$ in ${\cat W}$ with $ft^{'} = gt^{'}$.
\end{enumerate}
Then:
\begin{enumerate}
\item
Each map in $Hom_{{\bf ho}{\cat M}} (A, B)$ can be written as a right fraction $fs^{-1}$
\begin{center}
$\xymatrix{ A & A^{'} \ar[l]_s^\sim \ar[r]^f & B
  }$
\end{center}
with $s$ a weak equivalence. 
\item
Two fractions $fs^{-1}$, $gt^{-1}$ are equal in ${\bf ho}{\cat M}$ if and only if there exist weak equivalences $s^{'}, t^{'}$ as in the diagram below
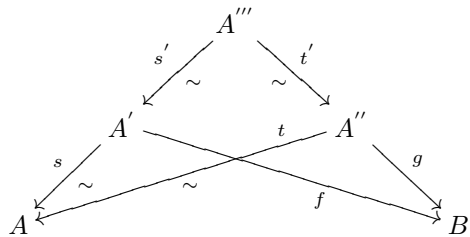
\begin{equation}
\label{eqn:descrhomotopy2}
\xymatrix{ 
    & & A^{'''} \ar[dl]^\sim_{s^{'}} \ar[dr]_\sim^{t^{'}} & & \\
    & A^{'} \ar[dl]^\sim_s \ar[drrr]_(.6)f & &  A^{''} \ar[dlll]^\sim_(.2)t \ar[dr]^g & \\
    A & & & & B
  }
\end{equation}
so that $ss^{'} = tt^{'}$ and $fs^{'} = gt^{'}$.
\end{enumerate}
If furthermore weak equivalences are right cancellable, in the sense that for any pair of maps $f, g \colon A \ra B$ and weak equivalence $h \colon A^{'} \ra A$ with $fh = gh$ we have $f = g$, then
\begin{enumerate}
\item[(3)]
Two maps $f, g \colon A \ra B$ are equal in ${\bf ho}{\cat M}$ if and only if $f = g$.
\end{enumerate}
\end{thm}

We only need to supply a 

\begin{proof}[Proof of \theoremref{thm::gabrielzisman}]
We construct a category ${\cat C}$ with objects $Ob {\cat M}$, and maps defined in terms of fractions as explained below.

Fix two objects $A$ and $B$. Consider the set of fractions $s^{-1}f$ as in (1). Denote $\sim$ the relation defined by (2) on the set of fractions $s^{-1}f$ from $A$ to $B$. The relation $\sim$ is clearly reflexive and symmetric. 

To see that the relation is transitive, assume ${s_1}^{-1}f_1 \sim {s_2}^{-1}f_2 \sim {s_3}^{-1}f_3$. We get a commutative diagram
\begin{center}
$\xymatrix{ 
    B^{''}_1 \ar[r]^{u_1} & B^{'''} & B^{''}_2 \ar[l]_{u_2} \\
    B^{'}_1 \ar[u]^{t_1} & B^{'}_2 \ar[ul]_{t_2} \ar[ur]^{t_3} & B^{'}_3 \ar[u]^{t_4} \\
    A \ar[u]^{f_1} \ar[ur]^(.4){f_2} \ar[urr]_(.3){f_3} && B \ar[u]_{s_3} \ar[ul]_(.4){s_2} \ar[ull]^(.3){s_1}
  }$
\end{center}
where the weak equivalences $t_1, t_2$ exist since ${s_1}^{-1}f_1 \sim {s_2}^{-1}f_2$, the weak equivalences $t_3, t_4$ exist since ${s_2}^{-1}f_2 \sim {s_3}^{-1}f_3$ and the weak equivalences $u_1, u_2$ exist from the hypothesis (b) applied to $t_2, t_3$. The compositions $u_1t_1, u_2t_4$ satisfy $u_1t_1f_1 = u_2 t_4 f_3$ and $u_1t_1s_1 = u_2 t_4 s_3$ which shows that ${s_1}^{-1}f_1 \sim {s_3}^{-1}f_3$, and we have proved that $\sim$ is transitive.

We let ${\cat C}(A, B)$ denote the set of fractions from $A$ to $B$ modulo the equivalence relation $\sim$. Given three objects $A, B, C$, we define composition ${\cat C}(A, B) \times {\cat C}(B, C) \ra {\cat C}(A, C)$ as follows. Given fractions $s^{-1}f$, $t^{-1}g$
\begin{center}
$\xymatrix{ 
    & & C^{''} & & \\
    & B^{'} \ar[ur]^{g^{'}} & &  C^{'} \ar[ul]^\sim_{t^{'}} & \\
    A \ar[ur]^f & & B \ar[lu]_s^\sim \ar[ru]_g & & C \ar[ul]^\sim_t
  }$
\end{center}
we use hypothesis (\ref{eqn:gabrielzismanhypothesis1}) to construct an object $C^{''}$, a map $g^{'}$ and a weak equivalence $t^{'}$ such that $g^{'}s = t^{'}g$, and then we define $t^{-1}g \comp s^{-1}f$ as $(t^{'}t)^{-1}(g^{'}f)$. The proof that the definition of composition does not depend on the choices involved uses hypotheses (b) and (c) (we leave this verification to the reader). Given an object $A$, the fraction $(1_A)^{-1}1_A$ is an identity element for the composition.

Define a functor $F \colon {\cat M} \ra {\cat C}$, by $F(A) = A$ and by sending $f \colon A \ra B$ to the fraction $F(f) = (1_B)^{-1}f$. It is not hard to see that $F$ is compatible with composition, and that if $s \colon A \ra B$ is a weak equivalence in ${\cat M}$ then $F(s)$ has $s^{-1}1_B$ as an inverse. 

Since $F$ sends weak equivalences to isomorphisms, it descends to a functor $\overline{F} \colon {\bf ho}{\cat M} \ra {\cat C}$, and it is straightforward to check that any other functor ${\bf ho}{\cat M} \ra {\cat C}^{'}$ factors uniquely through $\overline{F}$. It follows that the category ${\cat C}$ we constructed is equivalent to ${\bf ho}{\cat M}$, which ends the proof of (1) and (2).

If $f, g \colon A \ra B$ are equal in ${\bf ho}{\cat M}$, then by (2) there exists a weak equivalence $h \colon B \ra B^{'}$ such that $hf = hg$. If weak equivalences are left cancellable, then $f = g$, and this proves (3).
\end{proof}

Let us show that given a precofibration category ${\cat M}$, the category $\pi_l {\cat M}_{cof}$ satisfies the hypotheses of \theoremref{thm::gabrielzisman} up to homotopy.

\index{calculus of fractions!homotopic calculus of fractions}
\begin{thm}
\label{thm::homotopycalcoffrations}
\mbox{}
\begin{enumerate}
\item
In a precofibration category
\begin{enumerate}
\item
Any full diagram with cofibrant objects and weak equivalence $a$
\begin{center}
$\xymatrix{
  A \ar[r] \ar[d]_a^\sim & B \ar@{-->}[d]^b_\sim \\
  A^{'} \ar@{-->}[r] & B^{'}
}$
\end{center}
extends to a (strictly) homotopy commutative diagram with $b$ a weak equivalence and $B^{'}$ cofibrant
\item
For any $f,g \colon A \ra B$ with $A, B$ cofibrant
\begin{enumerate}
\item
If there is a weak equivalence $a \colon A^{'} \ra A$ with $A^{'}$ cofibrant such that $fa \simeq_l ga$, then $f \simeq_l g$
\item
If there is a weak equivalence $b \colon B \ra B^{'}$ with $B^{'}$ cofibrant such that $bf \simeq_l bg$, then $f \simeq_l g$
\end{enumerate}
\end{enumerate}
\item
In a prefibration category
\begin{enumerate}
\item
Any full diagram with fibrant objects and weak equivalence $a$
\begin{center}
$\xymatrix{
  B^{'} \ar@{-->}[r] \ar@{-->}[d]_b^\sim & A^{'} \ar[d]^a_\sim  \\
  B \ar[r] & A
}$
\end{center}
extends to a (strictly) homotopy commutative diagram with $b$ a weak equivalence and $B^{'}$ fibrant
\item
For any $f,g \colon A \ra B$ with $A, B$ fibrant
\begin{enumerate}
\item
If there is a weak equivalence $a \colon A^{'} \ra A$ with $A^{'}$ fibrant such that $fa \simeq_r ga$, then $f \simeq_r g$
\item
If there is a weak equivalence $b \colon B \ra B^{'}$ with $B^{'}$ fibrant such that $bf \simeq_r bg$, then $f \simeq_r g$
\end{enumerate}
\end{enumerate}
\end{enumerate}
\end{thm}

\begin{proof}
To prove (1) (a), denote $f \colon A \ra B$ and let $IA$ be a cylinder of $A$. The diagram
\begin{center}
$\xymatrix{
  A \ar[r]^f \ar[d]_a^\sim & B \ar[d]^{b} \\
  A^{'} \ar[r] & A^{'} \Sum_A IA \Sum_A B
}$
\end{center}
is strictly homotopy commutative. It remains to show that $b$ is a weak equivalence. 

Denote $F \colon IA \ra IA \Sum_A B$ the map induced by the $1^{st}$ component of the sum. In the pushout diagram
\begin{center}
$ \xymatrix {
        A \Sum A \ar[r]^{1 \Sum f} \ar@{>->}[d]_{i_0 + i_1} &
	A \Sum B \ar@{>->}[d] & \\
	IA \ar[r]^-F &
	IA \Sum_A B
	}$
\end{center}
the right vertical map is a cofibration and therefore $Fi_0 \colon A \ra IA \Sum_A B$ is a cofibration.

The map $b$ factors as $B \ra IA \Sum_A B \ra A^{'} \Sum_A IA \Sum_A B$. The first factor is a pushout of $i_1 \colon A \ra IA$, therefore a trivial cofibration. The second factor is a pushout of the weak equivalence $a$ by the cofibration $Fi_0$, therefore a weak equivalence by excision.

To prove (1) (b) (i), pick (\lemmaref{lem:existencerelcyl}) relative cylinders $IA^{'}, IA$ over $a$. By \lemmaref{lem::choosecyl}, there exists a homotopy $fa \simeq_l ga$ through $IA^{'}$ and through a trivial cofibration $b^{''}$. We get a commutative diagram
\begin{center}
$\xymatrix{
  A^{'} \Sum A^{'} \ar[r]^{a \Sum a}_\sim \ar@{>->}[d]_{i^{'}_0 + i^{'}_1} & A \Sum A \ar@{>->}[d]_{i_0 + i_1} \ar[r]^{f + g} & B \ar@{>->}[d]_\sim^{b^{''}} \\
  IA^{'} \ar[r]^-{h_1}_-\sim & IA^{'} \Sum_{A^{'} \Sum A^{'}} {A \Sum A} \ar[r]^-{h_2} \ar@{>->}[d]^j_\sim & B^{''} \ar@{>->}[d]^{b^{'}}_\sim \\
  & IA \ar[r]^H & B^{'}
}$
\end{center}
where the map $j$ is a cofibration because $IA^{'}, IA$ are relative cylinders. But $j$ is actualy a trivial cofibration. To see that, notice that since $a$ is a weak equivalence and $A, A^{'}$ are cofibrant, $a \Sum a$ is also a weak equivalence and by excision so is $h_1$. The map $IA^{'} \ra IA$ is a weak equivalence since $a$ is, and by the 2 out of 3 Axiom the map $j$ is a weak equivalence.

We define $B^{'}$ as the pushout of $j$ by $h_2$. The map $b^{'}$ is therefore a trivial cofibration. We let $b = b^{'}b^{''}$, and $IA, H, B^{'}, 1_{B^{'}}$ defines a homotopy $bf \simeq_l bg$.

Let us now prove (1) (b) (ii). Pick a homotopy $bf \simeq_l bg$ going through the cylinder $IA$ and through the trivial cofibration $b^{'}$. In the diagram below
\begin{center}
$\xymatrix{
    A \Sum A \ar[r]^-{f+g} \ar@{>->}[d]_{i_0 + i_1} & B \ar[rr]^-{b}_\sim \ar@{>->}[d]^{b_1} & & B^{'} \ar@{>->}[d]_\sim^{b^{'}} \\
    IA \ar[r]^{h_1} & B_1 \ar@{>->}[r]^{h_2} & B_2 \ar[r]^{h_3}_\sim & B^{''}
}$
\end{center}
construct $b_1$ as the pushout of $i_0 + i_1$ and $h_2, h_3$ as the factorization of $B_1 \ra B^{''}$ as a cofibration followed by a weak equivalence. Notice that $h_2b_1$ is a trivial cofibration, and we have constructed a homotopy $f \simeq_l g$ with homotopy map $h_2h_1$.

The proof of (2) is dual and is omitted.

\end{proof}

The Gabriel-Zisman left calculus of fractions applies therefore to the case of $\pi_l {\cat M}_{cof}$, if ${\cat M}$ is a precofibration category.

\begin{thm}[Brown's homotopy calculus of fractions, \cite{Brown}]
\label{thm::descrhomot}
\mbox{}
\begin{enumerate}
\item
Let ${\cat M}$ be a precofibration category, and $A, B$ be two cofibrant objects.
\begin{enumerate}
\item
Each map in $Hom_{{\bf ho}{\cat M}} (A, B)$ can be written as a left fraction $s^{-1}f$
\begin{center}
$\xymatrix{ A \ar[r]^f & B^{'} & B \ar[l]_-s^-\sim
  }$
\end{center}
with $s$ a weak equivalence and $B^{'}$ cofibrant. 
\item
Two such fractions $s^{-1}f$, $t^{-1}g$ are equal in ${\bf ho}{\cat M}$ if and only if there exist weak equivalences $s^{'}, t^{'}$ as in the diagram (\ref{eqn:descrhomotopy1}) with $B^{'''}$ cofibrant so that $s^{'}s \simeq_l t^{'}t$ and $s^{'}f \simeq_l t^{'}g$.
\item
Two maps $f, g \colon A \ra B$ are equal in ${\bf ho}{\cat M}$ if and only if they are homotopic $f \simeq_l g$
\end{enumerate}
\item
Let ${\cat M}$ be a prefibration category, and $A, B$ be two fibrant objects.
\begin{enumerate}
\item
Each map in $Hom_{{\bf ho}{\cat M}} (A, B)$ can be written as a right fraction $fs^{-1}$
\begin{center}
$\xymatrix{ A & A^{'} \ar[l]_s^\sim \ar[r]^f & B
  }$
\end{center}
with $s$ a weak equivalence and $A^{'}$ fibrant. 
\item
Two such fractions $fs^{-1}$, $gt^{-1}$ are equal in ${\bf ho}{\cat M}$ if and only if there exist weak equivalences $s^{'}, t^{'}$ as in the diagram (\ref{eqn:descrhomotopy2}) with $A^{'''}$ fibrant so that $ss^{'} \simeq_r tt^{'}$ and $fs^{'} \simeq_r gt^{'}$.
\item
Two maps $f, g \colon A \ra B$ are equal in ${\bf ho}{\cat M}$ if and only if $f \simeq_r g$.
\end{enumerate}
\end{enumerate}
\end{thm}

\begin{proof}
This is a consequence of \theoremref{thm::catcofobjequiv}, \theoremref{thm::gabrielzisman} and \theoremref{thm::homotopycalcoffrations}. 
\end{proof}

The proof of \theoremref{thm:homotmaps} is at this point complete, since \theoremref{thm::descrhomot} was its last prerequisite. From this point on we can freely write $\simeq$ instead of $\simeq_l$ and $\simeq_r$.

We can also prove a version \theoremref{thm::descrhomot} that describes ${\bf ho}{\cat M}$ in terms of fractions $s^{-1}f$ with $f, f+s$ cofibrations and $s$ a trivial cofibration: 
\begin{thm}
\label{thm::descrhomotcof}
\mbox{}
\begin{enumerate}
\item
Let ${\cat M}$ be a precofibration category, and $A, B$ be two cofibrant objects.
\begin{enumerate}
\item
Each map in $Hom_{{\bf ho}{\cat M}} (A, B)$ can be written as a left fraction $s^{-1}f$
\begin{center}
$\xymatrix{ A \ar@{>->}[r]^f & B^{'} & B \ar@{>->}[l]_-s^-\sim
  }$
\end{center}
with $f$, $f+s$ cofibrations and $s$ a trivial cofibration.
\item
Two fractions as in \theoremref{thm::descrhomot} (1) (a) $s^{-1}f$, $t^{-1}g$ with $s, t$ trivial cofibrations are equal in ${\bf ho}{\cat M}$ if and only if there exist trivial cofibrations $s^{'}, t^{'}$ as in the diagram (\ref{eqn:descrhomotopy1}) such that $s^{'}s = t^{'}t$ and $s^{'}f \simeq t^{'}g$.
\end{enumerate}
\item
Let ${\cat M}$ be a prefibration category, and $A, B$ be two fibrant objects.
\begin{enumerate}
\item
Each map in $Hom_{{\bf ho}{\cat M}} (A, B)$ can be written as a right fraction $fs^{-1}$
\begin{center}
$\xymatrix{ A & A^{'} \ar@{->>}[l]_s^\sim \ar@{->>}[r]^f & B
  }$
\end{center}
with $f$, $(f, s)$ fibrations and $s$ a trivial fibration
\item
Two fractions as in \theoremref{thm::descrhomot} (2) (a) $fs^{-1}$, $gt^{-1}$ with $s, t$ trivial fibrations are equal in ${\bf ho}{\cat M}$ if and only if there exist trivial fibrations $s^{'}, t^{'}$ as in the diagram (\ref{eqn:descrhomotopy2}) such that $ss^{'} = tt^{'}$ and $fs^{'} \simeq gt^{'}$.
\end{enumerate}
\end{enumerate}
\end{thm}

\begin{proof}
We only prove (1). Denote $\sim$ the equivalence relation defined by \theoremref{thm::descrhomot} (1) (b).

To prove (a) it suffices to show that any fraction $s^{-1}f$ with $s$ a weak equivalence is $\sim$ equivalent to a fraction $t^{-1}g$, with $g, g+t$ cofibrations and $t$ a trivial cofibration. Construct the commutative diagram
\begin{center}
$\xymatrix{
    & B^{'} & \\
    A \ar[ur]^f \ar@{>->}[dr]_{i_A} & B^{''} \ar[u]^\sim_{v} & B \ar[ul]_s^\sim \ar@{>->}[dl]^{i_B} \\
    & A \Sum B \ar@{>->}[u]_u &
  }$
\end{center}
where $vu$ is the factorization of $f + s$ as a cofibration followed by a weak equivalence. Define $g = ui_A$ and $t = ui_b$. The maps $g, g+t$ are cofibrations, the map $t$ is a trivial cofibration and the trivial map $v$ yields the desired $\sim$ equivalence between the fractions $s^{-1}f$ and $t^{-1}g$.

To prove (b), in the diagram below
\begin{center}
$\xymatrix{ 
    & & B^{'''} & & \\
    & B^{'} \ar[ur]_\sim^{s^{'}} & &  B^{''} \ar[ul]^\sim_{t^{'}} & \\
    A \ar[ur]^f \ar[urrr]^(.8)g & & & & B \ar[ulll]^(.4)s_(.8)\sim \ar[ul]^\sim_t
  }$
\end{center}
construct $s^{'}, t^{'}$ as the pushouts of $t, s$. We therefore have $s^{'}s = t^{'}t$. Since $s^{'}f$ and $t^{'}g$ are equal in ${\bf ho}{\cat M}$, we also get $s^{'}f \simeq t^{'}g$.
\end{proof}

Going back to the example of topological spaces $Top$ in \sectionref{sec:exampletopspaces}, recall that we have defined weak equivalences in $Top$ to be the homotopy equivalences given by the classic homotopy relation $\simeq$ in $Top$. We can now show that the definition of classic homotopy $\simeq$ in $Top$ is consistent with \definitionref{defn:homotopicmaps}. 

\begin{prop}
\label{prop:homotopyconsistence}
Two maps $f_0, f_1 \colon A \ra B$ in $Top$ have the same image in ${\bf ho}(Top)$ iff there exists a homotopy $h \colon A \times I \ra B$ which equals $f_k$ when restricted to $A \times {k}$, for $k = 0, 1$.
\end{prop}

\begin{proof}
If a homotopy $h$ exists, clearly $f_0, f_1$ have the same image in ${\bf ho}(Top)$. 

To prove the converse, we observe first that $\xymatrix{A \Sum A \ar[r]^{i_0 + i_1} & IA \ar[r]^p_\sim & A}$ is a cylinder with respect to the Hurewicz cofibration structure in $Top$. Using \theoremref{thm:homotmaps} and \lemmaref{lem::choosecyl}, we construct a left homotopy from $f$ to $g$ through the cylinder $IA$, i.e. a commutative diagram
\[
\xymatrix{
    A \Sum A \ar[r]^-{f + g} \ar@{>->}[d]_{i_0 + i_1} & B \ar@{>->}[d]^b_{\sim} \\
    IA \ar[r]^H & B^{'}
    }
\]
where $b$ is a trivial Hurewicz cofibration. By \lemmaref{lem:trivhurcof}, $b$ admits a retract $r$, and $h = rH$ is the desired homotopy.
\end{proof}

\chapter{Applications of the homotopy calculus of fractions}
\label{chap:homcalsappl}

Let us recapitulate what we have done so far. In \sectionref{sec:generalhomotopicmaps}, for a category pair $({\cat M}, {\cat W})$ we said that two maps $f \simeq g \colon A \ra B$ are homotopic if they have the same image in ${\bf ho}{\cat M}$. We have denoted $\pi {\cat M} = {\cat M} /_\simeq$, and $[A,B] = Hom_{\pi {\cat M}}(A, B)$. By definition therefore, the induced functor $\pi{\cat M} \ra {\bf ho}{\cat M}$ is faithful.

For a precofibration category $({\cat M}, {\cat W}, {\cat Cof})$, we have proved Brown's \theoremref{thm:homotmaps}, which says that if $A, B$ are cofibrant then $f \simeq g$ iff $f \simeq_l g$. The left homotopy relation $\simeq_l$ was defined in \sectionref{chap:homotopicmaps}. We have also proved Anderson's \theoremref{thm::catcofobjequiv} saying that inclusion induces an equivalence of categories ${\bf ho}{\cat M}_{cof} \cong {\bf ho}{\cat M}$. 

Finally, for $({\cat M}, {\cat W}, {\cat Cof})$ we have proved Brown's homotopy calculus of fractions \theoremref{thm::descrhomot}, which says that if $A, B$ are cofibrant, then any map in $Hom_{{\bf ho}{\cat M}}(A, B)$ is represented by a left fraction $s^{-1}f$
\[
\xymatrix{ A \ar[r]^f & B^{'} & B \ar[l]_-s^-\sim
  }
\]
with $s$ a weak equivalence and $B^{'}$ cofibrant. Furthermore, we can choose $f$ to be a cofibration and $s$ to be a trivial cofibration. 

Two such fractions $s^{-1}f$, $t^{-1}g$ are equal in ${\bf ho}{\cat M}$ if and only if there exist weak equivalences $s^{'}, t^{'}$ as in the diagram (\ref{eqn:descrhomotopy1}) with $B^{'''}$ cofibrant so that $s^{'}s \simeq t^{'}t$ and $s^{'}f \simeq t^{'}g$.

\bigskip

As an application, in this chapter we show that if ${\cat M}_k$ is a small set of precofibration categories, then the functor ${\bf ho}({\times} {\cat M}_k) \ra {\times}{\bf ho}{\cat M}_k$ is an {\it isomorphism} of categories. 

We also prove that any cofibration category $({\cat M}$, ${\cat W}$, ${\cat Cof})$ has saturated weak equivalences ${\cat W} = \overline{\cat W}$. A precofibration category $({\cat M}$, ${\cat W}$, ${\cat Cof})$ does not necessarily have saturated weak equivalences, but $({\cat M}, \overline{\cat W}, {\cat Cof})$ is again a precofibration category.

We then use Brown's homotopy calculus of fractions to give a number of sufficient conditions for ${\bf ho}{\cat M}$ in order to be ${\cat U}$-locally small, if ${\cat M}$ is ${\cat U}$-locally small.

\section{Products of cofibration categories}
\label{chap:prodcofcat}
Here is an application of homotopy calculus of fractions. If $({\cat M}_k, {\cat W}_k)$ for $k \eps K$ is a set of categories with weak equivalences, one can form the product $({\times}_{k \eps K} {\cat M}_k, {\times}_{k \eps K} {\cat W}_k)$ and its homotopy category denoted ${\bf ho}({\times} {\cat M}_k)$. Denote $p_k \colon {\times}_{k \eps K} {\cat M}_k \ra {\cat M}_k$ the projection. The components $({\bf ho}p_k)_{k \eps K}$ define a functor
\begin{center}
$P :{\bf ho}({\times} {\cat M}_k) \ra {\times}{\bf ho}{\cat M}_k$
\end{center}
If each category ${\cat M}_k$ carries a (pre)cofibration category structure $({\cat M}_k, {\cat W}_k, {\cat Cof}_k)$, then $(\times{\cat M}_k$, $\times{\cat W}_k$, $\times{\cat Cof}_k)$ defines the product (pre)cofibration category structure on ${\times} {\cat M}_k$. Dually, if each ${\cat M}_k$ carries a (pre)fibration category structure, then ${\times} {\cat M}_k$ carries a product (pre)fibration category structure.

Suppose that ${\cat M}_k$ are precofibration categories, and that $A = (A_k)_{k \eps K}$ is a cofibrant object of $\times {\cat M}_k$. Any factorization $A \Sum A$ $\ra$ $IA$ $\ra$ $A$ defines a cylinder in $\times {\cat M}_k$ iff each component $A_k \Sum A_k$ $\ra$ $(IA)_k$ $\ra$ $A_k$ is a cylinder in ${\cat M}_k$. 

If $B = (B_k)_{k \eps K}$ is a second cofibrant object and $f, g \colon A \ra B$ is a pair of maps in $\times {\cat M}_k$, then any homotopy $f \simeq g$ induces componentwise homotopies $f_k \simeq g_k$ in ${\cat M}_k$. Conversely, any set of homotopies $f_k \simeq g_k$ induces a homotopy $f \simeq g$.

\begin{thm}
\label{thm:prodhomotopycat}
If ${\cat M}_k$ for $k \eps K$ are each precofibration categories, or are each prefibration categories, then the functor
\begin{center}
$P \colon {\bf ho}({\times} {\cat M}_k) \ra {\times}{\bf ho}{\cat M}_k$
\end{center}
is an isomorphism of categories.
\end{thm}

\begin{proof}
Assume that each ${\cat M}_k$ is a precofibration category (the proof for prefibration categories is dual). Our functor $P$ is a bijection on objects, and we'd like to show that it is also fully faithful.

By \theoremref{thm::catcofobjequiv} we have equivalences of categories ${\bf ho}{\cat M}_k \cong {\bf ho}({\cat M}_k)_{cof}$ and ${\bf ho}(\times{\cat M}_k) \cong {\bf ho}(\times({\cat M}_k)_{cof})$. It suffices therefore in our proof to assume that ${\cat M}_k = ({\cat M}_k)_{cof}$ for all $k$. 

To prove fullness of $P$, let $A_k$ and $B_k$ be objects of ${\cat M}_k$ for $k \eps K$. Denote $A = (A_k)_{k \eps K}$, $B = (B_k)_{k \eps K}$ the corresponding objects of $\times{\cat M}_k$. Any map $\phi \colon A \ra B$ in $\times{\bf ho}{\cat M}_k$ can be expressed on components, using \theoremref{thm::descrhomot} (a) for each ${\cat M}_k$, as a left fraction $\phi_k = s_k^{-1}f_k$, with weak equivalences $s_k$. The map $\phi$ therefore is the image via $P$ of $(s_k)^{-1}(f_k)$.

To prove faithfulness of $P$, suppose that $\phi, \psi \colon A \ra B$ are maps in ${\bf ho}({\times} {\cat M}_k)$ which have the same image via $P$. Using \theoremref{thm::descrhomot} (a) applied to ${\times} {\cat M}_k$, we can write $\phi$ as $(s_k)^{-1}(f_k)$ and $\psi$ as $(t_k)^{-1}(g_k)$, with $s_k, t_k$ weak equivalences. From \theoremref{thm::descrhomot} (b) applied to each ${\cat M}_k$, we can find weak equivalences $s^{'}_k, t^{'}_k$ such that we have componentwise homotopies $s^{'}_ks_k \simeq t^{'}_kt_k$ and $s^{'}_kf_k \simeq t^{'}_kg_k$. The componentwise homotopies induce homotopies $s^{'}s \simeq t^{'}t$ and $s^{'}f \simeq t^{'}g$ in ${\times} {\cat M}_k$, so $\phi = \psi$ and we have shown that $P$ is faithful.
\end{proof}

\section{Saturation}
\label{sec:saturation}
Given a category with weak equivalences $({\cat M}, {\cat W})$, recall that $\overline{\cat W}$ denotes the saturation of ${\cat W}$, i.e. the class of maps of ${\cat M}$ that become isomorphisms in ${\bf ho}{\cat M}$.

\begin{lem}[Cisinski]
\label{lem::cofcatequivinhomotcat}
\mbox{}
\begin{enumerate}
\item
Suppose that $({\cat M}$, ${\cat W}$, ${\cat Cof})$ is a precofibration category, and that $f:A \ra B$ is a map with $A, B$ cofibrant.
\begin{enumerate}
\item
$f$ has a left inverse in ${\bf ho}{\cat M}$ if and only if there exists a cofibration $f^{'}:B \ra B^{'}$ such that $f^{'}f$ is a weak equivalence.
\item
$f$ is an isomorphism in ${\bf ho}{\cat M}$ if and only if there exist cofibrations $f^{'}:B \ra B^{'}, \, f^{''}:B^{'} \ra B^{''}$ such that $f^{'}f, f^{''}f^{'}$ are weak equivalences.
\end{enumerate}
\item
Suppose that $({\cat M}$, ${\cat W}$, ${\cat Fib})$ is a prefibration category, and that $f:A \ra B$ is a map with $A, B$ fibrant.
\begin{enumerate}
\item
$f$ has a right inverse in ${\bf ho}{\cat M}$ if and only if there exists a fibration $f^{'}:A^{'} \ra A$ such that $ff^{'}$ is a weak equivalence.
\item
$f$ is an isomorphism in ${\bf ho}{\cat M}$ if and only if there exist fibrations $f^{'}:A^{'} \ra A, \, f^{''}:A^{''} \ra A^{'}$ such that $ff^{'}, f^{'}f^{''}$ are weak equivalences.
\end{enumerate}
\end{enumerate}
\end{lem}

\begin{proof}
We only prove (1). The implications (a) $(\Leftarrow)$, (b) $(\Leftarrow)$ are clear. 

To prove (a) $(\Rightarrow)$, using \theoremref{thm::descrhomotcof} write the left inverse of $f$ in ${\bf ho}{\cat M}$ as a left fraction $s^{-1}f^{'}$ with $s$ a weak equivalence with cofibrant codomain. We get $1 = s^{-1}f^{'}f$ in ${\bf ho}{\cat M}$, therefore $s = f^{'}f$ in ${\bf ho}{\cat M}$ which means $s \simeq f^{'}f$. Since $s$ is a weak equivalence, $f^{'}f$ must be a weak equivalence.

Part (b) $(\Rightarrow)$ is a corollary of (a) applied first to the map $f$ to construct $f^{'}$ then to the map $f^{'}$ to construct $f^{''}$.
\end{proof}

\begin{lem}
\label{lem::cofcatissat}
\mbox{}
\begin{enumerate}
\item
If a precofibration category $({\cat M}$, ${\cat W}$, ${\cat Cof})$ satisfies CF6, then it has saturated weak equivalences ${\cat W} = \overline{\cat W}$.
\item
If a prefibration category $({\cat M}$, ${\cat W}$, ${\cat Fib})$ satisfies F6, then it has saturated weak equivalences ${\cat W} = \overline{\cat W}$.
\end{enumerate}
\end{lem}

\begin{proof}
We only prove (1). It suffices to show that any cofibration $A_0 \rightarrowtail A_1$ in $\overline{\cat W}$ is also in ${\cat W}$. Using \lemmaref{lem::cofcatequivinhomotcat}, we construct a sequence of cofibrations $A_0 \rightarrowtail A_1 \rightarrowtail A_2 ...$ with $A_n \rightarrowtail A_{n+2}$ a trivial cofibration, for all $n \ge 0$. From CF6, we see that $A_0 \rightarrowtail \colim^n A_n$ and $A_1 \rightarrowtail \colim^n A_n$ are trivial cofibrations, and we conclude from the 2 out of 3 axiom that $A_0 \rightarrowtail A_1$ is a trivial cofibration.
\end{proof}

As a consequence we have:
\begin{thm}
\label{thm::cofcatissat}
\mbox{}
\begin{enumerate}
\item
Any cofibration category $({\cat M}$, ${\cat W}$, ${\cat Cof})$ has saturated weak equivalences ${\cat W} = \overline{\cat W}$.
\item
Any fibration category $({\cat M}$, ${\cat W}$, ${\cat Fib})$ has saturated weak equivalences ${\cat W} = \overline{\cat W}$. $\square$
\end{enumerate}
\end{thm}

Since any Quillen model category is an ABC model category, we also have:
\begin{thm}
\label{thm:quillensaturatedweq}
Any Quillen model category $({\cat M}$, ${\cat W}$, ${\cat Cof}$, ${\cat Fib})$ has saturated weak equivalences ${\cat W}$ = $\overline{\cat W}$. $\square$
\end{thm}

In preparation for \theoremref{thm:cofcatissat2} below, we will recall two definitions. Suppose that $({\cat M}, {\cat W})$ is a category with a class of weak equivalences.

\begin{defn}[2 out of 6 axiom]
\label{defn:2outof6}
\index{axiom!2 out of 6 axiom}
${\cat W}$ satisfies the {\it 2 out of 6} property with respect to ${\cat M}$ if any sequence of composable maps $\overset{f}{\ra} \overset{g}{\ra} \overset{h}{\ra}$ in ${\cat M}$ for which the two compositions $gf$, $hg$ are in ${\cat W}$, the four maps $f$, $g$, $h$, $hgf$ are also in ${\cat W}$.
\end{defn}

The 2 out of 6 axiom is stronger than the 2 out of 3 axiom - this can be seen taking $f$, $g$ or $h$ to be identity maps in \definitionref{defn:2outof6}.

\begin{defn}[Weak saturation]
\label{defn:weaksaturation}
\index{weak saturation}
${\cat W}$ is {\it weakly saturated} with respect to ${\cat M}$ if:
\begin{description}
\item[WS1]
Every identity map is in ${\cat W}$
\item[WS2]
${\cat W}$ satisfies the 2 out of 3 axiom
\item[WS3]
If two maps $i \colon A \ra B$, $r \colon B \ra A$ in ${\cat M}$ satisfy $ri = 1_B$ and $ir \eps {\cat W}$ then $i, r \eps {\cat W}$
\end{description}
\end{defn}

If ${\cat W}$ satisfies 2 out of 6 and WS1, it is weakly saturated, for if $i, r$ are maps as in WS3 then the sequence $A \overset{i}{\ra} B \overset{r}{\ra} A \overset{i}{\ra} B$ has the 2 out of 6 property so $i, r \eps {\cat W}$.

If ${\cat W}$ satisfies 2 out of 3, WS1, and is closed under retracts, then it is also weakly saturated, for if $i, r$ are maps as in WS3 then we can exhibit $r$ as a retract of $ir$
\[
\xymatrix{
  B \ar[r]^1 \ar[d]_r & B \ar[r]^{1} \ar[d]^{ir} & B \ar[d]^r \\
  A \ar[r]^i & B \ar[r]^r & A
}
\]
and therefore $r, i \eps {\cat W}$.

Here is a characterization of precofibration and prefibration categories with saturated weak equivalences. 

\begin{thm}
\label{thm:cofcatissat2}
Suppose that $({\cat M}$, ${\cat W})$ admits either a precofibration or prefibration category structure. Then the following are equivalent:
\begin{enumerate}
\item
${\cat W}$ is weakly saturated
\item
${\cat W}$ satisfies the 2 out of 6 axiom
\item
${\cat W}$ is closed under retracts
\item
${\cat W}$ is saturated
\end{enumerate}
\end{thm}

\begin{proof}
We treat only the precofibration category case. For us ${\cat W}$ satisfies 2 out of 3 (this is CF2) and includes all isomorphisms (by CF1). In particular, ${\cat W}$ satisfies WS1.

Under these conditions, we have seen that (2) $\Ra$ (1), (3) $\Ra$ (1). We clearly have (4) $\Ra$ (1), (2), (3). It remains to show that (1) $\Ra$ (4).

Suppose that $({\cat M}$, ${\cat W}$, ${\cat Cof})$ is a precofibration category, with ${\cat W}$ weakly saturated. It suffices to show that any $\overline{\cat W}$-trivial cofibration $f \colon A_0 \rightarrowtail A_1$, with $A_0$ cofibrant, satisfies $f \eps {\cat W}$.

The over category $({\cat M} \downarrow A_1)$ carries an induced precofibration category by \propositionref{prop:overcofcat}. We apply \lemmaref{lem::cofcatequivinhomotcat} to the map 
\[
\xymatrix{
  A_0 \ar@{>->}[rr]^f \ar[dr]_f && A_1 \ar[dl]^{1_{A_1}} \\
  & A_1 &
}
\]
in $({\cat M} \downarrow A_1)$, and we construct a cofibration $g \colon A_1 \rightarrowtail A_2$ and a map $h \colon A_2 \ra A_1$, with $gf \eps {\cat W}$ and $hg = 1_{A_{2}}$. We observe that $ghg = g$, from which $gh$, $1_{A_2}$ have the same image in ${\bf ho}{\cat M}$, and therefore $gh \simeq 1_{A_2}$. Since $1_{A_2} \eps {\cat W}$, we have $gh \eps {\cat W}$. By WS3 applied to $g, h$, these two maps are in ${\cat W}$. By 2 out of 3 now $f \eps {\cat W}$.
\end{proof}

\begin{thm}[Cisinski]
\label{thm::cofcatsaturation}
\mbox{}
\begin{enumerate}
\item
If $({\cat M}$, ${\cat W}$, ${\cat Cof})$ is a precofibration category (resp. a CF1-CF5 cofibration category), then so is $({\cat M}$, $\overline{\cat W}$, ${\cat Cof})$. 
\item
If $({\cat M}$, ${\cat W}$, ${\cat Fib})$ is a prefibration category (resp. a F1-F5 cofibration category), then so is $({\cat M}$, $\overline{\cat W}$, ${\cat Fib})$. 
\end{enumerate}
\end{thm}

\begin{proof}
We only prove (1). Suppose that $({\cat M}$, ${\cat W}$, ${\cat Cof})$ is a precofibration category. 

(i) {\it The axioms CF1, CF2, CF3 (1), CF4} for $({\cat M}$, $\overline{\cat W}$, ${\cat Cof})$ are clearly satisfied.

(ii) {\it The axiom CF3 (2)} for $({\cat M}$, $\overline{\cat W}$, ${\cat Cof})$. Given a solid diagram in ${\cat M}$, with $A$, $C$ cofibrant and $i$ a $\overline{\cat W}$-trivial cofibration,
\begin{center}
$ \xymatrix {
	A \ar[r] \ar@{>->}[d]_i &
	C \ar@{-->}[d]^j \\
	B \ar@{-->}[r] &
	D
	}$
\end{center}
then by \lemmaref{lem::cofcatequivinhomotcat} there exist cofibrations $i^{'}, i^{''}$ such that $i^{'}i, i^{''}i^{'}$ are ${\cat W}$-trivial cofibrations. Denote $j^{'}, j^{''}$ the pushouts of the cofibrations $i^{'}, i^{''}$. We get that $j^{'}j, j^{''}j^{'}$ are ${\cat W}$-trivial cofibrations, and thereore $j$ is a $\overline{\cat W}$-trivial cofibration.

Assume now that $({\cat M}$, ${\cat W}$, ${\cat Cof})$ satisfies CF5.

(iii) {\it The axiom CF5 (1)} for $({\cat M}$, $\overline{\cat W}$, ${\cat Cof})$ is clearly satisfied.

(iv) {\it The axiom CF5 (2)}. Suppose that $f_i \colon A_i \ra B_i$ for $i \eps I$ is a set of $\overline{\cat W}$-trivial cofibrations with $A_i$ cofibrant. The map $\Sum f_i$ is a cofibration by axiom CF5 (1). By \lemmaref{lem::cofcatequivinhomotcat}, there exist cofibrations $f_i^{'} \colon B_i \ra B_i^{'}$, $f_i^{''} \colon B_i^{'} \ra B_i^{''}$ such that $f_i^{'}f_i$, $f_i^{''}f_i^{'}$ are ${\cat W}$-trivial cofibrations. It follows that $\Sum f_i^{'}f_i$, $\Sum f_i^{''}f_i^{'}$ are ${\cat W}$-trivial cofibrations, therefore $\Sum f_i$ is a $\overline{\cat W}$-trivial cofibration.
\end{proof}

We also note the following:

\begin{prop}
\label{prop:rightliftpropisweakequiv}
\mbox{}
\begin{enumerate}
\item
If $({\cat M}$, ${\cat W}$, ${\cat Cof})$ is a precofibration category, then all maps with the right lifting property with respect to all cofibrations are in $\overline{\cat W}$.
\item
If $({\cat M}$, ${\cat W}$, ${\cat Fib})$ is a prefibration category, then all maps with the left lifting property with respect to all fibrations are in $\overline{\cat W}$.
\end{enumerate}
\end{prop}

\begin{proof}
We only prove (1). Suppose that a map $p \colon C \ra D$ has the right lifting property with respect to all cofibrations. We construct a commutative square
\[
\xymatrix{
  A \ar[r]_{\sim}^{r} \ar@{>->}[d]_i & C \ar[d]^p \\
  B \ar[r]_\sim^{r^{'}} \ar[ru]^{s} & D
}
\]
where $A$ is a cofibrant replacement of $C$, and $r^{'}i$ is a factorization of $pr$ as a cofibration followed by a weak equivalence. A lift $s$ exists since $p$ has the right lifting property with respect to $i$, and applying the 2 out of 6 property to $p$, $s$, $i$ shows that $p \eps \overline{\cat W}$.
\end{proof}

Let us adapt this discussion to the case of left proper cofibration categories, and show that if $({\cat M}$, ${\cat W}$, ${\cat Cof})$ is a left proper cofibration category then so is $({\cat M}$, $\overline{\cat W}$, ${\cat Cof})$.

\begin{lem}
\label{lem::cofcatequivinhomotcatproper}
\mbox{}
\begin{enumerate}
\item
Suppose that $({\cat M}$, ${\cat W}$, ${\cat Cof})$ is a left proper precofibration category, and that $f:A \ra B$ is a map
\begin{enumerate}
\item
$f$ has a left inverse in ${\bf ho}{\cat M}$ if and only if there exists a left proper map $f^{'}:B \ra B^{'}$ such that $f^{'}f$ is a weak equivalence.
\item
$f$ is an isomorphism in ${\bf ho}{\cat M}$ if and only if there exist left proper maps $f^{'}:B \ra B^{'}, \, f^{''}:B^{'} \ra B^{''}$ such that $f^{'}f, f^{''}f^{'}$ are weak equivalences.
\end{enumerate}
\item
Suppose that $({\cat M}$, ${\cat W}$, ${\cat Fib})$ is a right proper prefibration category, and that $f:A \ra B$ is a map.
\begin{enumerate}
\item
$f$ has a right inverse in ${\bf ho}{\cat M}$ if and only if there exists a right proper map $f^{'}:A^{'} \ra A$ such that $ff^{'}$ is a weak equivalence.
\item
$f$ is an isomorphism in ${\bf ho}{\cat M}$ if and only if there exist right proper maps $f^{'}:A^{'} \ra A, \, f^{''}:A^{''} \ra A^{'}$ such that $ff^{'}, f^{'}f^{''}$ are weak equivalences.
\end{enumerate}
\end{enumerate}
\end{lem}

\begin{proof}
We only prove (1). The implications (a) $(\Leftarrow)$, (b) $(\Leftarrow)$ are trivial.
 
Let us prove (a) $(\Rightarrow)$. Construct the diagram
\begin{center}
$ \xymatrix {
    & B \ar[r]^{f^{'}} & B^{'} \\
    A \ar[r]^{f_2} \ar[ru]^f & B_2 \ar[u]_{s_2} \ar[r]^{f^{'}_2} & B^{'}_2 \ar[u] \\
    A_1 \ar@{>->}[r]^{f_1} \ar[u]^{r_1}_\sim & B_1 \ar@{>->}[r]^{f^{'}_1} \ar[u]_{s_1} & B^{'}_1 \ar[u]
  }$
\end{center}
as follows:
\begin{enumerate}
\item
$A_1$ is a cofibrant replacement of $A$, $f_1$ is a cofibrant replacement of $fr_1$. 
\item
It follows that $f_1$ has a left inverse in ${\bf ho}{\cat M}$. We use \lemmaref{lem::cofcatequivinhomotcat} to construct a cofibration $f^{'}_1$ such that $f^{'}_1f_1$ is a ${\cat W}$-weak equivalence.
\item
The maps $f_2$, resp. $f^{'}_2$ and $f^{'}$ are pushouts of $f_1$, resp. $f^{'}_1$. These pushouts can be constructed because ${\cat M}$ is left proper.
\end{enumerate}
The maps $r_1$ and $s_2s_1$ are ${\cat W}$-weak equivalences, and all horizontal maps are left proper. It follows that all vertical maps are ${\cat W}$-weak equivalences. Since $f^{'}_1f_1$ is a ${\cat W}$-weak equivalence, so is $f^{'}f$

Part (b) is proved applying (a) first to the map $f$ to construct $f^{'}$ and then a second time to the map $f^{'}$ to construct $f^{''}$.
\end{proof}

Using the previous lemma, the statement below is immediate:

\begin{prop}
\label{prop:leftpropermapsaturated}
\mbox{}
\begin{enumerate}
\item
In a left proper precofibration category $({\cat M}$, ${\cat W}$, ${\cat Cof})$, all the ${\cat W}$-left proper maps are $\overline{\cat W}$-left proper
\item
In a right proper prefibration category $({\cat M}$, ${\cat W}$, ${\cat Fib})$, all the ${\cat W}$-right proper maps are $\overline{\cat W}$-right proper
\end{enumerate}
\end{prop}

\begin{proof}
We only prove (1). Suppose that $f \colon A \ra B$ is a ${\cat W}$-left proper map. In the diagram with full maps
\begin{center}
$\xymatrix {
    A \ar[r] \ar[d]_f & C_1 \ar[r]^r \ar@{-->}[d] & C_2 \ar@{-->}[r]^{r_1} \ar@{-->}[d] & C_3 \ar@{-->}[r]^{r_2} \ar@{-->}[d] & C_4 \ar@{-->}[d] \\
    B \ar@{-->}[r] & D_1 \ar@{-->}[r]^{r^{'}} & D_2 \ar@{-->}[r]^{r^{'}_1} & D_3 \ar@{-->}[r]^{r^{'}_2} & D_4
  }$
\end{center}
suppose that $r$ is a $\overline{\cat W}$-weak equivalence. Use \lemmaref{lem::cofcatequivinhomotcatproper} to construct the maps $r_1$, $r_2$ such that $r_1r$ and $r_2r_1$ are ${\cat W}$-weak equivalences. Denote $r^{'}$, $r^{'}_1$, $r^{'}_2$ the pushouts along $f$. The maps $r^{'}_1r^{'}$ and $r^{'}_2r^{'}_1$ are ${\cat W}$-weak equivalences, therefore $r^{'}$ is a $\overline{\cat W}$-weak equivalence.
\end{proof}

We can now show the following result.

\begin{thm}
\label{thm::cofcatsaturationproper}
\mbox{}
\begin{enumerate}
\item
If $({\cat M}$, ${\cat W}$, ${\cat Cof})$ is a left proper precofibration category (resp. a left proper CF1-CF5 cofibration category), then so is $({\cat M}$, $\overline{\cat W}$, ${\cat Cof})$. 
\item
If $({\cat M}$, ${\cat W}$, ${\cat Fib})$ is a right proper prefibration category (resp. a right proper F1-F5 cofibration category), then so is $({\cat M}$, $\overline{\cat W}$, ${\cat Fib})$. 
\end{enumerate}
\end{thm}

\begin{proof}
Consequence of \theoremref{thm::cofcatsaturation} and \propositionref{prop:leftpropermapsaturated}.
\end{proof}

\section{Local smallness of ${\bf ho}{\cat M}$}
\label{sec:pim}

In this section, we will give a number of sufficient conditions for ${\bf ho}{\cat M}$ to be locally ${\cat U}$-small, if ${\cat M}$ is locally ${\cat U}$-small. We also show that a Quillen model category always has saturated weak equivalences.

\begin{defn}
\label{defn:fibrmodel}
\mbox{}
\begin{enumerate}
\item
Suppose that ${\cat M}$ is a precofibration category. An object $C$ is a {\it fibrant model} if for any cofibrant object $A$, map $f$ and trivial cofibration $i$ 
\[
\xymatrix{
  A \ar[rr]^f \ar@{>->}[d]_i^\sim && C \\
  B \ar@{-->}[rru]_h &&
}
\]
a lift $h$ exists making the diagram commutative. 
\item
Suppose that ${\cat M}$ is a prefibration category. An object $B$ is a {\it cofibrant model} if for any fibrant object $D$, map $g$ and trivial fibration $p$ 
\[
\xymatrix{
  && C \ar@{>->}[d]^p_\sim \\
  B \ar@{-->}[rru]^h \ar[rr]_g && D
}
\]
a lift $h$ exists making the diagram commutative. 
\end{enumerate}
\end{defn}

In (1), if ${\cat M}$ admits a terminal object $\terminal$, an object $C$ is a fibrant model iff $C \ra \terminal$ has the right lifting property with respect to all trivial cofibrations with cofibrant domain. In (2), if ${\cat M}$ admits an initial object $\initial$, an object $B$ is a fibrant model iff $\initial \ra B$ has the left lifting property with respect to all trivial fibrations with fibrant domain. 

\begin{prop}
\label{prop:smallhomotopycat}
Suppose that either
\begin{enumerate}
\item
${\cat M}$ is a precofibration category, with $A$ a cofibrant object and $B$ a fibrant model, or
\item
${\cat M}$ is a prefibration category, with $A$ is a cofibrant model and $B$ a fibrant object.
\end{enumerate}

Then the induced function $[A, B]$ $\ra$ $Hom_{{\bf ho}{\cat M}}(A, B)$ is bijective.
\end{prop}

\begin{proof}
This argument is due to Denis-Charles Cisinski \cite{Cisinski2}. The function $[A, B]$ $\ra$ $Hom_{{\bf ho}{\cat M}}(A, B)$ is injective for any $A, B$, since the functor $\pi{\cat M} \ra {\bf ho}{\cat M}$ is faithful. In case (1), any map $\phi$ in $Hom_{{\bf ho}{\cat M}}(A, B)$ can be represented by a zig-zag
\[
\xymatrix{ A \ar@{>->}[r]^f & B^{''} & B^{'} \ar@{>->}[l]_s^\sim \ar[r]^t_\sim & B}
\]
where a cofibrant object $B^{'}$ and a weak equivalence $t$ are constructed by CF4, a trivial cofibration $s$ and a cofibration $f$ by \theoremref{thm::descrhomotcof}. Since $B$ is a fibrant model, a lift $h \colon B^{''} \ra B$ exists with $hs = t$, and  $\phi$ is represented by $hf$.
\end{proof}

\begin{prop}
\label{prop:smallhomotopycat2}
Suppose that ${\cat M}$ is either
\begin{enumerate}
\item
A precofibration category, and each cofibrant object $A$ admits a weak equivalence map $A \overset{\sim}{\ra} A^{'}$ with $A^{'}$ a fibrant model, or
\item
A prefibration category, and each fibrant object $A$ admits a weak equivalence map $A^{'} \overset{\sim}{\ra} A$ with $A^{'}$ a cofibrant model.
\end{enumerate}

Then if ${\cat M}$ is ${\cat U}$-locally small, so is ${\bf ho}{\cat M}$.
\end{prop}

\begin{proof}
This follows from \propositionref{prop:smallhomotopycat}.
\end{proof}

Here is a result similar in spirit to \propositionref{prop:smallhomotopycat}.
\begin{prop}
\label{prop:smallhomotopycatbaues}
Suppose that either
\begin{enumerate}
\item
\label{prop:smallhomotopycat:baues1}
${\cat M}$ is a precofibration category satisfying the Baues axiom BCF6, and $A, B$ are two cofibrant objects.
\item
${\cat M}$ is a prefibration category satisfying the Baues axiom BF6 (dual to BCF6), and $A, B$ are two fibrant objects.
\end{enumerate}

Then the induced function $[A, B]$ $\ra$ $Hom_{{\bf ho}{\cat M}}(A, B)$ is bijective.
\end{prop}

\begin{proof}
Each map $\phi$ in $Hom_{{\bf ho}{\cat M}}(A, B)$ is represented by a left fraction $s^{-1}f$
\[
\xymatrix{ A \ar@{>->}[r]^f & B^{'} & B \ar@{>->}[l]_-s^-\sim
  }
\]
with $f$ a cofibration and $s$ a trivial cofibration. We now use BCF6 to construct a trivial cofibration $s^{'} \colon B^{'} \ra C$, with $C$ satisfying the property that each trivial cofibration to $C$ admits a left inverse. In particular, $s^{'}s$ has a left inverse $t$. Then $\phi$ is represented by $ts^{'}f$, which completes the proof.
\end{proof}

\begin{cor}
\label{cor:smallhomotopycat}
Suppose that ${\cat M}$ is either:
\begin{enumerate}
\item
A Baues cofibration category with an initial object
\item
A Baues fibration category with a terminal object
\item
An ABC premodel category for which trivial cofibrations with cofibrant domain have the left lifting property with respect to fibrations with fibrant codomain
\item
An ABC premodel category for which cofibrations with cofibrant domain have the left lifting property with respect to trivial fibrations with fibrant codomain
\item
A Quillen model category
\end{enumerate}

Then if ${\cat M}$ is ${\cat U}$-locally small, so is ${\bf ho}{\cat M}$.
\end{cor}

\begin{proof}
Parts (1), (2) follow from \propositionref{prop:smallhomotopycatbaues}, and (3)-(5) from \propositionref{prop:smallhomotopycat2}. 

Hans Baues in fact shows \cite{Baues1} that any ${\cat U}$-locally small Baues cofibration  (or fibration) category has an ${\cat U}$-locally small homotopy category, without requiring the existence of an initial (resp. terminal) object.
\end{proof}

\section{Homotopic maps in a Quillen model category}
We now describe Quillen's approach to homotopic maps, for a Quillen model category ${\cat M}$. The definition of $\simeq_{ql}$ below is similar to that of $\simeq_l$ defined in \sectionref{chap:homotopicmaps}.

Suppose that $A, B \eps {\cat M}$ are two objects, with $A$ cofibrant and $B$ fibrant. Two map $f, g \colon A \ra B$ are Quillen left homotopic $f \simeq_{ql} g$ if there exists a cylinder $IA$ and homotopy map $H$ making commutative the diagram
\[
\xymatrix{
  A \Sum A \ar[r]^{f + g} \ar[d]_{i_0 + i_1} & B  \\
  IA \ar[ru]^H & 
}
\].

Dually, $f$ and $g$ are Quillen right homotopic $f \simeq_{qr} g$ if there exists a path object $B^I$ and homotopy map $H$ making commutative the diagram
\[
\xymatrix{
  & B^I \ar[d]^{(p_0, p_1)} \\
  A \ar[r]^{(f,g)} \ar[ru]^H & B \times B   \\
}
\].

\begin{lem}
\label{lem:descrhomotopyquillen}
The following are equivalent:
\begin{enumerate}
\item
$f \simeq_{ql} g$
\item
$f \simeq_{qr} g$
\item
$f \simeq g$
\end{enumerate}
\end{lem}

\begin{proof}
We only prove (1) $\Leftrightarrow$ (3). Clearly (1) $\Ra$ (3).

Suppose $f \simeq g$. Pick a trivial cofibration $r \colon B^{'} \ra B$ with $B^{'}$ cofibrant by M5, and by M4 lifts $f^{'}, g^{'} \colon A \ra B^{'}$ with $rf^{'} = f$ and $rg^{'} = g$. 
\[
\xymatrix{
  A \Sum A \ar[r]^{f^{'} + g^{'}} \ar[d]_{i_0 + i_1} & B^{'} \ar@{>->}[d]_{b^{'}}^\sim \ar@{->>}[r]^r_\sim & B  \\
  IA \ar[r]^{H^{'}} & B^{''} \ar[ru]_h &   
}
\]
We have $f^{'} \simeq g^{'}$, and \theoremref{thm:homotmaps} yields $f^{'} \simeq_l g^{'}$ through a cylinder $IA$, a trivial cofibration $b^{'}$ and a homotopy map $H^{'}$. Since $B$ is fibrant, by M4 there exists a lift $h$ and now $H = hH^{'}$ provides the desired homotopy $f \simeq_{ql} g$.
\end{proof}

\section{A characterization of proper Quillen model categories}

For an object $A$ of a category ${\cat C}$, recall that the under category $(A \downarrow {\cat C})$ has:
\begin{enumerate}
\item
as objects, pairs $(B, f)$ of an object $B \eps {\cat C}$ and map $f \colon A \ra B$
\item
as maps $(B_1, f_1) \ra (B_2, f_2)$, the maps $g \colon B_1 \ra B_2$ such that $gf_1 = f_2$.
\end{enumerate}

The over category $({\cat C} \downarrow A)$ is defined by duality as $(A \downarrow {\cat C}^{op})^{op}$.

\begin{thm}[D.-C. Cisinski]
\label{thm:charactproperness}
\mbox{}
\begin{enumerate}
\item
If ${\cat M}$ is a left proper ABC precofibration category, then any weak equivalence $f \colon A \ra B$ induces an equivalence of categories ${\bf Ho}(B \downarrow {\cat M}) \ra {\bf Ho}(A \downarrow {\cat M})$.
\item
If ${\cat M}$ is a right proper ABC prefibration category, then any weak equivalence $f \colon A \ra B$ induces an equivalence of categories ${\bf Ho}({\cat M} \downarrow A) \ra {\bf Ho}({\cat M} \downarrow B)$.
\end{enumerate}
\end{thm}

\begin{proof}
We only prove (1). Denote $({\cat M}, {\cat W}, {\cat Cof})$ a left proper precofibration structure on ${\cat M}$. Using \theoremref{thm:leftproperstruct}, $({\cat M}, {\cat W}, {\cat PrCof})$ is also a cofibration structure. 

We notice that the latter cofibration structure induces a cofibration structure on $(A \downarrow {\cat M})$, where a map $g \colon (B_1, f_1) \ra (B_2, f_2)$ is a weak equivalence, resp. cofibration if $g \eps {\cat W}$, resp. $g \eps {\cat PrCof}$. We construct the following Quillen partial adjunction:
\begin{equation}
\label{eqn:thm:charactproperness}
\xymatrix{
    (A \downarrow {\cat M})_{\cat PrCof} \ar[r]^-{v_1} \ar@{_{(}->}[d]_{t_1} & (B \downarrow {\cat M}) \\
    (A \downarrow {\cat M}) & (B \downarrow {\cat M}) \ar[l]_-{v_2} \ar[u]_{t_2}^= &
    }
\end{equation}

We have denoted $(A \downarrow {\cat M})_{\cat PrCof}$ the full subcategory of  $(A \downarrow {\cat M})$ with objects the left proper maps $A \ra X$. The functor $t_1$ is inclusion, $t_2$ is the identity, $v_1$ is pushout along $f$, $v_2$ is restriction along $f$. The functor $t_1$ is a cofibrant approximation, therefore a left approximation. The functor $t_2$ is a (trivial) right approximation.A map $v_1t_1A^{'} \ra t_2B^{'}$ is a weak equivalence iff its adjoint $t_1A^{'} \ra v_2t_2B^{'}$ is a weak equivalence. The result now falls out from \theoremref{thm:generalexistencetotalderivedadjoint2}.

\end{proof}

A converse of \theoremref{thm:charactproperness} holds for Quillen model categories.

\begin{thm}[C. Rezk]
\label{thm:charactpropernessquillen}
Suppose that $({\cat M}$, ${\cat W}$, ${\cat Cof}$, ${\cat Fib})$ is a Quillen model category.
\begin{enumerate}
\item
${\cat M}$ is left proper iff any weak equivalence $f \colon A \ra B$ induces an equivalence of categories ${\bf Ho}(B \downarrow {\cat M}) \ra {\bf Ho}(A \downarrow {\cat M})$.
\item
${\cat M}$ is right proper iff any weak equivalence $f \colon A \ra B$ induces an equivalence of categories ${\bf Ho}({\cat M} \downarrow A) \ra {\bf Ho}({\cat M} \downarrow B)$.
\end{enumerate}
\end{thm}

\begin{proof}
We only prove (1). The {\it only if} part is proved by \theoremref{thm:charactproperness}. 

For the {\it if} part, in diagram \eqref{eqn:thm:charactproperness} the pair ${\bf L} v_1 \dashv {\bf ho} v_2$ is adjoint from \theoremref{thm:generalexistencetotalderivedadjoint2}. From our hypothesis ${\bf ho} v_2$ is an equivalence, therefore its left adjoint ${\bf L} v_1$ is an equivalence as well.

${\cat M}$ is a Quillen model category (has more structure than that of a mere cofibration category), which allows us to show that $v_1$ sends weak equivalences to weak equivalences. 

Indeed, $(A \downarrow {\cat M})_{cof}$ has a Quillen model category induced from ${\cat M}$, with all objects cofibrant. By Brown's factorization Lemma, any weak equivalence $u$ in $(A \downarrow {\cat M})_{cof}$ factors as a trivial cofibration followed by a left inverse to a trivial cofibration. Since trivial cofibrations push along $f$ to trivial cofibrations, it follows that $u$ pushes along $f$ to a weak equivalence.

So if $g \colon A \ra C$ is a cofibration, ${\bf L}v_1(g)$ is computed by $v_1g$. The adjunction unit $C \ra {\bf ho} v_2 \comp {\bf L}v_1(g)$ must be an isomorphism, therefore in the pushout
\[
\xymatrix{
  A \ar@{>->}[r]^g \ar[d]_f^\sim & C \ar[d]^{f^{'}} \\
  B \ar@{>->}[r]^{g^{'}} & D
}
\]
the map $f^{'}$ must be an isomorphism in ${\bf ho}{\cat M}$. Since ${\cat M}$ is a Quillen model category, by \theoremref{thm:quillensaturatedweq} we therefore have $f^{'} \eps {\cat W}$, which shows that ${\cat M}$ is left proper.
\end{proof}

As a consequence, in Quillen model categories left properness can be formulated without the use of cofibrations. If $({\cat M}$, ${\cat W}$, ${\cat Cof}$, ${\cat Fib})$ is a left proper Quillen model category, then it is left proper with respect to any other Quillen model category structure $({\cat M}$, ${\cat W}$, ${\cat Cof}^{'}$, ${\cat Fib}^{'})$.
\chapter{Review of category theory}
\label{chap:constrcat}

This chapter recalls classic constructions of category theory. We discuss over and under categories, inverse image categories, and diagram categories. The construction of limits and colimits is recalled. We discuss cofinal functors, the Grothendieck construction and some of its basic properties.

\section{Basic definitions and notations}
\label{sec:basicdefnotations}

\subsection{Initial and terminal categories}
\label{sec:initterminalcat}
The initial object in a category is denoted $\initial$, and the terminal object $\terminal$. The initial category (with an empty set of objects) is denoted $\initialcat$, and the terminal category (with one object and one identity map) is denoted $\terminalcat$.\index{$\initial$, $\terminal$, $\initialcat$, $\terminalcat$}

For a category ${\cat D}$ we denote $e_d \colon \terminalcat \ra {\cat D}$\index{$e_d$, $p_{\cat D}$} the functor that embeds $\terminalcat$ as the object $d \eps {\cat D}$, and $p_{\cat D} \colon {\cat D} \ra \terminalcat$ the terminal category projection.

\subsection{Inverse image category}
\label{sec:catinverseimage}
\index{category!inverse image category}
For a functor $u \colon {\cat A} \ra {\cat B}$ and an object $b$ of ${\cat B}$, the {\it inverse image} of $b$ is the subcategory $u^{-1}b$ of ${\cat A}$ consisting of objects $a$ with $ua = b$ and maps $f \colon a \ra a^{'}$ with $u f = 1_b$.

\subsection{Categories of diagrams}
\label{sec:catdiag}
If ${\cat D}$ and ${\cat M}$ are categories, the category of functors ${\cat D} \ra {\cat M}$ (or ${\cat D}$-diagrams of ${\cat M}$) is denoted ${\cat M}^{\cat D}$. 

A functor $u \colon {\cat D}_1 \ra {\cat D}_2$ induces a functor of diagram categories denoted $u^* \colon {\cat M}^{{\cat D}_2} \ra {\cat M}^{{\cat D}_1}$. If two functors $u, v$ are composable, then $(uv)^* = v^*u^*$. If two functors $u \dashv v$ are adjoint, then $v^* \dashv u^*$ are adjoint.

If ${\cat C}$ is a class of maps of ${\cat M}$ (for example the weak equivalences or the cofibrations in a cofibration category), we denote ${\cat C}^{\cat D}$ the class of maps $f$ of ${\cat M}^{\cat D}$ such that $f_d \eps C$ for any object $d \eps {\cat D}$.

We'd like to point out that we have defined the weak equivalences ${\cat W}$ and the cofibrations ${\cat Cof}$ of a cofibration category $({\cat M}$, ${\cat W}$, ${\cat Cof})$ as {\it classes of maps} rather than as subcategories. For a cofibration category, ${\cat W}^{\cat D}$ (and ${\cat Cof}^{\cat D}$) therefore denote the class of diagram maps $f \colon X \ra Y$ in ${\cat M}^{\cat D}$ such that each $f_d$ is a weak equivalence, resp. a cofibration. If $w{\cat M}$ denotes the subcategory of ${\cat M}$ generated by ${\cat W}$, then $(w{\cat M})^{\cat D}$ is {\it not} the same thing as ${\cat W}^{\cat D}$.

\section{Limits and colimits}
\label{sec:relativelimcolim}
Assume ${\cat D}_1$ and ${\cat D}_2$ are categories, $u \colon {\cat D}_1 \ra {\cat D}_2$ is a functor and ${\cat M}$ is a category. A colimit functor along $u$, denoted $\colim^u \colon {\cat M}^{{\cat D}_1} \ra {\cat M}^{{\cat D}_2}$, is by definition a left adjoint of $u^* \colon {\cat M}^{{\cat D}_2} \ra {\cat M}^{{\cat D}_1}$. A limit functor along $u$, denoted $\lim^u \colon {\cat M}^{{\cat D}_1} \ra {\cat M}^{{\cat D}_2}$, is by definition a right adjoint of $u^*$. $\colim^u$ is also called the left Kan extension along $u$, and $\lim^u$ is called the right Kan extension along $u$. 

Being a left adjoint, if $\colim^u$ exists then it is unique up to a unique natural isomorphism. As a right adjoint, if $\lim^u$ exists it is unique up to a unique natural isomorphism. 

We will also consider the case when $\colim^uX$ (resp. $\lim^uX$) exists for some, but not all objects $X$ of ${\cat M}^{{\cat D}_1}$. Each of $\colim^uX$ and $\lim^uX$ are defined by an universal property, and if $\colim^uX$ or $\lim^uX$ exist, then they are unique up to unique isomorphism.

If ${\cat D}_1 = {\cat D}$ and ${\cat D}_2$ is the point category, $\colim^u$ is the well-known 'absolute' colimit of ${\cat D}$-diagrams, denoted $\colim^{\cat D}$ (or simply $\colim$ if there is no confusion). Dually, in this case $\lim^u$ is the 'absolute' limit $\lim^{\cat D}$. If ${\cat D}$ has a terminal object $\final$, then $\colim^{\cat D} X$ always exists and the natural map $X_{\final} \ra \colim^{\cat D} X$ is an isomorphism. If ${\cat D}$ has an initial object $\initial$ then $\lim^{\cat D} X$ always exists and the natural map $\lim^{\cat D} X \ra X_{\initial}$ is an isomorphism.

Next lemma presents a base change formula for relative (co)limits. In the case of colimits, the lemma states that the relative colimit along a functor $u \colon {\cat D}_1 \ra {\cat D}_2$ can be pointwisely computed in terms of absolute colimits over $(u \downarrow d_2)$ for all objects $d_2 \eps {\cat D}_2$.

\begin{lem}
\label{lem::computerelcolim}
Suppose that $u \colon {\cat D}_1 \ra {\cat D}_2$ is a functor, and suppose that $X$ is a ${\cat D}_1$ diagram in a category ${\cat M}$.
\begin{enumerate}
\item
If $\colim^{(u \downarrow d_2)} X$ exists for all $d_2 \eps {\cat D}_2$, then $\colim^u X$ exists and 
\begin{center}
${\colim}^{(u \downarrow d_2)} X \cong ({\colim}^u X)_{d_2}$
\end{center}
\item
If $\lim^{(d_2 \downarrow u)} X$ exists for all $d_2 \eps {\cat D}_2$, then $\lim^u X$ exists and 
\begin{center}
$({\lim}^u X)_{d_2} \cong {\lim}^{(d_2 \downarrow u)} X$
\end{center}
\end{enumerate}
\end{lem}

\begin{proof}
We only prove (1). We show that $({\colim}^{(u \downarrow d_2)} X)_{d_2 \eps {\cat D}_2}$ satisfies the universal property that defines $\colim^u X$. Maps from $X$ to a diagram $u^* Y$ can be identified with maps $X_{d_1} \ra Y_{ud_1}$ that make the diagram below commutative for all maps $f \colon d_1 \ra d_1^{'}$
\begin{center}
$\xymatrix{
    X_{d_1} \ar[r] \ar[d]_{X_f} &  Y_{ud_1} \ar[d]^{Y_{uf}} \\
    X_{d_1^{'}} \ar[r] & Y_{ud_1^{'}}
    }$
\end{center}
They can be further identified with maps $X_{d_1} \ra Y_{d_2}$, defined for all $\phi \colon ud_1 \ra d_2$, that make the diagram below commutative 
\begin{center}
$\xymatrix{
    X_{d_1} \ar[r] \ar[d]_{X_f} &  Y_{d_2} \ar[d]^{Y_g} \\
    X_{d_1^{'}} \ar[r] & Y_{d_2^{'}}
    }$
\end{center}
for all maps $f, g$ that satisfy $\phi^{'} uf = g \phi$. They can  finally be identified with maps ${\colim}^{(u \downarrow d_2)} X \ra Y_{d_2}$ compatible in $d_2$, and from the universal property of the colimit we see that ${\colim}^u X$ exists and ${\colim}^{(u \downarrow d_2)} X \cong ({\colim}^u X)_{d_2}$.
\end{proof}

In particular, if ${\cat M}$ is closed under absolute colimits, it is also closed under relative colimits. By definition, if ${\cat M}$ is closed under small colimits we say that ${\cat M}$ is {\it cocomplete}. Dually, if ${\cat M}$ is closed under absolute limits, it is also closed under relative limits. If ${\cat M}$ is closed under small limits we say by definition that ${\cat M}$ is {\it complete}.

The following lemma presents another well known result - that the composition of (co)limits is the (co)limit of the composition.

\begin{lem}
\label{lem::colimcomposition}
Suppose that $u \colon {\cat D}_1 \ra {\cat D}_2$ and $v \colon {\cat D}_2 \ra {\cat D}_3$ are two functors, and suppose that $X$ is a ${\cat D}_1$ diagram in a category ${\cat M}$.
\begin{enumerate}
\item
Assume that $\colim^u X$ exists. If either $\colim^v \colim^u X$ or \\ $\colim^{vu} X$ exist, then they both exist and are canonically isomorphic.
\item
Assume that $\lim^u X$ exists. If either $\lim^v \lim^u X$ or $\lim^{vu} X$ exist, then they both exist and are canonically isomorphic.
\end{enumerate}
\end{lem}

\begin{proof}
We only prove (1). Since $\colim^u X$ exists, we have a bijection of sets natural in $Y \eps {\cat M}^{{\cat D}_3}$
\begin{center}
$Hom(X, u^*v^* Y) \cong Hom(\colim^u X, v^* Y)$ 
\end{center}
If $\colim^{vu} X$ exists, we also have a natural bijection
\begin{center}
$Hom(X, u^*v^* Y) \cong Hom(\colim^{vu} X, Y)$
\end{center}
therefore $\colim^{vu} X$ satisfies the universal property of $\colim^v \colim^u X$. If on the other hand $\colim^v \colim^u X$ exists, we have a natural bijection
\begin{center}
$Hom(\colim^u X, v^* Y) \cong Hom(\colim^v\colim^u X, Y)$
\end{center}
and $\colim^v \colim^u X$ satisfies the universal property of $\colim^{vu} X$. 
\end{proof}

In particular, if ${\cat M}$ is cocomplete then $\colim^{vu}$ and $\colim^v \colim^u$ are naturally isomorphic. Dually, if ${\cat M}$ is complete then $\lim^{vu}$ and $\lim^v \lim^u$ are naturally isomorphic.

\section{Simplicial sets}
\label{sec:simplindex}
\index{simplicial set}
This text assumes familiarity with simplicial sets, and the reader may refer for example to \cite{Goerss-Jardine} or \cite{Hovey} for a treatment of the standard theory of simplicial sets. 

We will denote by ${\bf n}$ the poset $\{0, 1, ... , n\}$, for $n \ge 0$, with the natural order. The cosimplicial indexing category $\Delta$\index{$\Delta$} is the category with objects ${\bf 0}, {\bf 1}, ..., {\bf n}, ...$ and maps the order-preserving maps ${\bf n}_1 \ra {\bf n}_2$. If ${\cat C}$ is a category, a simplicial object in ${\cat C}$ is then a functor $\Delta^{op} \ra {\cat C}$, and a cosimplicial object in ${\cat C}$ is a functor $\Delta \ra {\cat C}$.


Taking in particular simplicial objects in $Sets$ we get the category of simplicial sets, denoted $sSets$. The representable functor ${\Delta^{op}}(-, {\bf n})$ determines a simplicial set denoted $\Delta[n]$, called the standard $n$-simplex (for $n \ge 0$). For any map ${\bf m} \ra {\bf n}$ we have a simplicial set map $\Delta[m] \ra \Delta[n]$, thus $\Delta[-]$ determines a cosimplicial object in $sSets$.

In particular, we have simplicial set inclusions $d^0, d^1 \colon \Delta[0] \ra \Delta[1]$, which induce respectively functors $i_0, i_1 \colon X\sdot \cong X\sdot \times \Delta[0] \ra  X\sdot \times \Delta[1]$. We say that two maps $f, g \colon X\sdot \ra Y\sdot$ are simplcially homotopic if there exists a simplicial set map $h \colon X\sdot \times \Delta[1] \ra Y\sdot$ such that $f = hi_0, g = hi_1$. 





In general, the simplicial homotopy relation is not reflexive nor transitive.


\section{The nerve of a category}
\label{sec:catnerve}
Given a category ${\cat C}$, its nerve (or classifying space) $B{\cat C}$ is a simplicial set having as $n$-simplices the composable strings $A_0 \ra A_1 \ra ... \ra A_n$ of maps in ${\cat C}$. The poset ${\bf n}$ can be viewed as a category with objects $\{0, 1, ..., n\}$ and maps $n_1 \ra n_2$ for all $n_1 \le n_2$. $B{\cat C}_n$ can then be identified with the set of functors ${\bf n} \ra {\cat C}$, and we define the maps $({\bf n_1} \ra {\cat C}) \lra ({\bf n_2} \ra {\cat C})$ to be the commutative diagrams
\begin{center}
$ \xymatrix {
        {\bf n_2} \ar[rr] \ar[dr] & & {\bf n_1} \ar[dl] \\
	& {\cat C} &
	}$
\end{center}
The nerve functor $B{\cat C}$ commutes with arbitrary limits (has a left adjoint).

Denote $I$ the category ${\bf 1}$, with two objects and only one non-identity map. Its nerve $BI$ is isomorphic to $\Delta[1]$.

A functor $f \colon {\cat C}_1 \ra {\cat C}_2$ induces a map of simplicial sets $Bf \colon B{\cat C}_1 \ra B{\cat C}_2$. A natural map between two functors $h \colon f \Ra g$ induces a functor $H \colon {\cat C}_1 \times I \ra {\cat C}_2$, therefore a simplicial homotopy $BH \colon B{\cat C}_1 \times \Delta[1] \ra B{\cat C}_2$ between $Bf$ and $Bg$.

\bigskip

\section{Cofinal functors}
\label{sec:coffunctors1}
This section is a short primer on cofinal and homotopy cofinal functors. We only prove the minimal set of properties that we will need for the rest of the text. For more information on cofinal and homotopy cofinal functors, please refer to Hirschhorn \cite{Hirschhorn}.

\index{functor!left, right cofinal}
\index{functor!homotopy left, right cofinal}
\begin{defn}
\label{defn:coffunct}
A functor $u \colon {\cat D}_1 \ra {\cat D}_2$ beween small categories is
\begin{enumerate}
\item
{\it left cofinal} if for any object $d_2$ of ${\cat D}_2$ the space $B(u \downarrow d_2)$ is connected and non-empty,
\item
{\it homotopy left cofinal} if for any object $d_2$ of ${\cat D}_2$ the space $B(u \downarrow d_2)$ is contractible,
\item
{\it right cofinal} if for any object $d_2$ of ${\cat D}_2$ the space $B(d_2 \downarrow u)$ is connected and non-empty
\item
{\it homotopy right cofinal} if for any object $d_2$ of ${\cat D}_2$ the space $B(d_2 \downarrow u)$ is contractible
\end{enumerate}
\end{defn}

A homotopy left (resp. right) cofinal functor is in particular left (resp. right) cofinal. 

\begin{lem}
\label{lem::adjfunctcofinalfunct}
If $u : {\cat D}_1 \rightleftarrows {\cat D}_2 : v$ is an adjoint functor pair, then $u$ is homotopy left cofinal and $v$ is homotopy right cofinal.
\end{lem}
\begin{proof}
For any object $d_2$ of ${\cat D}_2$, the under category $(u \downarrow d_2)$ is isomorphic to $({\cat D}_1 \downarrow vd_2)$, and the latter has $id \colon vd_2 \ra vd_2$ as a terminal object. $B(u \downarrow d_2)$ is therefore contractible for any $d_2$, so $u$ is homotopy left cofinal.

A dual proof shows that $v$ is homotopy right cofinal.
\end{proof}

The following lemma gives a characterization of left (resp. right) cofinal functors in terms of their preservation of limits (resp. colimits).

\begin{prop}
\label{prop:cofinallimitconv}
Suppose that $u \colon {\cat D}_1 \ra {\cat D}_2$ is a functor of small categories. Then:
\begin{enumerate}
\item
$u$ is right cofinal iff for any cocomplete category ${\cat M}$ and any diagram $X \eps {\cat M}^{{\cat D}_2}$ the natural map $\colim^{{\cat D}_1}u^*X \ra \colim^{{\cat D}_2}X$ is an isomorphism in ${\cat M}$.
\item
$u$ is left cofinal iff for any complete category ${\cat M}$ and any diagram $X \eps {\cat M}^{{\cat D}_2}$ the natural map $\lim^{{\cat D}_2}X \ra \lim^{{\cat D}_1}u^*X$ is an isomorphism in ${\cat M}$.
\end{enumerate}
\end{prop}

Note the subtle positional difference of $u$ between the formulas of \lemmaref{lem::computerelcolim} (1) and \propositionref{prop:cofinallimitconv} (2). For a diagram $Y \eps {\cat M}^{{\cat D}_1}$ and an object $d_2 \eps {\cat D}_2$ we have that $(\colim^u Y)_{d_2} \cong \colim^{(u \downarrow d_2)}Y$. But for a diagram $X \eps {\cat M}^{{\cat D}_2}$, the natural map $\colim^{{\cat D}_1}u^*X \ra \colim^{{\cat D}_2}X$ is an isomorphism iff $(d_2 \downarrow u)$ is connected and non-empty for all objects $d_2 \eps {\cat D}_2$.

\begin{proof}[Proof of \propositionref{prop:cofinallimitconv}]
\mbox{}

(1) ($\La$). Let ${\cat M}$ be ${\cat Set}$ and $X$ be the ${\cat D}_2$-diagram ${{\cat D}_2}(d_2, -)$ for $d_2 \eps {\cat D}_2$. We notice that $\colim^{{\cat D}_2}X \cong \colim^{d^{'}_2 \eps {\cat D}_2} {{\cat D}_2}(d_2, d^{'}_2)$ is the one-point set, and $\colim^{{\cat D}_1}u^*X \cong \colim^{d_1 \eps {\cat D}_1} {{\cat D}_2}(d_2, ud_1)$ is the set of connected components of $B(d_2 \downarrow u)$, which must therefore be isomorphic with the one-point set. The conclusion is now proved.

(1) ($\Ra$). This will be a consequence of \propositionref{prop:cofinallimit} below. 

The proof of (2) is dual to the proof of (1).
\end{proof}

We next prove something a bit stronger than the right to left implication of \propositionref{prop:cofinallimitconv} (1).

\begin{prop}
\label{prop:cofinallimit}
Suppose that $u \colon {\cat D}_1 \ra {\cat D}_2$ is a functor, that ${\cat M}$ is a category and that $X$ an object of ${\cat M}^{{\cat D}_2}$.
\begin{enumerate}
\item
Assume that the functor $u$ is right cofinal. If either $\colim^{{\cat D}_1}u^*X$ or $\colim^{{\cat D}_2}X$ exist, then they both exist and the natural map $\colim^{{\cat D}_1}u^*X \ra \colim^{{\cat D}_2}X$ is an isomorphism.
\item
Assume that the functor $u$ is left cofinal. If either $\lim^{{\cat D}_1}u^*X$ or $\lim^{{\cat D}_2}X$ exist, then they both exist and the natural map $\lim^{{\cat D}_2}X \ra \lim^{{\cat D}_1}u^*X$ is an isomorphism.
\end{enumerate}
\end{prop}

\begin{proof}
We will prove (1); the statement (2) follows from duality. We denote $p_{{\cat D}_i} \colon {\cat D}_i \ra \terminalcat$, $i = 1, 2$ the terminal category projections.

We start by constructing a natural isomorphism in $X \eps {\cat M}^{{\cat D}_2}, Y \eps {\cat M}$
\begin{equation}
\label{eqn:cofinalcolimitiso}
{{\cat M}^{{\cat D}_1}}(u^*X, p_{{\cat D}_1}^*Y) \cong {{\cat M}^{{\cat D}_2}}(X, p_{{\cat D}_2}^*Y)
\end{equation}
From right to left, if $f \colon X \ra p_{{\cat D}_2}^*Y \eps {\cat M}^{{\cat D}_2}$ then we define $F(f) \colon u^*X \ra p_{{\cat D}_1}^*Y \eps {\cat M}^{{\cat D}_1}$ on component $d_1 \eps {\cat D}_1$ as $F(f)_{d_1} = f_{ud_1} \colon X_{ud_1} \ra Y$.

From left to right, if $f \colon u^*X \ra p_{{\cat D}_1}^*Y \eps {\cat M}^{{\cat D}_1}$ then we define $G(f) \colon X \ra p_{{\cat D}_2}^*Y \eps {\cat M}^{{\cat D}_2}$ on component $d_2 \eps {\cat D}_2$ as follows. Since $(d_2 \downarrow u)$ is non-empty we can pick a map $\alpha \colon d_2 \ra ud_1$, and we define $G(f)_{d_2} \colon X_{d_2} \ra X_{ud_1} \ra Y$ to be the composition of $X_\alpha$ with $f_{d1}$. Since $(d_2 \downarrow u)$ is connected, the map $G(f)$ does not depend on the choice involved, $G(f)$ is indeed a diagram map from $X$ to $p_{{\cat D}_2}^*Y$ and furthermore $F, G$ are inverses of each other.

If we view both terms of \ref{eqn:cofinalcolimitiso} as functors of $Y$, then $\colim^{{\cat D}_1}u^*X$ exists iff the left side of \ref{eqn:cofinalcolimitiso} is a representable functor of $Y$ iff the right side of \ref{eqn:cofinalcolimitiso} is a representable functor of $Y$ iff $\colim^{{\cat D}_2}X$ exists. Under these conditions it also follows that the natural map $\colim^{{\cat D}_1}u^*X \ra \colim^{{\cat D}_2}X$ is an isomorphism.
\end{proof}

As a consequence of \propositionref{prop:cofinallimitconv}, it is not hard to see that left (resp. right) cofinal functors are stable under composition. 

\section{The Grothendieck construction}
\label{sec:grothendieckconstruction}
\index{Grothendieck construction}
Denote ${\cat Cat}$ the category of small categories and functors. If ${\cat D}$ is a small category and $H$ is a functor ${\cat D}$ $\ra$ ${\cat Cat}$, the Grothendieck construction of $H$, denoted $\int_{\cat D} H$, is the category for which:
\begin{enumerate}
\item
objects are pairs $(d, x)$ with $d$ an object of ${\cat D}$ and $x$ an object of the category $H(d)$
\item 
maps $(d_1, x_1) \ra (d_2, x_2)$ are pairs of maps $(f, \phi)$ with $f \colon d_1 \ra d_2$ and $\phi \colon H(f)x_1 \ra x_2$
\item 
the composition of maps $(g, \psi)(f, \phi)$ is the map $(gf,  \psi \comp H(g) \phi)$
\item 
the identity of $(d, x)$ is $(1_d, 1_x)$
\end{enumerate}

The Grothendieck construction comes with a projection functor $p \colon \int_{\cat D} H \ra {\cat D}$, defined by $(d, x) \ra d$. If $d$ is an object of ${\cat D}$, then the inverse image category $p^{-1}d$ may be identified with the category $H(d)$. 

Dually if $H$ is a functor ${\cat D}^{op} \ra {\cat Cat}$, the contravariant Grothendieck construction of $H$ is defined as $\int^{\cat D} H = (\int_{{\cat D}^{op}} {H^{op}})^{op}$, and comes with a projection functor $p \colon \int^{\cat D} H \ra {\cat D}$.

\begin{prop}
\label{prop:grothendieckconstructionadjoint}
\mbox{}
\begin{enumerate}
\item
Suppose that $H \colon {\cat D} \ra {\cat Cat}$ is a functor from a small category to the category of small categories, with projection functor $p \colon \int_{\cat D} H \ra {\cat D}$. Then for any object $d \eps {\cat D}$, the inclusion $p^{-1}d \ra (p \downarrow d)$ admits a left adjoint.
\item
Suppose that $H \colon {\cat D}^{op} \ra {\cat Cat}$ is a functor from a small category to the category of small categories, with projection functor $p \colon \int^{\cat D} H \ra {\cat D}$. Then for any object $d \eps {\cat D}$, the inclusion $p^{-1}d \ra (d \downarrow p)$ admits a right adjoint.
\end{enumerate}
\end{prop}

\begin{proof}
To prove (1), denote $i_d \colon p^{-1}d \ra (p \downarrow d)$ the inclusion. A left adjoint $j_d \colon (p \downarrow d) \ra p^{-1}d$ of $i_d$ can be constructed such that:
\begin{enumerate}
\item[(a)]
$j_d$ sends the object $((d_1, x_1), f_1 \colon d_1 \ra d)$ of $(p \downarrow d)$ to the object $H(f_1)x_1$ of $H(d) \cong p^{-1}d$
\item[(b)]
For objects $((d_1, x_1), f_1 \colon d_1 \ra d)$ and $((d_2, x_2), f_2 \colon d_2 \ra d)$ of $(p \downarrow d)$, $j_d$ sends a map $(f, \phi) \colon (d_1, x_1) \ra (d_2, x_2)$ with $f_2f = f_1$ to a map $H(f_1)x_1 \ra H(f_2)x_2$ in $H(d)$ defined by $H(f_2)\phi$
\item[(c)]
The adjunction counit is $id \colon j_di_d \Ra 1_{p^{-1}d}$
\item[(d)]
The adjunction unit $1_{(p \downarrow d)} \Ra i_dj_d$ maps 
\begin{center}
$((d_1, x_1), f_1 \colon d_1 \ra d) \Ra ((d, H(f_1)x_1), 1_d \colon d \ra d)$
\end{center}
via the map $(d_1, x_1) \ra (d, H(f_1)x_1)$ in $\int_{\cat D} H$ defined by $(f_1, 1_{H(f_1)x_1})$.
\end{enumerate}
The proof of (2) is dual.
\end{proof}

As a consequence we can prove a Fubini-type formula for the computation of (co)limits indexed by Grothendieck constructions.

\begin{prop}
\label{prop:grothendieckconstrcolims}
\mbox{}
\begin{enumerate}
\item
Suppose that $H \colon {\cat D} \ra {\cat Cat}$ is a functor from a small category to the category of small categories. Let ${\cat M}$ be a category, and $X$ be a $\int_{\cat D} H$-diagram of ${\cat M}$. Assume that the inner colimit of the right side of the equation below exists for all $d \eps {\cat D}$. If either the left side or the outer right side colimits exist, then both exist and we have a natural isomorphism
\begin{center}
$\colim^{\int_{\cat D} H}X \cong \colim^{d \eps {\cat D}} \colim^{H(d)} X$
\end{center}
\item
Suppose that $H \colon {\cat D}^{op} \ra {\cat Cat}$ is a functor from a small category to the category of small categories. Let ${\cat M}$ be a category, and $X$ be a $\int^{\cat D} H$-diagram of ${\cat M}$. Assume that the inner limit of the right side of the equation below exists for all $d \eps {\cat D}$. If either the left side or the outer right side limits exist, then both exist and we have a natural isomorphism
\begin{center}
$\lim^{\int^{\cat D} H}X \cong \lim^{d \eps {\cat D}} \lim^{H(d)^{op}} X$
\end{center}
\end{enumerate}
\end{prop}

\begin{proof}
We only prove (1) - statement (2) uses a dual proof.

For any object $d$ of ${\cat D}$, the inclusion $p^{-1}d \ra (p \downarrow d)$ has a left adjoint, therefore by \lemmaref{lem::adjfunctcofinalfunct} it is homotopy right cofinal, and in particular right cofinal. We identify $p^{-1}d \cong H(d)$. Since $\colim^{H(d)}X$ exists for all $d$, by \propositionref{prop:cofinallimit} $\colim^{(p \downarrow d)}X \cong \colim^{H(d)}X$ exists for all $d$, and by \lemmaref{lem::computerelcolim} $\colim^p X$ exists. By \lemmaref{lem::colimcomposition}, if either $\colim^{d \eps {\cat D}} \colim^{H(d)} X \cong \colim^{\cat D}\colim^p X$ or $\colim^{\int_{\cat D} H}X$ exist, then they both exist and are canonically isomorphic.
\end{proof}

If $H \colon {\cat D}_1$ $\ra$ ${\cat Cat}$ is a constant functor with value ${\cat D}_2$, then the Grothendieck construction $\int_{\cat D} H$ is {\it isomorphic} to the product of categories ${\cat D}_1 \times {\cat D}_2$. In this case, \propositionref{prop:grothendieckconstrcolims} yields

\begin{cor}
\label{cor:grothendieckconstrcolims}
Suppose that ${\cat D}_1$, ${\cat D}_2$ are two small categories and suppose that $X$ is a ${\cat D}_1 \times {\cat D}_2$ diagram in a category ${\cat M}$.
\begin{enumerate}
\item
Assume that the inner colimit of the right side of the equation below exists for all $d_1 \eps {\cat D}_1$. If either the left side or the outer right side colimits exist, then both exist and we have a natural isomorphism
\begin{center}
$\colim^{{\cat D}_1 \times {\cat D}_2}X \cong \colim^{d_1 \eps {\cat D}_1} \colim^{\{d_1\} \times {\cat D}_2} X$
\end{center}
\item
Assume that the inner limit of the right side of the equation below exists for all $d_1 \eps {\cat D}_1$. If either the left side or the outer right side limits exist, then both exist and we have a natural isomorphism
\begin{center}
$\lim^{{\cat D}_1 \times {\cat D}_2}X \cong \lim^{d_1 \eps {\cat D}_1} \lim^{\{d_1\} \times {\cat D}_2} X$ $\square$
\end{center}
\end{enumerate}
\end{cor}

\begin{rem}
\label{rem:grothendieckconstrcolims}
A typical application of \corollaryref{cor:grothendieckconstrcolims} (1) will be for the case when $\colim^{\{d_1\} \times {\cat D}_2} X$ and $\colim^{{\cat D}_1 \times \{d_2\}} X$ exist for all objects $d_1 \eps {\cat D}_1$ and $d_2 \eps {\cat D}_2$. Then in the equation below if either the left or the right outer side colimits exist, they both exist and we have a natural isomorphism
\begin{center}
$\colim^{d_1 \eps {\cat D}_1} \colim^{\{d_1\} \times {\cat D}_2} X \cong \colim^{d_2 \eps {\cat D}_2} \colim^{{\cat D}_1 \times \{d_2\}} X$ 
\end{center}
\end{rem}

\chapter{Homotopy colimits in a cofibration category}
\label{chap:limcolimcofinal}

In this chapter we prove the existence of homotopy colimits in cofibration categories, and dually the existence of homotopy limits in fibration categories. The homotopy colimit should be thought of as the total left derived functor of the colimit - at least if the base cofibration category is cocomplete. The actual construction of the homotopy colimit proceeds in two steps - first, for diagrams of direct categories, and second, reducing the general case to the case of diagrams of direct categories. This essentially follows Anderson's original argument \cite{Anderson1}, simplified by Cisinski \cite{Cisinski2}, \cite{Cisinski4}.

Throughout this chapter, we will assume the entire set of cofibration category axioms CF1-CF6. An exercise left for the reader is to see which of the results can be reformulated and proved within the restricted cofibration category axioms CF1-CF5, or even within the precofibration category axioms CF1-CF4.

We show that given a cofibration category ${\cat M}$ and a small direct category ${\cat D}$, the diagram category ${\cat M}^{\cat D}$ with pointwise weak equivalences ${\cat W}^{\cat D}$ admits two cofibration category structures - the Reedy $({\cat M}^{\cat D}, {\cat W}^{\cat D}, {\cat Cof}^{\cat D}_{Reedy})$ and the pointwise structure $({\cat M}^{\cat D}, {\cat W}^{\cat D}, {\cat Cof}^{\cat D})$. A general small category ${\cat D}$ only yields a pointwise cofibration structure $({\cat M}^{\cat D}, {\cat W}^{\cat D}, {\cat Cof}^{\cat D})$, but ${\cat M}^{\cat D}$ admits a certain cofibrant approximation functor that allows us to reduce the construction of homotopy colimits in ${\cat M}^{\cat D}$ to the construction of homotopy colimits indexed by direct diagrams in ${\cat M}$.

We then proceed to prove a number of properties of homotopy limits and of ${\bf ho}{\cat M}^{\cat D}$ that will allow us to show in \chapterref{chap:derivators} that that homotopy colimits in a cofibration category satisfy the axioms of a left Heller derivator. Dually, homotopy limits in a fibration category satisfy the axioms of a right Heller derivator.

\section{Direct and inverse categories}
\label{sec:catdirinv}
If we start with a small {\it direct} category ${\cat D}$, then homotopy colimits in ${\cat M}^{\cat D}$ for a cofibration category ${\cat M}$ can be constructed using just cofibrant replacements, small sums and pushouts. Dually, if ${\cat D}$ is a small {\it inverse} category and ${\cat M}$ is a fibration category, homotopy limits in ${\cat M}^{\cat D}$ can be constructed using fibrant replacements, small products and pullbacks.

Here are a few definitions that we will need.

\index{category!direct, inverse}
\begin{defn}[Direct and inverse categories]
\label{def:directinvcat}
\mbox{}
Let ${\cat D}$ be a category. A non-negative degree function on its objects is a function $deg \colon Ob {\cat D} \ra {\mathbb Z}_+$.
\begin{enumerate}
\item
The category ${\cat D}$ is {\it direct} if there is a non-negative degree function $deg$ on the objects of ${\cat D}$ such that any non-identity map $d_1 \ra d_2$ satisfies $deg(d_1) < deg(d_2)$.
\item
The category ${\cat D}$ is {\it inverse}  if there is a non-negative degree function $deg$ on the objects of ${\cat D}$ such that any non-identity map $d_1 \ra d_2$ satisfies $deg(d_1) > deg(d_2)$.
\end{enumerate}
\end{defn}

\index{latching and matching!category}
\begin{defn}[Latching and matching categories]
\label{def:latchingmatchingcat}
Let $d$ be an object of a category ${\cat D}$.
\begin{enumerate}
\item
If ${\cat D}$ is direct, the {\it latching category} $\partial({\cat D} \downarrow d)$ is the full subcategory of the over category $({\cat D} \downarrow d)$ consisting of all objects except the identity of $d$. 
\item
If ${\cat D}$ is inverse, the {\it matching category} $\partial(d \downarrow {\cat D})$ is the full subcategory of the under category $(d \downarrow {\cat D})$ consisting of all objects except the identity of $d$.
\end{enumerate}
\end{defn}

If ${\cat D}$ is a direct category, then $\partial({\cat D} \downarrow d)$ is also direct with $deg(d^{'} \ra d) = deg (d^{'})$. All the objects of $\partial({\cat D} \downarrow d)$ have therefore degree $< deg(d)$. For any map $f \colon d^{'} \ra d$, forgetting codomains of maps yields a canonical isomorphism of categories 
\begin{equation}
\label{eqn:latchlatch}
\partial(\partial({\cat D} \downarrow d) \downarrow f) \cong \partial ({\cat D} \downarrow d^{'})
\end{equation}
between the latching category of $\partial({\cat D} \downarrow d)$ at $f$ and the latching category of ${\cat D}$ at $d^{'}$.

Dually, if the category ${\cat D}$ is inverse, then $\partial(d \downarrow {\cat D})$ is an inverse category whose objects have all degree $< deg(d)$. For any map $f \colon d \ra d^{'}$, forgetting domains of maps yields an isomorphism of matching categories
\begin{equation}
\label{eqn:matchmatch}
\partial( f \downarrow \partial(d \downarrow {\cat D})) \cong \partial (d^{'} \downarrow {\cat D})
\end{equation}

\index{latching and matching!object}
\begin{defn}[Latching and matching objects]
\label{def:latchingmatchingobject}
Suppose that ${\cat M}$ is a category, that $X$ is a ${\cat D}$ -diagram of ${\cat M}$ and that $d$ is an object of ${\cat D}$.
\begin{enumerate}
\item
If ${\cat D}$ is direct, the {\it latching object} of $X$ at $d$, if it exists, is by definition 
\begin{center}
$LX_d = \colim^{\partial({\cat D} \downarrow d)}X$
\end{center}
\item
If ${\cat D}$ is inverse, the {\it matching object} of $X$ at $d$, if it exists, is by definition 
\begin{center}
$MX_d = \lim^{\partial(d \downarrow {\cat D})}X$
\end{center}
\end{enumerate}
\end{defn}

The latching object $LX_d$, if it exists, is therefore defined only up to an unique isomorphism. If the category ${\cat M}$ is cocomplete, then $LX_d$ always exists.

If ${\cat D}$ is direct or inverse and $n \eps {\mathbb Z}_+$ then we will denote ${\cat D}^{<n}, {\cat D}^{\leq n}$ the full subcategories with objects of degree $< n$ respectively $\leq n$. If ${\cat D}$ is direct, then $\partial({\cat D} \downarrow d) \cong$ $\partial({\cat D}^{\leq deg(d)} \downarrow d) \cong$ $(\partial({\cat D} \downarrow d))^{< deg(d)}$. 

Suppose that $X$ is a ${\cat D}$-diagram in ${\cat M}$, with ${\cat D}$ direct, and that $f \colon d^{'} \ra d$ is a map in ${\cat D}$. Using \eqref{eqn:latchlatch}, if either $L(X|_{\partial({\cat D} \downarrow d)})_f$ or $LX_{d^{'}}$ exist, then they both exist and are canonically isomorphic.

Latching objects are related to Kan extensions in the following sense. Denote $\delta_d \colon {\cat D}^{<deg(d)} \ra {\cat D}$ the inclusion functor. If ${\cat M}$ is a cocomplete category and $X$ is a ${\cat D}$ -diagram of ${\cat M}$, then the latching object $LX_d$ is isomorphic to $(\colim^{\delta_d} X)_d$.

Dually, if the category ${\cat D}$ is inverse, then $\partial(d \downarrow {\cat D})$ is an inverse category with $deg(d \ra d^{'}) = deg (d^{'})$, and $\partial(d \downarrow {\cat D}) \cong \partial(d \downarrow {\cat D}^{\leq deg(d)}) \cong (\partial(d \downarrow {\cat D}))^{< deg(d)}$. If $X$ is a ${\cat D}$ -diagram of ${\cat M}$ and if $f \colon d \ra d^{'}$ is a map in ${\cat D}$, if either matching space $M(X|_{\partial(d \downarrow {\cat D})})_f$ or $MX_{d^{'}}$ exist then they both exist and are canonically isomorphic. If ${\cat M}$ is furthermore a complete category, then $MX_d \cong (\lim^{\delta_d} X)_d$.

\section{Reedy and pointwise cofibration structures for direct diagrams}
\label{sec:reedycofstruct}

Gived a cofibration category ${\cat M}$ and a small direct category ${\cat D}$, we define two cofibration category structures on ${\cat M}^{\cat D}$ - the Reedy and the pointwise cofibration structure. The pointwise cofibration structure has pointwise weak equivalences and pointwise cofibrations. The Reedy cofibration structure also has pointwise weak equivalences, but has a more restrictive set of cofibrations - these are called the Reedy cofibrations. In particular, Reedy cofibrations are pointwise cofibrations.

The axioms for the pointwise cofibration structure are easily verified except for the factorization axiom. To prove the factorization axiom, we will construct the Reedy cofibration structure, prove all the axioms (including the factorization axiom) for the Reedy cofibration structure, and obtain as a corollary the factorization axiom for the pointwise cofibration structure.

Dually, if ${\cat M}$ is a fibration category and ${\cat D}$ is a small inverse category then we define two fibration category structures on ${\cat M}^{\cat D}$ - the pointwise and the Reedy fibration category structures. 

We start with the definition of pointwise weak equivalences, pointwise cofibrations and pointwise fibrations:

\index{pointwise!(co)fibration}
\index{pointwise!weak equivalence}
\begin{defn}
\label{def:pointwisecof}
Let ${\cat D}$ be a small category. Given a cofibration (respectively fibration) category ${\cat M}$, a map of ${\cat D}$-diagrams $X \ra Y$ in ${\cat M}^{\cat D}$ is a pointwise weak equivalence (resp. a pointwise cofibration, resp. a pointwise fibration) if all maps $X_d \ra Y_d$ are weak equivalences (resp. cofibrations, resp. fibrations) for all objects $d$ of ${\cat D}$.
\end{defn}

A diagram $X$ is therefore pointwise cofibrant (resp. pointwise fibrant) if all $X_d$ are cofibrant (resp. fibrant) for all objects $d$ of ${\cat D}$.

The category of pointwise weak equivalences is then just ${\cat W}^{\cat D}$. The category of pointwise cofibrations is ${\cat Cof}^{\cat D}$. The category of pointwise fibrations is ${\cat Fib}^{\cat D}$.

We continue with the definition of Reedy cofibrations and Reedy fibrations:

\index{Reedy (co)fibration}
\begin{defn}
\label{def:reedycof}
\mbox{}
\begin{enumerate}
\item
Let ${\cat M}$ be a cofibration category and ${\cat D}$ be a small direct category.
\begin{enumerate}
\item
A ${\cat D}$-diagram $X$ of ${\cat M}$ is called {\it Reedy cofibrant} if for any object $d$ of ${\cat D}$, the latching object $LX_d$ exists and is cofibrant, and the natural map $i_d \colon LX_d \ra X_d$ is a cofibration.
\item
A map of Reedy cofibrant ${\cat D}$-diagrams $f \colon X \ra Y$ is called a {\it Reedy cofibration} if for any object $d$ of ${\cat D}$, the natural map $X_d \Sum_{LX_d} LY_d \ra Y_d$ is a cofibrantion. (Notice that the pushout $X_d \Sum_{LX_d} LY_d$ always exists if $X, Y$ are Reedy cofibrant because of the pushout axiom.) The class of Reedy cofibrations will be denoted ${\cat Cof}^{\cat D}_{Reedy}$.
\end{enumerate}
\item
Let ${\cat M}$ be a fibration category and ${\cat D}$ be a small inverse category.
\begin{enumerate}
\item
A ${\cat D}$-diagram $X$ of ${\cat M}$ is called {\it Reedy fibrant} if for any object $d$ of ${\cat D}$, the matching object $MX_d$ exists and is fibrant, and the natural map $p_d \colon X_d \ra MX_d$ is a fibration.
\item
A map of Reedy fibrant ${\cat D}$-diagrams $f \colon X \ra Y$ is called a {\it Reedy fibration} if for any object $d$ of ${\cat D}$, the natural map $X_d \ra MX_d \times_{MY_d} Y_d$ is a fibrantion. The class of Reedy cofibrations will be denoted ${\cat Fib}^{\cat D}_{Reedy}$.
\end{enumerate}
\end{enumerate}
\end{defn}

We'd like to stress that a Reedy cofibration $X \ra Y$ has by definition a Reedy cofibrant domain $X$. This is required because our cofibration categories are not necessarily cocomplete. Dually, a Reedy fibration $X \ra Y$ has by definition a Reedy fibrant codomain $Y$. In this regard, our definition of Reedy (co)fibrations is different than the usual one in Quillen model categories, which are by definition complete and cocomplete, where Reedy cofibrations (resp. fibrations) are allowed to have non-Reedy cofibrant domain (resp. non-Reedy fibrant codomain).

A pointwise cofibration $X \ra Y$ is not required to have a pointwise cofibrant domain $X$.

If ${\cat M}$ is a cofibration category and ${\cat D}$ is a small direct category, notice that the constant initial-object ${\cat D}$-diagram $c\initial$ is Reedy cofibrant (as well as pointwise cofibrant). A ${\cat D}$-diagram $X$ is Reedy cofibrant iff the map $c\initial \ra X$ is a Reedy cofibration. 

If $f \colon X \ra Y$ is a Reedy cofibration in ${\cat M}^{\cat D}$ for ${\cat D}$ direct, then on account of \eqref{eqn:latchlatch} for any object $d \eps {\cat D}$ the restriction of $f$ as a diagram map over $\partial ({\cat D} \downarrow d)$ is also a Reedy cofibration. 

\begin{prop}
\label{prop::reedycofformcat}
\mbox{}
\begin{enumerate}
\item
Let $({\cat M}$, ${\cat W}$, ${\cat Cof})$ be a cofibration category and ${\cat D}$ be a small direct category. Then the Reedy cofibrations ${\cat Cof}^{\cat D}_{Reedy}$ are stable under composition, and include the isomorphisms with a Reedy cofibrant domain.
\item
Let $({\cat M}$, ${\cat W}$, ${\cat Fib})$ be a fibration category and ${\cat D}$ be a small inverse category. Then the Reedy fibrations ${\cat Fib}^{\cat D}_{Reedy}$ are stable under compositions, and include the isomorphisms with a Reedy fibrant codomain.
\end{enumerate}
\end{prop}

\begin{proof}
We only prove (1). If $X \ra Y \ra Z$ is a composition of Reedy cofibrations in ${\cat M}^{\cat D}$, then for any object $d$ of ${\cat D}$ we factor $X_d \Sum_{LX_d} LZ_d \ra Z_d$ as the composition 
\begin{center}
$X_d \Sum_{LX_d} LZ_d \ra Y_d \Sum_{LY_d} LZ_d \ra Z_d$
\end{center}
The second map is a cofibration since $Y \ra Z$ is a Reedy cofibration. The first map is a cofibration as well, as a pushout of the cofibration $X_d \Sum_{LX_d} LY_d \ra Y_d$ along $X_d \Sum_{LX_d} LY_d \ra X_d \Sum_{LX_d} LZ_d$, where $X_d \Sum_{LX_d} LZ_d$ is a cofibrant object. We deduce that $X \ra Z$ is a Reedy cofibration.

If $X \ra Y$ is an isomorphism in ${\cat M}^{\cat D}$ with $X$ Reedy cofibrant, then $Y$ is Reedy cofibrant as well. The latching maps $X_d \Sum_{LX_d} LY_d \ra Y_d$ are isomorphisms with cofibrant domain therefore cofibrations, so $X\ra Y$ is a Reedy cofibration.
\end{proof}

\begin{thm}[Reedy and pointwise (co)fibration structures]
\label{thm:reedycofibstruct}
\mbox{}
\begin{enumerate}
\item
If $({\cat M}$, ${\cat W}$, ${\cat Cof})$ is a cofibration category and ${\cat D}$ is a small direct category, then
\begin{enumerate}
\item
$({\cat M}^{\cat D}, {\cat W}^{\cat D}, {\cat Cof}^{\cat D}_{Reedy})$ is a cofibration category - called the Reedy cofibration structure on ${\cat M}^{\cat D}$.
\item
$({\cat M}^{\cat D}, {\cat W}^{\cat D}, {\cat Cof}^{\cat D})$ is a cofibration category - called the pointwise cofibration structure on ${\cat M}^{\cat D}$.
\end{enumerate}
\item
If $({\cat M}$, ${\cat W}$, ${\cat Fib})$ is a fibration category and ${\cat D}$ is a small inverse category, then
\begin{enumerate}
\item
$({\cat M}^{\cat D}, {\cat W}^{\cat D}, {\cat Fib}^{\cat D}_{Reedy})$ is a fibration category - called the Reedy fibration structure on ${\cat M}^{\cat D}$.
\item
$({\cat M}^{\cat D}, {\cat W}^{\cat D}, {\cat Fib}^{\cat D})$ is a fibration category - called the pointwise fibration structure on ${\cat M}^{\cat D}$.
\end{enumerate}
\end{enumerate}
\end{thm}

The proof of this theorem is deferred until the next section, after we work out some basic properties of colimits and limits. 

\section{Colimits in direct categories (the absolute case)}
\label{sec:colimdirectabsolute}

The Lemmas \ref{lem:colimpushout}, \ref{lem:colimpushoutset} below provide the essential inductive step required for building colimits of Reedy cofibrant diagrams in a direct category. 

The colimits of any two weakly equivalent Reedy cofibrant diagrams are weakly equivalent, cf. \theoremref{thm:colimdirect} below. In that sense, the colimit of a Reedy cofibrant diagram actually computes its {\it homotopy} colimit. In view of this, the meaning of \lemmaref{lem:colimpushoutset} below is that homotopy colimits over a direct category can be constructed, after a Reedy cofibrant replacement, as iterated pushouts of small sums of latching object maps.

If $\{d_k | k \eps K\}$ is a set of objects of a category ${\cat D}$, denote ${\cat D} \backslash \{d_k | k \eps K \}$ the maximal full subcategory of a category ${\cat D}$ without the objects $d_k$. Our first lemma applies to the case of a direct (resp. inverse) category ${\cat D}$ and an object $d$ of ${\cat D}$ such that ${\cat D} \backslash \{d\} \ra {\cat D}$ is an open (resp. closed) embedding.

For example, if ${\cat D}$ is a direct category with all objects of degree $\leq n$ and $d$ is an object of ${\cat D}$ of degree $n$ then ${\cat D} \backslash \{d\} \ra {\cat D}$ is an open embedding. Dually, if ${\cat D}$ an inverse category with all objects of degree $\leq n$ and $d$ is an object of ${\cat D}$ of degree $n$ then ${\cat D} \backslash \{d\} \ra {\cat D}$ is a closed embedding.

\begin{lem}
\label{lem:colimpushout}
\mbox{}
\begin{enumerate}
\item
Let ${\cat M}$ be a cofibration category. Let ${\cat D}$ be a direct category, $d$ an object of ${\cat D}$ such that ${\cat D} \backslash \{d\} \ra {\cat D}$ is an open embedding, and $X$ be a ${\cat D}$-diagram of ${\cat M}$. Assume that $LX_d$ and $\colim^{{\cat D} \backslash \{d\}} X$ exist and are cofibrant, and that $i_{d} \colon LX_d \ra X_d$ is a cofibration. Denote $f_d$ the colimit map induced by ${\cat D} \backslash \{d\} \ra {\cat D}$.
\begin{center}
$ \xymatrix {
        LX_d \ar[r]^-{f_d} \ar@{>->}[d]_{i_d} & 
	\colim^{{\cat D} \backslash \{d\}} X \ar@{>-->}[d]^{j_d} \\
	X_d \ar@{-->}[r] & 
	\colim^{\cat D} X
	}$
\end{center}
Then the pushout of ($i_d, f_d$) is isomorphic to $\colim^{\cat D} X$. 

In particular $\colim^{\cat D} X$ exists and is cofibrant, and $j_d$ is a cofibration.
\item
Let ${\cat M}$ be a fibration category. Let ${\cat D}$ be an inverse category, $d$ an object of ${\cat D}$ such that ${\cat D} \backslash \{d\} \ra {\cat D}$ is a closed embedding, and $X$ be a ${\cat D}$-diagram of ${\cat M}$. Assume that $MX_d$ and $\lim^{{\cat D} \backslash \{d\}} X$ exist and are fibrant, and that $p_{d} \colon X_d \ra MX_d$ is a fibration. Denote $g_d$ the limit map induced by ${\cat D} \backslash \{d\} \ra {\cat D}$.
\begin{center}
$ \xymatrix {
	\lim^{\cat D} X \ar@{-->}[r] \ar@{-->>}[d]_{q_d} &
	X_d \ar@{->>}[d]^{p_d} \\
	\lim^{{\cat D} \backslash \{d\}} \ar[r]_-{g_d} & 
        MX_d
	}$
\end{center}
Then the pullback of ($p_d, g_d$) is isomorphic to $\lim^{\cat D} X$.

In particular $\lim^{\cat D} X$ exists and is fibrant, and $q_d$ is a fibration.
\end{enumerate}

\end{lem}

\begin{proof}
We will prove (1) - the proof of (2) is dual.

The pushout of (1) exists from axiom CF3. The latching object $LX_d$ is by definition $\colim^{\partial({\cat D} \downarrow d)}X$. Since ${\cat D} \backslash \{d\} \ra {\cat D}$ is an open embedding, the pushout of ($i_d, f_d$) satisfies the universal property that defines $\colim^{\cat D} X$. The map $j_d$ is a cofibration as the pushout of $i_d$, and therefore $\colim^{\cat D} X$ is cofibrant.
\end{proof}

Here is a slightly more general statement of the previous lemma:

\begin{lem}
\label{lem:colimpushoutset}
\mbox{}
\begin{enumerate}
\item
Let ${\cat M}$ be a cofibration category. Let ${\cat D}$ be a direct category, $\{d_k | k \eps K\}$ a set object of ${\cat D}$ of the same degree such that ${\cat D} \backslash \{d_k | k \eps K\} \ra {\cat D}$ is an open embedding, and $X$ be a ${\cat D}$-diagram of ${\cat M}$. Assume that $LX_{d_k}$ for all $k$ and $\colim^{{\cat D} \backslash \{d_k | k \eps K\}} X$ exist and are cofibrant, and that $i_{k} \colon LX_{d_k} \ra X_{d_k}$ is a cofibration for all $k$.  Denote $f_K$ the colimit map induced by ${\cat D} \backslash \{d_k | k \eps K\} \ra {\cat D}$.
\begin{center}
$ \xymatrix {
        \Sum LX_{d_k} \ar[r]^-{f_K} \ar@{>->}[d]_{\Sum i_k} & 
	\colim^{{\cat D} \backslash \{d_k | k \eps K\}} X \ar@{>-->}[d]^{j_K} \\
	\Sum X_{d_k} \ar@{-->}[r] & 
	\colim^{\cat D} X
	}$
\end{center}
Then the pushout of ($\Sum i_{k}, f_K$) is isomorphic to $\colim^{\cat D} X$. 

In particular $\colim^{\cat D} X$ exists and is cofibrant, and $j_K$ is a cofibration.
\item
Let ${\cat M}$ be a fibration category. Let ${\cat D}$ be an inverse category, $\{d_k | k \eps K\}$ a set object of ${\cat D}$ of the same degree such that ${\cat D} \backslash \{d_k | k \eps K\} \ra {\cat D}$ is a closed embedding, and $X$ be a ${\cat D}$-diagram of ${\cat M}$. Assume that $MX_{d_k}$ and $\lim^{{\cat D} \backslash \{d_k | k \eps K\}} X$ exist and are fibrant, and that $p_{k} \colon X_{d_k} \ra MX_{d_k}$ is a fibration.  Denote $g_K$ the limit map induced by ${\cat D} \backslash \{d_k | k \eps K\} \ra {\cat D}$.
\begin{center}
$ \xymatrix {
	\lim^{\cat D} X \ar@{-->}[r] \ar@{-->>}[d]_{q_K} &
	\times X_{d_k} \ar@{->>}[d]^{\times p_{k}} \\
	\lim^{{\cat D} \backslash \{d_k | k \eps K\}} X \ar[r]_-{g_K} & 
        \times MX_{d_k}
	}$
\end{center}
Then the pullback of ($\times p_k, g_K$) is isomorphic to $\lim^{\cat D} X$.

In particular $\lim^{\cat D} X$ exists and is fibrant, and $q_K$ is a fibration.
\end{enumerate}
\end{lem}

\begin{proof}
We only prove (1). By axiom CF5, the map $\Sum i_k$ is a cofibration, and the pushout of ($\Sum i_{k}, f_K$) exists. Since ${\cat D} \backslash \{d_k | k \eps K\} \ra {\cat D}$ is an open embedding, the pushout satisfies the universal property that defines $\colim^{\cat D} X$.
\end{proof}

We will also need the following two lemmas.

\begin{lem}
\label{lem::compospushout}
\mbox{}
\begin{enumerate}
\item
Let ${\cat M}$ be a cofibration category, and consider a commutative diagram
\begin{center}
$\xymatrix{
    A_0 \ar@{>->}[r]^{a_{01}} \ar[d]_{f_0} & A_1 \ar@{>->}[r]^{a_{12}} \ar[d]_{f_1} & A_2 \ar[d]_{f_2} \\
    B_0 \ar[r]^{b_{01}} & B_1 \ar[r]^{b_{12}} & B_2
    }$
\end{center}
with $a_{01}, a_{12}$ cofibrations and $A_0, B_0, B_1, B_2$ cofibrant.
\begin{enumerate}
\item
If $B_0 \Sum_{A_0} A_1 \ra B_1$ and $B_1 \Sum_{A_1} A_2 \ra B_2$ are cofibrations, then so is $B_0 \Sum_{A_0} A_2 \ra B_2$
\item
If $B_0 \Sum_{A_0} A_1 \ra B_1$ and $B_1 \Sum_{A_1} A_2 \ra B_2$ are weak equivalences, then so is $B_0 \Sum_{A_0} A_2 \ra B_2$
\end{enumerate}
\item
Let ${\cat M}$ be a fibration category, and consider a commutative diagram
\begin{center}
$\xymatrix{
    A_2 \ar[r]^{a_{21}} \ar[d]_{f_2} & A_1 \ar[r]^{a_{10}} \ar[d]_{f_1} & A_0 \ar[d]_{f_0} \\
    B_2 \ar@{->>}[r]^{b_{21}} & B_1 \ar@{->>}[r]^{b_{10}} & B_0
    }$
\end{center}
with $b_{10}, b_{21}$ fibrations and $A_0, A_1, A_2, B_0$ fibrant.
\begin{enumerate}
\item
If $A_1 \ra B_1 \times_{B_0} A_0$, $A_2 \ra B_2 \times_{B_1} A_1$ are fibrations, then so is $A_2 \ra B_2 \times_{B_0} A_0$
\item
If $A_1 \ra B_1 \times_{B_0} A_0$, $A_2 \ra B_2 \times_{B_1} A_1$ are weak equivalences, then so is $A_2 \ra B_2 \times_{B_0} A_0$
\end{enumerate}
\end{enumerate}
\end{lem}

\begin{proof}
We only prove (1). The pushouts $B_0 \Sum_{A_0} A_1$, $B_1 \Sum_{A_1} A_2$ and $B_0 \Sum_{A_0} A_2$ exist since $a_{01}$ and $a_{12}$ are cofibrations and $A_0,, B_0, B_1, B_2$ are cofibrant. In the diagram
\begin{center}
$\xymatrix{
    A_0 \ar@{>->}[r]^{a_{01}} \ar[d]_{f_0} & A_1 \ar@{>->}[r]^{a_{12}} \ar[d] & A_2 \ar[d] \\
    B_0 \ar[rd]_{b_{01}} \ar@{>->}[r] & B_0 \Sum_{A_0} A_1 \ar@{>->}[r]^-b \ar[d]_a & B_0 \Sum_{A_0} A_2 \ar[d]^c \\
    & B_1 \ar[rd]_{b_{12}} \ar@{>->}[r] & B_1 \Sum_{A_1} A_2 \ar[d]^d \\
    & & B_2 \\
    }$
\end{center}
the map $b$ is a cofibration, as a pushout of $a_{12}$.

If $a$ and $d$ are cofibrations, then $c$ is a cofibration as a pushout of $a$, therefore $dc$ is a cofibration. This proves part (a).

If $a$ and $d$ are weak equivalences, then by excision $c$ is a weak eqivalence, and therefore $dc$ is a weak equivalence. This proves part (b).
\end{proof}

\begin{lem}
\label{lem::compospushout2}
\mbox{}
\begin{enumerate}
\item
Let ${\cat M}$ be a cofibration category, and consider a map of countable direct sequences of cofibrations
\begin{center}
$ \xymatrix {
        A_0 \ar@{>->}[r]^{a_{0}} \ar[d]_{f_0} &
        A_1 \ar@{>->}[r]^{a_{1}} \ar[d]_{f_1} &
	... \ar@{>->}[r] &
        A_n \ar[d]_{f_n} \ar@{>->}[r]^{a_n} & 
        ... & \\
        B_0 \ar@{>->}[r]^{b_{0}} &
        B_1 \ar@{>->}[r]^{b_{1}} &
	... \ar@{>->}[r] &
        B_n \ar@{>->}[r]^{b_n} & 
        ... &
	}$
\end{center}
with $A_0$, $B_0$ cofibrant and $a_{n}, b_{n}$ cofibrations for $n \ge 0$. 
\begin{enumerate}
\item
If all maps $B_{n-1} \Sum_{A_{n-1}} A_{n} \ra B_{n}$ are cofibrations, then 
\begin{center}
$B_0 \Sum_{A_0} \colim A_n \ra \colim B_n$ 
\end{center}
is a cofibration.
\item
If all maps $B_{n-1} \Sum_{A_{n-1}} A_{n} \ra B_{n}$ are weak equivalences, then 
\begin{center}
$B_0 \Sum_{A_0} \colim A_n \ra \colim B_n$ 
\end{center}
is a weak equivalence.
\end{enumerate}
\item
Let ${\cat M}$ be a fibration category and consider a map of countable inverse sequences of fibrations
\begin{center}
$ \xymatrix {
    & 
    ... \ar@{->>}[r]^{a_n} & 
    A_n \ar@{->>}[r] \ar[d]^{f_n} & 
    ... \ar@{->>}[r]^{a_{1}} & \
    A_1 \ar@{->>}[r]^{a_{0}} \ar[d]^{f_1} &
    A_0 \ar[d]^{f_0} \\
    &
    ... \ar@{->>}[r]^{b_n} &
    B_n \ar@{->>}[r] &
    ... \ar@{->>}[r]^{b_{1}} &
    B_1 \ar@{->>}[r]^{b_{0}} &
    B_0
  }$
\end{center}
with $A_0$, $B_0$ fibrant and $a_{n}, b_{n}$ fibrations for $n \ge 0$.
\begin{enumerate}
\item
If all maps $A_{n} \ra B_{n} \times_{B_{n-1}} A_{n-1}$ are fibrations, then 
\begin{center}
$\lim A_n \ra \lim B_n \times_{B_0} A_0$
\end{center}
is a fibration.
\item
If all maps $A_{n} \ra B_{n} \times_{B_{n-1}} A_{n-1}$ are weak equivalences, then 
\begin{center}
$\lim A_n \ra \lim B_n \times_{B_0} A_0$
\end{center}
is a weak equivalence.
\end{enumerate}
\end{enumerate}
\end{lem}

\begin{proof}
We only prove(1). Denote $A = \colim A_n$ and $B = \colim B_n$.

The map $B_0 \Sum_{A_0} A \ra B$ factors as the composition of the direct sequence of maps $B_{n-1} \Sum_{A_{n-1}} A \ra B_{n} \Sum_{A_{n}} A$, and each map in the sequence is a pushout of $B_{n-1} \Sum_{A_{n-1}} A_{n} \ra B_{n}$ along the cofibration $B_{n-1} \Sum_{A_{n-1}} A_{n} \ra B_{n-1} \Sum_{A_{n-1}} A$.

\begin{center}
$\xymatrix{
    A_0 \ar@{>->}[r]^{a_{0}} \ar[d]_{f_0} & A_1 \ar@{>->}[r]^{a_{1}} \ar[d] & ... \ar[r] & A \ar[d] \\
    B_0 \ar[rd]_{b_{0}} \ar@{>->}[r] & B_0 \Sum_{A_0} A_1 \ar@{>->}[r] \ar[d] & ... \ar[r] & B_0 \Sum_{A_0} A \ar[d] \\
    & B_1 \ar@{>->}[r] \ar@{>->}[rd]_{b_1} & ... \ar[r] & B_1 \Sum_{A_1} A \ar[d] \\
    & & ... \ar[rd] & ... \ar[d] \\
    &&& B
    }$
\end{center}

For part (a), each map $B_{n-1} \Sum_{A_{n-1}} A_{n} \ra B_{n}$ is a cofibration, therefore by CF3 so is its pushout $B_{n-1} \Sum_{A_{n-1}} A \ra B_{n} \Sum_{A_{n}} A$, and by CF6 so is the sequence composition $B_0 \Sum_{A_0} A \ra B$.

For part (b), consider the map of direct sequences of cofibrations $\phi_n \colon B_0 \Sum_{A_0} A_n \ra B_n$. Each map $\phi_n$ is a weak equivalence using excision, and the result folows from \lemmaref{lem:transcompequiv}.
\end{proof}

We are ready to state the following result.
\begin{thm}
\label{thm:colimdirect}
\mbox{}
\begin{enumerate}
\item
Let ${\cat M}$ be a cofibration category and ${\cat D}$ be a small direct category. Then:
\begin{enumerate}
\item
If $X$ is Reedy cofibrant in ${\cat M}^{\cat D}$, then $\colim^{\cat D} X$ exists and is cofibrant in ${\cat M}$
\item
If $f \colon X \ra Y$ is a Reedy cofibration in ${\cat M}^{\cat D}$, then $\colim^{\cat D} f$ is a cofibration in ${\cat M}$
\item
If $f \colon X \ra Y$ is a pointwise weak equivalence of Reedy cofibrant objects in ${\cat M}^{\cat D}$, then $\colim^{\cat D} f$ is a weak equivalence in ${\cat M}$.
\end{enumerate}
\item
Let ${\cat M}$ be a fibration category and ${\cat D}$ be a small inverse category. Then:
\begin{enumerate}
\item
If $X$ is Reedy fibrant in ${\cat M}^{\cat D}$, then $\lim^{\cat D} X$ exists and is fibrant in ${\cat M}$
\item
If $f \colon X \ra Y$ is a Reedy fibration in ${\cat M}^{\cat D}$, then $\lim^{\cat D} f$ is a fibration in ${\cat M}$
\item
If $f \colon X \ra Y$ is a pointwise weak equivalence of Reedy fibrant objects in ${\cat M}^{\cat D}$, then $\lim^{\cat D} f$ is a weak equivalence in ${\cat M}$.
\end{enumerate}
\end{enumerate}
\end{thm}

\begin{proof}
We only prove (1). Denote ${\cat D}_n = {\cat D}^{\le n}$, and ${\cat D}_{-1} = \emptyset$.

For (1) (a), an inductive argument using \lemmaref{lem:colimpushoutset} shows that each $\colim^{{\cat D}_n} X$ exists and that all maps  $\colim^{{\cat D}_{n-1}} X \ra \colim^{{\cat D}_n} X$ are cofibrations. Using CF6 (1), we see that $\colim^{{\cat D}} X$ exists and is cofibrant.

We now prove (1) (b). In view of \lemmaref{lem::compospushout2} (a), it suffices to show for any Reedy cofibrantion $f \colon X \ra Y$ in ${\cat M}^{{\cat D}}$ that
\begin{equation}
\label{eqn:colimreedycof3}
\colim\nolimits^{{\cat D}_{n}} X \Sum_{\colim^{{\cat D}_{n-1}} X} \colim\nolimits^{{\cat D}_{n-1}} Y \ra \colim\nolimits^{{\cat D}_{n}} Y
\end{equation}
is a cofibration in ${\cat M}$ for any $n$. Observe that the sum of the left term of (\ref{eqn:colimreedycof3}) exists and is cofibrant from CF3, because $\colim^{{\cat D}_{n-1}} X \ra \colim^{{\cat D}_n} X$ is a cofibration with cofibrant domain, and $\colim\nolimits^{{\cat D}_{n-1}} Y$ is cofibrant as proved in (1) (a). 

We use induction on $n$. Assume that (\ref{eqn:colimreedycof3}) is a cofibration for indices $< n$, for any choice of ${\cat D}$. From \lemmaref{lem::compospushout} (a), the map 
\begin{equation}
\label{eqn:colimreedycof5}
\colim\nolimits^{{\cat D}_{n-1}} X \ra \colim\nolimits^{{\cat D}_{n-1}} Y
\end{equation}
is a cofibration, for any choice of ${\cat D}$.

Using the notation of \lemmaref{lem:colimpushoutset}, we denote $\{d_k | k \eps K\}$ the set of objects of ${\cat D}$ of degree $n$. Consider the diagram below
\begin{equation}
\label{eqn:colimreedycof6}
\xymatrix {
        \;\;\;\; \Sum LX_{d_k} \;\;\;\; \ar[rr] \ar[dd]_{u_1} \ar@{>->}[rd]_{\Sum i_k(X)} &  & \colim\nolimits^{{\cat D}_{n-1}} X \ar'[d][dd]^{u_3} \ar@{>->}[rd] & \\
        & \;\;\;\; \Sum X_{d_k}  \;\;\;\; \ar[rr] \ar[dd]^(.3){u_2} & & \colim\nolimits^{{\cat D}_{n}} X \ar[dd]^{u_4} \\
        \;\;\;\; \Sum LY_{d_k} \;\;\;\; \ar'[r][rr] \ar@{>->}[rd]_{\Sum i_k(Y)} & & \colim\nolimits^{{\cat D}_{n-1}} Y \ar@{>->}[rd] & \\
        & \;\;\;\; \Sum Y_{d_k} \;\;\;\; \ar[rr] & & \colim\nolimits^{{\cat D}_{n}} Y
	}
\end{equation}
 
The maps $\Sum i_k(X), \Sum i_k(Y)$ are cofibrations because $X, Y$ are Reedy cofibrant. The top and bottom faces are pushouts by \lemmaref{lem:colimpushoutset}. 

The vertical map $u_3$ is a cofibration from statement (\ref{eqn:colimreedycof5}). 

The vertical map $u_1$ is a cofibration also from statement (\ref{eqn:colimreedycof5}). This is because the latching space at $d_k$ is a colimit over the direct category $\partial ({\cat D} \downarrow d_k)$. All objects of $\partial ({\cat D} \downarrow d_k)$ have degree $< n$, and the restriction of $f \colon X \ra Y$ to $\partial ({\cat D} \downarrow d_k)$ is Reedy cofibrant, so from (\ref{eqn:colimreedycof5}) we get that each $\colim^{\partial ({\cat D} \downarrow d_k)} f$ is a cofibration in ${\cat M}$.

At last, the map $(\Sum X_{d_k}) \Sum_{\Sum LX_{d_k}} (\Sum LY_{d_k}) \ra \Sum Y_{d_k}$ is a cofibration because $f$ is a Reedy cofibration.

The hypothesis of the Gluing Lemma \ref{lem:gluinglemma1} (1) (a) then applies, and we can conclude that the map \eqref{eqn:colimreedycof3} is a cofibration. The induction step is completed, and with it (1) (b) is proved.

The proof of (1) (c) follows the same exact steps as the proof of (1) (b) - only we use 
\begin{enumerate}
\item[-]
The Gluing Lemma \ref{lem:gluinglemma1} (b) instead of the Gluing Lemma \ref{lem:gluinglemma1} (a)
\item[-]
\lemmaref{lem::compospushout} (b) instead of (a)
\item[-]
\lemmaref{lem::compospushout2} (b) instead of (a)
\end{enumerate}
The reader is invited to verify the remaining details.
\end{proof}

In order to go back and complete the missing proof of the previous section, we need to state the following corollary - whose proof is implicit in the proof of \theoremref{thm:colimdirectemb}.
\begin{cor}
\label{cor:colimdirectcorollary}
\mbox{}
\begin{enumerate}
\item
Let ${\cat M}$ be a cofibration category and ${\cat D}$ be a small direct category. If $X \ra Y$ is a Reedy cofibrantion in ${\cat M}^{\cat D}$, then both $LX_d \ra LY_d$ and $X_d \ra Y_d$ are cofibrations in ${\cat M}$ for any object $d$ of ${\cat D}$.
\item
Let ${\cat M}$ be a fibration category and ${\cat D}$ be a small inverse category. If $X \ra Y$ is a Reedy fibrantion in ${\cat M}^{\cat D}$, then both $MX_d \ra MY_d$ and $X_d \ra Y_d$ are fibrations in ${\cat M}$ for any object $d$ of ${\cat D}$.
\end{enumerate}
\end{cor}
\begin{proof}
We only prove (1). The latching category $\partial({\cat D} \downarrow d)$ is a direct category, and the restriction of $X \ra Y$ to $\partial({\cat D} \downarrow d)$ is Reedy cofibrant. It follows from \theoremref{thm:colimdirect} (1) (b) that $LX_d \ra LY_d$ is a cofibrantion. The map $X_d \ra Y_d$ factors as $X_d \ra X_d \Sum_{LX_d} LY_d \ra Y_d$; the second factor is a cofibration since $X \ra Y$ is a Reedy cofibration, and the first factor is a pushout of $LX_d \ra LY_d$, therefore a cofibration. It follows that $X_d \ra Y_d$ is a cofibration.
\end{proof}

In particular, we have proved that Reedy cofibrations are pointwise cofibrations. Dually, Reedy fibrations are pointwise fibrations. We can finally go back and provide a

\begin{proof}[Proof of \theoremref{thm:reedycofibstruct}]
Given a cofibration category $({\cat M}$, ${\cat W}$, ${\cat Cof})$ and ${\cat D}$ a small direct category, we want to show that $({\cat M}^{\cat D}, {\cat W}^{\cat D}, {\cat Cof}^{\cat D}_{Reedy})$ forms a cofibration category. Recall that any Reedy cofibration $X \ra Y$ has $X$ {\it Reedy cofibrant} by definition.

(i) {\it Axiom CF1} is verified in view of \lemmaref{prop::reedycofformcat}.

(ii) {\it Axiom CF2} is easily verified.

(iii) {\it Pushout axiom CF3 (1)}. Let $f \colon X \ra Y$ be a Reedy cofibration, and $g \colon X \ra Z$ be a map with $Z$ Reedy cofibrant. The pushout $Y \Sum_X Z$ exists since each $X_d \ra Y_d$ is a cofibration by \corollaryref{cor:colimdirectcorollary}, and each $Z_d$ is cofibrant. Denote $T = Y \Sum_X Z$, and $f^{'}_d$ the cofibration $LX_d \ra LY_d$

Consider the diagram below:
\begin{center}
$\xymatrix {
        LX_d \ar[rr] \ar@{>->}[dd]_{u_1} \ar@{>->}[rd]_{f^{'}_d} &  & LZ_d \ar@{>->}'[d][dd]^{u_3} \ar@{>->}[rd] & \\
        & LY_d \ar[rr] \ar@{>->}[dd]^(.3){u_2} & & LT_d \ar[dd]^{u_4} \\
        X_d \ar'[r][rr] \ar@{>->}[rd]_{f_d} & & Z_d \ar@{>->}[rd] & \\
        & Y_d \ar[rr] & & T_d
	}$
\end{center}

From the universal property of $LT_d$, it follows that $LT_d$ exists and is isomorphic to $LY_d \Sum_{LX_d} LZ_d$, so the top face is a pushout. The bottom face is a pushout as well, and the hypothesis of the Gluing Lemma \ref{lem:gluinglemma1} (1) (a) applies: namely, $u_1, u_3$ are cofibrations since $X, Z$ are Reedy cofibrant, $f_d, f^{'}_d$ are cofibrations by \corollaryref{cor:colimdirectcorollary}, and $X_d \Sum_{LX_d} LY_d \ra Y_d$ is a cofibration since $f$ is a Reedy cofibration .

It follows from the Gluing Lemma that $u_4$ is a cofibration and $LT_d \Sum_{LZ_d} Z_d \ra T_d$ is a cofibration, which proves that $T$ is Reedy cofibrant and that $Z \ra T$ is a Reedy cofibration. This completes the proof that Reedy cofibrations are stable under pushout. 

(iv) {\it Axiom CF3 (2)} follows from a pointwise application of CF3 (2) in ${\cat M}$, since Reedy cofibrations are pointwise cofibrations by \corollaryref{cor:colimdirectcorollary}.

(v) {\it Factorization axiom CF4}. Let $f \colon X \ra Y$ be a map with $X$ Reedy cofibrant. We construct a factorization $f=rf^{'}$ with $r$ a weak equivalence and $f^{'}$ a Reedy cofibration.
\begin{center}
$f \colon \xymatrix {X \ar@{>->}[r]^{f^{'}} & Y^{'} \ar[r]^r_{\sim} & Y}$
\end{center}

We employ induction on the degree $n$. Assume $Y^{'}, f^{'}, r$ constructed in degree $<n$. Let $d$ be an object of ${\cat D}$ of degree $n$, and let us define $Y^{'}_d, f^{'}_d, r_d$. Define $LY^{'}_d$ as $\colim^{\partial({\cat D} \downarrow d)} Y^{'}$, which exists and is cofibrant by \theoremref{thm:colimdirect} applied to $\partial({\cat D} \downarrow d)$. In the diagram below
\begin{center}
$\xymatrix {
    LX_d \ar@{>->}[r] \ar@{>->}[d] & LY^{'}_d \ar@{>->}[d] & & \\
    X_d \ar@{>->}[r]^-{g_d} & X_d \Sum_{LX_d} LY^{'}_d \ar@{>->}[r]^-{h_d} & Y^{'}_d \ar[r]^{r_d}_{\sim} & Y_d
  }$
\end{center}
the pushout exists by axiom CF3, and $Y^{'}_d, h_d, r_d$ is defined as the CF4 factorization of $X_d \Sum_{LX_d} LY^{'}_d \ra Y_d$. Define $f^{'}_d = h_dg_d$, and the inductive step is complete.

(vi) {\it The axiom CF5.} 
Suppose that $f_i \colon X_i \ra Y_i$, $i \eps I$ is a set of Reedy cofibrations. The objects $X_i$ are Reedy cofibrant, therefore pointwise cofibrant by \corollaryref{cor:colimdirectcorollary}. The maps $f_i$ are in particular pointwise cofibrations by \corollaryref{cor:colimdirectcorollary}, and a pointwise application of CF5 shows that $\Sum X_i$, $\Sum Y_i$ exist and $\Sum (f_i)$ is a pointwise cofibration. Furthermore, latching spaces commute with direct sums, from which one easily sees that $\Sum f_i$ is actually a Reedy cofibration. If each $f_i$ is a trivial Reedy cofibration, a pointwise application of CF5 yields that $\Sum f_i$ is a weak equivalence.

(vii) {\it The axiom CF6.} 
Consider a countable direct sequence of Reedy cofibrations
\begin{center}
$ \xymatrix {
        X_0 \ar@{>->}[r]^{a_{0}} &
        X_1 \ar@{>->}[r]^{a_{1}} &
        X_2 \ar@{>->}[r]^{a_{2}} &
	... \\
	}$
\end{center}
The object $X_0$ is Reedy cofibrant, therefore pointwise cofibrant by \corollaryref{cor:colimdirectcorollary}. The maps $a_n$ are in particular pointwise cofibrations by \corollaryref{cor:colimdirectcorollary}, and a pointwise application of CF6 shows that $\colim X_n$ exists and $X_0 \ra \colim X_n$ is a pointwise cofibration. We'd like to show that $X_0 \ra \colim X_n$ is a Reedy cofibration.

For any object $d$ of ${\cat D}$, in the diagram below
\begin{center}
$ \xymatrix {
        L(X_0)_d \ar@{>->}[r] \ar@{>->}[d] &
        L(X_1)_d \ar@{>->}[r] \ar@{>->}[d] &
        L(X_2)_d \ar@{>->}[r] \ar@{>->}[d] &
	... \\
        (X_0)_d \ar@{>->}[r] &
        (X_1)_d \ar@{>->}[r] &
        (X_2)_d \ar@{>->}[r] &
	... \\
	}$
\end{center}
the vertical maps are cofibrations in ${\cat M}$ because the diagrams $X_n$ are Reedy cofibrant. By \corollaryref{cor:colimdirectcorollary} the top and bottom horizontal maps are cofibrations. A pointwise application of CF6 shows $\colim^n L(X_n)_d$ exists, and one easily sees that it satisfies the universal property that defines $L(\colim^n (X_n))_d$.

Furthermore, each map $(X_{n-1})_d \Sum_{L(X_{n-1})_d} {L(X_{n})_d} \ra (X_{n})_d$ is a cofibration, since $a_{n-1}$ is a Reedy cofibration. From \lemmaref{lem::compospushout2} we see that each map $(X_{0})_d \Sum_{L(X_{0})_d} {L(\colim X_{n})_d} \ra (\colim X_{n})_d$ is a cofibration, which implies that $\colim X_n$ is Reedy cofibrant and that $X_0 \ra \colim X_n$ is a Reedy cofibration. This proves CF6 (1).

If additionally all $a_n$ are trivial Reedy cofibrations, a pointwise application of CF6 (2) shows that $X_0 \ra \colim X_n$ is a pointwise weak equivalence, therefore a trivial Reedy cofibration. 

We have completed the proof that $({\cat M}^{\cat D}, {\cat W}^{\cat D}, {\cat Cof}^{\cat D}_{Reedy})$ is a cofibration category. Let's prove now that $({\cat M}^{\cat D}, {\cat W}^{\cat D}, {\cat Cof}^{\cat D})$ is also a cofibration category.

(i) {\it Axioms CF1-CF3, CF5-CF6} are trivially verified.

(ii) {\it The factorization axiom CF4.} Let $f \colon X \ra Y$ be a map of ${\cat D}$-diagrams with $X$ pointwise cofibrant. Consider the commutative diagram
\begin{center}
$\xymatrix{ 
    & & Y \\
    X \ar[r]_{f^{'}} \ar[urr]^f & Y^{'} \ar[ur]^(.3)r & \\
    X_1 \ar[u]_\sim^a \ar[r]_{f_1} & Y_1^{'} \ar[u]^b_\sim \ar[uur]_{r_1}^\sim &
  }$
\end{center}
where $X_1$ is a Reedy cofibrant replacement of $X$, $r_1f_1$ is a factorization of $fa$ as a Reedy cofibration $f_1$ followed by a weak equivalence $r_1$, and $Y^{'} = X \Sum_{X_1} Y_1^{'}$. The map $f_1$ is in particular a pointwise cofibration. Its pushout $f_1$ is therefore a pointwise cofibration, and by excision the map $b$ is a weak equivalence, so $r$ is a weak equivalence by the 2 out of 3 axiom. We have thus constructed a factorization $f = rf^{'}$ as a pointwise cofibration $f^{'}$ followed by a pointwise weak equivalence $r$.

The proof of \ref{thm:reedycofibstruct} part (1) is now complete, and part (2) is proved by duality.
\end{proof}

As a corollary of \theoremref{thm:reedycofibstruct}, we can construct the Reedy and the pointwise cofibration category structures on {\it restricted} small direct diagrams in a cofibration category. 

\index{restricted diagrams}
\begin{defn}
\label{defn::restricteddiagrams}
If $({\cat M}, {\cat W})$ is a category with weak equivalences and $({\cat D_1}, {\cat D_2})$ is a category pair, a ${\cat D_1}$ diagram $X$ is called {\it restricted} with respect to ${\cat D_2}$ if for any map $d \ra d^{'}$ of ${\cat D}_2$ the map $X_{d} \ra X_{d^{'}}$ is a weak equivalence. 

We will denote ${\cat M}^{({\cat D_1}, {\cat D_2})}$\index{${\cat M}^{({\cat D_1}, {\cat D_2})}$} the full subcategory of ${\cat D}_2$-restricted diagrams in ${\cat M}^{{\cat D_1}}$.
\end{defn}

With this definition we have

\begin{thm}
\label{thm:reedycofibstruct2}
\mbox{}
\begin{enumerate}
\item
If $({\cat M}$, ${\cat W}$, ${\cat Cof})$ is a cofibration category and $({\cat D}_1, {\cat D}_2)$ is a small direct category pair, then
\begin{enumerate}
\item
$({\cat M}^{({\cat D}_1, {\cat D}_2)} , {\cat W}^{{\cat D}_1}, {\cat Cof}^{{\cat D}_1}_{Reedy} \cap {\cat M}^{({\cat D}_1, {\cat D}_2)})$ is a cofibration category - called the ${\cat D}_2$-restricted Reedy cofibration structure on ${\cat M}^{({\cat D}_1, {\cat D}_2)}$.
\item
$({\cat M}^{({\cat D}_1, {\cat D}_2)} , {\cat W}^{{\cat D}_1}, {\cat Cof}^{({\cat D}_1, {\cat D}_2)})$ is a cofibration category - called the ${\cat D}_2$-restricted pointwise cofibration structure on ${\cat M}^{({\cat D}_1, {\cat D}_2)}$.
\end{enumerate}
\item
If $({\cat M}$, ${\cat W}$, ${\cat Fib})$ is a fibration category and $({\cat D}_1, {\cat D}_2)$ is a small inverse category pair, then
\begin{enumerate}
\item
$({\cat M}^{({\cat D}_1, {\cat D}_2)} , {\cat W}^{{\cat D}_1}, {\cat Fib}^{{\cat D}_1}_{Reedy} \cap {\cat M}^{({\cat D}_1, {\cat D}_2)})$ is a fibration category - called the ${\cat D}_2$-restricted Reedy fibration structure.on ${\cat M}^{({\cat D}_1, {\cat D}_2)}$
\item
$({\cat M}^{({\cat D}_1, {\cat D}_2)} , {\cat W}^{{\cat D}_1}, {\cat Fib}^{({\cat D}_1, {\cat D}_2)})$ is a fibration category - called the ${\cat D}_2$- restricted pointwise fibration structure on ${\cat M}^{({\cat D}_1, {\cat D}_2)}$.
\end{enumerate}
\end{enumerate}
\end{thm}

\begin{proof}
We only prove part (1) - part (2) is dual. 

(i) {\it Axioms CF1 and CF2} are clearly verified for both the pointwise and the Reedy restricted cofibration structures. 

(ii) {\it The pushout axiom CF3 (1)}. Given a pointwise cofibration $i$ and a map $f$ with $X, Y, Z$ pointwise cofibrant in ${\cat M}^{({\cat D}_1, {\cat D}_2)}$
\begin{center}
$ \xymatrix {
	X \ar[r]^f \ar@{>->}[d]_i &
	Z \ar@{>-->}[d]^j \\
	Y \ar@{-->}[r]_g &
	T
	}$
\end{center}
then the pushout $j$ of $i$ exists in ${\cat M}^{{\cat D_1}}$, and $j$ is a pointwise cofibration. For any map $d \ra d^{'}$ of ${\cat D}_2$ using the Gluing Lemma \ref{lem:gluinglemma1} applied to the diagram
\begin{center}
$ \xymatrix {
        X_d \ar[rr]^{f_{d}} \ar[dd]_{\sim} \ar@{>->}[rd]_{i_{d}} &  & Z_d \ar'[d][dd]^{\sim} \ar@{>->}[rd] & \\
        & Y_d \ar[rr] \ar[dd]^(.3){\sim} & & T_d \ar[dd] \\
        X_{d^{'}} \ar'[r][rr]_{f_{d^{'}}} \ar@{>->}[rd]_{i_{d^{'}}} & & Z_{d^{'}} \ar@{>->}[rd] & \\
        & Y_{d^{'}} \ar[rr] & & T_{d^{'}}
	}$
\end{center}
it follows that $T_{d} \ra T_{d^{'}}$ is an equivalence, therefore $T$ is a ${\cat D}_2$- restricted diagram.

Furthermore, $j$ is a Reedy cofibration if $i$ is one and $X$, $Z$ are Reedy cofibrant, by \theoremref{thm:reedycofibstruct}. This shows that the pushout axiom is satisfied for both the pointwise and the Reedy restricted cofibration structures.

(iii) {\it The axiom CF3 (2)} is clearly verified for both the pointwise and the Reedy restricted cofibration structures. 

(iv) {\it The factorization axiom CF4}. Let $f \colon X \ra Y$ be a map in ${\cat M}^{({\cat D}_1, {\cat D}_2)}$. If $X$ is a pointwise (resp. Reedy) cofibrant diagram, by \theoremref{thm:reedycofibstruct} $f$ factors as $f = rf^{'}$ with $f^{'} \colon X \ra Y^{'}$ a pointwise (resp. Reedy) cofibration and $r \colon Y^{'} \ra Y$ a pointwise weak equivalence. In both cases $Y^{'}$ is restricted, and CF4 is satisfied for both the pointwise and the Reedy restricted cofibration structures.

(v) {\it Axiom CF5} for both the restricted pointwise and restricted Reedy cofibration structures follows from the fact that if $X_i$, $i \eps I$ is a set of restricted pointwise (resp. Reedy) cofibrant diagrams, then $\Sum X_i$ is a restricted pointwise (resp. Reedy) cofibrant diagram by \lemmaref{lem:smallcofsums}.

(vi) {\it Axiom CF6}. Given a countable direct sequence of pointwise (resp. Reedy) cofibrations with $X_0$ pointwise (resp. Reedy) cofibrant
\begin{center}
$ \xymatrix {
        X_0 \ar@{>->}[r]^{a_{01}} &
        X_1 \ar@{>->}[r]^{a_{12}} &
        X_2 \ar@{>->}[r]^{a_{23}} &
	... \\
	}$
\end{center}
the colimit $\colim X_n$ exists and is pointwise (resp. Reedy) cofibrant. If all $X_n$ are restricted, from \lemmaref{lem:transcompequiv} $\colim X_n$ is restricted. The axiom CF6 now follows for both the restricted pointwise and restricted Reedy cofibration structures.
\end{proof}

\section{Colimits in direct categories (the relative case)}
\label{sec:colimdirectrelative}

The contents of this section is not used elsewhere in this text, and may be skipped at a first reading. 

Recall that a functor $u \colon {\cat D} \ra {\cat D}^{'}$ is an {\it open embedding} (or a {\it crible}) \index{open embedding, crible} if $u$ is a full embedding with the property that any map $f \colon d^{'} \ra ud$ with $d \eps {\cat D}, d^{'} \eps {\cat D}^{'}$ has $d^{'}$ (and $f$) in the image of $u$. Dually, $u$ is a {\it closed embedding} (or a {\it cocrible}) \index{closed embedding, cocrible} if $u^{op}$ is an open embedding, meaning that $u$ is a full embedding and any map $f \colon ud \ra d^{'}$ with $d \eps {\cat D}, d^{'} \eps {\cat D}^{'}$ has $d^{'}$ (and $f$) in the image of $u$.

Open and closed embedding functors are stable under compositions. 

If ${\cat D}$ is a direct category, note that ${\cat D}^{\le n} \ra {\cat D}^{\le n+1}$ and ${\cat D}^{\le n} \ra {\cat D}$ are open embeddings. If ${\cat D}$ is an inverse category, then ${\cat D}^{\le n} \ra {\cat D}^{\le n+1}$ and ${\cat D}^{\le n} \ra {\cat D}$ are closed embeddings.

Recall that a functor $u \colon {\cat D} \ra {\cat D}^{'}$ is an {\it open embedding} (or a {\it crible}) \index{open embedding, crible} if $u$ is a full embedding with the property that any map $f \colon d^{'} \ra ud$ with $d \eps {\cat D}, d^{'} \eps {\cat D}^{'}$ has $d^{'}$ (and $f$) in the image of $u$. Dually, $u$ is a {\it closed embedding} (or a {\it cocrible}) \index{closed embedding, cocrible} if $u^{op}$ is an open embedding, meaning that $u$ is a full embedding and any map $f \colon ud \ra d^{'}$ with $d \eps {\cat D}, d^{'} \eps {\cat D}^{'}$ has $d^{'}$ (and $f$) in the image of $u$.

Open and closed embedding functors are stable under compositions. 

If ${\cat D}$ is a direct category, note that ${\cat D}^{\le n} \ra {\cat D}^{\le n+1}$ and ${\cat D}^{\le n} \ra {\cat D}$ are open embeddings. If ${\cat D}$ is an inverse category, then ${\cat D}^{\le n} \ra {\cat D}^{\le n+1}$ and ${\cat D}^{\le n} \ra {\cat D}$ are closed embeddings.

Here is a statement that is slightly more general than \theoremref{thm:colimdirect}.

\begin{thm}
\label{thm:colimdirectemb}
\mbox{}
\begin{enumerate}
\item
Let ${\cat M}$ be a cofibration category and ${\cat D}^{'} \ra {\cat D}$ be an open embedding of small direct categories. Then:
\begin{enumerate}
\item
If $X$ is Reedy cofibrant in ${\cat M}^{\cat D}$, then $\colim^{{\cat D}^{'}} X, \colim^{\cat D} X$ exist and $\colim^{{\cat D}^{'}} X \ra \colim^{\cat D} X$ is a cofibration in ${\cat M}$
\item
If $f \colon X \ra Y$ is a Reedy cofibration in ${\cat M}^{\cat D}$, then 
\begin{center}
$\colim^{\cat D} X \Sum_{\colim^{{\cat D}^{'}} X} \colim^{{\cat D}^{'}} Y \ra \colim^{\cat D} Y$ 
\end{center}
is a cofibration in ${\cat M}$
\item
If $f \colon X \ra Y$ is a pointwise weak equivalence between Reedy cofibrant diagrams in ${\cat M}^{\cat D}$, then 
\begin{center}
$\colim^{\cat D} X \Sum_{\colim^{{\cat D}^{'}} X} \colim^{{\cat D}^{'}} Y \ra \colim^{\cat D} Y$ 
\end{center}
is a weak equivalence in ${\cat M}$
\end{enumerate}
\item
Let ${\cat M}$ be a fibration category and  ${\cat D}^{'} \ra {\cat D}$ be a closed embedding of small inverse categories. Then:
\begin{enumerate}
\item
If $X$ is Reedy fibrant in ${\cat M}^{\cat D}$, then $\lim^{{\cat D}^{'}} X, \lim^{\cat D} X$ exist and $\lim^{{\cat D}} X \ra \lim^{{\cat D}^{'}} X$ is a fibration in ${\cat M}$
\item
If $f \colon X \ra Y$ is a Reedy fibration in ${\cat M}^{\cat D}$, then 
\begin{center}
$\lim^{{\cat D}} X \ra \lim^{{\cat D}^{'}} X \times_{\lim^{{\cat D}^{'}} Y} \lim^{{\cat D}} Y $ 
\end{center}
is a fibration in ${\cat M}$
\item
If $f \colon X \ra Y$ is a pointwise weak equivalence between Reedy fibrant diagrams in ${\cat M}^{\cat D}$, then 
\begin{center}
$\lim^{{\cat D}} X \ra \lim^{{\cat D}^{'}} X \times_{\lim^{{\cat D}^{'}} Y} \lim^{{\cat D}} Y $ 
\end{center}
is a weak equivalence in ${\cat M}$
\end{enumerate}
\end{enumerate}
\end{thm}

If we take ${\cat D}^{'}$ to be empty in \theoremref{thm:colimdirectemb}, then we get the statement of \theoremref{thm:colimdirect}.

\begin{proof}[Proof of \theoremref{thm:colimdirectemb}]
Part (2) is dual to part (1), so we only need to prove part (1). The proof, as expected, retracks that of \theoremref{thm:colimdirect}.

We will use on ${\cat D}^{'}$ the degree induced from ${\cat D}$. Denote ${\cat D}_{-1} = {\cat D}^{'}$, and let ${\cat D}_n$ be the full subcategory of ${\cat D}$ having as objects all the objects of ${\cat D}^{'}$ and all the objects of degree $\le n$ of ${\cat D}$. Notice that all inclusions ${\cat D}_{n-1} \ra {\cat D}_n$ are open embeddings. Also, all inclusions ${\cat D}^{' < n} \ra {\cat D}^{' \le n}$ are open embeddings.

Let us prove (1) (a). From \theoremref{thm:colimdirect} we know that $\colim^{{\cat D}^{'}} X, \colim^{\cat D} X$ exist and are cofibrant. An inductive argument using \lemmaref{lem:colimpushoutset} shows that each $\colim^{{\cat D}_n} X$ exists and that all maps $\colim^{{\cat D}_{n-1}} X \ra  \colim^{{\cat D}_n} X$ are cofibrations. Using CF6 (1), we see that $\colim^{{\cat D}^{'}}X \ra \colim^{{\cat D}} X$ is a cofibration.

We now prove (1) (b). In view of \lemmaref{lem::compospushout2} (a), it suffices to show that
\begin{equation}
\label{eqn:colimreedycof3bis}
\colim\nolimits^{{\cat D}_{n}} X \Sum_{\colim^{{\cat D}_{n-1}} X} \colim\nolimits^{{\cat D}_{n-1}} Y \ra \colim\nolimits^{{\cat D}_{n}} Y
\end{equation}
is a cofibration in ${\cat M}$ for any $n$. Observe that the sum of the left term of (\ref{eqn:colimreedycof3bis}) exists and is cofibrant because of (1) (a) and CF3. 

As in \lemmaref{lem:colimpushoutset}, we denoted $\{d_k | k \eps K\}$ the set of objects of ${\cat D} \backslash {\cat D}^{'}$ of degree $n$. We now look back at the cubic diagram \eqref{eqn:colimreedycof6}:

The maps $\Sum i_k(X), \Sum i_k(Y)$ are cofibrations because $X, Y$ are Reedy cofibrant. The top and bottom faces are pushouts by \lemmaref{lem:colimpushoutset}. 

The vertical maps $u_1, u_3$ are cofibrations by \theoremref{thm:colimdirect}.

The map $(\Sum X_{d_k}) \Sum_{\Sum LX_{d_k}} (\Sum LY_{d_k}) \ra \Sum Y_{d_k}$ is a cofibration because $f$ is a Reedy cofibration.

From the Gluing Lemma \ref{lem:gluinglemma1} (1) (a), we conclude that the map \eqref{eqn:colimreedycof3bis} is a cofibration.

The proof of (1) (c) follows the same exact steps as the proof of (1) (b) - only we use as for \theoremref{thm:colimdirect} (1) (c):
\begin{enumerate}
\item[-]
The Gluing Lemma \ref{lem:gluinglemma1} (b) instead of the Gluing Lemma \ref{lem:gluinglemma1} (a)
\item[-]
\lemmaref{lem::compospushout} (b) instead of (a)
\item[-]
\lemmaref{lem::compospushout2} (b) instead of (a)
\end{enumerate}
\end{proof}

We have shown that given a cofibration category ${\cat M}$ and a small direct category ${\cat D}$, then $\colim^{\cat D}$ carries Reedy cofibrations (resp. weak equivalences between Reedy cofibrant diagrams) in ${\cat M}^{\cat D}$ to cofibrations (resp. weak equivalences between cofibrant objects) in ${\cat M}$. This result is extended below to the case of relative colimits from a small direct category to an arbitrary small category (\theoremref{thm::colimfromdirectcat}) and between two small direct categories (\theoremref{thm::colimpairdirectcat}). 

\begin{thm}
\label{thm::colimfromdirectcat}
\mbox{}
\begin{enumerate}
\item
Let ${\cat M}$ be a cofibration category and $u \colon {\cat D}_1 \ra {\cat D}_2$ be a functor of small categories with ${\cat D}_1$ direct.
\begin{enumerate}
\item
If $X$ is Reedy cofibrant in ${\cat M}^{{\cat D}_1}$, then $\colim^u X$ exists and is pointwise cofibrant in ${\cat M}^{{\cat D}_2}$
\item
If $f \colon X \ra Y$ is a Reedy cofibration in ${\cat M}^{{\cat D}_1}$, then $\colim^u f$ is a pointwise cofibration in ${\cat M}^{{\cat D}_2}$
\item
If $f \colon X \ra Y$ is a pointwise weak equivalence of Reedy cofibrant objects in ${\cat M}^{{\cat D}_1}$, then $\colim^u f$ is a pointwise weak equivalence in ${\cat M}^{{\cat D}_2}$.
\end{enumerate}
\item
Let ${\cat M}$ be a fibration category and $u \colon {\cat D}_1 \ra {\cat D}_2$ be a functor of small categories  with ${\cat D}_1$ inverse.
\begin{enumerate}
\item
If $X$ is Reedy fibrant in ${\cat M}^{{\cat D}_1}$, then $\lim^u X$ exists and is pointwise fibrant in ${\cat M}^{{\cat D}_2}$
\item
If $f \colon X \ra Y$ is a Reedy fibration in ${\cat M}^{{\cat D}_1}$, then $\colim^u f$ is a pointwise fibration in ${\cat M}^{{\cat D}_2}$
\item
If $f \colon X \ra Y$ is a pointwise weak equivalence of Reedy fibrant objects in ${\cat M}^{{\cat D}_1}$, then $\lim^u f$ is a pointwise weak equivalence in ${\cat M}^{{\cat D}_2}$.
\end{enumerate}
\end{enumerate}
\end{thm}

\begin{proof}
We only prove statement (1) - statement (2) follows from duality.

Let us prove (a). To prove that $\colim^u X$ exists, cf. \lemmaref{lem::computerelcolim} it suffices to show for any object $d_2 \eps {\cat D}_2$ that $\colim^{(u \downarrow d_2)} X$ exists. But the over category $(u \downarrow d_2)$ is direct and the restriction of $X$ to $(u \downarrow d_2)$ is Reedy cofibrant, therefore by \theoremref{thm:colimdirect} $\colim^{(u \downarrow d_2)} X$ exists and is cofibrant in ${\cat M}$. It follows that $\colim^u X$ exists, and since $\colim^{(u \downarrow d_2)} X \cong (\colim^u X)_{d_2}$ we have that $\colim^u X$ is pointwise cofibrant in ${\cat M}^{{\cat D}_2}$.

We now prove (b). If $f \colon X \ra Y$ is a Reedy cofibration in ${\cat M}^{{\cat D}_1}$, then the restriction of $f$ to $(u \downarrow d_2)$ is Reedy cofibrant for any object $d_2 \eps {\cat D}_2$, therefore by \theoremref{thm:colimdirect} $\colim^{(u \downarrow d_2)} f$ is a cofibration in ${\cat M}$. Since $(\colim^u f)_{d_2} \cong \colim^{(u \downarrow d_2)} f$ by the naturality of the isomorphism in \theoremref{thm:colimdirect}, it follows that $\colim^u f$ is pointwise cofibrant in ${\cat M}^{{\cat D}_2}$.

To prove (c), assume that $f \colon X \ra Y$ is a pointwise weak equivalence between Reedy cofibrant diagrams in ${\cat M}^{{\cat D}_1}$. The restrictions of $X$ and $Y$ to $(u \downarrow d_2)$ are Reedy cofibrant for any object $d_2 \eps {\cat D}_2$, and by \theoremref{thm:colimdirect} (1) (c) the map $\colim^{(u \downarrow d_2)} X \ra \colim^{(u \downarrow d_2)} Y$ is a weak equivalence. In conclusion, the map $\colim^u X \ra \colim^u Y$ is a pointwise weak equivalence in ${\cat M}^{{\cat D}_2}$.
\end{proof}

\begin{thm}
\label{thm::colimpairdirectcat}
\mbox{}
\begin{enumerate}
\item
Let ${\cat M}$ be a cofibration category and $u \colon {\cat D}_1 \ra {\cat D}_2$ be a functor of small direct categories.
\begin{enumerate}
\item
If $X$ is Reedy cofibrant in ${\cat M}^{{\cat D}_1}$, then $\colim^u X$ exists and is Reedy cofibrant in ${\cat M}^{{\cat D}_2}$
\item
If $f \colon X \ra Y$ is a Reedy cofibration in ${\cat M}^{{\cat D}_1}$, then $\colim^u f$ is a Reedy cofibration in ${\cat M}^{{\cat D}_2}$
\item
If $f \colon X \ra Y$ is a pointwise weak equivalence of Reedy cofibrant objects in ${\cat M}^{{\cat D}_1}$, then $\colim^u f$ is a pointwise weak equivalence in ${\cat M}^{{\cat D}_2}$.
\end{enumerate}
\item
Let ${\cat M}$ be a fibration category and $u \colon {\cat D}_1 \ra {\cat D}_2$ be a functor of small inverse categories.
\begin{enumerate}
\item
If $X$ is Reedy fibrant in ${\cat M}^{{\cat D}_1}$, then $\lim^u X$ exists and is Reedy fibrant in ${\cat M}^{{\cat D}_2}$
\item
If $f \colon X \ra Y$ is a Reedy fibration in ${\cat M}^{{\cat D}_1}$, then $\lim^u f$ is a Reedy fibration in ${\cat M}^{{\cat D}_2}$
\item
If $f \colon X \ra Y$ is a pointwise weak equivalence of Reedy fibrant objects in ${\cat M}^{{\cat D}_1}$, then $\lim^u f$ is a pointwise weak equivalence in ${\cat M}^{{\cat D}_2}$.
\end{enumerate}
\end{enumerate}
\end{thm}

\begin{proof}
We only prove statement (1) - statement (2) follows from duality.

Let us prove (a). We know from \theoremref{thm::colimfromdirectcat} (1) (a) that $\colim^u X$ exists and is pointwise cofibrant in ${\cat M}^{{\cat D}_2}$, and we would like to show that $\colim^u X$ is Reedy cofibrant. For that, fix an object $d_2$ of ${\cat D}_2$, and let us try to identify the latching object $\colim^{\partial({\cat D}_2 \downarrow d_2)} \colim^u X$.

The functor $H \colon {\cat D}_2 \ra {\cat Cat}$, $Hd_2^{'} = (u \downarrow d_2^{'})$ restricts to a functor $\partial({\cat D}_2 \downarrow d_2) \ra {\cat Cat}$. The Grothendieck construction $\int_{\partial({\cat D}_2 \downarrow d_2)} {H}$ has as objects 4-tuples $(d_1^{'}, d_2^{'}, d_2, ud_1^{'} \ra d_2^{'} \ra d_2)$ of objects $d_1^{'} \eps {\cat D}_1, d_2^{'}, d_2 \eps {\cat D}_2$ and maps $ud_1^{'} \ra d_2^{'} \ra d_2$ such that $d_2^{'} \ra d_2$ is a non-identity map. 

Denote $\partial(u \downarrow d_2)$ the full subcategory of the over category $(u \downarrow d_2)$ consisting of triples $(d_1^{'}, d_2, ud_1^{'} \ra d_2)$ with a non-identity map $ud_1^{'} \ra d_2$. We have an adjoint pair of functors $F : \partial(u \downarrow d_2) \rightleftarrows \int_{\partial({\cat D}_2 \downarrow d_2)} {H} : G$, defined as follows:
\begin{center}
$F(d_1^{'}, d_2, ud_1^{'} \ra d_2) = (d_1^{'}, ud_1^{'}, d_2, ud_1^{'} \ra ud_1^{'} \ra d_2)$
\end{center}
\begin{center}
$G(d_1^{'}, d_2^{'}, d_2, ud_1^{'} \ra d_2^{'} \ra d_2) = (d_1^{'}, d_2, ud_1^{'} \ra d_2)$
\end{center}
\begin{center}
$id \colon 1_{\partial(u \downarrow d_2)} \Ra GF$
\end{center}
\begin{center}
$FG(d_1^{'}, d_2^{'}, d_2, ud_1^{'} \ra d_2^{'} \ra d_2) = (d_1^{'}, ud_1^{'}, d_2, ud_1^{'} \ra ud_1^{'} \ra d_2) \Ra (d_1^{'}, d_2^{'}, d_2, ud_1^{'} \ra d_2^{'} \ra d_2)$ given on $2^{nd}$ component by $ud_1^{'} \ra d_2^{'}$.
\end{center}
The category $\partial(u \downarrow d_2)$ is direct, and the restriction of $X$ to $\partial(u \downarrow d_2)$ is Reedy cofibrant, therefore $\colim^{\partial(u \downarrow d_2)} X$ exists and is cofibrant. 

$G$ is a right adjoint functor, therefore right cofinal, and by \propositionref{prop:cofinallimit} we have that $\colim^{\int_{\partial({\cat D}_2 \downarrow d_2)} {H}} X$ exists and is $\cong \colim^{\partial(u \downarrow d_2)} X$. From \propositionref{prop:grothendieckconstrcolims}, $\colim^{(d_2^{'} \ra d_2) \eps \partial({\cat D}_2 \downarrow d_2)}$ $\colim^{(u \downarrow d^{'}_2)} X$ exists and is $\cong \colim^{\int_{\partial({\cat D}_2 \downarrow d_2)} {H}} X \cong \colim^{\partial(u \downarrow d_2)} X$. In conclusion, the latching object $\colim^{\partial({\cat D}_2 \downarrow d_2)} \colim^u X$ exists and is $\cong \colim^{\partial(u \downarrow d_2)} X$. In particular the latching object is cofibrant.

The inclusion functor $\partial(u \downarrow d_2) \ra (u \downarrow d_2)$ is an open embedding, and from \theoremref{thm:colimdirectemb} (1) (a) the latching map $\colim^{\partial(u \downarrow d_2)} X \ra \colim^{(u \downarrow d_2)} X$ is a cofibration, therefore $\colim^u X$ is Reedy cofibrant in ${\cat M}^{{\cat D}_2}$. The proof of statement (a) of our theorem is complete.

Let us prove (b). If $f \colon X \ra Y$ is a Reedy cofibration, by \theoremref{thm:colimdirectemb} (1) (b) the map
\begin{center}
$\colim^{(u \downarrow d_2)} X \Sum_{\colim^{\partial(u \downarrow d_2)} X} \colim^{\partial(u \downarrow d_2)} Y \ra \colim^{(u \downarrow d_2)} Y$
\end{center}
is a cofibration, therefore $\colim^u X \ra \colim^u Y$ is a Reedy cofibration in ${\cat M}^{{\cat D}_2}$.

Part (c) has been already proved as \theoremref{thm::colimfromdirectcat} (1) (c).
\end{proof}

\section{Colimits in arbitrary categories}
\label{sec:colimarbitrarycat}

Denote $\Delta^{'}$\index{$\Delta^{'}$, $\Delta^{'}{\cat D}$, $\Delta^{'op}{\cat D}$} the subcategory of the cosimplicial indexing category $\Delta$, with same objects as $\Delta$ and with maps the order-preserving {\it injective} maps ${\bf n}_1 \ra {\bf n}_2$. 

If ${\cat D}$ is a category, we define $\Delta^{'}{\cat D}$ to be the category with objects the functors ${\bf n} \ra {\cat D}$, and with maps $({\bf n_1} \ra {\cat D}) \lra ({\bf n_2} \ra {\cat D})$ the commutative diagrams
\begin{center}
$ \xymatrix {
        {\bf n_1} \ar[rr]^f \ar[dr] & & {\bf n_2} \ar[dl] \\
	& {\cat D} &
	}$
\end{center}
where $f$ is injective and order-preserving.

The category $\Delta^{'}{\cat D}$ is direct, and comes equipped with a {\it terminal projection} functor $p_t \colon \Delta^{'}{\cat D} \ra {\cat D}$\index{$p_i, p_t$} that sends ${\bf n} \ra {\cat D}$ to the image of the terminal object $n$ of the poset ${\bf n}$.

The opposite category of $\Delta^{'}{\cat D}$ is denoted $\Delta^{'op}{\cat D}$. It is an inverse category, and comes equipped with an {\it initial projection} functor $p_i \colon \Delta^{'op}{\cat D} \ra {\cat D}$ that sends ${\bf n} \ra {\cat D}$ to the image of the initial object $0$ of the poset ${\bf n}$.

\begin{lem}
\label{lem:deltaprimeovercat}
If $u \colon {\cat D}_1 \ra {\cat D}_2$ is a functor and $d_2 \eps {\cat D}_2$ is an object
\begin{enumerate}
\item
There is a natural isomorphism $(up_t \downarrow d_2) \cong \Delta^{'}(u \downarrow d_2)$
\item
There is a natural isomorphism $(d_2 \downarrow up_i) \cong \Delta^{'op}(d_2 \downarrow u)$
\end{enumerate}
\end{lem}

\begin{proof}
Left to the reader.
\end{proof}

\begin{lem}
\label{lem:deltaprimelcofinal}
If ${\cat D}$ is a category and $d \eps {\cat D}$ is an object
\begin{enumerate}
\item
\begin{enumerate}
\item
The category $p_t^{-1}d$ has an initial object and hence has a contractible nerve.
\item
The category $(p_t \downarrow d)$ has a contractible nerve.
\item
The inclusion $p_t^{-1}d \ra (p_t \downarrow d)$ is homotopy right cofinal.
\end{enumerate}
\item
\begin{enumerate}
\item
The category $p_i^{-1}d$ has a terminal object and hence has a contractible nerve.
\item
The category $(d \downarrow p_i)$ has a contractible nerve.
\item
The inclusion $p_i^{-1}d \ra (d \downarrow p_i)$ is homotopy left cofinal.
\end{enumerate}
\end{enumerate}
\end{lem}

\begin{proof}
We denote an object ${\bf n} \ra {\cat D}$ of $\Delta^{'}{\cat D}$ as $(d_0 \ra ... \ra d_n)$, where $d_i$ is the image of $i \eps {\bf n}$ and $d_i \ra d_{i+1}$ is the image of $i \ra i+1$. Each map $i \colon {\bf k} \ra {\bf n}$ determines a map $(d_{i0} \ra ... \ra d_{ik}) \ra (d_0 \ra ... \ra d_n)$ in the category $\Delta^{'}{\cat D}$.

We denote an object of $(p_t \downarrow d)$ as $(d_0 \ra ... \ra d_n) \ra d$, where $d_n \ra d$ is the map $p_t (d_0 \ra ... \ra d_n) \ra d$.

The category $p_t^{-1}d$ has the initial object $(d) \eps \Delta^{'}{\cat D}$, which proves part (1) (a). For (1) (b), a contraction of the nerve of $(p_t \downarrow d)$ is defined by
\begin{center}
$\xymatrix{
    (p_t \downarrow d) \rruppertwocell^{1_d}{\alpha} \rrlowertwocell_{c_d}{^\beta} \ar[rr]^(.3){a_d} && (p_t \downarrow d)
  }$
\end{center}

In this diagram, $c_d$ is the constant functor that takes as value the object $(d) \overset{1_d}{\ra} d$. The functor $a_d$ sends an object $(d_0 \ra ... \ra d_n) \overset{f}{\ra} d$ to $(d_0 \ra ... \ra d_n \overset{f}{\ra} d) \overset{1_d}{\ra} d$, and a map defined by $i \colon {\bf k} \ra {\bf n}$ to its extension $\overline{i} \colon {\bf k+1} \ra {\bf n+1}$ with $\overline{i}(k+1) = n+1$.

On an object $(d_0 \ra ... \ra d_n) \ra d$, the natural map $\alpha$ is given by the map $i \colon {\bf n} \ra {\bf n+1}$, $ik = k$ and the natural map $\beta$ is given by the map $i \colon {\bf 0} \ra {\bf n+1}$, $i0 = n+1$.

Let us prove (1) (c). Denote $i_d \colon p_t^{-1}d \ra (p_t \downarrow d)$ the full inclusion functor. Note that the images of $a_d$ and $c_d$ are inside $i_d(p_t^{-1}d)$. Let $x$ be an object of $(p_t \downarrow d)$ of the form $(d_0 \ra ... \ra d_n) \overset{f}{\ra} d$. 

If $d_n \overset{f}{\ra} d$ is the identity map, then $x$ is in the image of $i_d$ therefore $(x, x \overset{1_x}{\ra} x)$ is an initial object of $(x \downarrow i_d)$. If $d_n \overset{f}{\ra} d$ is not the identity map, then $(a_d x, x \ra a_d x)$ defined by the map $i \colon {\bf n} \ra {\bf n+1}$, $ik = k$ is an initial object in $(x \downarrow i_d)$. In both cases, $(x \downarrow i_d)$ is contractible therefore $i_d$ is homotopy right cofinal.

The statements of part (2) follow from duality.
\end{proof}

The category $p_t^{-1}d$ is direct for any object $d \eps {\cat D}$, since it is a subcategory of the direct category $\Delta^{'}{\cat D}$. Dually, the category $p_i^{-1}d$ is inverse for any object $d \eps {\cat D}$.

\begin{lem}
\label{lem:restrictedreedycofibrant}
\mbox{}
\begin{enumerate}
\item
Let ${\cat M}$ be a cofibration category and ${\cat D}$ be a small category. Then for any Reedy cofibrant diagram $X \eps {\cat M}^{\Delta^{'}{\cat D}}$ and any object $d \eps {\cat D}$, the restriction $X|_{p_t^{-1}d}$ is Reedy cofibrant in ${\cat M}^{p_t^{-1}d}$.
\item
Let ${\cat M}$ be a fibration category and ${\cat D}$ be a small category. Then for any Reedy fibrant diagram $X \eps {\cat M}^{\Delta^{'op}{\cat D}}$ and any object $d \eps {\cat D}$, the restriction $X|_{p_i^{-1}d}$ is Reedy fibrant in ${\cat M}^{p_i^{-1}d}$.
\end{enumerate}
\end{lem}

\begin{proof}
We only prove (1). For $d \eps {\cat D}$, fix an object $\underline{d} = (d_0 \ra ... \ra d_n)$ $\eps \Delta^{'}{\cat D}$ with $d_n = d$. We need some notations. Assume that 
\[
\underline{i} = \{i_1, ..., i_u \}, \,\,\,\,\,\, \underline{j} = \{j_1, ..., j_v \}, \,\,\,\,\,\, \underline{k} = \{k_1, ..., k_w \} 
\]
is a partition of $\{0, ..., n \}$ into three (possibly empty) subsets with $u + v + w = n+1$. Denote $\underline{n}_{\underline{i}, \underline{\widehat{j}}}$ the full subcategory of $\Delta^{'}{\cat D}$. with objects $d_{l_0} \ra ... \ra d_{l_x}$ with $\underline{i} \subset \underline{l}$ and $\underline{j} \cap \underline{l} = \emptyset$, where $\underline{l} = \{l_0, .., l_x \}$. Although not apparent from the notation, the category $\underline{n}_{\underline{i}, \underline{\widehat{j}}}$ depends on the choice of $\underline{d} \eps \Delta^{'}{\cat D}$.

The category $\underline{n}_{\underline{i}, \underline{\widehat{j}}}$ is direct, and has a terminal object denoted $\underline{d}_{\underline{\widehat{j}}}$. Denote $\partial \underline{n}_{\underline{i}, \underline{\widehat{j}}}$ the maximal full subcategory of $\underline{n}_{\underline{i}, \underline{\widehat{j}}}$ without its terminal element.

$Claim.$ The restriction of $X$ to $\underline{n}_{\underline{i}, \underline{\widehat{j}}}$ is Reedy cofibrant.

Taking in particular $\underline{i} = \{ n \}$ and $\underline{j} = \emptyset$, the Claim implies that the Reedy condition is satisfied for $X|_{p_t^{-1}d}$ at $\underline{d}$. It remains to prove the Claim, and we will proceed by induction on $n$.

\begin{enumerate}
\item[-]
The Claim can be directly verified for $n = 0$. Assume that the Claim was proved for $n^{'} < n$, and let's prove it for $n \ge 1$.
\item[-]
If $\underline{j} \neq \emptyset$, the claim follows from the inductive hypothesis for smaller $n$. It remains to prove the Claim for $\underline{n}_{\underline{i}, \underline{\widehat{\emptyset}}}$.
\item[-]
We only need to prove the Reedy condition for $X$ at the terminal object $\underline{d}_{\underline{\widehat{\emptyset}}} = \underline{d}$ of $\underline{n}_{\underline{i}, \underline{\widehat{\emptyset}}}$. The Reedy condition at any other object of $\underline{n}_{\underline{i}, \underline{\widehat{\emptyset}}}$ follows from the inductive hypothesis for smaller $n$.
\item[-]
If $\underline{i} = \{i_1, ..., i_u \}$ is nonempty, in the diagram
\begin{equation}
\label{eqn:claimindstep}
\xymatrix{
    \colim^{\partial \underline{n}_{\{i_1, ..., i_u\} \backslash \{ i_s \}, \widehat{\{ i_s \}}}} X \ar[r] \ar@{>->}[d] & 
    \colim^{\partial \underline{n}_{\{i_1, ..., i_u\}, \widehat{\emptyset}}} X \ar@{>->}[d] \\
    X_{\underline{d}_{\widehat{\{ i_s \}}}} \ar[r] & \colim^{\partial \underline{n}_{\{i_1, ..., i_u\} \backslash \{ i_s \}, \widehat{\emptyset}}} X
    }
\end{equation}
the left vertical map is a cofibration with cofibrant domain, from the inductive hypothesis for smaller $n$. If the top right colimit in the diagram exists and is cofibrant, then the bottom right colimit exists, the diagram is a pushout and therefore the right vertical map is a cofibration.
\item[-]
The colimit $\colim^{\partial \underline{n}_{\{0, ..., n\}, \widehat{\emptyset}}} X$ exists and is the initial object of ${\cat M}$
\item[-]
The map $\colim^{\partial \underline{n}_{\emptyset, \widehat{\emptyset}}} X$ $\ra$ $X_{\underline{d}}$ is a cofibration with cofibrant domain since $X$ is Reedy cofibrant in $\Delta^{'} {\cat D}$. 
\item[-]
An iterated use of \equationref{eqn:claimindstep} shows that the colimit below exists and 
\[
\colim\nolimits^{\partial \underline{n}_{\{i_1, ..., i_u\} , \widehat{\emptyset}}} X \ra X_{\underline{d}}
\]
is a cofibration with cofibrant domain. This completes the proof of the Claim.
\end{enumerate}

\comment{
For any object $\underline{d} = (d_0 \ra ... \ra d_n = d)$ of $p_t^{-1}d$ denote $\underline{d^{'}}$ the object $(d_0 \ra ... \ra d_{n-1})$ of $\Delta^{'}{\cat D}$. Denote $LX_{\underline{d}}, LX_{\underline{d^{'}}}$ the latching objects of $X$ at $\underline{d}$ resp. $\underline{d^{'}}$ in ${\cat M}^{\Delta^{'}{\cat D}}$, and denote $L(X|_{p_t^{-1}d})_{\underline{d}}$ the latching space of the restriction $X|_{p_t^{-1}d}$ at $\underline{d}$ in ${\cat M}^{p_t^{-1}d}$.

The square in the diagram below is a pushout
\begin{center}
$\xymatrix{
    LX_{\underline{d^{'}}} \ar[r] \ar@{>->}[d] & L(X|_{p_t^{-1}d})_{\underline{d}} \ar@{>->}[d]_g & \\
    X_{\underline{d^{'}}} \ar[r] & LX_{\underline{d}} \ar@{>->}[r]^h & X_{\underline{d}}
    }$
\end{center}
The maps $f$ and $h$ are cofibrations since $X$ is Reedy cofibrant in ${\cat M}^{\Delta^{'}{\cat D}}$. The map $g$ is a cofibration as the pushout of $f$, therefore $hg$ is a cofibration. This shows that $X|_{p_t^{-1}d}$ is Reedy cofibrant.

Part (2) follows from duality.
}
\end{proof}

For a category ${\cat D}$, the subcategories $p_t^{-1}d \subset \Delta^{'}{\cat D}$ are disjoint for $d \eps {\cat D}$, and their union $\cup_{d \eps {\cat D}} \; p_t^{-1}d$ forms a category that we will denote $\Delta^{'}_{res}{\cat D}$\index{$\Delta^{'}_{res}{\cat D}$, $\Delta^{'op}_{res}{\cat D}$}. We will also denote $\Delta^{'op}_{res}{\cat D}$ the opposite of $\Delta^{'}_{res}{\cat D}$.

Given a cofibration category ${\cat M}$ we use the shorthand notation ${\cat M}^{\Delta^{'}{\cat D}}_{res}$\index{${\cat M}^{\Delta^{'}{\cat D}}_{res}$, ${\cat M}^{\Delta^{'op}{\cat D}}_{res}$} for the category  ${\cat M}^{({\Delta^{'}{\cat D}}, {\Delta^{'}_{res}{\cat D}})}$ of $\Delta^{'}_{res}{\cat D}$ restricted $\Delta^{'}{\cat D}$ diagrams in ${\cat M}$. ${\cat M}^{\Delta^{'}{\cat D}}_{res}$ is a full subcategory of ${\cat M}^{\Delta^{'}{\cat D}}$, and (\theoremref{thm:reedycofibstruct2}) it carries a restricted Reedy cofibration structure as well as a restricted pointwise cofibration structure. Furthermore, the functor $p_t^* \colon {\cat M}^{\cat D} \ra {\cat M}^{\Delta^{'}{\cat D}}$ has its image inside ${\cat M}^{\Delta^{'}{\cat D}}_{res}$.

Dually, for a fibration category ${\cat M}$ we denote ${\cat M}^{\Delta^{'op}{\cat D}}_{res}$ the category ${\cat M}^{({\Delta^{'op}{\cat D}}, {\Delta^{'op}_{res}{\cat D}})}$ of restricted diagrams. ${\cat M}^{\Delta^{'op}{\cat D}}_{res}$ carries a restricted Reedy fibration structure as well as a restricted pointwise fibration structure. The functor $p_i^* \colon {\cat M}^{\cat D} \ra {\cat M}^{\Delta^{'op}{\cat D}}$ has its image inside ${\cat M}^{\Delta^{'op}{\cat D}}_{res}$.

\begin{prop}
\label{prop:restrictedproperty}
\mbox{}
\begin{enumerate}
\item
Let ${\cat M}$ be a cofibration category and ${\cat D}$ be a small category. Then for every restricted Reedy cofibrant diagrams $X, X^{'} \eps {\cat M}^{\Delta^{'}{\cat D}}_{res}$ and every diagram $Y \eps {\cat M}^{{\cat D}}$
\begin{enumerate}
\item
A map $X \ra p_t^* Y$ is a pointwise weak equivalence iff \\ its adjoint $\colim^{p_t} X$ $\ra$ $Y$ is a pointwise weak equivalence.
\item
A map $X \ra X^{'}$ is a pointwise weak equivalence iff \\ the map $\colim^{p_t} X$ $\ra$ $\colim^{p_t}X^{'}$ is a pointwise weak equivalence.
\end{enumerate}
\item
Let ${\cat M}$ be a fibration category and ${\cat D}$ be a small category. Then for every restricted Reedy fibrant diagrams $X, X^{'} \eps {\cat M}^{\Delta^{'op}{\cat D}}_{res}$ and every diagram $Y \eps {\cat M}^{{\cat D}}$
\begin{enumerate}
\item
A map $p_i^* Y \ra X$ is a pointwise weak equivalence iff \\ its adjoint $Y$ $\ra$ $\lim^{p_i} X$ is a pointwise weak equivalence.
\item
A map $X \ra X^{'}$ is a pointwise weak equivalence iff \\ the map $\lim^{p_t} X$ $\ra$ $\lim^{p_t}X^{'}$ is a pointwise weak equivalence.
\end{enumerate}
\end{enumerate}
\end{prop}

\begin{proof}
We will prove (1), and (2) will follow from duality. Let $d$ be an object of ${\cat D}$. The categories $p_t^{-1}d$ and $(p_t \downarrow d)$ are direct, and since $X$ is Reedy cofibrant so are its restrictions to ${p_t^{-1}d}$ (by \lemmaref{lem:restrictedreedycofibrant}) and to $(p_t \downarrow d)$.

The colimits $\colim^{p_t^{-1}d}X$, $\colim^{(p_t \downarrow d)} X$ therefore exist. Since $i_d \colon p_t^{-1}d \ra (p_t \downarrow d)$ is right cofinal, we have $\colim^{p_t^{-1}d}X \cong \colim^{(p_t \downarrow d)} X$. We also conclude that $\colim^{p_t} X$ exists and $\colim^{(p_t \downarrow d)} X \cong (\colim^{p_t} X)_d$.

The diagram $X$ is restricted and $p_t^{-1}d$ has an initial object that we will denote $e(d)$. We get a pointwise weak equivalence $cX_{e(d)} \ra X|_{p_t^{-1}d}$ in ${\cat M}^{p_t^{-1}d}$ from the constant diagram to the restriction of $X$. But the diagram $cX_{e(d)}$ is Reedy cofibrant since $e(d)$ is initial in $p_t^{-1}d$. The map $cX_{e(d)} \ra X|_{p_t^{-1}d}$ is a pointwise weak equivalence between Reedy cofibrant diagrams in ${\cat M}^{p_t^{-1}d}$, therefore $X_{e(d)} \cong \colim^{p_t^{-1}d}cX_{e(d)} \ra \colim^{p_t^{-1}d}X$ is a weak equivalence.

In summary, we have showed that the composition 
\begin{center}
$X_{e(d)} \ra \colim^{p_t^{-1}d}X \cong \colim^{(p_t \downarrow d)} X \cong (\colim^{p_t} X)_d$
\end{center}
is a weak equivalence. We can now complete the proof of \propositionref{prop:restrictedproperty}.

To prove (a), the map $\colim^{p_t}X \ra Y$ is a pointwise weak equivalence iff the map $X_{e(d)} \ra Y_d$ is a weak equivalence for all objects $d$ of ${\cat D}$. Since $X$ is restricted, this last statement is true iff the map $X \ra p_t^*Y$ is a pointwise weak equivalence.

To prove (b), the map $\colim^{p_t}X \ra \colim^{p_t}X^{'}$ is a pointwise weak equivalence iff the map $X_{e(d)} \ra X^{'}_{e(d)}$ is a weak equivalence for all objects $d$ of ${\cat D}$. Since $X, X^{'}$ are restricted, this last statement is true iff the map $X \ra X^{'}$ is a pointwise weak equivalence.
\end{proof}

As an application, we can now prove that the category of diagrams in a cofibration (resp. fibration) category admits a pointwise cofibration (resp. fibration) structure. This result is due to Cisinski. We should note that the statement below is not true if we replace ``cofibration category'' with ``Quillen model category''.

\begin{thm}[Pointwise (co)fibration structure]
\label{thm:generalpointwisecofstructure}
\mbox{}
\begin{enumerate}
\item
If $({\cat M}$, ${\cat W}$, ${\cat Cof})$ is a cofibration category and ${\cat D}$ is a small category then $({\cat M}^{\cat D}$, ${\cat W}^{\cat D}$, ${\cat Cof}^{\cat D})$ is a cofibration category.
\item
If $({\cat M}$, ${\cat W}$, ${\cat Fib})$ is a fibration category and ${\cat D}$ is a small category then $({\cat M}^{\cat D}$, ${\cat W}^{\cat D}$, ${\cat Fib}^{\cat D})$ is a fibration category.
\end{enumerate}
\end{thm}

\begin{proof}
To prove (1), axioms CF1-CF3 and CF5-CF6 are easily verified. To prove axiom CF4, we replay an argument found in the proof of \theoremref{thm:reedycofibstruct} (1) (b). Let $f \colon X \ra Y$ be a map of ${\cat D}$-diagrams with $X$ pointwise cofibrant. Let $a \colon X_1 \ra p^*_t X$ be a Reedy cofibrant replacement in ${\cat M}^{\Delta^{'}{\cat D}}_{res}$. We factor $X_1 \ra p^*_t Y$ as a Reedy cofibration $f_1$ followed by a pointwise weak equivalence $r_1$
\begin{center}
$\xymatrix {
    X_1 \ar@{>->}[r]^{f_1} & Y_1 \ar[r]^{r_1}_\sim & p^*_tY
  } $
\end{center}

We then construct a commutative diagram
\begin{center}
$\xymatrix{ 
    &&& Y \\
    X \ar@{>->}[rr]_{f^{'}} \ar[urrr]^f && Y^{'} \ar[ur]^(.3)r & \\
    \colim^{p_t} X_1 \ar[u]_\sim^{a^{'}} \ar@{>->}[rr]_{\colim^{p_t} f_1} && \colim^{p_t} Y_1 \ar[u]^{b^{'}}_\sim \ar[uur]_{r_1^{'}}^\sim &
  }$
\end{center}
In this diagram $a^{'}$ resp. $r^{'}_1$ are the adjoints of $a$ resp. $r_1$, therefore by \propositionref{prop:restrictedproperty} (1) (a) $a^{'}$ and $r^{'}_1$ are weak equivalences. Since $f_1$ is a Reedy cofibration, $\colim^{p_t} f_1$ is a pointwise cofibration, and we construct $f^{'}$ as the pushout of $\colim^{p_t} f_1$. It follows that $f^{'}$ is a pointwise cofibration. By pointwise excision, $b^{'}$ and therefore r are pointwise weak equivalences. The factorization $f = rf^{'}$ is the desired decomposition of $f$ as a pointwise cofibration followed by a weak equivalence, and CF4 is proved.

The proof of (2) is dual.
\end{proof}

As an immediate corollary, we can show that for a small category pair $({\cat D}_1, {\cat D}_2)$, the category of reduced diagrams ${\cat M}^{({\cat D}_1, {\cat D}_2)}$ carries a pointwise cofibration structure if ${\cat M}$ is a cofibration category.

\begin{thm}[Reduced pointwise (co)fibration structure]
\label{thm:generalpointwisecofibstruct2}
\mbox{}
\begin{enumerate}
\item
If $({\cat M}$, ${\cat W}$, ${\cat Cof})$ is a cofibration category and $({\cat D}_1, {\cat D}_2)$ is a small category pair, then $({\cat M}^{({\cat D}_1, {\cat D}_2)} , {\cat W}^{{\cat D}_1}, {\cat Cof}^{({\cat D}_1, {\cat D}_2)})$ is a cofibration category - called the ${\cat D}_2$-restricted pointwise cofibration structure on ${\cat M}^{({\cat D}_1, {\cat D}_2)}$.
\item
If $({\cat M}$, ${\cat W}$, ${\cat Fib})$ is a fibration category and $({\cat D}_1, {\cat D}_2)$ is a small category pair, then $({\cat M}^{({\cat D}_1, {\cat D}_2)} , {\cat W}^{{\cat D}_1}, {\cat Fib}^{({\cat D}_1, {\cat D}_2)})$ is a fibration category - called the ${\cat D}_2$- restricted pointwise fibration structure on ${\cat M}^{({\cat D}_1, {\cat D}_2)}$.
\end{enumerate}
\end{thm}

\begin{proof}
Entirely similar to that of \theoremref{thm:reedycofibstruct2}.
\end{proof}

Denote ${\cat M}^{\Delta^{'}{\cat D}}_{res, rcof}$\index{${\cat M}^{\Delta^{'}{\cat D}}_{res, rcof}$, ${\cat M}^{\Delta^{'op}{\cat D}}_{res, rfib}$} the full subcategory of ${\cat M}^{\Delta^{'}{\cat D}}_{res}$ of restricted Reedy cofibrant diagrams for a cofibration category ${\cat M}$. The next proposition states that restricted Reedy cofibrant diagrams in ${\cat M}^{\Delta^{'}{\cat D}}$ form a cofibrant approximation of ${\cat M}^{\cat D}$. 

Dually for a fibration category ${\cat M}$ denote ${\cat M}^{\Delta^{'op}{\cat D}}_{res, rfib}$ the full subcategory of ${\cat M}^{\Delta^{'op}{\cat D}}_{res}$ of restricted Reedy fibrant diagrams. We show that restricted Reedy fibrant diagrams in ${\cat M}^{\Delta^{'op}{\cat D}}$ form a fibrant approximation of ${\cat M}^{\cat D}$.

\begin{prop}
\label{prop:generalpointwisecofstructureapprox}
Let ${\cat D}$ be a small category.
\begin{enumerate}
\item
If ${\cat M}$ is a cofibration category then $\colim^{p_t} \colon {\cat M}^{\Delta^{'}{\cat D}}_{res, rcof} \ra {\cat M}^{\cat D}$ is a cofibrant approximation for the pointwise cofibration structure on ${\cat M}^{\cat D}$.
\item
If ${\cat M}$ is a fibration category then $\lim^{p_i} \colon {\cat M}^{\Delta^{'op}{\cat D}}_{res, rfib} \ra {\cat M}^{\cat D}$ is a fibrant approximation for the pointwise fibration structure on ${\cat M}^{\cat D}$.
\end{enumerate}
\end{prop}

\begin{proof}
We only prove (1), since the proof of (2) is dual.

The functor $\colim^{p_t}$ sends Reedy cofibrations to pointwise cofibrations by \theoremref{thm:colimdirect} and preserves the initial object, which proves CFA1. CFA2 is a consequence of \propositionref{prop:restrictedproperty} (1) (b). CFA3 is a consequence of \remarkref{rem:grothendieckconstrcolims}.

For CFA4, let $f \colon \colim^{p_t}X \ra Y$ be a map in ${\cat M}^{\cat D}$ with $X \eps {\cat M}^{\Delta^{'}{\cat D}}_{res, rcof}$ restricted and Reedy cofibrant. Factor its adjoint $f^{'} \colon X \ra p_t^*Y$ as a Reedy cofibration $f_1$ followed by a pointwise weak equivalence $r_1$
\begin{center}
$\xymatrix {
    X \ar@{>->}[r]^{f_1} & Y^{'} \ar[r]^{r_1}_\sim & p^*_tY
  } $
\end{center}
We get the following CFA4 factorization of $f$
\begin{center}
$\xymatrix {
    \colim^{p_t}X \ar@{>->}[rr]^{\colim^{p_t} f_1} && \colim^{p_t} Y^{'} \ar[r]^-{r_1^{'}}_-\sim & Y
  } $
\end{center}
where $r_1^{'}$ is the adjoint of $r_1$, therefore a weak equivalence.
\end{proof}

Here is another application of \propositionref{prop:restrictedproperty}:

\begin{thm}
\label{thm:pointwisecofstructureequivalence}
Let ${\cat D}$ be a small category.
\begin{enumerate}
\item
If ${\cat M}$ is a cofibration category then 
\begin{center}
${\bf ho}p_t^{*} \colon {\bf ho}{\cat M}^{\cat D} \ra {\bf ho}{\cat M}^{\Delta^{'}{\cat D}}_{res}$ 
\end{center}
is an equivalence of categories.
\item
If ${\cat M}$ is a fibration category then 
\begin{center}
${\bf ho}p_i^{*} \colon {\bf ho}{\cat M}^{\cat D} \ra {\bf ho}{\cat M}^{\Delta^{'op}{\cat D}}_{res}$ 
\end{center}
is an equivalence of categories.
\end{enumerate}
\end{thm}

\begin{proof}
Let us prove (1). We will apply the Abstract Partial Quillen Adjunction \theoremref{thm:generalexistencetotalderivedadjoint2} to
\begin{center}
$\xymatrix{
    {\cat M}^{\Delta^{'}{\cat D}}_{res, rcof} \ar[rr]^{v_1 = \colim^{p_t}} \ar@{_{(}->}[d]_{t_1} && {\cat M}^{\cat D} \\
    {\cat M}^{\Delta^{'}{\cat D}}_{res} && {\cat M}^{\cat D} \ar[ll]_{v_2 = p_t^*} \ar[u]_{t_2 = 1_{{\cat M}^{\cat D}}}
    }$
\end{center}
The functor $t_1$ is the full inclusion of ${\cat M}^{\Delta^{'}{\cat D}}_{res, rcof}$ in ${\cat M}^{\Delta^{'}{\cat D}}_{res}$. We have that $v_1, v_2$ is an abstract Quillen partially equivalent pair with respect to $t_1, t_2$. Indeed:
\begin{enumerate}
\item
The functor pair $v_1, v_2$ is partially adjoint with respect to $t_1, t_2$.
\item
$t_1$ is a cofibrant approximation of the cofibration category ${\cat M}^{\Delta^{'}{\cat D}}_{res}$ with the Reedy reduced structure. In particular $t_1$ is a left approximation. The functor $t_2 = 1_{{\cat M}^{\cat D}}$ is a right approximation. 
\item
$v_1$ preserves weak equivalences from \theoremref{thm::colimfromdirectcat}, and so does $v_2 = p_t^*$.
\item
\propositionref{prop:restrictedproperty} (1) (a) states that for any objects $X \eps {\cat M}^{\Delta^{'}{\cat D}}_{res, rcof}$, $Y \eps {\cat M}^{\cat D}$, a map $v_1X \ra t_2Y$ is a weak equivalence iff its partial adjoint $t_1X \ra v_2Y$ is a weak equivalence
\end{enumerate}

In conclusion we have a pair of equivalences of categories
\begin{center}
$\xymatrix {{\bf ho}{\cat M}^{\Delta^{'}{\cat D}}_{res} \ar@<2pt>[rrrr]^-{{\bf ho}(v_1) \, {\bf s}_1} &&&& {\bf ho}{\cat M}^{{\cat D}} \ar@<2pt>[llll]^-{{\bf ho}(v_2) \, {\bf s}_2 \cong {\bf ho}(p^*_t)}}$
\end{center}
where ${\bf s}_1$ is a quasi-inverse of ${\bf ho}t_1$ and ${\bf s}_2$ is a quasi-inverse of ${\bf ho}t_2$, and therefore ${\bf ho}(v_2) \, {\bf s}_2$ is naturally isomorphic to ${\bf ho}(p^*_t)$. This proves that ${\bf ho}(p^*_t)$ is an equivalence of categories.

The proof of part (2) is dual.
\end{proof}

\section{Homotopy colimits}
\label{sec:hocolimcolimdirect}
Suppose that $({\cat M}, {\cat W})$ is a category with weak equivalences and suppose that $u \colon {\cat D}_1 \ra {\cat D}_2$ is a functor of small categories. We denote $\gamma_{{\cat D}_i} \colon {\cat M}^{{\cat D}_i} \ra {\bf ho} {\cat M}^{{\cat D}_i}$, $i = 1, 2$ the localization functors.

We define ${\cat M}^{{\cat D}_1}_{{\colim}^u}$\index{${\cat M}^{{\cat D}_1}_{{\colim}^u}$, ${\cat M}^{{\cat D}_1}_{{\lim}^u}$} to be the full subcategory of ${\cat M}^{{\cat D}_1}$ of ${\cat D}_1$ diagrams $X$ with the property that $\colim^uX$ exists in ${\cat M}^{{\cat D}_2}$. Denote $i_{\colim^u} \colon {\cat M}^{{\cat D}_1}_{{\colim}^u} \ra {\cat M}^{{\cat D}_1}$ the inclusion. In general ${\cat M}^{{\cat D}_1}_{{\colim}^u}$ may be empty, but if ${\cat M}$ is cocomplete then ${\cat M}^{{\cat D}_1}_{{\colim}^u} = {\cat M}^{{\cat D}_1}$. Let ${\cat W}^{{\cat D}_1}_{{\colim}^u}$ be the class of pointwise weak equivalences of ${\cat M}^{{\cat D}_1}_{{\colim}^u}$. 

Dually, let ${\cat M}^{{\cat D}_1}_{{\lim}^u}$ be the full subcategory of ${\cat M}^{{\cat D}_1}$ of ${\cat D}_1$ diagrams $X$ with the property that $\lim^uX$ exists in ${\cat M}^{{\cat D}_2}$, let $i_{\lim^u} \colon {\cat M}^{{\cat D}_1}_{{\lim}^u} \ra {\cat M}^{{\cat D}_1}$ denote the inclusion, and let ${\cat W}^{{\cat D}_1}_{{\lim}^u}$ be the class of pointwise weak equivalences of ${\cat M}^{{\cat D}_1}_{{\lim}^u}$.

\index{homotopy (co)limit}
\begin{defn}
\label{defn:homotopycolim}
\mbox{}
\begin{enumerate}
\item
The {\it homotopy colimit} of $u$, if it exists, is the left Kan extension of $\gamma_{{\cat D}_2}\colim^u$ along $\gamma_{{\cat D}_1}i_{\colim^u}$
\begin{center}
$\xymatrix {
    {\cat M}^{{\cat D}_1}_{{\colim}^u} \ar[rr]^{\colim^u} \ar[d]_{\gamma_{_{{\cat D}_1}}i_{\colim^u}} && {\cat M}^{{\cat D}_2} \ar[d]^{\gamma_{{\cat D}_2}} _(.5)\,="b" \\ 
    {\bf ho} {\cat M}^{{\cat D}_1} \ar[rr]_{{\bf L} \colim^u} ^(.5)\,="a" & & {\bf ho} {\cat M}^{{\cat D}_2} \ulltwocell<\omit>{<0>\eps_u \;}
  }$
\end{center}
and is denoted for simplicity $({\bf L}\colim^u, \eps_u)$\index{${\bf L}\colim^u$, ${\bf R}\lim^u$} instead of the more complete notation $({\bf L}_{_{\gamma_{_{{\cat D}_1}}i_{\colim^u}}}\colim^u, \eps_u)$
\item
The {\it homotopy limit} of $u$, if it exists, is the right Kan extension of $\gamma_{{\cat D}_2}\lim^u$ along $\gamma_{{\cat D}_1}i_{\lim^u}$
\begin{center}
$\xymatrix {
    {\cat M}^{{\cat D}_1}_{{\lim}^u} \ar[rr]^{\lim^u} \ar[d]_{\gamma_{_{{\cat D}_1}}i_{\lim^u}} \drrtwocell<\omit>{<0>\; \nu_u} && {\cat M}^{{\cat D}_2} \ar[d]^{\gamma_{{\cat D}_2}} \\ 
    {\bf ho} {\cat M}^{{\cat D}_1} \ar[rr]_{{\bf R} \lim^u} && {\bf ho} {\cat M}^{{\cat D}_2}
  }$
\end{center}
and is denoted for simplicity $({\bf R}\lim^u, \nu_u)$
\end{enumerate}
\end{defn}

If $({\cat M}$, ${\cat W}$, ${\cat Cof})$ is a cofibration category, then we also define the class ${\cat Cof}^{{\cat D}_1}_{{\colim}^u}$ of pointwise cofibrations $f \colon X \ra Y$ in ${\cat M}^{{\cat D}_1}_{{\colim}^u}$ with the property that $\colim^u f \colon \colim^u X \ra \colim^u Y$ is well defined and pointwise cofibrant in ${\cat M}^{{\cat D}_2}$.

Dually, if $({\cat M}$, ${\cat W}$, ${\cat Fib})$ is a fibration category, then we define the class ${\cat Fib}^{{\cat D}_1}_{{\lim}^u}$ of pointwise fibration maps $f$ in ${\cat M}^{{\cat D}_1}_{{\lim}^u}$ with the property that $\lim^u f$ is well defined and pointwise fibrant in ${\cat M}^{{\cat D}_2}$.

\begin{lem}
\label{lem:homotcatcolim}
Let $u \colon {\cat D}_1 \ra {\cat D}_2$ be a functor of small categories.
\begin{enumerate}
\item
If ${\cat M}$ is a cofibration category, then 
\begin{enumerate}
\item
$({\cat M}^{{\cat D}_1}_{{\colim}^u}, {\cat W}^{{\cat D}_1}_{{\colim}^u}, {\cat Cof}^{{\cat D}_1}_{{\colim}^u})$ is a cofibration category
\item
$\colim^{p_t} \colon {\cat M}^{\Delta^{'}{\cat D}_1}_{res, rcof} \ra {\cat M}^{{\cat D}_1}_{\colim^u}$ is a cofibrant approximation
\item
${\bf ho}{\cat M}^{{\cat D}_1}_{{\colim}^u} \ra {\bf ho}{\cat M}^{{\cat D}_1}$ is an equivalence of categories.
\end{enumerate}
\item
If ${\cat M}$ is a fibration category, then 
\begin{enumerate}
\item
$({\cat M}^{{\cat D}_1}_{{\lim}^u}, {\cat W}^{{\cat D}_1}_{{\lim}^u}, {\cat Fib}^{{\cat D}_1}_{{\lim}^u})$ is a cofibration category
\item
$\lim^{p_i} \colon {\cat M}^{\Delta^{'op}{\cat D}_1}_{res, rfib} \ra {\cat M}^{{\cat D}_1}_{\lim^u}$ is a fibrant approximation
\item
${\bf ho}{\cat M}^{{\cat D}_1}_{{\lim}^u} \ra {\bf ho}{\cat M}^{{\cat D}_1}$ is an equivalence of categories.
\end{enumerate}
\end{enumerate}
\end{lem}

\begin{proof}
We only prove (1), since the proof of (2) is dual.

For an object $X \eps {\cat M}^{\Delta^{'}{\cat D}_1}_{res, rcof}$, we have that the colimits $\colim^{p_t}X$, $\colim^{up_t}X$ exist and are pointwise cofibrant by \theoremref{thm::colimfromdirectcat} since $X$ is Reedy cofibrant in ${\cat M}^{\Delta^{'}{\cat D}_1}$. Since $\colim^{up_t}X \cong \colim^{u}\colim^{p_t}X$, we conclude that the functor $\colim^{p_t} \colon {\cat M}^{\Delta^{'}{\cat D}_1}_{res, rcof} \ra {\cat M}^{{\cat D}_1}$ has its image inside ${\cat M}^{{\cat D}_1}_{\colim^u}$.

Let us prove (a). Axioms CF1-CF2 and CF5-CF6 are easily verified for ${\cat M}^{{\cat D}_1}_{{\colim}^u}$. The pushout axiom CF3 (1) follows from the fact that if 
\begin{center}
$ \xymatrix {
	X \ar[r] \ar@{>->}[d]_i &
	Z \ar@{>-->}[d]^j \\
	Y \ar@{-->}[r] &
	T
	}$
\end{center}
is a pushout in ${\cat M}^{{\cat D}_1}$ with $X, Y, Z$ cofibrant in ${\cat M}^{{\cat D}_1}_{{\colim}^u}$ and $i$ a cofibration in ${\cat M}^{{\cat D}_1}_{{\colim}^u}$, then $\colim^u X, \colim^u Y, \colim^u Z$ are pointwise cofibrant and $\colim^u i$ is a pointwise cofibration in ${\cat M}^{{\cat D}_2}$, therefore by \remarkref{rem:grothendieckconstrcolims} we have that $\colim^uT$ exists and is the pushout in ${\cat M}^{{\cat D}_2}$ of
\begin{center}
$ \xymatrix {
	\colim^u X \ar[r] \ar@{>->}[d]_{\colim^u i} &
	\colim^u Z \ar@{>-->}[d]^{\colim^u j} \\
	\colim^u Y \ar@{-->}[r] &
	\colim^u T
	}$
\end{center}
The axiom CF3 (2) follows from a pointwise application of CF3 (2) in ${\cat M}$.

Let us prove the factorization axiom CF4. We repeat the argument in the proof of \theoremref{thm:generalpointwisecofstructure}.

Let $f \colon X \ra Y$ be a map in ${\cat M}^{{\cat D}_1}_{{\colim}^u}$ with $X$ cofibrant. Let $a \colon X_1 \ra p_t^*X$ be a Reedy cofibrant replacement in ${\cat M}^{\Delta^{'}{\cat D}_1}_{red}$. We factor $X_1 \ra p^*_t Y$ as a Reedy cofibration $f_1$ followed by a pointwise weak equivalence $r_1$
\begin{center}
$\xymatrix {
    X_1 \ar@{>->}[r]^{f_1} & Y_1 \ar[r]^{r_1}_\sim & p^*_tY
  } $
\end{center}

We then construct a commutative diagram
\begin{center}
$\xymatrix{ 
    &&& Y \\
    X \ar@{>->}[rr]_{f^{'}} \ar[urrr]^f && Y^{'} \ar[ur]^(.3)r & \\
    \colim^{p_t} X_1 \ar[u]_\sim^{a^{'}} \ar@{>->}[rr]_{\colim^{p_t} f_1} && \colim^{p_t} Y_1 \ar[u]^{b^{'}}_\sim \ar[uur]_{r_1^{'}}^\sim &
  }$
\end{center}
In this diagram $a^{'}$ resp. $r^{'}_1$ are the adjoints of $a$ resp. $r_1$, therefore by \propositionref{prop:restrictedproperty} (1) (a) $a^{'}$ and $r^{'}_1$ are weak equivalences. Since $f_1$ is a Reedy cofibration, $\colim^{p_t} f_1$ is a cofibration in ${\cat M}^{{\cat D}_1}_{{\colim}^u}$, and we construct $f^{'}$ as the pushout of $\colim^{p_t} f_1$. It follows that $f^{'}$ is a cofibration of ${\cat M}^{{\cat D}_1}_{{\colim}^u}$. By pointwise excision, $b^{'}$ and therefore r are pointwise weak equivalences. The factorization $f = rf^{'}$ is the desired decomposition of $f$ as a pointwise cofibration followed by a weak equivalence in ${\cat M}^{{\cat D}_1}_{{\colim}^u}$, and CF4 is proved.

Part (b) follows from \propositionref{prop:generalpointwisecofstructureapprox} and the fact that $\colim^{p_t} \colon {\cat M}^{\Delta^{'}{\cat D}_1}_{res, rcof} \ra {\cat M}^{{\cat D}_1}$ has its image inside ${\cat M}^{{\cat D}_1}_{{\colim}^u}$. 

To prove part (c), since both functors $\colim^{p_t} \colon {\cat M}^{\Delta^{'}{\cat D}_1}_{res, rcof} \ra {\cat M}^{{\cat D}_1}$ and its corestriction $\colim^{p_t} \colon {\cat M}^{\Delta^{'}{\cat D}_1}_{res, rcof} \ra {\cat M}^{{\cat D}_1}_{{\colim}^u}$ are cofibrant approximations it follows (\theoremref{thm:existenceleftderivedequiv}) that both induced functors ${\bf ho}{\cat M}^{\Delta^{'}{\cat D}_1}_{res, rcof} \ra {\bf ho}{\cat M}^{{\cat D}_1}$ and ${\bf ho}{\cat M}^{\Delta^{'}{\cat D}_1}_{res, rcof} \ra {\bf ho}{\cat M}^{{\cat D}_1}_{{\colim}^u}$ are equivalences of categories. The functor ${\bf ho}{\cat M}^{{\cat D}_1}_{{\colim}^u} \ra {\bf ho}{\cat M}^{{\cat D}_1}$ is therefore an equivalence of categories.
\end{proof}

We now state the main result of this section. 
\begin{thm}[Existence of homotopy (co)limits]
\label{thm:genholim}
\mbox{}
\begin{enumerate}
\item
Let ${\cat M}$ be a cofibration category and $u \colon {\cat D}_1 \ra {\cat D}_2$ be a functor of small categories. Then the homotopy colimit $({\bf L}\colim^u, \eps_u)$ exists and 
\begin{center}
${\bf L}\colim^u \colon {\bf ho} {\cat M}^{{\cat D}_1} \rightleftarrows {\bf ho} {\cat M}^{{\cat D}_2} \colon {\bf ho}u^{*}$ 
\end{center}
forms a naturally adjoint pair.
\item
Let ${\cat M}$ be a fibration category and $u \colon {\cat D}_1 \ra {\cat D}_2$ be a functor of small categories. Then the homotopy limit $({\bf R}\lim^u, \nu_u)$ exists and 
\begin{center}
${\bf ho}u^{*} \colon {\bf ho} {\cat M}^{{\cat D}_2} \rightleftarrows {\bf ho} {\cat M}^{{\cat D}_1} \colon {\bf R}\lim^u$ 
\end{center}
forms a naturally adjoint pair.
\end{enumerate}
\end{thm}

\begin{proof}
Parts (1) and (2) are dual, and we will only prove part (1).

From \lemmaref{lem:homotcatcolim} and \corollaryref{cor:neededbygenhomlim}, to prove the existence of the left Kan extension $({\bf L}\colim^u, \eps_u)$ it suffices to prove the existence of the total left derived of $\colim^u \colon {\cat M}^{{\cat D}_1}_{{\colim}^u} \ra {\cat M}^{{\cat D}_2}$. But the latter is a consequence of \theoremref{thm:existenceleftderivedapprox} applied to the cofibrant approximation $\colim^{p_t} \colon {\cat M}^{\Delta^{'}{\cat D}_1}_{res, rcof} \ra {\cat M}^{{\cat D}_1}_{{\colim}^u}$.

To prove that ${\bf L}\colim^u \dashv {\bf ho}u^{*}$ forms a naturally adjoint pair, we will apply the Abstract Quillen Partial Adjunction \theoremref{thm:generalexistencetotalderivedadjoint2} to
\begin{center}
$\xymatrix{
    {\cat M}^{\Delta^{'}{\cat D}_1}_{res, rcof} \ar[rr]^{v_1 = \colim^{up_t}} \ar[d]_{t_1 = \colim^{p_t}} && {\cat M}^{{\cat D}_2} \\
    {\cat M}^{{\cat D}_1} && {\cat M}^{{\cat D}_2} \ar[ll]_{v_2 = u^*} \ar[u]_{t_2 = 1_{{\cat M}^{{\cat D}_2}}}
    }$
\end{center}

We have that $v_1, v_2$ is an abstract Quillen partially adjoint pair with respect to $t_1, t_2$:
\begin{enumerate}
\item
The functor pair $v_1, v_2$ is partially adjoint with respect to $t_1, t_2$.
\item
$t_1$ is a cofibrant approximation of the cofibration category ${\cat M}^{{\cat D}_1}$ by \propositionref{prop:generalpointwisecofstructureapprox}. In particular $t_1$ is a left approximation. The functor $t_2 = 1_{{\cat M}^{{\cat D}_2}}$ is a right approximation. 
\item
$v_1$ preserves weak equivalences by \theoremref{thm::colimfromdirectcat}, and so does $v_2 = p_t^*$.
\end{enumerate}

In conclusion we have a naturally adjoint pair
\begin{center}
$\xymatrix {{\bf ho}{\cat M}^{{\cat D}_1} \ar@<2pt>[rrrr]^-{{\bf ho}(v1) \, {\bf s}_1 \cong {\bf L}\colim^u} &&&& {\bf ho}{\cat M}^{{\cat D}_2} \ar@<2pt>[llll]^-{{\bf ho}(v_2) \, {\bf s}_2 \cong {\bf ho}(u^*)}}$
\end{center}
where ${\bf s}_1$ is a quasi-inverse of ${\bf ho}t_1$ and ${\bf s}_2$ is a quasi-inverse of ${\bf ho}t_2$.
\end{proof}

As a consequence of \theoremref{thm:genholim} we can verify that

\begin{cor}
\label{cor:composhomotcolimits}
Suppose that $u \colon {\cat D}_1 \ra {\cat D}_2$ and $v \colon {\cat D}_2 \ra {\cat D}_3$ are two functors of small categories.
\begin{enumerate}
\item
If ${\cat M}$ is a cofibration category, then ${\bf L}\colim^{vu} \cong {\bf L}\colim^{v}{\bf L}\colim^{u}$.
\item
If ${\cat M}$ is a fibration category, then ${\bf R}\lim^{vu} \cong {\bf R}\lim^{v}{\bf R}\lim^{u}$.
\end{enumerate}
\end{cor}

\begin{proof}
This is a consequence of the adjunction property of the homotopy (co)limit and the fact that ${\bf ho}(vu)^* \cong {\bf ho}u^* {\bf ho}v^* $.
\end{proof}

Consider a small diagram
\begin{equation}
\label{eqn:2catsquare4}
\xymatrix {
    {\cat D}_1 \ar[r]^u \ar[d]_f & 
    {\cat D}_2 \ar[d]^g _\,="b" \\
    {\cat D}_3 \ar[r]_v ^\,="a" & 
    {\cat D}_4 \ultwocell<\omit>{<0>\phi}
  }
\end{equation}
If ${\cat M}$ is a cofibration category, from the adjunction property of the homotopy colimit we get a natural map denoted
\begin{center}
$\phi_{{\bf L}\colim} \colon {\bf L}\colim^u {\bf ho}f^* \Ra {\bf ho}g^* {\bf L}\colim^v$
\end{center}
and dually if ${\cat M}$ is a fibration category we get a natural map denoted
\begin{center}
$\phi_{{\bf R}\colim} \colon {\bf ho}v^* {\bf R}\lim^g  \Ra {\bf R}\lim^f {\bf ho}u^* $
\end{center}

Suppose that $u \colon {\cat D}_1 \ra {\cat D}_2$ is a functor of small categories. For any object $d_2 \eps {\cat D_2}$, the {\it standard over 2-category diagram} of $u$ at $d_2$ is defined as
\begin{equation}
\label{eqn:stdover2cat}
\xymatrix {
    (u \downarrow d_2) \ar[r]^-{p_{(u \downarrow d_2)}} \ar[d]_{i_{u, d_2}} & 
    \terminalcat \ar[d]^{e_{d_2}} _\,="b" \\
    {\cat D}_1 \ar[r]_u ^\,="a" & 
    {\cat D}_2 \ultwocell<\omit>{<0>{{\phi}_{u, d_2} \;\;\;\;\;\;}}
  }
\end{equation}

In this diagram, the functor $p_{{\cat D}} \colon {\cat D} \ra \terminalcat$ denotes the terminal category projection. The functor $i_{u, d_2}$ sends an object $(d_1, f: ud_1 \ra d_2)$ to $d_1 \eps {\cat D}_1$ and a map $(d_1, f) \ra (d^{'}_1, f^{'})$ to the component map $d_1 \ra d^{'}_1$. The functor $e_{d_2}$ embeds the terminal category $\terminalcat$ as the object $d_2 \eps {\cat D}_2$, and for an object $(d_1 \eps {\cat D}_1, f \colon ud_1 \ra d_2)$ of $(u \downarrow d_2)$ the natural map $\phi_{u, d_2}(d_1, ud_1 \ra d_2)$ is $f \colon ud_1 \ra d_2$. If ${\cat M}$ is a cofibration category, we obtain a natural map
\begin{equation}
\label{eqn:homotopycolimbasechange}
{\bf L}{\colim}^{(u \downarrow d_2)} X \Ra ({\bf L}{\colim}^u X)_{d_2}
\end{equation}

Dually, the {\it standard under 2-category diagram} of $u$ at $d_2$ is
\begin{equation}
\label{eqn:stdunder2cat}
\xymatrix {
    (d_2 \downarrow u) \ar[r]^-{i_{d_2, u}} \ar[d]_-{p_{(d_2 \downarrow u)}} & 
    {\cat D}_1 \ar[d]^u _\,="b" \\
    \terminalcat \ar[r]_{e_{d_2}} ^\,="a" & 
    {\cat D}_2 \ultwocell<\omit>{<0>{{\phi}_{d_2, u} \;\;\;\;\;\;}}
  }
\end{equation}
If ${\cat M}$ is a fibration category, we obtain a natural map
\begin{equation}
\label{eqn:homotopylimbasechange}
({\bf R}{\lim}^u X)_{d_2} \Ra {\bf R}{\lim}^{(d_2 \downarrow u)} X
\end{equation}

The next theorem proves a base change formula for homotopy (co)limits. This lemma is a homotopy colimit analogue to the well known base change formula for ordinary colimits \lemmaref{lem::computerelcolim}. 

\begin{thm}[Base change property]
\label{thm:computerelhocolim}
\mbox{}
\begin{enumerate}
\item
If ${\cat M}$ is a cofibration category and $u \colon {\cat D}_1 \ra {\cat D}_2$ is a functor of small categories, then the natural map (\ref{eqn:homotopycolimbasechange}) induces an isomorphism
\begin{center}
${\bf L}{\colim}^{(u \downarrow d_2)} X \cong ({\bf L}{\colim}^u X)_{d_2}$
\end{center}
for objects $X \eps {\cat M}^{{\cat D}_1}$ and $d_2 \eps {\cat D}_2$.
\item
If ${\cat M}$ is a fibration category and $u \colon {\cat D}_1 \ra {\cat D}_2$ is a functor of small categories, then the natural map (\ref{eqn:homotopylimbasechange}) induces an isomorphism
\begin{center}
$({\bf R}{\lim}^u X)_{d_2} \cong {\bf R}{\lim}^{(d_2 \downarrow u)} X$
\end{center}
for objects $X \eps {\cat M}^{{\cat D}_1}$ and $d_2 \eps {\cat D}_2$.
\end{enumerate}
\end{thm}

\begin{proof}
We only prove (1). For a diagram $X \eps {\cat M}^{{\cat D}_1}$, pick a reduced Reedy cofibrant replacement $X^{'} \eps \Delta^{'}_{res}{\cat D}_1$ of $p_t^* X$. We have ${\colim}^{up_t} X^{'} \cong {\bf L}{\colim}^u X$. By \lemmaref{lem:deltaprimeovercat}, $(up_t \downarrow d_2) \cong \Delta^{'}(u \downarrow d_2)$. The restriction of $X^{'}$ to the direct category $\Delta^{'}(u \downarrow d_2)$ is a Reedy cofibrant replacement of the restriction of $X$ to $\Delta^{'}(u \downarrow d_2)$, so in ${\bf ho}{\cat M}$ we have isomorphisms ${\colim}^{(up_t \downarrow d_2)} X^{'}$ $\cong$ ${\colim}^{\Delta^{'}(u \downarrow d_2)} X^{'}$ $\cong$ ${\bf L}{\colim}^{(u \downarrow d_2)} X$. Using \lemmaref{lem::computerelcolim}, the top map and therefore all maps in the commutative diagram 
\begin{center}
$\xymatrix{
    {{\colim}^{(up_t \downarrow d_2)} X^{'}}_{\,} \ar[r] \ar[d] & ({\colim}^{up_t} X^{'})_{d_2} \ar[d] \\
    {{\bf L}{\colim}^{(u \downarrow d_2)} X}_{\,} \ar[r] & ({\bf L}{\colim}^u X)_{d_2}
  }$
\end{center}
are isomorphisms, and the conclusion is proved.
\end{proof}

Suppose that $({\cat M}, {\cat W})$ is a category with weak equivalences, and that $u \colon {\cat D}_1 \ra {\cat D}_2$ is a functor of small categories. The next result describes a sufficient condition for $({\bf L}\colim^u$, $\eps_u)$ to exist - without requiring ${\cat M}$ to carry a cofibration category structure. 

In preparation, notice that the natural map (\ref{eqn:homotopycolimbasechange}) actually exists under the weaker assumption that $({\cat M}, {\cat W})$ is a category with weak equivalences, that ${\bf L}{\colim}^u$ exists and is a left adjoint of ${\bf ho}u^*$, and that ${\bf L}{\colim}^{(u \downarrow d_2)}$ exists and is a left adjoint of ${\bf ho}p^*_{(u \downarrow d_2)}$. A dual statement holds for the map (\ref{eqn:homotopylimbasechange}).

\begin{thm}
\label{thm:openembeddingdiagrams}
Suppose that $({\cat M}, {\cat W})$ is a pointed category with weak equivalences.
\begin{enumerate}
\item
If $u \colon {\cat D}_1 \ra {\cat D}_2$ is a small closed embedding functor, then
\begin{enumerate}
\item
$\colim^u$ and its left Kan extension $({\bf L}\colim^u, \eps_u)$ exist, and ${\bf L}\colim^u$ is a fully faithful left adjoint to ${\bf ho}u^*$
\item
For any object $d_2$ of ${\cat D}_2$, the functors ${\colim}^{(u \downarrow d_2)}$ and ${\bf L}{\colim}^{(u \downarrow d_2)}$ are well defined and the natural map (\ref{eqn:homotopycolimbasechange}) induces an isomorphism
\begin{center}
${\bf L}{\colim}^{(u \downarrow d_2)} X \cong ({\bf L}{\colim}^u X)_{d_2}$
\end{center}
for objects $X \eps {\cat M}^{{\cat D}_1}$ and $d_2 \eps {\cat D}_2$.
\end{enumerate}
\item
If $u \colon {\cat D}_1 \ra {\cat D}_2$ is a small open embedding functor, then
\begin{enumerate}
\item
$\lim^u$ and its right Kan extension $({\bf R}\lim^u, \nu_u)$ exist, and ${\bf R}\lim^u$ is a fully faithful left adjoint to ${\bf ho}u^*$
\item
For any object $d_2$ of ${\cat D}_2$, the functors ${\lim}^{(d_2 \downarrow u)}$ and ${\bf R}{\lim}^{(d_2 \downarrow u)}$ are well defined and the natural map (\ref{eqn:homotopylimbasechange}) induces an isomorphism
\begin{center}
$({\bf R}{\lim}^u X)_{d_2} \cong {\bf R}{\lim}^{(d_2 \downarrow u)} X$
\end{center}
for objects $X \eps {\cat M}^{{\cat D}_1}$ and $d_2 \eps {\cat D}_2$.
\end{enumerate}
\end{enumerate}
\end{thm}

\begin{proof}
We start with (1) (a). Denote $u_! \colon {\cat M}^{{\cat D}_1} \ra {\cat M}^{{\cat D}_2}$ the 'extension by zero' functor, that sends a diagram $X \eps {\cat M}^{{\cat D}_1}$ to the diagram given by $(u_!X)_{d_2} = X_{d_2}$ for $d_2 \eps u{\cat D}_1$ and $(u_!X)_{d_2} = \initial$ otherwise. The functor $u_!$ is a fully faithful left adjoint to $u^*$, and sends weak equivalences to weak equivalences. In particular, $\colim^u \cong u_1$ exists and is defined on the entire ${\cat M}^{{\cat D}_1}$, and ${\bf L}\colim^u$ as in \definitionref{defn:homotopycolim} exists and is isomorphic to ${\bf ho}u_!$. We apply the Abstract Quillen Adjunction \theoremref{thm:generalexistencetotalderivedadjoint} to 
\begin{center}
$\xymatrix {{\cat M}^{{\cat D}_2} \ar[r]^{t_1 = id} & {\cat M}^{{\cat D}_2} \ar@<2pt>[r]^{u_1 = u_!} & {\cat M}^{{\cat D}_1} \ar@<2pt>[l]^{u_2 = u^*} & {\cat M}^{{\cat D}_1} \ar[l]_{t_2 = id} }$
\end{center}
observing that its hypotheses (1)-(3) and (4l) apply. We deduce that ${\bf L}\colim^u \cong {\bf ho}u_!$ is a fully faithful right adjoint to ${\bf ho}u^*$. 

To prove the isomorphism (1) (b), observe that $(u \downarrow d_2)$ has $d_2$ as a terminal object if $d_2 \eps u{\cat D}_1$ and is empty otherwise, so ${\colim}^{(u \downarrow d_2)} X \cong X_{d_2}$ if $d_2 \eps u{\cat D}_1$ and $\cong \initial$ otherwise. Both functors ${\colim}^{(u \downarrow d_2)}$ and $\colim^u$ preserve weak equivalences, and we have adjoint pairs ${\bf L}{\colim}^{(u \downarrow d_2)} \dashv {\bf ho}p^*_{(u \downarrow d_2)}$ and ${\bf L}{\colim}^u \dashv {\bf ho}u^*$. The natural isomorphism ${\colim}^{(u \downarrow d_2)} X \cong (\colim^u X)_{d_2}$ yields the desired isomorphism (\ref{eqn:homotopycolimbasechange}).

The proof of part (2) is dual.
\end{proof}

The two lemmas below are part of the proof of \theoremref{thm:leftadjointoflcolim} below. We keep the notations used in \theoremref{thm:openembeddingdiagrams} and its proof, and introduce a few new ones. ${\cat M}$ denotes a pointed cofibration category, and $u \colon {\cat D}_1 \ra {\cat D}_2$ a small closed embedding functor. $v \colon {\cat D}_2 \backslash {\cat D}_1 \ra {\cat D}_2$ denotes the inclusion functor - it is an {\it open} embedding. We also denote $V = \Delta^{'}v \colon {\Delta^{'}({\cat D}_2 \backslash {\cat D}_1)} \ra {\Delta^{'}{\cat D}_2}$. 

We denote $u_! = \colim^u \colon {\cat M}^{{\cat D}_1} \ra {\cat M}^{{\cat D}_2}$ and $v_* = \lim^v \colon {\cat M}^{{\cat D}_2 \backslash {\cat D}_1} \ra {\cat M}^{{\cat D}_2}$  - these are the 'extension by zero' functors. We also denote $v_! = \colim^v$ and $V_! = \colim^V$, and we will keep in mind that they are defined only on a {\it full subcategory} of ${\cat M}^{{\cat D}_2 \backslash {\cat D}_1}$ resp. $ {\cat M}^{\Delta^{'}({\cat D}_2 \backslash {\cat D}_1)}$.

We denote $({\cat M}^{{\cat D}_2})_0$ the full subcategory of ${\cat M}^{{\cat D}_2}$ consisting of objects $X$ with the property that $v_!v^*X$ exists and is pointwise cofibrant, and the map $v_!v^*X \ra X$ is a pointwise cofibration. $\initial$ is a cofibrant object, and we define the functor $u^{'}_1 \colon ({\cat M}^{{\cat D}_2})_0 \ra {\cat M}^{{\cat D}_2}$ as the pushout
\begin{center}
$\xymatrix{
    v_!v^*X \ar@{>->}[r] \ar[d] & X \ar@{-->}[d] \\
    \initial \ar@{>-->}[r] & u_1^{'}X
  }$
\end{center}

\begin{lem}
\label{lem:leftadjointoflcolim1}
\mbox{}
\begin{enumerate}
\item
There exists a canonical partial adjunction
\begin{center}
$\xymatrix{
    ({\cat M}^{{\cat D}_2})_0 \ar[r]^{u^{'}_1} \ar@{_{(}->}[d] & {\cat M}^{{\cat D}_2} \\
    {\cat M}^{{\cat D}_2} & {\cat M}^{{\cat D}_2} \ar[l]^{u_!u^*} \ar[u]_{1_{{\cat M}^{{\cat D}_2}}}
  }$
\end{center}
\item
There exists a canonical partial adjunction
\begin{center}
$\xymatrix{
    ({\cat M}^{{\cat D}_2})_0 \ar[r]^{u^*u^{'}_1} \ar@{_{(}->}[d] & {\cat M}^{{\cat D}_1} \\
    {\cat M}^{{\cat D}_2} & {\cat M}^{{\cat D}_1} \ar[l]^{u_!} \ar[u]_{1_{{\cat M}^{{\cat D}_1}}}
  }$
\end{center}
\end{enumerate}
\end{lem}

\begin{proof}
Denote $\eps_X \colon v_!v^*X \rightarrowtail X$ the natural map defined for $X \eps ({\cat M}^{{\cat D}_2})_0$.

For any object $Z \eps {\cat M}^{{\cat D}_2}$, the diagram 
\begin{center}
$\xymatrix{
    u_!u^*Z \ar[r] \ar[d] & Z \ar[d]^{\nu_Z} \\
    \initial \ar[r] & v_*v^*Z
  }$
\end{center}
is both a pushout and a pullback. For $X \eps ({\cat M}^{{\cat D}_2})_0$, the maps $X \ra u_!u^*Z$ are in a 1-1 correspondence with maps $f: X \ra Z$ such that $\nu_Z f: X \ra v_*v^*Z$ is null, therefore in 1-1 correspondence with maps $f: X \ra Z$ such that $f\eps_X \colon v_!v^*X \ra Z$ is null, therefore in 1-1 correspondence with maps $u^{'}_1 X \ra Z$. This shows that we have a natural bijection $Hom (u^{'}_1 X, Z) \cong Hom (X, u_!u^*Z)$, which proves (1). 

For an object $Y \eps {\cat M}^{{\cat D}_1}$ and a map $u^*u^{'}_1 X \ra Y$, denote $\overline{Y}^{X} \eps {\cat M}^{{\cat D}_2}$ the object with $\overline{Y}^{X}_{d} = Y_d$ for $d \eps {\cat D}_1$ and $\overline{Y}^{X}_{d} = (u^{'}_1 X)_d$ otherwise. We see that $Hom (u^*u^{'}_1 X, Y)$ $\cong$ $Hom (u^{'}_1 X, \overline{Y}^{X})$ $\cong$ $Hom (X, u_!u^*\overline{Y}^{X})$ $\cong$ $Hom (X, u_!Y)$, which completes the proof of (2).
\end{proof}

\begin{lem}
\label{lem:leftadjointoflcolim2}
For any diagrams $Y, Y^{'} \eps {\cat M}^{\Delta^{'}{\cat D}_2}$ that are Reedy cofibrant, denote $X = \colim^{p_t} Y$ and $X^{'} = \colim^{p_t} Y^{'}$ in ${\cat M}^{{\cat D}_2}$. Then
\begin{enumerate}
\item
The colimit $v_!v^*X$ exists and is pointwise cofibrant, and $v_!v^*X \rightarrowtail X$ is a pointwise cofibration
\item
If $Y \ra Y^{'}$ is a pointwise weak equivalence, then so is $v_!v^*X \ra v_!v^*X^{'}$
\end{enumerate}
\end{lem}

\begin{proof}
The objects of ${\Delta^{'}{\cat D}_2}$ are all of the form 
\begin{center}
$\underline{d} = (d_0 \ra ... \ra d_i \ra d^{'}_{i+1} \ra ... \ra d^{'}_n)$
\end{center}
where $d_0, ..., d_i \eps {\cat D}_2 \backslash {\cat D}_1$ and $d^{'}_{i+1}, ..., d^{'}_n \eps {\cat D}_1$. As a consequence, given $Y \eps {\cat M}^{\Delta^{'}({\cat D}_2 \backslash {\cat D}_1)}$ by \lemmaref{lem::computerelcolim} we have $(V_!Y)_{\underline{d}} \cong Y_{d_0 \ra ... \ra d_i}$. The functor $V_!$ is thus defined on the entire ${\cat M}^{\Delta^{'}({\cat D}_2 \backslash {\cat D}_1)}$ - note that $v_!$ may not be defined on the entire ${\cat M}^{{\cat D}_2 \backslash {\cat D}_1}$. 

The latching map of $V_!Y$ at $\underline{d}$ is $LY_{\underline{d}} \ra Y_{\underline{d}}$ if $i = n$, and $Y_{d_0 \ra ... \ra d_i} \overset{id}{\ra}$ $Y_{d_0 \ra ... \ra d_i}$ if $i < n$. Based on this, we see that if $Y \eps {\cat M}^{\Delta^{'}{\cat D}_2}$ is Reedy cofibrant then $V^*Y$ is Reedy cofibrant and $V_!V^*Y \ra Y$ is a Reedy cofibration.

Pick $Y \eps {\cat M}^{\Delta^{'}{\cat D}_2}_{res, rcof}$ with $\colim^{p_t} Y = X$. Using \lemmaref{lem::computerelcolim} we see that $v^*X \cong \colim^{p_t} V^* Y$.
\begin{center}
$\xymatrix{
    {\cat M}^{\Delta^{'}{\cat D}_2 \backslash {\cat D}_1}_{rcof} \ar[r]^V \ar[d]_{\colim^{p_t}} & {\cat M}^{\Delta^{'}{\cat D}_2}_{rcof} \ar[d]^{\colim^{p_t}} \\
    {\cat M}^{{\cat D}_2 \backslash {\cat D}_1} \ar[r]^v & {\cat M}^{{\cat D}_2}
  }$
\end{center}
The colimits $\colim^{p_t} V^* Y \cong v^*X$ and $\colim^{p_t} V_!V^*Y$ exist (the latter since $V_!V^*Y$ is Reedy cofibrant). By \lemmaref{lem::colimcomposition} we have that $v_!\colim^{p_t} V^* Y \cong v_!v^*X$ exists and is $\cong \colim^{p_t}V_!V^*Y$. Applying $\colim^{p_t}$ to the Reedy cofibration $V_!V^*Y \ra Y$ yields $v_!v^*X \ra X$, which is a pointwise cofibration between pointwise cofibrant objects by \theoremref{thm::colimfromdirectcat}. This proves part (1).

If $Y \ra Y^{'}$ is a pointwise weak equivalence between Reedy cofibrant objects, then so is $V_!V^*Y \ra V_!V^*Y^{'}$, so by \theoremref{thm::colimfromdirectcat} the map $\colim^{p_t}V_!V^*Y \ra \colim^{p_t}V_!V^*Y^{'}$ is a weak equivalence. This proves part (2).
\end{proof}

\begin{thm}
\label{thm:leftadjointoflcolim}
\mbox{}
\begin{enumerate}
\item
If ${\cat M}$ is a pointed cofibration category and $u \colon {\cat D}_1 \ra {\cat D}_2$ is a small closed embedding functor, then ${\bf L}\colim^u$ admits a left adjoint.
\item
If ${\cat M}$ is a pointed fibration category and $u \colon {\cat D}_1 \ra {\cat D}_2$ is a small open embedding functor, then ${\bf R}\lim^u$ admits a right adjoint.
\end{enumerate}
\end{thm}

\begin{proof}
We only prove (1). We will apply the Abstract Partial Quillen Adjunction \theoremref{thm:generalexistencetotalderivedadjoint2} to
\begin{center}
$\xymatrix{
    {\cat M}^{\Delta^{'}{\cat D}_2}_{res, rcof} \ar[rrr]^-{v_1 = u^*u^{'}_1\colim^{p_t}} \ar[d]_{t_1 = \colim^{p_t}} &&& {\cat M}^{{\cat D}_1} \\
    {\cat M}^{{\cat D}_2} &&& {\cat M}^{{\cat D}_1} \ar[lll]_-{v_2 = u_!} \ar[u]_{t_2 = 1_{{\cat M}^{{\cat D}_1}}}
    }$
\end{center}
We have that $v_1, v_2$ is an abstract Quillen partially adjoint pair with respect to $t_1, t_2$:
\begin{enumerate}
\item
From \lemmaref{lem:leftadjointoflcolim2} (1) we have $Im \, \colim^{p_t} \subset ({\cat M}^{{\cat D}_2})_0$, so the functor $v_1$ is correctly defined. Using \lemmaref{lem:leftadjointoflcolim1} (2) we see that $v_1, v_2$ are partially adjoint with respect to $t_1, t_2$.
\item
The functor $t_1$ is a cofibrant approximation of ${\cat M}^{{\cat D}_2}$ by \propositionref{prop:generalpointwisecofstructureapprox}, therefore a left approximation. The functor $t_2 = 1_{{\cat M}^{{\cat D}_1}}$ is a right approximation, and $u_2t_2$ preserves weak equivalences.
\item
By \lemmaref{lem:leftadjointoflcolim2} (2), the functor $v_!v^*t_1$ preserves weak equivalences. But $u^{'}_1t_1Y$ is the pushout of $v_!v^*t_1Y \rightarrowtail t_1 Y$ by $v_!v^*t_1Y \ra 0$, and an application of the Gluing Lemma shows that $u^{'}_1t_1$ and therefore $v_1$ also preserve weak equivalences.
\item
The functor $v_2$ preserves weak equivalences.
\end{enumerate}

We conclude that ${\bf R}v_2 \cong {\bf ho}u_! \cong {\bf L}\colim^u$ admits a left adjoint.
\end{proof}


\section{The conservation property}
\label{sec:conservationproperty}

\index{functor!conservative functors}
Recall that a family of functors $u_i \colon {\cat A} \ra {\cat B}_i, i \eps I$ is {\it conservative} if for any map $f \eps {\cat A}$ with $u_if$ an isomorphism in ${\cat B}_i$ for all $i \eps I$ we have that $f$ is an isomorphism in ${\cat A}$. A family of functors $\{ u_i\}$ is conservative iff the functor $(u_i)_i \colon {\cat A} \ra \times_i {\cat B}_i$ is conservative.

\begin{thm}
\label{thm:homotcatconservative}
Let ${\cat D}$ be a small category, and suppose that ${\cat M}$ is either a cofibration or a fibration category. The projections $p_d \colon {\cat M}^{\cat D} \ra {\cat M} $ on the $d$ componenent for all objects $d \eps {\cat D}$ then induce a conservative family of functors ${\bf ho}(p_d) \colon {\bf ho}{\cat M}^{\cat D} \ra {\bf ho}{\cat M}$.
\end{thm}

\begin{proof}
Assume that ${\cat M}$ is a cofibration category (the proof for fibration categories is dual).

Let $\overline{f} \colon A \ra B$ be a map in ${\bf ho}{\cat M}^{\cat D}$ such that ${\bf ho}(p_d)\overline{f} \eps {\bf ho}{\cat M}$ are isomorphisms for all objects $d \eps {\cat D}$. We want to show that $\overline{f}$ is an isomorphism in ${\bf ho}{\cat M}^{\cat D}$.

Using the factorization axiom CF4 applied to the pointwise cofibration structure $({\cat M}^{\cat D}$, ${\cat W}^{\cat D}$, ${\cat Cof}^{\cat D})$, we may assume that $A, B$ are pointwise cofibrant. From \theoremref{thm::descrhomotcof} applied to $({\cat M}^{\cat D}$, ${\cat W}^{\cat D}$, ${\cat Cof}^{\cat D})$, we may further assume that $\overline{f}$ is the image of a pointwise cofibration $f$ of ${\cat M}^{\cat D}$.

Assume we proved our theorem for all {\it direct} small categories. The map $p_t^*f$ in ${\cat M}^{\Delta^{'}{\cat D}}$ satisfies the hypothesis of our theorem and $\Delta^{'}{\cat D}$ is direct. It follows that $p_t^*f$ is an isomorphism in ${\bf ho}{\cat M}^{\Delta^{'}{\cat D}}$, so by \lemmaref{lem::cofcatequivinhomotcat} there exist pointwise cofibrations $f^{'}, f^{''} \eps {\cat M}^{\Delta^{'}{\cat D}}$ such that $f^{'}p^{*}f, f^{''}f^{'}$ are pointwise weak equivalences in ${\cat M}^{\Delta^{'}{\cat D}}$. 

But $p^{*}f$ is $\Delta^{'}_{res}{\cat D}$ restricted, and therefore so are $f^{'}$ and $f^{''}$. By the same \lemmaref{lem::cofcatequivinhomotcat} $p^{*}f$ is an isomorphism in ${\bf ho}{\cat M}^{\Delta^{'}{\cat D}}_{res}$, and by \theoremref{thm:pointwisecofstructureequivalence} $f$ is an isomorphism in ${\bf ho}{\cat M}^{{\cat D}}$

It remains to prove our theorem in the case when ${\cat D}$ is direct.

As before, let $\overline{f} \colon A \ra B$ be a map in ${\bf ho}{\cat M}^{\cat D}$ such that ${\bf ho}(p_d)\overline{f} \eps {\bf ho}{\cat M}$ are isomorphisms for all objects $d \eps {\cat D}$, and we want to show that $\overline{f}$ is an isomorphism in ${\bf ho}{\cat M}^{\cat D}$. Repeating the previous argument applied to the Reedy cofibration structure $({\cat M}^{\cat D}, {\cat W}^{\cat D}, {\cat Cof}^{\cat D}_{Reedy})$, we may further assume that $\overline{f}$ is the image of a Reedy cofibration $f \colon A \ra B$.  

We will construct a Reedy cofibration $f^{'} \colon B \ra B^{'}$ such that $f^{'}f \eps {\cat W}^{\cat D}$. Once the construction is complete, we will be able to apply the same construction to $f^{'}$ and obtain a Reedy cofibration $f^{''} \colon B^{'} \ra B^{''}$ such that $f^{''}f^{'} \eps {\cat W}^{\cat D}$. As a consequence, it will follow that $f$ is an isomorphism in ${\bf ho}{\cat M}^{\cat D}$.

To summarize, given a Reedy cofibration $f \colon A \ra B$ with the property that $f_d$ is an isomorphism in ${\bf ho}{\cat M}$, it remains to construct a Reedy cofibration $f^{'} \colon B \ra B^{'}$ such that $f^{'}f \eps {\cat W}^{\cat D}$. We will construct $f^{'}$ by induction on degree.

For $n = 0$, ${\cat D}^{0}$ is discrete and the existence of $f^{' \, 0}$ follows from \lemmaref{lem::cofcatequivinhomotcat}. Assume now $f^{' \, < n} \colon B^{< n} \ra B^{' \, < n}$ constructed. 

For any object $d \eps {\cat D}^n$, we construct the following diagram:
\begin{center}
$ \xymatrix {
        L A_d \ar@{>->}[r]^-{i_d} \ar@{>->}[d]_{L(f^{'} f)_d}^\sim &
        A_d \ar@{>->}[r]^{f_d} \ar@{>->}[d]^\sim_{\beta_d} & B_d \ar@{>->}[r]^{\overline{f^{'}_d}} & \overline{B_d} \ar@{>->}[d]^\sim_{\overline{\overline{f^{'}_d}}} \\
        L B^{'}_d \ar@{>->}[r]^-{\gamma_d} & \overline{B^{'}_d} \ar@{>->}[rr]_\sim^{\delta_d} & & B^{'}_d
	}$
\end{center}
where:
\begin{enumerate}
\item
$LA_d$, $LB^{'}_d$ exist and are cofibrant because $A, B^{' \, < n}$ are Reedy cofibrant
\item
Since $f^{' \, <n}f^{<n}$ is a trivial Reedy cofibration, $L(f^{'}f)_d$ is a trivial cofibration
\item
The map $i_d$ is a cofibration since $A$ is Reedy cofibrant
\item
$\beta_d$ is constructed as the pushout of $L(f^{'}f)_d$, therefore a trivial cofibration. $\gamma_d$ is a pushout of $i_d$, therefore a cofibration.
\item
$\overline{f^{'}_d}$ is a cofibration constructed by \lemmaref{lem::cofcatequivinhomotcat} applied to $f_d$, so that $\overline{f^{'}_d} f_d$ is a trivial cofibration
\item
$\delta_d$ is constructed as the pushout of $\overline{f^{'}}_d f_d$, therefore a trivial cofibration. $\overline{\overline{f^{'}_d}}$ is a pushout of $\beta_d$, therefore a trivial cofibration.
\end{enumerate}

We define $f^{'}_d = \overline{\overline{f^{'}_d}} \, \overline{f^{'}_d}$, for all objects $d \eps {\cat D}^n$. The map $f^{' \, \le n}$ is a Reedy cofibration, and $(f^{'}f)^{\le n} \eps {\cat W}^{\le n}$. The inductive step is now complete, and the proof is finished.
\end{proof}

For each category with weak equivalences $({\cat M}, {\cat W})$ and small category ${\cat D}$, the functor ${\cat M}^{\cat D} \ra ({\bf ho}{\cat M})^{\cat D}$ induces a functor denoted $\mathrm{dgm}_{_{{\cat D}, {\cat M}}} \colon {\bf ho}({\cat M}^{\cat D}) \ra ({\bf ho}{\cat M})^{\cat D}$, or simply
\begin{equation}
\label{eqn:dgm}
\mathrm{dgm}_{_{{\cat D}}} \colon {\bf ho}({\cat M}^{\cat D}) \ra ({\bf ho}{\cat M})^{\cat D}
\end{equation}
when $({\cat M}, {\cat W})$ is inferred from the context.

\begin{cor}
\label{cor:homotcatconservative}
If ${\cat M}$ is either a cofibration or a fibration category, and if ${\cat D}$ is a small category, then the functor $\mathrm{dgm}_{_{\cat D}}$ of (\ref{eqn:dgm}) is a conservative functor.
\end{cor}

\begin{proof}
Consequence of the fact that a map $f$ of ${\cat D}$ diagrams is an isomorphism iff each $f_d$ for $d \eps {\cat D}$ is an isomorphism.
\end{proof}

\section{Realizing diagrams}
\label{sec:realizingdiagrams}

Start with a category with weak equivalences $({\cat M}, {\cat W})$ and a small category ${\cat D}$. Given a ${\cat D}$-diagram $X$ in ${\bf ho}{\cat M}$, we would like to know under what conditions the diagram $X$ is isomorphic to the image of a diagram $X^{'} \eps {\bf ho}({\cat M}^{\cat D})$ under the functor $\mathrm{dgm}_{_{\cat D}}$ of (\ref{eqn:dgm}). If such a diagram $X^{'}$ exists, it is called a {\it realization} of the diagram $X$ in the homotopy category ${\bf ho}({\cat M}^{\cat D})$.

We will show that if ${\cat M}$ is a cofibration category and ${\cat D}$ is a small, direct, free category then any diagram $X \eps ({\bf ho}{\cat M})^{\cat D}$ admits a realization $X^{'} \eps {\bf ho}({\cat M}^{\cat D})$. Furthermore any two such realizations $X^{'}, X^{''}$ of $X$ are non-canonically isomorphic in ${\bf ho}({\cat M}^{\cat D})$. Denis-Charles Cisinski \cite{Cisinski2} proves this result for {\it finite} direct, free categories, but his techniques extend to our situation.

Dually, if ${\cat M}$ is a fibration category and ${\cat D}$ is a small, inverse, free category then any diagram $X \eps ({\bf ho}{\cat M})^{\cat D}$ admits a realization $X^{'} \eps {\bf ho}({\cat M}^{\cat D})$, and $X^{'}$ is unique up to a non-canonical isomorphism in ${\bf ho}({\cat M}^{\cat D})$. 

Let us recall the definition of a small free category. A {\it directed graph} ${\cat G}$ consists of a set of vertices ${\cat G}_0$, a set of arrows ${\cat G}_1$ along with two functions $s, t \colon {\cat G}_1 \ra {\cat G}_0$ giving each arrow a source respectively a target. A {\it map} of directed graphs $u \colon {\cat G} \ra {\cat G}^{'}$ consists of two functions $u_i \colon {\cat G}_i \ra {\cat G}^{'}_i$ for $i = 0, 1$ that commute with the source and destination maps. We denote ${\cat Graph}$ the category of directed graphs. 

There is a functor $F \colon {\cat Cat} \ra {\cat Graph}$ which sends a small category ${\cat D}$ to the underlying graph $F{\cat D}$ - with the objects of ${\cat D}$ as vertices and the maps of ${\cat D}$ as arrows, forgetting the composition of maps. The functor $F$ has a left adjoint $G \colon {\cat Graph} \ra {\cat Cat}$, which sends a graph ${\cat G}$ to the small category $G{\cat G}$ with objects ${\cat G}_0$ and with non-identity maps between $x, x^{'} \eps {\cat G}_0$ defined as all formal compositions $f_nf_{n-1}...f_0$ of arrows $f_i$, for $i \ge 0$, such that $sf_0 = x, tf_n = x^{'}$ and $tf_i = sf_{i+1}$ for $0 \le i < n$.

A small category ${\cat D}$ is {\it free} if it is {\it isomorphic} to $G{\cat G}$ for a graph ${\cat G}$. The generators of a category are its undecomposable maps, i.e. maps that cannot be written as compositions of two non-identity maps. For a small category ${\cat D}$, one can form the graph ${\cat G}$ of undecomposable maps, with the objects of ${\cat D}$ as vertices and the undecomposable maps of ${\cat D}$ as arrows. Since ${\cat G}$ is a subgraph of $F{\cat D}$, there is a functor $G{\cat G} \ra {\cat D}$ which is an {\it isomorphism} iff the category ${\cat D}$ is free.



A vertex $x_0$ of a graph ${\cat G}$ has {\it uniformly bounded ascending chains} if there is a positive integer $k$ such that any sequence of arrows $x_0 \ra x_1 \ra x_2 ...$ has at most length $k$. The vertex $x_0$ has {\it uniformly bounded descending chains} if there is a positive integer $k$ such that any sequence of arrows $... \ra x_2 \ra x_1 \ra x_0$ has at most length $k$. We leave the proof of the following lemma to the reader:

\begin{lem}
\label{lem:descrfinitefreecat}
\mbox{}
\begin{enumerate}
\item
The following statements are equivalent for a small category ${\cat D}$.
\begin{enumerate}
\item
${\cat D}$ is direct and free
\item
${\cat D}$ is free, and all vertices of its graph of undecomposable maps have uniformly bounded descending chains (in particular, its graph of undecomposable maps has no loops)
\item
${\cat D}$ is direct and for any $d \eps {\cat D}$ the latching category $\partial({\cat D} \downarrow d)$ is a disjoint sum of categories with a terminal object.
\end{enumerate}
\item
The following statements are equivalent for a small category ${\cat D}$.
\begin{enumerate}
\item
${\cat D}$ is inverse and free
\item
${\cat D}$ is free, and all vertices of its graph of undecomposable maps have uniformly bounded ascending chains
\item
${\cat D}$ is inverse and for any $d \eps {\cat D}$ the matching category $\partial(d \downarrow {\cat D})$ is a disjoint sum of categories with an initial object. 
\end{enumerate}
\item
A finite, free category is both direct and inverse. $\square$
\end{enumerate}

\end{lem}

For example, the categories $0 \dbra 1$, $\,\,1 \la 0 \ra 2$ and $1 \ra 0 \la 2$ (given by their subgraph of non-identity maps) are at the same time direct, inverse and free. The category $\mathbf{N}$ of nonnegative integers $0 \ra 1 \ra ...$ is free and direct, and its opposite $\mathbf{N}^{op}$ is free and inverse. The category given by the commutative diagram
\[
\xymatrix{
  0 \ar[r] \ar[d] & 1 \ar[d] \\
  2 \ar[r] & 3
}
\]
is finite, direct and inverse but not free.

Let us now investigate the more trivial problem of realizing ${\cat D}$-diagrams for a {\it discrete} small category ${\cat D}$.

\begin{lem}
\label{lem:cofcatsumsofdiagrams}
Suppose that ${\cat D}_k$, $k \eps K$ is a set of small categories. If ${\cat M}$ is either a cofibration or a fibration category, then the functor
\begin{center}
$\xymatrix {
    {\bf ho}{\cat M}^{{\Sum}_{k \eps K} {\cat D}_k}
    \ar[rr]^-{{\bf ho}(i_k^*)_{k \eps K}} &&
    {\times}_{k \eps K} {\bf ho}{\cat M}^{{\cat D}_k}
  }$
\end{center}
is an isomorphism of categories, where $i_{k} \colon {\cat D}_{k} \ra {\Sum}_{k^{'} \eps K} {\cat D}_{k^{'}}$ denotes the component inclusion for $k \eps K$.
\end{lem}

\begin{proof}
Consequence of \theoremref{thm:prodhomotopycat} for ${\cat M}_k = {\cat M}^{{\cat D}_k}$, which is a cofibration category by \theoremref{thm:generalpointwisecofstructure}.
\end{proof}

We can apply the existence of the homotopy (co)limit functor for the particular case of discrete diagrams to prove:

\begin{lem}
\label{lem:cofcathomotopysum}
\mbox{}
\begin{enumerate}
\item
If ${\cat M}$ is a cofibration category, then ${\bf ho}{\cat M}$ admits all (small) sums of objects. The functor ${\cat M}_{cof} \ra {\bf ho}{\cat M}$ commutes with sums of objects.
\item
If ${\cat M}$ is a fibration category, then ${\bf ho}{\cat M}$ admits all (small) products of objects. The functor ${\cat M}_{fib} \ra {\bf ho}{\cat M}$ commutes with products of objects.
\end{enumerate}
\end{lem}

\begin{proof}
We only prove (1). Given a set of objects $X_k, k \eps K$ of a cofibration category ${\cat M}$, if we view $K$ as a discrete category then by \theoremref{thm:prodhomotopycat} ${\bf L}\colim^{K} X_k$ satisfies the universal property of the sum of $X_k$ in ${\bf ho}{\cat M}$. If all $X_k$ are cofibrant, then $\Sum_{k \eps K} X_k$ exists in ${\cat M}$ and computes ${\bf L}\colim^{K} X_k$, which proves the second part.
\end{proof}

We now turn to the problem of realizing ${\cat D}$-diagrams in a cofibration category for a small direct category ${\cat D}$.

\begin{lem}
\label{lem:sufficientfinitefull}
Suppose that either ${\cat M}$ is a cofibration category and ${\cat D}$ is a small direct category, or ${\cat M}$ is a fibration category and ${\cat D}$ is a small inverse category.

Suppose that $X, Y$ are objects of ${\cat M}^{\cat D}$. If for each $n \ge 0$ there exists 
\[
f_n \eps Hom_{{\bf ho}{{\cat M}^{{\cat D}^{\le n}}}} (X^{\le n}, Y^{\le n})
\]
 such that $f_{n}$ restricts to $f_{n-1}$ for all $n > 0$, then there exists $f \eps Hom_{{\bf ho}{{\cat M}^{{\cat D}}}} (X, Y)$ such that $f$ restricts to $f_n$ for all $n \ge 0$.
\end{lem}

\begin{proof}
Assume that ${\cat M}$ is a cofibration category and that ${\cat D}$ is a small direct category. (The proof using the alternative hypothesis is dual).

We may assume that $X, Y$ are Reedy cofibrant. We fix a sequence of cylinders with respect to the Reedy cofibration structure $I^nX$ and $I^nY$, and denote $I^{\infty}X = \colim (X \overset{i_0}{\ra} IX \overset{i_0}{\ra} I^2X ...)$ and $I^{\infty}Y = \colim (Y \overset{i_0}{\ra} IY \overset{i_0}{\ra} I^2Y ...)$. Denote $j_{n} \colon X \ra I^{n}X$ and $k_{n} \colon Y \ra I^{n}Y$ the trivial cofibrations given by iterated compositions of $i_0$. Using axiom CF6, the maps $j_{\infty} = \colim j_n \colon X \ra I^{\infty}X$ and $k_{\infty} = \colim k_n \colon Y \ra I^{\infty}Y$ are trivial cofibrations.

By induction on $n$, we will construct a factorization in ${\cat D}^{\le n}$
\begin{center}
$\xymatrix{
    I^nX^{\le n} \ar[r]^-{a_n} & 
    Z_n &
    I^nY^{\le n} \ar[l]^-{\sim}_-{b_n}
  }$
\end{center}
with $Z_n$ Reedy cofibrant and $b_n$ a weak equivalence in ${\cat D}^{\le n}$, such that $b_n^{-1}a_n$ has the homotopy type of $f_n$, and a trivial Reedy cofibration $c_{n-1}$ in ${\cat D}^{< n}$ that make the next diagram commutative
\begin{center}
$\xymatrix{
    I^{n-1}X^{< n} \ar[r]^-{a_{n-1}} \ar@{>->}[d]_{i_0}^\sim  & 
    Z_{n-1} \ar@{>->}[d]_-{c_{n-1}}^-\sim & 
    I^{n-1}Y^{< n} \ar[l]_-{b_{n-1}}^-\sim \ar@{>->}[d]_{i_0}^\sim \\
    I^{n}X^{< n} \ar[r]_-{a^{< n}_n} & 
    Z^{< n}_{n} &
    I^{n}Y^{< n} \ar[l]^-{b^{< n}_n}_-\sim
  }$
\end{center}

Once the inductive construction is complete, we get a factorization
\begin{center}
$\xymatrix{
    X \ar@{>->}[r]^-{j_{\infty}}_-\sim & 
    I^{\infty}X \ar[r]^-{a} & 
    Z &
    I^{\infty}Y \ar[l]^-{\sim}_-{b} &
    Y \ar@{>->}[l]^-{\sim}_-{k_{\infty}}
  }$
\end{center}
defined by $a = \colim^n a_n$ and $b = \colim^n b_n$. Using axiom CF6 and \lemmaref{lem:transcompequiv}, $Z$ is Reedy cofibrant and $b$ is a weak equivalence. The map $k_{\infty}^{-1}b^{-1}aj_{\infty}$ in ${\bf ho}{{\cat M}^{{\cat D}}}$ restricts to $f_n$ for all $n \ge 0$.

To complete the proof, let us perform the inductive construction of $Z_n, a_n, b_n, c_{n-1}$ for $n \ge 0$. The initial step construction of $Z_0, a_0, b_0$ folows from \lemmaref{lem:cofcatsumsofdiagrams}. Assume that $Z_{n-1}, a_{n-1}, b_{n-1}$ have been constructed.

From \theoremref{thm::descrhomotcof}, there exists a factorization in ${\cat D}^{\le n}$
\begin{center}
$\xymatrix{
    I^{n-1}X^{\le n} \ar@{>->}[r]^-{a^{'}_n} & 
    Z^{'}_n &
    I^{n-1}Y^{\le n} \ar@{>->}[l]^-{\sim}_-{b^{'}_n}
  }$
\end{center}
with $a^{'}_n$ a Reedy cofibration and $b^{'}_n$ a trivial Reedy cofibration in ${\cat D}^{\le n}$, such that $b_n^{'-1}a_n^{'}$ has the homotopy type of $f_n$ in ${\bf ho}{{\cat M}^{{\cat D}^{\le n}}}$. Using the inductive hypothesis and \theoremref{thm::descrhomot}, there exist a Reedy cofibrant diagram $Z^{'' < n}_{n}$ and weak equivalences $\alpha, \beta$ in ${\cat D}^{< n}$ that make the next diagram {\it homotopy} commutative
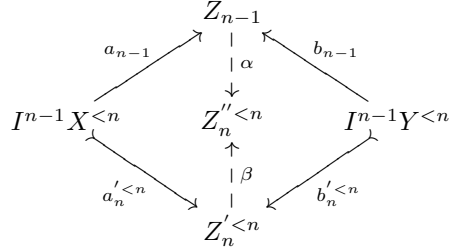
\begin{center}
$\xymatrix{
    & 
    Z_{n-1} \ar@{-->}[d]^-{\alpha} & 
    \\
    I^{n-1}X^{< n} \ar[ru]^-{a_{n-1}} \ar@{>->}[rd]_-{a^{' < n}_n} & 
    Z^{'' < n}_{n}  & 
    I^{n-1}Y^{< n} \ar[lu]_-{b_{n-1}} \ar@{>->}[ld]^-{b^{' < n}_n} \\
    & 
    Z^{' < n}_{n} \ar@{-->}[u]_-{\beta} &
  }$
\end{center}
We may assume that $\alpha, \beta$ are trivial Reedy cofibrations - if they are not, we may replace them with $\alpha^{'}, \beta^{'}$ defined by a cofibrant replacement
\begin{center}
$\alpha + \beta \colon Z_{n-1} \Sum Z^{' < n}_{n} \overset{\alpha^{'} + \beta^{'}}{\rightarrowtail} Z^{''' < n}_{n} \overset{\sim}{\ra} Z^{'' < n}_{n}$
\end{center}

Using \lemmaref{lem::choosecyl} twice, we can embed our homotopy commutative diagram in a diagram commutative on the nose with the maps $c_{n-1}, \gamma$ trivial Reedy cofibrations
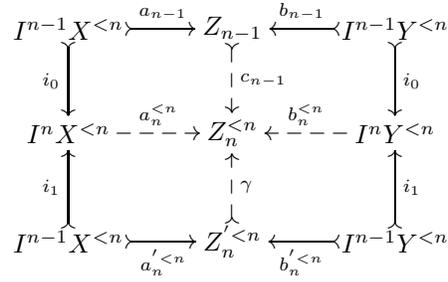
\begin{center}
$\xymatrix{
    I^{n-1}X^{< n} \ar@{>->}[r]^-{a_{n-1}} \ar@{>->}[d]_{i_0} & 
    Z_{n-1} \ar@{>-->}[d]^-{c_{n-1}} & 
    I^{n-1}Y^{< n} \ar@{>->}[l]_-{b_{n-1}} \ar@{>->}[d]^{i_0} \\
    I^{n}X^{< n} \ar@{-->}[r]^-{a^{< n}_n}  & 
    Z^{< n}_{n} & 
    I^{n}Y^{< n}  \ar@{-->}[l]_-{b^{< n}_n} \\
    I^{n-1}X^{< n} \ar@{>->}[r]_-{a^{' < n}_n} \ar@{>->}[u]^{i_1} & 
    Z^{' < n}_{n} \ar@{>-->}[u]_-{\gamma} &
    I^{n-1}Y^{< n}  \ar@{>->}[l]^-{b^{' < n}_n} \ar@{>->}[u]_{i_1}
  }$
\end{center}
This diagram defines a Reedy cofibrant object $Z_n^{< n}$, a map $a_n^{< n}$, a weak equivalence $b_n^{< n}$ and the trivial Reedy cofibration $c_{n-1}$. It remains to extend $Z_n^{< n}, a_n^{< n}, b_n^{< n}$ in degree $n$ such that $b_n^{-1}a_n = f_n$ in ${\bf ho}{{\cat M}^{{\cat D}^{\le n}}}$.

For each object $d \eps {\cat D}^{n}$, we construct $Z^{'''}_d$ as the colimit of the diagram of solid maps
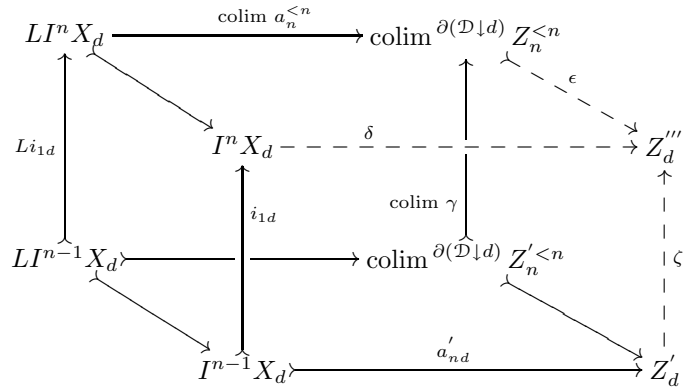
\begin{center}
$ \xymatrix {
        LI^{n}X_d \ar[rr]^{\colim a_n^{< n}} \ar@{>->}[rd]_{} &  & \colim^{\partial ({\cat D} \downarrow d)}Z^{< n}_{n} \ar@{>-->}[rd]^{\epsilon} & \\
        & I^{n}X_d \ar@{-->}[rr]^(.3){\delta}  & & Z^{'''}_d  \\
        LI^{n-1}X_d \ar@{>->}'[r][rr]_(.7){} \ar@{>->}[rd]_{} \ar@{>->}[uu]^{Li_{1d}} & & \colim^{\partial ({\cat D} \downarrow d)}Z^{' < n}_{n} \ar@{>->}[rd] \ar@{>->}'[u]^{\colim \gamma}[uu] & \\
        & I^{n-1}X_d \ar@{>->}[rr]^{a^{'}_{nd}} \ar@{>->}[uu]_(.7){i_{1d}} & & Z^{'}_d \ar@{-->}[uu]_{\zeta}
	}$
\end{center}
The bottom and left faces are Reedy cofibrant, and the maps $Li_{1d}, i_{1d}$ and $\colim \gamma$ are trivial Reedy cofibrations. A quick computation shows that the colimit $Z^{'''}_d$ actually exists, that $\eps$ is a Reedy cofibration and that the map $\zeta$ is a weak equivalence. In consequence $\delta$ has the homotopy type of $f_{nd}$. 

We also construct $Z^{''''}_d$ as the colimit of the diagram of solid maps
\begin{center}
$ \xymatrix {
        \colim^{\partial ({\cat D} \downarrow d)}Z^{< n}_{n} \ar@{>-->}[rd]_{\mu} &  & LI^{n}Y_d \ar@{>->}[rd] \ar[ll]_{\colim b_n^{< n}} & \\
        & Z^{''''}_d  & & I^{n}Y_d \ar@{-->}[ll]^(.7){\nu} \\
        \colim^{\partial ({\cat D} \downarrow d)}Z^{' < n}_{n}  \ar@{>->}[rd]_{} \ar@{>->}[uu]^{\colim \gamma} & & LI^{n-1}Y_d \ar@{>->}[rd] \ar@{>->}'[u]^{Li_{1d}}[uu] \ar@{>->}'[l]^(.8){\colim b_{n-1}}[ll] & \\
        & Z^{'}_d \ar@{-->}[uu]_(.7){\theta} & & I^{n-1}Y_d \ar@{>->}[uu]_{i_{1d}} \ar@{>->}[ll]_{b_{n-1,d}} 
	}$
\end{center}
The bottom and right faces are Reedy cofibrant, the maps $Li_{1d}$, $i_{1d}$, $\colim \gamma$, $\colim b_{n-1}$ and $b_{n-1,d}$ are trivial Reedy cofibrations and the map $\colim b_n^{< n}$ is a weak equivalence. By computation one shows that the colimit $Z^{''''}_d$ actually exists, that $\mu$ is a Reedy cofibration and that $\nu, \theta$ are weak equivalences. We now construct $Z_d$ as the pushout
\begin{center}
$\xymatrix{
    \colim^{\partial ({\cat D} \downarrow d)}Z^{< n}_{n} \ar@{>->}[r]^-{\mu} \ar@{>->}[d]_{\eps} & Z^{''''}_d \ar@{>-->}[d]_{{\eps}^{'}} \\
    Z^{'''}_d \ar@{>-->}[r]^{{\mu}^{'}} & Z_d
  }$
\end{center}
From the Gluing Lemma applied to the diagram
\begin{center}
$ \xymatrix {
        \colim^{\partial ({\cat D} \downarrow d)}Z^{< n}_{n} \ar@{>->}[rr]_{\mu} \ar@{>->}[rd]_{\eps} &  & \;\; Z^{''''}_d \;\; \ar@{>->}[rd]^{\eps^{'}} & \\
        & \;\; Z^{'''}_d \;\; \ar@{>->}[rr]^(.3){\mu^{'}}  & & \;\; Z_d \;\;  \\
        \colim^{\partial ({\cat D} \downarrow d)}Z^{' < n}_{n} \ar[uu]^{\colim \gamma} \ar@{>->}'[r][rr] \ar@{>->}[rd] && \;\; Z^{'}_d \;\; \ar[rd]^{id} \ar'[u]^{\theta}[uu] & \\
        & \;\; Z^{'}_d \;\; \ar[rr]^{id} \ar[uu]_(.7){\zeta} & & \;\; Z^{'}_d \;\; \ar@{>->}[uu]_{\iota}
	}$
\end{center}
where the top and bottom faces are pushouts and the vertical maps $\colim \gamma$, $\zeta$, $\theta$ are weak equivalences we see that $\iota$ and therefore ${\mu}^{'}, {\eps}^{'}$ are weak equivalences.

We define $a_{nd} = {\mu}^{'}\delta$, $b_{nd} = {\eps}^{'}\nu$. Repeating our construction for any $d \eps {\cat D}^n$ allows us to define $Z_n, a_n$ and $b_n$.

For all objects $d \eps {\cat D}^n$, the latching map $\eps^{'}\mu$ is a cofibration, therefore $Z_n$ is Reedy cofibrant. The map $b_n$ is a weak equivalence because $b_n^{< n}$ is a weak equivalence and for all objects $d \eps {\cat D}^n$ the maps ${\eps}^{'}, \nu$ are weak equivalences. The map $b_n^{-1}a_n$ has the homotopy type of $f_n \eps Hom_{{\bf ho}{\cat M}^{{\cat D}^{\le n}}}(X^{\le n}, Y^{\le n})$ because $b_n^{' -1}a_n^{'}$ has the homotopy type of $f_n$ and for all $d \eps {\cat D}^n$ the maps $\zeta$, $\theta$, $\nu$, $\mu^{'}$ and $\eps^{'}$ are weak equivalences.
\end{proof}

We now turn to the main theorem of this section.

\begin{thm}[Cisinski]
\label{thm:freedirectdiagramrealization}
\mbox{}
\begin{enumerate}
\item
If ${\cat M}$ is a cofibration category and ${\cat D}$ is a small, direct and free category, then the functor $\mathrm{dgm}_{_{\cat D}} \colon {\bf ho}({\cat M}^{\cat D}) \ra ({\bf ho}{\cat M})^{\cat D}$ is full and essentially surjective.
\item
If ${\cat M}$ is a fibration category and ${\cat D}$ is small, inverse and free then the functor $\mathrm{dgm}_{_{\cat D}} \colon$ ${\bf ho}({\cat M}^{\cat D})$ $\ra$ $({\bf ho}{\cat M})^{\cat D}$ is full and essentially surjective.
\end{enumerate}
\end{thm}

\begin{proof}
We only prove part (1). The category ${\cat M}^{\cat D}$ is endowed with a pointwise cofibration structure. Since ${\bf ho}({\cat M}^{\cat D}) \cong {\bf ho}({\cat M}^{\cat D}_{cof})$ and ${\bf ho}{\cat M} \cong {\bf ho}{\cat M}_{cof}$, we may assume for simplicity that ${\cat M} = {\cat M}_{cof}$.

Let us first show that $\mathrm{dgm}_{_{\cat D}}$ is essentially surjective. For an object $X$ of $({\bf ho}{\cat M})^{\cat D}$, we will construct by induction on degree a Reedy cofibrant diagram $X^{'} \eps {\cat M}^{\cat D}$ and an isomorphism $f \colon \mathrm{dgm}_{_{\cat D}}(X^{'}) \overset{\cong}{\ra} X$. The inductive hypothesis is that $X^{' \le n}$ is Reedy cofibrant and that $f^{\le n} \colon \mathrm{dgm}_{{\cat D}^{\le n}}(X^{' \le n}) \overset{\cong}{\ra} X^{\le n}$ is an isomorphism.

The initial step $n = 0$ follows from \lemmaref{lem:cofcatsumsofdiagrams}, since ${\cat D}^{\le 0}$ is discrete. Assume that $X^{'<n}, f^{< n}$ have been constructed, and let's try to extend them over each object $d \eps {\cat D}^n$. 

By \lemmaref{lem:descrfinitefreecat}, the latching category $\partial({\cat D} \downarrow d)$ is a disjoint sum of categories with terminal objects denoted $d_i \ra d$. We therefore have $\colim^{\partial({\cat D} \downarrow d)} X^{'} \cong \Sum_i X^{'}_{d_i}$, and $\Sum_i X^{'}_{d_i}$ computes the sum of $X^{'}_{d_i}$ also in ${\bf ho}{\cat M}$. We can thus construct a map $\colim^{\partial({\cat D} \downarrow d)} X^{'} \ra X_d$ in ${\bf ho}{\cat M}$, compatible with $f^{< n}$. This map yields a factorization $\colim^{\partial({\cat D} \downarrow d)} X^{'} \rightarrowtail X^{'}_d \overset{\sim}{\leftarrowtail} X_d$ in ${\cat M}$, since we assumed that ${\cat M} = {\cat M}_{cof}$. Define $f_d$ as the inverse of $X^{'}_d \overset{\sim}{\leftarrowtail} X_d$. Repeating the construction of $X^{'}_d, f_d$ for each $d \eps {\cat D}^n$ yields the desired extension $X^{' \le n}, f^{\le n}$. The map $f \colon \mathrm{dgm}_{_{\cat D}}(X^{'}) \overset{\cong}{\ra} X$ we constructed is a degreewise isomorphism, therefore an isomorphism.

We have shown that $\mathrm{dgm}_{_{\cat D}}$ is essentially surjective, let us now show that it is full. Using \lemmaref{lem:sufficientfinitefull}, it suffices to show that for any Reedy cofibrant diagrams $X^{'}, Y^{'} \eps {\cat M}^{\cat D}$ and map $f \colon \mathrm{dgm}_{{{\cat D}}}(X^{'}) \ra \mathrm{dgm}_{{{\cat D}}}(Y^{'})$, we can construct a set of maps $f^{'}_n \eps Hom_{{\bf ho}{{\cat M}^{{\cat D}^{\le n}}}} (X^{\le n}, Y^{\le n})$ for $n \ge 0$ such that $f^{'}_{n}$ restricts to $f^{'}_{n-1}$ and $f$ restricts to $\mathrm{dgm}_{{\cat D}^{\le n}}(f^{'}_n)$.

We will construct such a map $f^{'}_n$ by induction on $n$. The initial step map $f^{'}_0$ exists as a consequence of \lemmaref{lem:cofcatsumsofdiagrams}. Assume that $f^{'}_{n-1}$ has been constructed.

Using \theoremref{thm::descrhomotcof} we can construct a factorization in ${\cat D}^{< n}$
\begin{center}
$\xymatrix{
    X^{' < n} \ar@{>->}[r]^-{a_{n-1}} & 
    Z^{'}_{n-1} &
    Y^{' < n} \ar@{>->}[l]^-{\sim}_-{b_{n-1}}
  }$
\end{center}
with $a_{n-1}$ a Reedy cofibration and $b_{n-1}$ a trivial Reedy cofibration in ${\cat D}^{< n}$, such that $b_{n-1}^{-1}a_{n-1} = f^{'}_{n-1}$.

For any object $d \eps {\cat D}^n$, the latching category $\partial({\cat D} \downarrow d)$ is a disjoint sum of categories with terminal objects, and the functor ${\cat M} = {\cat M}_{cof} \ra {\bf ho}{\cat M}$ commutes with sums of objects. As a consequence, the diagram 
\begin{center}
$\xymatrix{
    LX^{'}_d \ar@{>->}[rrr]^-{\colim^{\partial({\cat D} \downarrow d)} a_{n-1}} \ar[d] &&&
    \colim^{\partial({\cat D} \downarrow d)} Z^{'}_{n-1} &&&
    LY^{'}_d \ar@{>->}[lll]^-{\sim}_-{\colim^{\partial({\cat D} \downarrow d)} b_{n-1}} \ar[d] \\
    X^{'}_d \ar[rrrrrr]^{f_d} &&&&&& Y^{'}_d 
  }$
\end{center}
commutes in ${\bf ho}{\cat M}$. All maps of this diagram are maps of ${\cat M}$, with the exception of $f_d$ which is a map of ${\bf ho}{\cat M}$. Since ${\cat M}$ admits a homotopy calculus of left fractions, a quick computation shows that this diagram can be embedded in a homotopy commutative diagram with all maps in ${\cat M}$
\begin{center}
$\xymatrix{
    LX^{'}_d \ar@{>->}[rrr]^-{\colim^{\partial({\cat D} \downarrow d)} a_{n-1}} \ar[d] &&&
    \colim^{\partial({\cat D} \downarrow d)} Z^{'}_{n-1} \ar@{>-->}[d] &&&
    LY^{'}_d \ar@{>->}[lll]^-{\sim}_-{\colim^{\partial({\cat D} \downarrow d)} b_{n-1}} \ar[d] \\
    X^{'}_d \ar@{>-->}[rrr] &&& {S_d} &&& Y^{'}_d \ar@{>-->}[lll]^\sim
  }$
\end{center}
with the composition in ${\bf ho}{\cat M}$ of the bottom edge having the homotopy type of $f_d$.

Since the cylinders are with respect to the Reedy cofibration structure, notice that $LIX^{'}_d$, $LIY^{'}_d$ are cylinders of $(LX^{'})_d$ and $(LY^{'})_d$. Using \lemmaref{lem::choosecyl} twice, we can embed our homotopy commutative diagram into a diagram commutative on the nose
\begin{center}
$\xymatrix{
    LIX^{'}_d \ar@{-->}[ddrr]_{t_d} & LX^{'}_d \ar@{>->}[l]^-{\sim}_-{i_0} \ar@{>->}[r] &
    \colim^{\partial({\cat D} \downarrow d)} Z^{'}_{n-1} \ar@{>->}[d] & LY^{'}_d \ar@{>->}[l]^-{\sim} \ar@{>->}[r]_-{\sim}^-{i_0} & LIY^{'}_d \ar@{-->}[ddll]^{t^{'}_d} \\
    LX^{'}_d \ar@{>->}[u]_-{\sim}^-{i_1} \ar[d] && {S_d} \ar@{>-->}[d]^\sim && LY^{'}_d \ar@{>->}[u]^-{\sim}_-{i_1} \ar[d] \\
    X^{'}_d \ar@{>-->}[rr] && {T_d} && Y^{'}_d \ar@{>-->}[ll]^\sim
  }$
\end{center}
where the bottom edge has the homotopy type of $f_d$. 

Denote $z_d$ the map $\colim^{\partial({\cat D} \downarrow d)} Z^{'}_{n-1} \ra T_d$. We repeat the construction above for all $d \eps {\cat D}^n$. The maps $t_d, t_di_0, z_d, t^{'}_di_0, t^{'}_d$ define Reedy cofibrant extensions $\overline{LIX^{' \le n}}$, $\overline{LX^{' \le n}}$, $\overline{Z^{'}_{n-1}}$, $\overline{LY^{' \le n}}$, $\overline{LIY^{' \le n}}$ over ${\cat D}^{\le n}$ and a zig-zag of maps 
\begin{center}
$X^{' \le n} \ra \overline{LIX^{' \le n}} \overset{\sim}{\la} \overline{LX^{' \le n}} \ra \overline{Z^{'}_{n-1}} \overset{\sim}{\la} \overline{LY^{' \le n}} \ra \overline{LIY^{' \le n}} \overset{\sim}{\la} Y^{' \le n}$
\end{center}
which defines the desired $f^{'}_n$.
\end{proof}

Suppose that ${\cat C}$ is a category and $u \colon {\cat D}_1 \ra {\cat D}_2$ is a functor. Let $X \eps {\cat C}^{\cat D}_1$ be a diagram. A {\it weak colimit}\index{weak (co)limit} of $X$ along $u$, if it exists, is a diagram $\wcolim^u X \eps {\cat C}^{\cat D}_2$ along with a map $\eps \colon X \ra u^* \wcolim^u X$ with the property that for any other diagram $Y \eps {\cat D}_2$ and map $f \colon X \ra u^*Y$ there exists a (not necessarily unique) map $g \colon \wcolim^u X \ra Y$ such that $f = u^* g \comp \eps$. A weak colimit $(\wcolim^u X \eps)$, if it exists, is not necessarily unique. If the base category ${\cat C}$ is cocomplete, then the colimits along $u$ are weak colimits. Weak limits are defined in a dual fashion, and are denoted $\wlim^u X$.

The homotopy category of a cofibration category is not cocomplete, in general. As a benefit of \theoremref{thm:freedirectdiagramrealization}, the homotopy category of a cofibration category admits weak colimits indexed by diagrams that satisfy the conclusion of \theoremref{thm:freedirectdiagramrealization}. Furthermore, there is a way to make these weak colimits {\it unique} up to a {\it non-canonical} isomorphism. This is an idea we learned from Alex Heller, \cite{Heller}.

Suppose that ${\cat M}$ is a cofibration category and that $u \colon {\cat D}_1 \ra {\cat D}_2$ is a functor of small categories where ${\cat D}_1$ has the property that $\mathrm{dgm}_{{\cat D}_1}$ is full and essentially surjective. Let $X \eps ({\bf ho}{\cat M})^{{\cat D}_1}$ be a diagram in the homotopy category. There exists $X^{'} \eps {\bf ho}({\cat M}^{{\cat D}_1})$ with $\mathrm{dgm}_{{\cat D}_1}X^{'} \cong X$, using the essential surjectivity of $\mathrm{dgm}_{{\cat D}_1}$. Denote $\eps^{'} \colon X^{'} \ra u^* {\bf L}\colim^u X^{'}$ the canonical map. Since $\mathrm{dgm}_{{\cat D}_1}$ is full, we see that 
\[
(\mathrm{dgm}_{{\cat D}_1}({\bf L}\colim\nolimits^u X^{'}), \mathrm{dgm}_{{\cat D}_1}\eps^{'}) 
\] 
is a weak colimit of $X$ along $u$. 

We denote ${\bf W}\colim^u X$ the weak colimit that we constructed. If $X^{''} \eps {\bf ho}({\cat M}^{{\cat D}_1})$ also satisfies $\mathrm{dgm}_{{\cat D}_1}X^{''} \cong X$, since $\mathrm{dgm}_{{\cat D}_1}$ is full we can construct a {\it non-canonical} map $f \colon X^{'} \ra X^{''}$ with $\mathrm{dgm}_{{\cat D}_1}f$ the isomorphism $\mathrm{dgm}_{{\cat D}_1}X^{'} \cong \mathrm{dgm}_{{\cat D}_1}X^{''}$. By \corollaryref{cor:homotcatconservative}, $\mathrm{dgm}_{{\cat D}_1}$ is conservative therefore $f$ is an isomorphism. This shows that ${\bf W}\colim^u X$ as constructed is unique {\it up to a non-canonical isomorphism}.

The weak colimit ${\bf W}\colim^u X$ will be called a {\it privileged weak colimit}, following Alex Heller's terminology. Dually, if ${\cat M}$ is a fibration category and $u \colon {\cat D}_1 \ra {\cat D}_2$ is a functor of small categories with $\mathrm{dgm}_{{\cat D}_1}$ full and essentially surjective, we can construct privileged weak limits of $X \eps ({\bf ho}{\cat M})^{{\cat D}_1}$ along $u$, which are denoted ${\bf W}\lim^u X$. 
\begin{thm}
\label{thm:existenceweakcolimits}
\mbox{}
\begin{enumerate}
\item
Suppose that ${\cat M}$ is a cofibration category and that $u \colon {\cat D}_1 \ra {\cat D}_2$ is a functor of small categories with ${\cat D}_1$ direct and free. Then the privileged weak colimits ${\bf W}\colim^u X$ exist for $X \eps ({\bf ho}{\cat M})^{{\cat D}_1}$, and are unique up to a non-canonical isomorphism.
\item
Suppose that ${\cat M}$ is a fibration category and that $u \colon {\cat D}_1 \ra {\cat D}_2$ is a functor of small categories with ${\cat D}_1$ inverse and free. Then the privileged weak limits ${\bf W}\lim^u X$ exist for $X \eps ({\bf ho}{\cat M})^{{\cat D}_1}$, and are unique up to a non-canonical isomorphism.
\end{enumerate}
\end{thm}

\begin{proof}
Consequence of the fact that by \theoremref{thm:freedirectdiagramrealization}, for both (1) and (2) the functor $\mathrm{dgm}_{{\cat D}_1}$ is full and essentially surjective.
\end{proof}

\chapter{Derivators}
\label{chap:derivators}

In this chapter, we collect all our previous results and interpret them to say that a cofibration category has a canonically associated {\it left Heller derivator}, and dually a fibration category has a canonically associated {\it right Heller derivator}. 

Derivators are best thought of as an axiomatization of homotopy (co)limits. They were introduced by Alexandre Grothendieck in his manuscripts \cite{Grothendieck1}, \cite{Grothendieck2}. Don Anderson \cite{Anderson2} and Alex Heller \cite{Heller} present alternative definitions of derivators. For an elementary introduction to the theory of derivators, the reader is invited to refer to Georges Maltsiniotis \cite{Maltsiniotis1}.

Grothendieck's original definition had derivators as contravariant on both 1- and 2- cells; we will refer to this type of derivators as {\it Grothendieck derivators}. Heller's definition had derivators as contravariant on 1- cells but covariant on 2-cells; we will refer to them as {\it Heller derivators}, although the axioms we use are set up a little bit differently than in \cite{Heller}. The axioms are arranged so that Grothendieck and Heller derivators correspond to each other by duality.

We will recall the basic definitions we need. Fix a pair of universes ${\cat U} \subset {\cat U}^{'}$. We denote $2{\cat Cat}$ the 2-category of ${\cat U}$-small categories, and $2CAT$ the 2-category of ${\cat U}^{'}$-small categories. We denote $\initialcat$ the initial category (having an empty set of objects), and $\terminalcat$ the terminal category (having one object and the identity of that object as the only map). We write $p_{\cat D} \colon {\cat D} \ra \terminalcat$ for the functor to the terminal category and $e_{d} \colon \terminalcat \ra {\cat D}$ for the functor that embeds the terminal category $\terminalcat$ as the object $d \eps {\cat D}$. 

For a 2-category ${\cat C}$, we denote ${\cat C}^{1-op}$ (${\cat C}^{2-op}$, resp. ${\cat C}^{1,2-op}$) the 2-category obtained reversing the direction of the 1-cells (resp. 2-cells, resp. both 1- and 2-cells) of ${\cat C}$.

We will denote by ${\cat Dia}$ any 2-full 2-subcategory of $2{\cat Cat}$, with the property that its objects (viewed as small categories) include the set of finite direct categories, and are stable under the following category operations:
\begin{enumerate}
\item
small disjoint sums of categories
\item
finite products of categories
\item
stable under taking overcategories and undercategories, i.e. if $f \colon {\cat D}_1$ $\ra$ ${\cat D}_2$ is a functor with ${\cat D}_1, {\cat D}_2 \eps {\cat Dia}$ then $(f \downarrow d_2)$, $(d_2 \downarrow f)$ $\eps {\cat Dia}$ for any object $d_2 \eps {\cat D_2}$
\end{enumerate}
In particular from (1) and (2), we will assume that $\initialcat$, $\terminalcat$ $\eps {\cat Dia}$.

For example, we can take ${\cat Dia}$ to be
\begin{enumerate}
\item
the entire $2{\cat Cat}$, or
\item
the 2-full subcategory $2{\cat FinDirCat}$ whose objects are small disjoint sums of finite direct categories
\end{enumerate}

\section{Prederivators}
\label{sec:preder}

\begin{defn}
\label{defn:hellerpreder}
A {\it Heller prederivator}\index{prederivator} of domain ${\cat Dia}$ is a pseudo 2-functor ${\mathbb D} \colon {\cat Dia}^{1-op} \ra  2CAT$. 
\end{defn}

For any morphism $u \colon {\cat D}_1 \ra {\cat D}_2$ in ${\cat Dia}$ we denote $u^* = {\mathbb D}u \colon {\mathbb D}({\cat D}_2) \ra {\mathbb D}({\cat D}_1)$, and for any natural map $\alpha : u \Ra v$ in ${\cat Dia}$ we denote $\alpha^* = {\mathbb D}\alpha : u^* \ra v^*$.

A Heller prederivator is {\it strict} if it is strict as a 2-functor.

For example, any category with weak equivalences $({\cat M}, {\cat W})$ gives rise to a strict Heller prederivator ${\mathbb D}_{({\cat M}, {\cat W})}$ of domain $2{\cat Cat}$ defined by ${\mathbb D}_{({\cat M}, {\cat W})} ({\cat D}) = {\bf ho}{\cat M}^{{\cat D}}$, where the homotopy category is taken with respect to pointwise weak equivalences ${\cat W}^{\cat D}$. Any natural map $\alpha \colon u \Ra v$ yields a natural map $\alpha^* \colon u^* \Ra v^*$
\begin{center}
  $\xymatrix{
    {\cat D}_1 \ar@/^1pc/[r]^u _\,="a" \ar@/_1pc/[r]_v ^\,="b" & {\cat D_2} && {\bf ho}{\cat M}^{{\cat D}_1}  && {\bf ho}{\cat M}^{{\cat D}_2} \ar@/^1pc/[ll]^{v^*} _\,="c" \ar@/_1pc/[ll]_{u^*} ^\,="d"
    \ar@{=>} "a";"b" ^{\alpha}
    \ar@{=>} "d";"c" ^{\alpha^*}
    }$
\end{center}
where $\alpha^*$ is defined on components as $X_{\alpha d_1}: X_{vd_1} \ra X_{ud_1}$, for diagrams $X \eps {\cat M}^{{\cat D}_1}$.

\section{Derivators}
\label{sec:der}
Suppose that $u \colon {\cat D}_1 \ra {\cat D}_2$ is a functor in ${\cat Dia}$.  A left adjoint for $u^* \colon {\mathbb D}({\cat D}_2) \ra {\mathbb D}({\cat D}_1)$, if it exists, is denoted $u_{!} \colon {\mathbb D}({\cat D}_1) \ra {\mathbb D}({\cat D}_2)$. A right adjoint for $u^*$, if it exists, is denoted $u_{*} \colon {\mathbb D}({\cat D}_1) \ra {\mathbb D}({\cat D}_2)$. If the functors $u_!$ or $u_*$ exist, they are only defined up to unique isomorphism.

Suppose we have a diagram in ${\cat Dia}$ with $\phi \colon vf \Ra gu$ a natural map
\begin{equation}
\label{eqn:2catsquareder}
\xymatrix {
    {\cat D}_1 \ar[r]^u \ar[d]_f & 
    {\cat D}_2 \ar[d]^g _\,="b" \\
    {\cat D}_3 \ar[r]_v ^\,="a" & 
    {\cat D}_4 \ultwocell<\omit>{<0>\phi}
  }
\end{equation}
We apply the 2-functor ${\mathbb D}$ and we get a natural map $\phi^* \colon f^*v^* \Ra u^*g^*$. If the left adjoints $u_!$ and $v_!$ exist, the natural map $\phi^*$ yields by adjunction the {\it cobase change} natural map denoted
\begin{equation}
\label{eqn:2catsquaredermap}
\phi_! \colon u_! f^* \Ra g^* v_!
\end{equation}
If $\phi_!$ is an isomorphism for a choice of left adjoints $u_!, v_!$, then $\phi_!$ is an isomorphism for any such choice of right adjoints $u_!, v_!$. Dually, if we assume that the right adjoints $f^*$ and $g^*$ exist, the {\it base change} morphism associated to $\phi$ is
\begin{equation}
\label{eqn:2catsquaredermapop}
\phi_* \colon v^* g_* \Ra f_* u^*
\end{equation}
and if $\phi_*$ is an isomorphism for a choice of $f_*, g_*$ then $\phi_*$ is an isomorphism for any such choice.

Given a functor $u \colon {\cat D}_1 \ra {\cat D_2}$ and an object $d_2 \eps {\cat D_2}$, in the standard over 2-category diagram of $u$ at $d_2$ (see (\ref{eqn:stdover2cat}))
\begin{equation}
\label{eqn:stdover2catbasechangedgm}
\xymatrix {
    (u \downarrow d_2) \ar[r]^-{p_{(u \downarrow d_2)}} \ar[d]_{i_{u, d_2}} & 
    \terminalcat \ar[d]^{e_{d_2}} _\,="b" \\
    {\cat D}_1 \ar[r]_u ^\,="a" & 
    {\cat D}_2 \ultwocell<\omit>{<0>_{{\phi}_{u, d_2} \;\;\;\;\;\;}}
  }
\end{equation}

In \eqref{eqn:stdover2catbasechangedgm}, if the left adjoints $u_!$ and $p_{(u \downarrow d_2)!}$ exist, the associated cobase change morphism is denoted
\begin{equation}
\label{eqn:stdover2catbasechange}
({\phi}_{u, d_2})_! \colon p_{(u \downarrow d_2)!}i_{u, d_2}^* \Ra e_{d_2}^*u_!
\end{equation}

Dually, for an object $d_2 \eps {\cat D_2}$, in the standard under 2-category diagram of $u$ at $d_2$ (see (\ref{eqn:stdunder2cat}))
\begin{equation}
\label{eqn:stdover2catcobasechangedgm}
\xymatrix {
    (d_2 \downarrow u) \ar[r]^-{i_{d_2, u}} \ar[d]_{p_{(d_2 \downarrow u)}}& 
    {\cat D}_1 \ar[d]^u _\,="b" \\
    \terminalcat \ar[r]_{e_{d_2}} ^\,="a" & 
    {\cat D}_2 \ultwocell<\omit>{<0>_{{\phi}_{d_2, u} \;\;\;\;\;\;}}
  }
\end{equation}

In \eqref{eqn:stdover2catcobasechangedgm}, if the right adjoints $u^*$ and $p_{(d_2 \downarrow u)}^*$ exist, the associated base change morphism is denoted
\begin{equation}
\label{eqn:stdunder2catcobasechange}
({\phi}_{d_2, u})^* \colon e_{d_2}^*u_* \Ra p_{(d_2 \downarrow u)*}i_{d_2, u}^*
\end{equation}

\index{derivator}
\index{derivator!Heller derivator}
\begin{defn}[Heller derivators]
\label{defn:hellerderivator}
\mbox{} \\
A Heller prederivator ${\mathbb D} \colon {\cat Dia}^{1-op} \ra  2CAT$ is a {\it left Heller derivator} if it satisfies the following axioms:
\begin{description}
\item[Der1]
For any set of small categories ${\cat D}_k$, $k \eps K$, the functor
\begin{center}
$\xymatrix {
    {\mathbb D}( {\Sum}_{k \eps K} {\cat D}_k) 
    \ar[rr]^-{(i_k^*)_{k \eps K}} &&
    {\times}_{k \eps K} {\mathbb D}({\cat D}_k)
  }$
\end{center}
is an equivalence of categories, where $i_{k} \colon {\cat D}_{k} \ra {\Sum}_{k^{'} \eps K} {\cat D}_{k^{'}}$ denotes the component inclusions for $k \eps K$.

\item[Der2]
For any ${\cat D}$ in ${\cat Dia}$, the family of functors $e_{d}^* \colon {\mathbb D}({\cat D}) \ra {\mathbb D}(\terminalcat)$ for all objects $d \eps {\cat D}$ is conservative.

\item[Der3l]
For any $u \colon {\cat D}_1 \ra {\cat D}_2$ in ${\cat Dia}$, the functor $u^* \colon {\mathbb D}({\cat D}_2) \ra {\mathbb D}({\cat D}_1)$ admits a left adjoint $u_! \colon {\mathbb D}({\cat D}_1) \ra {\mathbb D}({\cat D}_2)$

\item[Der4l]
For any $u \colon {\cat D}_1 \ra {\cat D}_2$ in ${\cat Dia}$ and any object $d_2 \eps {\cat D}_2$, the cobase change morphism $({\phi}_{u, d_2})_! \colon p_{(u \downarrow d_2)!}i_{u, d_2}^* \Ra e_{d_2}^*u_!$ of (\ref{eqn:stdover2catbasechange}) associated to the standard over 2-category diagram (\ref{eqn:stdover2cat}) is an isomorphism.
\end{description}

${\mathbb D}$ is a {\it right Heller derivator} if it satisfies axioms Der1 and Der2 above and:
\begin{description}
\item[Der3r]
For any $u \colon {\cat D}_1 \ra {\cat D}_2$ in ${\cat Dia}$, the functor $u^* \colon {\mathbb D}({\cat D}_2) \ra {\mathbb D}({\cat D}_1)$ admits a right adjoint $u_* \colon {\mathbb D}({\cat D}_1) \ra {\mathbb D}({\cat D}_2)$

\item[Der4r]
For any $u \colon {\cat D}_1 \ra {\cat D}_2$ in ${\cat Dia}$ and any object $d_2 \eps {\cat D}_2$, the base change morphism $({\phi}_{d_2, u})^* \colon e_{d_2}^*u_* \Ra p_{(d_2 \downarrow u)+}i_{d_2, u}^*$ of (\ref{eqn:stdunder2catcobasechange}) associated to the standard under 2-category diagram (\ref{eqn:stdunder2cat}) is an isomorphism.
\end{description}

${\mathbb D}$ is a {\it Heller derivator} if it is both a left and a right Heller derivator. 
\end{defn}

\index{derivator}
\index{derivator!Grothendieck derivator}
\begin{defn}[Grothendieck derivators]
\label{defn:grothderivator}
A 2-functor ${\mathbb D} \colon {\cat Dia}^{1,2-op} \ra  2CAT$ is a left (right, two-sided) {\it Grothendieck derivator} if the functor ${\cat Dia}^{1-op} \ra  2CAT$, ${\cat D} \mapsto {\mathbb D}({\cat D}^{op})$ is a right (left, two-sided) Heller derivator.
\end{defn}

Notice that for a left Grothendieck derivator the map $u^* = {\mathbb D}(u)$ by definition admits a {\it right} adjoint $u_*$, and for a right Grothendieck derivator ${\mathbb D}(u)$ admits a {\it left} adjoint $u_!$. This observation should be kept in mind when reading \cite{Grothendieck1}, \cite{Grothendieck2}, \cite{Maltsiniotis1}, \cite{Cisinski2}, \cite{Cisinski3}, \cite{Cisinski4}. 

We have allowed for pseudo-functors ${\mathbb D}$ rather than just strict 2-functors in the definition of Heller, resp. Grotherndieck derivators. Grothendieck and Heller had originally just used strict 2-functors ${\mathbb D}$ in the derivator definition, but their arguments apply just as well for pseudo-functors.

The axioms for a Grothendieck derivator that we propose are stronger than the ones originally introduced by Grothendieck in that we require Der1 to hold for arbitrary small sums of categories rather than just finite sums of categories. The axioms for a Heller derivator that we propose are weaker than the axioms of a Heller {\it homotopy theory}, as defined in \cite{Heller}, in that our Heller derivators don't satisfy Heller's axiom H2 (same as our Der5 below).

We will use the term {\it derivator} to refer to Heller derivators, when no confusion is possible.

For ${\cat D}_1$ and ${\cat D}_2$ in ${\cat Dia}$, we define the functor
\begin{center}
$\mathrm{dgm}_{_{{\cat D}_1, {\cat D}_2}} \colon \mathbb{D}({\cat D}_1 \times {\cat D}_2) \ra \mathbb{D}({\cat D}_2)^{{\cat D}_1}$ \\
$(\mathrm{dgm}_{_{{\cat D}_1, {\cat D}_2}} X)_{d_1} = i_{d_1}^* X$ \\
$(\mathrm{dgm}_{_{{\cat D}_1, {\cat D}_2}} X)_{f} = i_{f}^* X$
\end{center}
for any object $d_1$ and map $f \colon d_1 \ra d^{'}_1$ of ${\cat D}_1$, where $i_{d_1} \colon {\cat D}_2 \ra {\cat D}_1 \times {\cat D}_2$ is the functor $d_2 \mapsto (d_1, d_2)$ and $i_f \colon i_{d_1} \Ra i_{d^{'}_1}$ is the natural map given on components by $(f, 1_{{\cat D}_1})$.

\index{derivator!strong}
\begin{defn}[Strong Heller derivators]
\label{defn:hellerstrongderivator}
\mbox{} \\
A left (resp. right, two-sided) Heller derivator ${\mathbb D} \colon {\cat Dia}^{1-op} \ra  2CAT$ is {\it strong} if it satisfies the additional axiom Der5
\begin{description}
\item[Der5]
For any finite, free category ${\cat D}_1$ and category ${\cat D}_2$ in ${\cat Dia}$, the functor $\mathrm{dgm}_{_{{\cat D}_1, {\cat D}_2}}$ is full and essentially surjective.
\end{description}
\end{defn}

Note that in the literature the axiom Der5 is sometimes stated in the weaker form
\begin{description}
\item[Der5w]
For any category ${\cat D}_2$ in ${\cat Dia}$ the functor $\mathrm{dgm}_{_{I, {\cat D}_2}}$ is full and essentially surjective, where $I$ denotes the category $0 \ra 1$, with two objects and one non-identity map.
\end{description}

\comment{
\index{derivator!pointed}
\begin{defn}[Pointed Heller derivators]
\label{defn:hellerpointedderivator}
\mbox{} \\
A left (resp. right, two-sided) Heller derivator ${\mathbb D} \colon {\cat Dia}^{1-op} \ra  2CAT$ is {\it pointed} if it satisfies the additional axiom Der6l (resp. Der6r, resp. Der6l and Der6r).
\begin{description}
\item[Der6l]
For any closed embedding $u \colon {\cat D}_1 \ra {\cat D}_2$ in ${\cat Dia}$, the functor $u_!$ admits a left adjoint denoted $u^?$
\item[Der6r]
For any open embedding $u \colon {\cat D}_1 \ra {\cat D}_2$ in ${\cat Dia}$, the functor $u_*$ admits a right adjoint denoted $u^!$
\end{description}
\end{defn}
} 

\begin{thm}
\label{thm:hellerderishtpytheory}
A Heller prederivator of domain $2{\cat Cat}$ is a {\it homotopy theory} in the sense of \cite{Heller} iff it is a strong, strict Heller derivator in the sense of \definitionref{defn:hellerderivator}.
\end{thm}

\begin{proof}
This is a consequence of the results of \cite{Heller}. 
\end{proof}

In fact, a careful reading of \cite{Heller} shows that most proofs stated there for one-sided homotopy theories are true for one-sided {\it strong} Heller derivators. However we suspect that the theory developed in \cite{Heller} can be reworked from the axioms of a one-sided Heller derivator alone. 

\section{Derivability of cofibration categories}
\label{sec:derivability}
Let us turn back to the Heller prederivator ${\mathbb D}_{({\cat M}, {\cat W})}$ associated to a category with weak equivalences $({\cat M}, {\cat W})$ and introduce the following

\index{derivable category}
\begin{defn}
\label{defn:derivablecat}
A category with weak equivalences $({\cat M}, {\cat W})$ is {\it left Heller derivable} (resp. strongly left Heller derivable, right Heller derivable, etc) over ${\cat Dia}$ if the prederivator ${\mathbb D}_{({\cat M}, {\cat W})}$ is a left (left strong, right etc.) Heller derivator over ${\cat Dia}$.
\end{defn}

A left Heller derivable (resp. strongly left Heller derivable, etc.) category by convention is left Heller derivable over the entire $2{\cat Cat}$. 

With the definitions of derivators and of derivable categories at hand, we can state

\begin{thm}[Cisinski]
\label{thm:cofcatderivator}
\mbox{}
\begin{enumerate}
\item
Any ABC cofibration category is strongly left Heller derivable.
\item
Any ABC fibration category is strongly right Heller derivable.
\item
Any ABC model category is strongly Heller derivable.
\end{enumerate}
\end{thm}

\begin{proof}
Part (2) is dual to (1) and part (3) is a consequence of (1) and (2). Therefore it suffices to prove part (1). 

Part (1) is proved by verifying the derivator axioms as follows:
\begin{enumerate}
\item[-]
Der1 is proved by \lemmaref{lem:cofcatsumsofdiagrams}
\item[-]
Der2 by \theoremref{thm:homotcatconservative}
\item[-]
Der3r by \theoremref{thm:genholim}
\item[-]
Der4r by \theoremref{thm:computerelhocolim}
\item[-]
Der5 in \theoremref{thm:generalpointwisecofstructure} and \theoremref{thm:freedirectdiagramrealization}
\end{enumerate}
\end{proof}

As a corollary of \theoremref{thm:cofcatderivator} and \propositionref{prop:relclosedmodcat} we obtain

\begin{thm}
\label{thm:modelcatderivator}
Any Quillen model category is strongly Heller derivable. $\square$
\end{thm}

In particular, all the results of \cite{Heller} apply to ABC model categories (and therefore to Quillen model categories).

We should mention that our terminology for derivable categories is different than the one used by Denis-Charles Cisinski. His paper \cite{Cisinski2} calls left derivable categories what we call F1-F4 Anderson-Brown-Cisinski fibration categories. The naming anomaly is explained by the fact that \cite{Cisinski2} is written in the language of Grothendieck derivators, and a left Grothendieck derivator corresponds to a right Heller derivator.

\comment{
To illustrate the power of \theoremref{thm:cofcatderivator}, we end with the following application.

\begin{thm}
\label{thm:homotcofinalhomotocolim}
Suppose that $u \colon {\cat D}_1 \ra {\cat D}_2$ is a functor of small categories. Then:
\begin{enumerate}
\item
$u$ is homotopy right cofinal iff for any cofibration category ${\cat M}$ and any diagram $X \eps {\cat M}^{{\cat D}_2}$ the natural map 
\begin{center}
${\bf L}\colim^{{\cat D}_1} \; {\bf ho}u^*X \ra {\bf L}\colim^{{\cat D}_2}X$
\end{center}
is an isomorphism in ${\bf ho}{\cat M}$.
\item
$u$ is homotopy left cofinal iff for any fibration category ${\cat M}$ and any diagram $X \eps {\cat M}^{{\cat D}_2}$ the natural map 
\begin{center}
${\bf R}\lim^{{\cat D}_2}X \ra {\bf R}\lim^{{\cat D}_1} \; {\bf ho}u^*X$ 
\end{center}
is an isomorphism in ${\bf ho}{\cat M}$.
\end{enumerate}
\end{thm}
} 

\section{Odds and ends}
\label{sec:oddsandends}
The author has not been able to complete his research project as thoroughly as he would have liked, but it may be worthwhile to at least state how this body of work may be improved and extended. First of all, it is proved by Denis-Charles Cisinski that 

\begin{thm}[\cite{Cisinski2}]
\label{thm:cofcatderivator2}
\mbox{}
\begin{enumerate}
\item
Any CF1-CF4 cofibration category is strongly left Heller derivable over $2{\cat FinDirCat}$. 
\item
Any F1-F4 fibration category is strongly right Heller derivable over $2{\cat FinDirCat}$.
\end{enumerate}
\end{thm}

Indeed, we have modeled our proofs in a very large part on the arguments of \cite{Cisinski2}.

Second, consider a weaker replacement of axiom CF6.

\begin{description}
\item[CF6']
For any sequence of weak equivalence maps $A_0 \la B_0 \ra A_1 \la B_1 \ra A_2 ...$ such that each $B_n$ is cofibrant, and such that $B_0 \ra A_0$ and all $B_n \Sum B_{n+1} \ra A_n$ are cofibrations, we have that 
\[
A_0 \ra \colim (A_0 \la B_0 \ra A_1 \la B_1 \ra A_2 ...)
\]
is a weak equivalence.
\end{description}

The colimit in CF6' can be shown to always exist using a modified version of \lemmaref{lem:colimpushout}. If  $A^{'}_0 \la B^{'}_0 \ra A^{'}_1 ...$ is another sequence of weak equivalence maps such that each $B_n$ is cofibrant, and such that $B^{'}_0 \ra A^{'}_0$ and all $B^{'}_n \Sum B^{'}_{n+1} \ra A^{'}_n$ are cofibrations, and if $a_n \colon A_n \ra A^{'}_n$, $b_n \colon B_n \ra B^{'}_n$ form a {\it diagram} map, then using a modified version of \theoremref{thm:colimdirect} one can show that the induced map
\[
\colim (A_0 \la B_0 \ra A_1 \la B_1 \ra A_2 ...) \ra \colim (A^{'}_0 \la B^{'}_0 \ra A^{'}_1 \la B^{'}_1 \ra A^{'}_2 ...)
\]
is a weak equivalence.

We conjecture that a precofibration category $({\cat M}, {\cat W}, {\cat Cof})$ satisfying CF5 and CF6' is strongly left Heller derivable. It is an open question actually if axiom CF6' is needed at all - so an even stronger conjecture would be that axioms CF1-CF5 imply strong left Heller derivability.

The dual axiom for fibration categories is 

\begin{description}
\item[F6']
For any sequence of weak equivalence maps $A_0 \ra B_0 \la A_1 \ra B_1 \la A_2 ...$ such that each $B_n$ is fibrant, and such that $A_0 \ra B_0$ and all $A_{n+1} \ra B_n \times B_{n+1}$ are fibrations, we have that 
\[
\lim (A_0 \ra B_0 \la A_1 \ra B_1 \la A_2 ...) \ra A_0
\]
is a weak equivalence.
\end{description}

The dual weak conjecture is that a prefibration category $({\cat M}, {\cat W}, {\cat Fib})$ satisfying F5 and F6' is strongly right Heller derivable. The dual strong conjecture is that a prefibration category satisfying F5 is strongly right Heller derivable.

\backmatter
\bibliographystyle{amsalpha}
\bibliography{bibliography}
\printindex

\end{document}